\documentclass[reqno]{amsart}
\usepackage[margin = 1in]{geometry}
\usepackage{amsmath, amssymb, amsthm, fancyhdr, verbatim, graphicx}
\usepackage[notref,notcite]{}
\usepackage{enumerate}
\usepackage[all]{xy}
\usepackage[usenames,dvipsnames]{xcolor}
\usepackage{mathrsfs}
\usepackage{tikz-cd}
\usetikzlibrary{decorations.pathmorphing}

\usepackage[hyperindex,breaklinks]{hyperref}

\usepackage[T1]{fontenc} 
\usepackage{framed}
\usepackage[OT2, T1]{fontenc}
\usepackage[titletoc]{appendix}
\usepackage{bbm} 
\usepackage{stmaryrd}

\usetikzlibrary{calc} 

\usepackage{adjustbox, marvosym, rotating}

\input xy
\xyoption{all}

\newcommand*\leftdash{\rotatebox[origin=c]{-45}{$\dabar@\dabar@\dabar@$}}
\newcommand*\rightdash{\rotatebox[origin=c]{45}{$\dabar@\dabar@\dabar@$}}

\usepackage{amssymb}

\numberwithin{equation}{subsection}

\DeclareSymbolFont{cyrletters}{OT2}{wncyr}{m}{n}
\DeclareMathSymbol{\Sha}{\mathalpha}{cyrletters}{"58}

\usepackage{color}
\newcommand{\tony}[1]{{\color{blue} \sf
    $\spadesuit\spadesuit\spadesuit$ TONY: [#1]}}
    
\newcommand{\Bao}[1]{{\color{red} \sf
    $\spadesuit\spadesuit\spadesuit$ Bao: [#1]}}

\newcommand\muSh{\mu\textup{Sh}}
\newcommand\musupp{\mrm{SS}}
\newcommand\muloc{\mu\mrm{Loc}}

\newcommand{\F}{\mathbf{F}}

\newcommand{\C}{\mathbf{C}}
\newcommand{\G}{\mathbf{G}}
\newcommand{\tr}[0]{\operatorname{tr}}
\newcommand{\wt}[1]{\widetilde{#1}}

\newcommand{\Q}{\mathbf{Q}}
\newcommand{\Z}{\mathbf{Z}}

\newcommand{\mf}[1]{\mathfrak{#1}}

\newcommand{\sgn}{\operatorname{sgn}}
\newcommand{\Gal}{\operatorname{Gal}}
\newcommand{\Mat}{\operatorname{Mat}}
\newcommand{\cl}{\overline}

\newcommand{\R}{\mathbf{R}}

\newcommand{\ul}[1]{{#1}}
\newcommand{\ol}[1]{\overline{#1}}
\newcommand{\wh}[1]{\widehat{#1}}
\newcommand{\PP}{\mathbf{P}}
\newcommand{\mbb}[1]{\mathbb{#1}}

\newcommand{\A}{\mathbf{A}}

\newcommand{\mbf}[1]{\mathbf{#1}}

\newcommand{\co}{\colon}
\newcommand{\mrm}[1]{\mathrm{#1}}

\newcommand{\bs}{\backslash}
\newcommand{\dtimes}{\stackrel{\mrm{L}}\times}

\newcommand{\bu}{\bullet}
\newcommand{\bup}{\bullet_p}
\newcommand{\inj}{\hookrightarrow}
\newcommand{\surj}{\twoheadrightarrow}
\newcommand{\sph}{\mbb{S}}
\newcommand{\mBM}{\mrm{H}^{\mrm{BM}}}
\newcommand{\mBMg}{\mrm{H}^{\mrm{BM}}}
\newcommand{\mBMT}{\mrm{H}^{\mrm{BM, \chT}}}
\newcommand{\mBMulT}{\rH_{\mrm{top}}^{\BM, \ul{\chT}}}
\newcommand{\fsp}{\mathfrak{sp}}

\newcommand{\dWext}{\wt{W}}
\newcommand{\dWaff}{W_{\aff}}

\newcommand{\tw}[1]{\langle #1 \rangle}
\newcommand{\D}{\mbf{D}}

\newcommand\cA{\mathcal{A}}
\newcommand\cB{\mathcal{B}}
\newcommand\cC{\mathcal{C}}
\newcommand\cD{\mathcal{D}}
\newcommand\cE{\mathcal{E}}
\newcommand\cF{\mathcal{F}}
\newcommand\cG{\mathcal{G}}

\newcommand\cL{\mathcal{L}}
\newcommand\cM{\mathcal{M}}
\newcommand\cN{\mathcal{N}}
\newcommand\cO{\mathcal{O}}

\newcommand\cR{\mathcal{R}}
\newcommand\cS{\mathcal{S}}

\newcommand\cU{\mathcal{U}}

\newcommand\cX{\mathcal{X}}
\newcommand\cY{\mathcal{Y}}
\newcommand\cZ{\mathcal{Z}}

\newcommand{\rH}{\mathrm{H}}
\newcommand{\rK}{\mathrm{K}}
\newcommand{\rT}{\mathrm{T}}

\newcommand{\rX}{\mathrm{X}}
\newcommand{\rY}{\mathrm{Y}}

\newcommand{\BB}{\mathbb{B}}
\newcommand{\DD}{\mathbb{D}}

\newcommand{\II}{\mathbb{I}}

\newcommand{\chB}{\check{B}}
\newcommand{\chG}{\check{G}}
\newcommand{\lsup}[2]{{}^{#1}\hspace{-0.6mm}#2}
\newcommand{\LG}{\lsup LG}
\newcommand{\chK}{\mbf{\check{K}}}
\newcommand{\chN}{\check{N}}
\newcommand{\chT}{\check{T}}
\newcommand{\ulchT}{\ul{\chT}}
\newcommand{\chI}{\mbf{\check{I}}}

\newcommand{\chcG}{\check{\cG}}

\newcommand{\chP}{\check{P}}
\newcommand{\chQ}{\check{Q}}

\newcommand{\chPhi}{\check{\Phi}}
\newcommand{\chDelta}{\check{\Delta}}

\newcommand{\ulC}{\ul{C}}
\newcommand{\ulG}{\ul{G}}
\newcommand{\ulT}{\ul{T}}
\newcommand{\ulW}{\ul{W}}
\newcommand{\ulX}{\ul{\rX}}
\newcommand{\ulY}{\ul{\rY}}

\newcommand{\chfb}{\check{\mf{b}}}
\newcommand{\chfg}{\check{\mf{g}}}

\newcommand{\chfn}{\check{\mf{n}}}
\newcommand{\chft}{\check{\mf{t}}}

\DeclareMathOperator{\GL}{GL}
\DeclareMathOperator{\GSp}{GSp}

\DeclareMathOperator{\Frob}{Frob}

\DeclareMathOperator{\N}{\mathbf{N}}

\DeclareMathOperator{\Rep}{Rep}

\DeclareMathOperator{\Spec}{Spec\,}

\DeclareMathOperator{\Lie}{Lie}

\DeclareMathOperator{\ad}{ad}

\DeclareMathOperator{\Res}{Res}
\DeclareMathOperator{\Frac}{Frac}

\DeclareMathOperator{\bupn}{Bun}

\DeclareMathOperator{\Id}{Id}
\DeclareMathOperator{\Ad}{Ad}
\DeclareMathOperator{\Gr}{Gr}

\DeclareMathOperator{\Sym}{Sym}

\DeclareMathOperator{\dom}{dom}
\DeclareMathOperator{\pt}{pt}

\DeclareMathOperator{\len}{len}

\DeclareMathOperator{\Fl}{Fl}

\DeclareMathOperator{\ev}{ev}
\DeclareMathOperator{\Weil}{Weil}

\DeclareMathOperator{\pr}{pr}
\DeclareMathOperator{\red}{red}

\DeclareMathOperator{\reg}{reg}

\DeclareMathOperator{\Fr}{Fr}

\DeclareMathOperator{\EG}{EG}
\DeclareMathOperator{\Ch}{CH}

\DeclareMathOperator{\Spf}{Spf}
\DeclareMathOperator{\Bun}{Bun}

\DeclareMathOperator{\JH}{JH}
\DeclareMathOperator{\aff}{aff}

\DeclareMathOperator{\Adm}{Adm}
\DeclareMathOperator{\Conv}{Conv}
\DeclareMathOperator{\St}{St}

\DeclareMathOperator{\Coh}{Coh}
\DeclareMathOperator{\Loc}{Loc}
\DeclareMathOperator{\BM}{BM}

\DeclareMathOperator{\transfer}{transfer}

\newcommand{\zhc}{\mf{Z}_{\mrm{HC}}}
\newcommand{\zp}{\mf{Z}_{\Fr}}

\newcommand{\colim}{\varinjlim}

\newcommand{\Ql}{\Q_\ell}
\newcommand{\Qll}[1]{\Q_{\ell, #1}}
\newcommand{\spc}{\mathfrak{sp}}

\newcommand{\Ug}{\cU\mf{g}}

\newcommand{\topBM}{\rH^{\BM}_{\mrm{top}}}
\newcommand{\topBMg}{\rH^{\BM}_{\mrm{top}}}
\newcommand{\dZ}{Z_d}
\newcommand{\topCh}{\Ch_{\mrm{top}}}

\newcommand{\topBMT}{\rH^{\BM, \chT}_{\mrm{top}}}
\newcommand{\topBMulT}{\rH^{\BM, \ul{\chT}}_{\mrm{top}}}
\newcommand{\topBMTg}{\rH^{\BM,\chT}_{\mrm{top}}}
\newcommand{\topom}{\omega_{\mrm{top}}}

\newcommand{\nocontentsline}[3]{}
\newcommand{\tocless}[2]{\bgroup\let\addcontentsline=\nocontentsline#1{#2}\egroup}

\newcommand\semis{\mrm{ss}}

\newcommand{\nv}{\text{nv}}

\newcommand{\bG}{\check{G}}
\newcommand{\bP}{\check{P}}
\newcommand{\bN}{\check{N}}
\newcommand{\bU}{\check{U}}
\newcommand{\bM}{\check{M}}

\newcommand{\bB}{\check{B}}
\newcommand{\bT}{\check{T}}

\newcommand{\Inf}{{\operatorname{inf}}}

\newcommand{\Img}{{\operatorname{Im}}}
\newcommand{\dlog}{{\operatorname{dlog}}}

\newtheorem{thm}{Theorem}[subsection]
\newtheorem{lemma}[thm]{Lemma}
\newtheorem{prop}[thm]{Proposition}
\newtheorem{cor}[thm]{Corollary}

\newtheorem{conj}[thm]{Conjecture}

\theoremstyle{remark}
\newtheorem{remark}[thm]{Remark} 
\newtheorem{reminder}[thm]{Reminder} 
\newtheorem{defn}[thm]{Definition}
\newtheorem{const}[thm]{Construction}
\newtheorem{example}[thm]{Example}

\newtheorem{warning}[thm]{Warning}

\newtheorem{situation}[thm]{Situation}
\newtheorem{notation}[thm]{Notation}
\newtheorem{hypothesis}[thm]{Hypothesis}

\makeatletter
\def\th@remark{%

  \thm@headfont{\bfseries}%
  \normalfont 
  \thm@preskip \thm@preskip 
  \thm@postskip\thm@preskip
}
\def\imod#1{\allowbreak\mkern5mu({\operator@
font mod}\,\,#1)}
\makeatother

\numberwithin{equation}{subsection}

\title[Mirror symmetry and the Breuil-M\'{e}zard Conjecture]{Mirror symmetry and the Breuil-M\'{e}zard Conjecture}

\author{Tony Feng and Bao Le Hung}

\begin{document}

\begin{abstract} 
The Breuil-M\'{e}zard Conjecture predicts the existence of hypothetical ``Breuil--M\'ezard cycles'' in the moduli space of mod $p$ Galois representations of $\Gal(\ol{\Q}_p/\Q_p)$ that should govern congruences between mod $p$ automorphic forms. We formulate a generalization of the Breuil--M\'ezard Conjecture for general unramified groups over $\Q_p$, and prove it for sufficiently generic parameters. More specifically, we construct Breuil-M\'{e}zard cycles for generic Serre weights, and verify that they satisfy the Breuil-M\'{e}zard Conjecture for all sufficiently generic tame types and small Hodge-Tate weights. 

Our results are uniform, and apply in particular to exceptional groups (which were previously inaccessible), while also improving significantly upon the previous state-of-the-art results in the previously accessible case of general linear groups. Our method is group-theoretic and purely local, and completely distinct from previous approaches to the Breuil--M\'ezard Conjecture. It leverages an analogy between the Breuil--M\'ezard Conjecture and homological mirror symmetry.
\end{abstract}

\maketitle

\tableofcontents

\section{Introduction}
This paper introduces a new approach to the \emph{Breuil-M\'{e}zard Conjecture}, in its refined form due to Emerton-Gee \cite{EG23}. The main features of our approach, as we see them, are the following.
\begin{itemize}
\item \emph{Its generality.} We will prove results towards the Breuil--M\'ezard Conjecture for general unramified groups over $\Q_p$, whereas previous results were confined to (unramified Weil restrictions of) general linear groups, and $\GSp_4$. We note that this \textbf{improves upon previous versions of this paper}, which handled only general linear groups (see \S \ref{ssec:comparison-versions} for discussion). Even in the case of general linear groups, our results improve substantially upon the previous state-of-the-art for rank larger than $2$. 

\item \emph{Its locality.} The Breuil--M\'ezard Conjecture is expected to be a numerical shadow of a conjectural $p$-adic Local Langlands correspondence, for which prior evidence has been limited to very specific groups, or dependent upon inaccessible conjectures about independence of global realization. By contrast, our arguments are both purely local and completely uniform in the group, giving an unprecedented type of evidence for the $p$-adic Local Langlands correspondence. 

\item \emph{Its connection to other themes within mathematics}. We will connect the Breuil--M\'ezard Conjecture to threads within geometric representation theory and homological mirror symmetry, giving a new entry in the old analogy between $\Q_p$ and topological surfaces. Indeed, we will use a particular (known) case of such homological mirror symmetry as a crucial step in our construction of Breuil--M\'ezard cycles for unramified groups. We remark that this connection extends beyond the present paper: in a sequel work, we will see that the Breuil--M\'ezard Conjecture for \emph{ramified} groups is related to other cases of homological mirror symmetry, which are still conjectural. 
\end{itemize}

 Before formulating the Breuil--M\'ezard Conjecture and stating our main results, we recall some motivation and context. Since much of the paper operates in realms outside number theory, the introduction will be aimed at a somewhat broader audience than usual.

\subsection{Motivation: Serre's Conjectures} The phenomenon of \emph{congruences between modular forms} underpins many facets of modern algebraic number theory. It is therefore natural to try to classify the possible congruences between modular forms. The mod $p$ reductions of modular forms are organized by representation-theoretic parameters called \emph{Serre weights}, which are irreducible representations of $\GL_2(\F_p)$ over $\ol{\F}_p$. The \emph{weight part of Serre's Conjecture} \cite{Ser87} addresses the question of when two mod $p$ eigenforms with different Serre weights can be congruent to each other. It was proved in the early 1990s as the culmination of work by many authors; see the introduction of \cite{GLS15} for references. 

There is a much more general notion of complex-valued automorphic form on a reductive group $G$ over a number field, wherein modular forms constitute the special case $G = \GL_2$ (over $\Q$), to which one would like to generalize Serre's Conjecture. It is subtler to formulate the notion of \emph{mod $p$} automorphic form for general reductive groups, but one option is to proceed as follows. By the Eichler-Shimura relation, mod $p$ modular forms can be interpreted as classes in the cohomology of modular curves, with coefficients in the local systems induced by Serre weights. More generally, let $G$ be a reductive group over a number field, with good reduction at a prime $p$. Then the irreducible representations of $G(\F_p)$ are the \emph{Serre weights} of $G$ (at $p$). Each Serre weight $\sigma$ induces a mod $p$ local system on the locally symmetric space associated to $G$, and it is natural to regard the Hecke eigensystems in the cohomology of this local system as the analogue of mod $p$ eigenforms on $G$ with weight $\sigma$. One can then ask to classify the possible congruences between different weights. This problem has come to also be called the ``weight part of Serre's conjectures'', and beyond $\GL_2$ it becomes much more complicated. At present, the only proposed answer that applies in general (even conjecturally) is itself contingent upon another conjecture, that we will describe next. 

\subsection{The Breuil-M\'{e}zard Conjecture} 

Let $G$ be an unramified reductive group over $\Q_p$. Certain technical conditions are imposed upon $G$ in the main text, not as a restriction of our method but because of genuine group-theoretic subtleties (e.g., pertaining to the $C$-group versus the $L$-group), which we ignore in this Introduction. We regard the dual group $\chG$ as a reductive group over some sufficiently large finite extension $\cO/\Z_p$, and denote by $\LG  = \chG \rtimes \Gal(\ol \Q_p/\Q_p)$ its $L$-group. 

\subsubsection{The Emerton-Gee stack} Emerton-Gee have constructed in \cite{EG23} a formal algebraic stack, which is roughly meant to be a moduli stack of $p$-adic representations of $\Gal(\ol{\Q}_p/\Q_p)$. Work of Lin \cite{Lin23a, Lin23b} and Min \cite{Min24} generalizes the Emerton--Gee stack to reductive groups.\footnote{We are not aware that the compatibility of their constructions is currently discussed in the literature, so let us clarify that in this paper, we will only invoke Lin's work, which applies only for tamely ramified groups.} We let $\cX^{\EG}_{\LG}$ be the Emerton--Gee stack corresponding to $\LG$, which morally parametrizes $p$-adic $L$-parameters $\Gal(\ol{\Q}_p/\Q_p) \rightarrow  \LG$. We regard $\cX^{\EG}_{\LG}$ as defined over $\Spf \cO$ for $\cO \subset E$ the ring of integers in a sufficiently large finite extension $E/\Q_p$, and write $\F$ for the residue field of $\cO$. The reduced substack $(\cX_{\LG}^{\EG})_{\red}$ is an equidimensional algebraic stack over $\F$. 

\subsubsection{Potentially crystalline substacks} Let $\lambda \in X^*(T)^+$ be a dominant weight and $\tau$ be an \emph{inertial parameter} for $G$, which is an admissible representation from the inertia group of $\Q_p$ to $\LG$ that extends to $\Weil(\ol{\Q}_p/\Q_p)$. Then there is a substack $\cX^{\lambda, \tau} \inj \cX^{\EG}$, whose $\ol{\Q}_p$ points correspond to Galois representations with potentially crystalline Hodge-Tate weights $\lambda$ and Weil-Deligne inertial parameter $\tau$ (the point is that these are the local conditions that one expects to see on the local Galois representations associated to algebraic automorphic forms with infinitesimal character $\lambda-\rho$\footnote{For a semisimple simply connected $G$ we would take $\rho$ to be the usual half sum of the positive roots. For a reductive group $G$ with simply connected derived subgroup, we may take any extension of the $\rho$ for its derived group. For $G = \GL_n$ we prefer the choice $\rho = (n-1, n-2, \ldots, 0)$ for reasons of convention.} and ``level'' $\tau$, by $p$-adic Hodge theory). The construction of this substack is due to Emerton--Gee for $G = \GL_n$, and was generalized to tame groups in \cite{Lin23b}.

\subsubsection{The Breuil--M\'ezard Conjecture} Suppose for the moment that $G = \GL_n$. Then to each inertial parameter $\tau$, the inertial Local Langlands correspondence associates an \emph{inertial type} $\sigma(\tau)$, which is a smooth finite-dimensional representation of $G(\Z_p)$ over a finite extension of $\Q_p$. Fix an algebraic closure $k$ of $\F_p$. For a finite-dimensional representation $R$ of $G(\Z_p)$ over a finite extension of $\Q_p$, we let $[\ol{R}] \in K_0(\Rep_k ( G(\Z_p)))$ be the class of its reduction modulo $p$; note that since the kernel of $G(\Z_p) \rightarrow G(\F_p)$ is pro-$p$, we may regard $[\ol{R}] \in K_0(\Rep_k ( G(\F_p))) \xrightarrow{\sim} K_0(\Rep_k(G(\Z_p)))$. Let $W(\lambda) \in \Rep_{\Q_p}(G(\Z_p))$ be the Weyl module of highest weight $\lambda$. 

For every $\lambda$ and $\tau$, we have 
\[
\dim \cX^{\lambda+\rho, \tau}_{\F} = \dim \cX^{\EG}_{\red} = n(n-1)/2,
\]
and we denote this common dimension by $d$. (For general $G$, $d$ is the dimension of the flag variety of $\chG$.) Let $\dZ(\cX^{\EG}_{\red})$ be the group of $d$-dimensional algebraic cycles in $\cX^{\EG}_{\red}$ with $\Q$-coefficients. (Since $d = \dim \cX^{\EG}_{\red}$, this is the same as $\Ch_d(\cX^{\EG}_{\red})$.)

\begin{conj}[Geometric Breuil-M\'{e}zard Conjecture]\label{conj: intro BM}
Suppose $G = \GL_n$. Then there is a map 
\[
\cZ \co K_0( \Rep_k(G(\F_p))) \rightarrow \dZ(\cX^{\EG}_{\red})
\]
such that for every $\lambda \in X^*(T)^+$ and every inertial parameter $\tau$, we have
\begin{equation}\label{eq: intro BM condition}
\cZ [\ol{W(\lambda) \otimes \sigma(\tau)}] = [\cX^{\lambda+\rho, \tau}_{\F}] \in \dZ(\cX^{\EG}_{\red}),
\end{equation}
where $[\cX^{\lambda+\rho, \tau}_{\F}]$ is the cycle class of $\cX_{\F}^{\lambda+\rho, \tau}$. 
\end{conj}

To connect Conjecture \ref{conj: intro BM} with more typical formulations of the (refined geometric) Breuil-M\'{e}zard Conjecture, we note that the group $K_0 (\Rep_k(G(\F_p)))$ is free abelian on the set $\{\sigma\}$ of Serre weights of $G$. Therefore, giving a map $K_0 (\Rep_k(G(\F_p))) \rightarrow \dZ(\cX^{\EG}_{\red})$ amounts to specifying cycles $\cZ(\sigma) \in \dZ(\cX^{\EG}_{\red})$ for each Serre weight $\sigma$, such that if 
\[
[\ol{W(\lambda) \otimes \sigma(\tau)}]  = \sum_\sigma n_\sigma(\lambda, \tau) [\sigma] \in K_0 (\Rep_k(G(\F_p))) 
\]
then 
\begin{equation}\label{eq: intro BM equations}
[\cX^{\lambda + \rho, \tau}_{\F}]  = \sum_\sigma n_\sigma(\lambda, \tau) \cZ(\sigma) \in \dZ(\cX^{\EG}_{\red}).
\end{equation}
This is the formulation of the Breuil--M\'ezard Conjecture as it appears in \cite[\S 8.2.2]{EG23}; we note that it goes significantly beyond the original conjectures of Breuil-M\'{e}zard from \cite{BM02}. The hypothetical cycles $\cZ(\sigma)$ in Conjecture \ref{conj: intro BM} are called \emph{Breuil-M\'{e}zard cycles}. Note that there are approximately $p^n$ Breuil-M\'{e}zard cycles while there are infinitely many possibilities each for $\lambda$ and $\tau$; therefore, the conjecture can be thought of as positing the existence of a solution to a massively overdetermined system of equations.

\subsubsection{Generalization to other groups}For general $G$, we cannot yet formulate an analogue of \eqref{eq: intro BM condition}, because we lack the Local Langlands correspondence for $G$, which supplies the assignment $\tau \rightarrow \sigma(\tau)$. Even the formulation of the (inertial) Local Langlands correspondence for general $G$ is arguably not completely clear. 

However, we know \emph{part} of the Local Langlands correspondence for general classes of groups. For example, if $G$ is unramified then we know what the ``tame part'' of the Local Langlands correspondence should look like. In particular, one has access to the map $\tau \mapsto \sigma(\tau)$ at least for tamely ramified $\tau$ \cite[Theorem D]{Lin23c}. In this vein, we are able to formulate a ``partial Breuil--M\'ezard Conjecture'', which should be the tame part of a hypothetical stronger Conjecture, whose very formulation is still conditional.

\begin{conj}[Geometric Breuil-M\'{e}zard Conjecture]\label{conj: intro BM G}
Let $G$ be an unramified group over $\Q_p$ satisfying \cite[Hypothesis 9.1.1]{GHS18}: $G_{\mrm{der}}$ is simply connected, $Z(G)$ is connected, and $G$ admits a ``local twisting element''. Let $d$ be the dimension of the flag variety of $\chG$.\footnote{It is expected that $d = \dim \cX^{\EG}_{\red}$. However, this is currently only known if $\chG$ has no simple factors of type E or F. } Then there is a map 
\[
\cZ \co K_0( \Rep_k(G(\F_p))) \rightarrow \dZ(\cX^{\EG}_{\red})
\]
such that for every $\lambda \in X^*(T)^+$ and every tame inertial parameter $\tau$, we have
\begin{equation}
\cZ [\ol{W(\lambda) \otimes \sigma(\tau)}] = [\cX^{\lambda+\rho, \tau}_{\F}] \in \dZ(\cX^{\EG}_{\red}),
\end{equation}
where $[\cX^{\lambda+\rho, \tau}_{\F}]$ is the cycle class of $\cX_{\F}^{\lambda+\rho, \tau}$. 
\end{conj}

\begin{remark}
Conjecture \ref{conj: intro BM G} already pins down the map $\cZ$ uniquely (if it exists at all), so that if a stronger Conjecture were eventually formulated, it would necessarily have the same map. 
\end{remark}

\begin{remark}
The hypothesis of the Conjecture is not a technical artifact of our arguments, but is rather a reflection of genuine subtleties in the formulation of Langlands duality, related (for example) to the necessity of using the $C$-group rather than the $L$-group in \cite{BG14}. In practice the hypothesis is harmless, as one can always make some central extension so that it is satisfied. We note that the hypothesis is already satisfied in all cases where the Breuil--M\'ezard has seen previous progress (e.g., $G = \GL_n$ or $G = \GSp_4$), and also in many cases which have never seen progress until this paper (e.g., $G$ an exceptional group).  
\end{remark}

\subsubsection{The weight part of Serre's Conjecture for general $G$} Returning to the thread of Serre's Conjectures, Gee-Kisin \cite{GK14} suggested using the Breuil-M\'{e}zard Conjecture to formulate the generalization of the weight part of Serre's Conjecture to general $G$. We will summarize this formulation as it is described in \cite[\S 6]{GHS18}. Suppose that the Breuil-M\'{e}zard cycles $\cZ(\sigma)$ are given. Let $\ol{r} \co \Gal(\ol{\Q}/\Q) \rightarrow \lsup LG(\ol{\F}_p)$ be an irreducible $L$-parameter. Let $\{\sigma\}_{\ol{r}}$ be the set of all Serre weights for which $\ol{r}|_{\Gal(\ol{\Q}_p/\Q_p)}$ lies on the Breuil--M\'ezard cycle $\cZ(\sigma)$. Then one expects that \emph{$\{\sigma\}_{\ol{r}}$ is the set of Serre weights in which the Hecke eigensystem associated to $\ol{r}$ occurs at some level which is good at $p$}. See \cite{GHS18} and \cite[\S 8]{EG23} for a more detailed discussion.

The picture just described is obviously meaningless without a definition of the Breuil-M\'{e}zard cycles $\cZ(\sigma)$, which would seem to require proving Conjecture \ref{conj: intro BM}, but in fact it is meaningful to construct candidate Breuil-M\'{e}zard cycles without proving the full conjecture. The point is that the cycles $\cZ(\sigma)$ are already uniquely determined by a finite number of equations of the form \eqref{eq: intro BM equations}, so one can know that a candidate construction of $\cZ(\sigma)$ is ``correct'' as long as it satisfies a sufficiently large subset of such equations. The candidate cycles can then be fed into the previous paragraph to give an unconditional formulation of the weight part of Serre's Conjecture for $\GL_n$. 

\subsection{Main results}
 The main results of this paper will follow along the lines just described. We will construct candidate Breuil-M\'{e}zard cycles $\cZ(\sigma)$ for ``sufficiently generic'' $\sigma$, and verify that they satisfy conditions \eqref{eq: intro BM condition} whenever $\lambda$ is ``small'' enough and $\tau$ is a ``sufficiently generic'' tame type. To make this more precise, we need to recall a bit of (modular) representation theory. 

\subsubsection{Modular representation theory}\label{sssec: mod rep theory} For $T \subset G$ the standard maximal torus, the irreducible algebraic representations of $G$ over $k$ are in bijection with the dominant weights $X^*(T)^+$, with $\lambda \in X^*(T)^+$ corresponding to the highest weight representation $L(\lambda)$ of $G$. 
\begin{itemize}
\item The \emph{$p$-restricted weights} $X^*_1(T)\subset X^*(T)$ consist of $\lambda$ such that $0 \leq \tw{\lambda, \alpha^\vee} < p$ for all simple roots $\alpha$. The simple representations (i.e., Serre weights) of $G(\F_p)$ over $k$ are in bijection with $X_1^*(T)$ modulo an equivalence relation on central characters that we suppress, with $\lambda \in X_1^*(T)$ corresponding to $F(\lambda) := L(\lambda)|_{G(\F_p)}$. 
\item  For $(w,\mu)\in W\times X^*(T)$ there is an explicitly constructed tame inertial parameter $\tau=\tau(w,\mu)$. This process gives a bijection between tame inertial parameters and equivalence classes of $(w,\mu)$. When $\tau=\tau(w,\mu)$, the corresponding $\sigma(\tau)$ can be taken to be a certain Deligne-Lusztig representation $R(w,\mu)$ of $G(\F_p)$. 
\end{itemize}

The affine hyperplanes defined by the condition $\tw{-, \alpha^\vee} \in p \Z$ divide $X^*(T)$ into alcoves, and  genericity will be measured by the ``distance'' to the walls of these alcoves. More precisely, we say that $\lambda \in X^*(T)$ is \emph{$m$-generic} (in the lowest alcove) if $m < |\tw{\lambda, \alpha^\vee}| < p-m$ for any root $\alpha$. Note that $m$-generic $\lambda$ only exist if $p>2m$. We say that \emph{$\tau$ is $m$-generic} if $\sigma(\tau) = R(w, \mu)$ where $\mu$ is $m$-generic. We say that \emph{$\sigma = F(\lambda)$ is $m$-generic} if $\lambda+\rho$ is $m$-generic. For $\lambda \in X^*(T)$, we define $h_\lambda$ to be the maximum of $\tw{\lambda, \alpha^\vee}$ among all roots $\alpha$.

\subsubsection{Main theorem} We can now state our main theorem, which proves Theorem \ref{conj: intro BM G} in the generic range of parameters. We have in mind that $p$ is large relative to $G$; for example, $p$ must be larger than four times the Coxeter number of $G$ for the results to be non-vacuous. 

\begin{thm}\label{thm: intro main} Let $G$ be an unramified reductive group over $\Q_p$ satisfying \cite[Hypothesis 9.1.1]{GHS18}: $G_{\mrm{der}}$ is simply connected, $Z(G)$ is connected, and $G$ admits a ``local twisting element''.

(1) There exists a collection of cycles $\cZ^{\EG}(\sigma) \in \dZ(\cX^{\EG}_{\red})$ such that for each $\lambda \in X^*(T)$ and each tame inertial parameter $\tau$ which is $2h_{\lambda+\rho}$-generic, \eqref{eq: intro BM equations} is satisfied. 

(2) If Conjecture \ref{conj: intro BM G} is true, then the ``true'' $\cZ(\sigma)$ agrees with the $\cZ^{\EG}(\sigma)$ from (1) if $\sigma$ is $6h_\rho$-generic. 
\end{thm}

\begin{remark}The proof of the Theorem gives a practical algorithm to compute the Breuil--M\'ezard cycles using equivariant cohomology, which has been implemented at least in low rank by the second author in joint work with Zhongyipan Lin. 
\end{remark}

\begin{remark}[Optimality of the constants] The constant $2h_{\lambda+\rho}$ from part (1) is almost optimal with our method. It is conceivable that it can be slightly reduced to around $h_{\lambda+\rho}$, as this is the natural barrier for the governing modular representation theory to behave in a stable (e.g. independent of $p$) way. 

 The constant $6h_\rho$ from part (2) is not optimal; the argument sacrifices optimality for simplicity by using ``trivial bounds'' where possible. We expect that the argument for (2) can probably be optimized, with some more effort, to give a constant of $3h_\rho$ instead of $6 h_\rho$. 
\end{remark}

\subsubsection{Brief remarks on the proof} The proof of Theorem \ref{thm: intro main} does not follow any existing approaches to the Breuil-M\'{e}zard Conjecture, which are either based on $p$-adic Local Langlands or automorphy lifting (hence have limited generality in the group aspect). Instead our argument is more geometric, and hearkens to the analogy between $\Q_p$ and a 2-manifold, which suggests in turn an analogy between the Emerton-Gee stack and the moduli space of local systems on a 2-manifold. The latter object has a natural \emph{symplectic} structure, and under this analogy the potentially crystalline loci correspond to Lagrangian subspaces. Kontsevich's philosophy of \emph{homological mirror symmetry} posits (roughly) that Lagrangian subspaces of a symplectic manifold can be assembled into a Fukaya category which should then admit a mirror description in terms of coherent sheaves on a mirror variety. In particular, it predicts that Lagrangians can be indexed by mirror data of a ``dual'' nature, which bears a loose similarity to the Breuil--M\'ezard Conjecture. 

In reality, the Emerton-Gee stack does not literally have a symplectic structure. However, the above metaphor can be substantiated by using $p$-adic Hodge theory to ``approximate'' the potentially crystalline substacks by explicit algebraic varieties (certain \emph{affine Springer fibers}), which really do comprise a Lagrangian skeleton of a certain symplectic space of Higgs bundles. Then constructions of Bezrukavnikov--Boixeda Alvarez--McBreen--Yun \cite{BBMY2} provide the necessary ``mirror symmetry'' input to prescribe Lagrangians using coherent sheaves on a mirror variety $A$. Finally, to connect this to the Breuil-M\'{e}zard Conjecture, we ``approximate'' the representation theory of $G(\F_p)$ by representations of $\mf{g}_{\F_p}$, which we then transform into coherent sheaves on $A$ using the modular localization theory of Bezrukavnikov-Mirkovic-Rumynin \cite{BMR08}. More details are given in \S \ref{ssec: overview} below.

\subsection{Comparison to other results} As has been mentioned already, one of the main features of our approach is its generality. Prior to Theorem \ref{thm: intro main}, the Breuil--M\'ezard Conjecture has only been seen progress in the cases $G = \GL_n$ and $G = \GSp_4$, because of genuine limitations in methods. Even when restricted to these groups, we obtain substantial improvements over previous results. We will comment briefly on prior work. 

\subsubsection{The case $G = \GL_2$}
For $G = \GL_2/\Q_p$, our understanding is that Conjecture \ref{conj: intro BM} has now been proven in full, thanks mostly to work of Kisin \cite{Kis09} and Pa\v{s}k\={u}nas \cite{Pas15}, which left out cases that were completed by Hu-Tang, Sander, and Tung. We refer to \cite[\S 8.5]{EG23} for references and more detailed descriptions. The proof is based on the $p$-adic Local Langlands correspondence for $\GL_2/\Q_p$, which has resisted generalization to other groups despite much effort. 

For $\GL_2$ over a finite extension $K/\Q_p$ and $\lambda = (0,0)$, the Breuil--M\'ezard Conjecture is proved by Gee-Kisin \cite{GK14} (at the level of deformation rings) and Caraiani-Emerton-Gee-Savitt \cite{CEGS} (in the geometric form). The proof is based on patching arguments, and therefore is limited to groups whose automorphic theory is sufficiently well-understood, which aside from special low rank examples just leaves $\GL_n$. 

We digress to comment on the difference between the nature of the problem for $\GL_2$ versus higher rank groups. In the case of $\GL_2$ (and only in this case) the Serre weights all lift to Weyl modules in characteristic 0. Therefore in this case (only), Conjecture \ref{conj: intro BM} already includes in its formulation the definition of $\cZ(F(\lambda))$: it must be $[\cX^{\lambda, \mrm{triv}}_{\F}]$. Furthermore, $\cZ(\sigma)$ turns out to be relatively simple: it is simply the irreducible component labeled by $\sigma$, for all $\sigma$ which are not Steinberg.

\subsubsection{The case $G = \GL_n$ for $n>2$} By contrast, the higher rank cases of the Breuil--M\'ezard Conjecture have a very different texture. We quote from \cite[\S 8.7]{EG23}: 
\begin{quote}``We expect the situation for $\GL_n$, $n > 2$, to be considerably more complicated than that for $\GL_2$. Experience to date suggests that the weight part of Serre's conjecture in high dimension is consistently more complicated than is anticipated, and so it seems unwise to engage in much speculation.'' 
\end{quote}
For example, it is expected that $\cZ(\sigma)$ can be reducible even for very generic $\sigma$ \cite[Remark 1.5.11]{LLLM22}. In fact, our construction relates $\cZ(\sigma)$ to the characteristic cycles of simple objects in a representation-theoretic category, and then experience in geometric representation theory suggests that their decompositions with respect to the $\cC_{\sigma}$ are rather subtle in higher dimension, no matter how generic $\sigma$ is.

 The only general result towards the Breuil--M\'ezard Conjecture in arbitrary rank, prior to the present paper, was	 the work of \cite{LLLM22}, which applies to $G = \GL_n$ and unramified $K/\Q_p$. Given a finite set $\Lambda = \{\lambda\}$ of dominant weights, they produce candidate Breuil--M\'ezard cycles $\cZ(\sigma)$ which satisfy \eqref{eq: intro BM equations} whenever $\tau$ is very generic. Here, ``very generic'' depends only on $\Lambda$ and $n$, but is of the nature that some parameters avoid some proper (possibly non-linear) subvarieties in an affine space (in particular, this is rather more stringent than familiar notions of genericity in representation theory, which are of the nature that some parameter avoids some union of hyperplanes). Worse, unlike in Theorem \ref{thm: intro main} \emph{it seems hard to quantify this condition effectively}, so that one cannot say how large $p$ needs to be for the statement to be non-vacuous. 

\begin{remark}A back-of-the-envelope estimate indicates that when $G = \GL_n$ \cite{LLLM22} verifies $O(p^n)$ of the equations \eqref{eq: intro BM equations} where the implicit constant is not effective, while Theorem \ref{thm: intro main} verifies $O(p^{2n})$ of the equations \eqref{eq: intro BM equations} where the implicit constant is effective. 
\end{remark}

The approach of \cite{LLLM22} is based on the patching method as in \cite{GK14}. It is therefore limited by global constraints on executing automorphy lifting, which obstruct it from generalizing to most other groups, except $G = \GSp_4$, which is treated by a similar method in \cite{Lee23}.

\subsection{Comparison to earlier versions}\label{ssec:comparison-versions}Earlier versions of this paper treated only the case $G = \Res_{K/\Q_p} \GL_n$ for an unramified extension $K/\Q_p$. This is the setting where the Breuil--M\'ezard Conjecture has predominantly been studied in the past, although Lee \cite{Lee23} has examined the case of $\GSp_4$, as mentioned above. One motivation for setting up the foundations more generally was the discovery of a ``Langlands functoriality principle for the Emerton--Gee stack'' in our sequel work, which flourishes most naturally in the general framework of reductive groups. 

The additions required to accommodate this generalization are mostly in Appendix \ref{app: EG stacks} (by the second author and Zhongyipan Lin), which is entirely new. Appendix \ref{app: B} has also been rewritten with improved results. Smaller adaptations are also diffused throughout the pre-Appendix body of the paper, but there the structure and content remains largely the same.

\subsection{Discussion of proof and overview of paper}\label{ssec: overview}

We will now summarize our approach to Theorem \ref{thm: intro main} and how it is distributed over the various sections of the paper. 

\subsubsection{Inspirations} Our work has several inspirations, especially the work of Breuil-Hellmann-Schraen, whose \cite[Theorem 1.9]{BHS19} might be regarded as a locally analytic analogue of the Breuil--M\'ezard Conjecture. In their case, the locally analytic Breuil-M\'{e}zard cycles are imported via local model diagrams from the characteristic cycles of $D$-modules that are connected to representation theory through Beilinson-Bernstein localization. Separately, Emerton had remarked to the authors that the geometry of $\cX^{\lambda, \tau}$ resembled that of a Lagrangian in a symplectic manifold. This led us to try to construct our Breuil-M\'{e}zard cycles by importing characteristic cycles from somewhere, although our importing process ends up being more difficult than in \cite{BHS19}, and the construction of the characteristic cycles is also much more involved.

\subsubsection{Construction of Breuil--M\'ezard cycles}

We will divide this discussion into two parts: the first is the construction of the Breuil-M\'{e}zard cycles $\cZ^{\EG}(\sigma)$, and the second is the verification that they satisfy the Breuil-M\'{e}zard relations \eqref{eq: intro BM equations}. 

For the first part, we want to construct (top-dimensional) cycles in $\cX^{\lambda, \tau}_{\F} \subset \cX^{\EG}_{\F}$ from representations. To begin, we construct algebraic models for $\cX^{\lambda, \tau}_{\F}$ inside spaces that we call $\rY_\gamma^{\varepsilon = 1}$, where $\gamma$ is cooked up from $\tau$ by a recipe that we omit for now. As the name suggests, these models can be deformed in a parameter $\varepsilon$, and the specialization $\rY_\gamma^{\varepsilon = 0}$ at $\varepsilon=0$ is an \emph{affine Springer fiber}, which has been previously studied in geometric representation theory. There is a ``microlocal support'' map from a $K$-group of \emph{graded Lie algebra representations} $K_0(\Rep_0( \Ug^0, T))$ towards the group of top-dimensional cycles $Z_d(\rY_\gamma^{\varepsilon = 0})$, which we will explain further below. The Breuil-M\'{e}zard cycles are the images of simple representations in $\Rep_0( \Ug^0, T) $ under this sequence of transformations, summarized in the following diagram. 
\begin{equation}\label{eq: intro BM cycle summary}
\begin{tikzcd}[row sep = tiny]
\substack{\text{Representation}\\\text{theory}} \ar[r, rightsquigarrow, "\text{Part 2}"] & \substack{\text{Affine} \\ \text{Springer fiber}} \ar[r, rightsquigarrow, "\text{Part 1}"] & \substack{\text{Models}} \ar[r, rightsquigarrow, "\text{Part 3}"] & \substack{\text{Emerton-Gee}\\ \text{stack}} \\
K_0(\Rep_0( \Ug^0, T)) \ar[r, rightsquigarrow, "{\substack{\text{microlocal}\\\text{support}}}"] &  \dZ(\rY_\gamma^{\varepsilon = 0}) \ar[r, rightsquigarrow, "{\text{deform $\varepsilon$}}"] &  \dZ(\rY_\gamma^{\varepsilon = 1})  \ar[r, rightsquigarrow, "\substack{\text{$p$-adic}\\ \text{Hodge} \\ \text{theory}}"] & \dZ(\cX^{\lambda, \tau}_{\F}) 
\end{tikzcd}
\end{equation}
We will now explain these steps in more detail, starting from the right side. 

\subsubsection{Local models via $p$-adic Hodge theory} Kisin \cite{Kis08} developed an approach to understanding the moduli of (potentially semistable) $\Gal(\ol{\Q}_p/\Q_p)$-representations using the \emph{$p$-adic Hodge theory} of \cite{Kis06}, which interprets Galois representations in terms of more manageable linear algebraic data called Breuil-Kisin modules. The passage from the models to the Emerton-Gee stack in \eqref{eq: intro BM cycle summary} builds on Kisin's ideas, as further developed in \cite{LLLM22}, \cite{companion}. The space $\rY_\gamma^{\varepsilon = 1}$ is the special fiber of an explicit algebraic $\Z_p$-scheme $\cX_\gamma^{\varepsilon  =1 }$ obtained by truncating the $p$-adic expansions of some equations in a moduli space of Breuil-Kisin modules. It is not clear a priori that this truncation provides a ``good'' approximation to $\cX^{\lambda, \tau}$. A key point is that our strategy asks for relatively little about the quality of this approximation. 

For comparison, we note that the strategy of \cite{LLLM22} requires a smooth local model for $\cX^{\lambda, \tau}$ which is moreover unibranch, and this is the source of their stringent (and hard to compute) genericity conditions. Our arguments need much less, which is part of why we are able to relax the genericity hypotheses. All we need is to find a subscheme $\rY_\gamma^{\varepsilon = 1}(\lambda)\subset \rY_\gamma^{\varepsilon = 1}$ which serves as a ``homological model'' in a very loose sense: (part of) its top homology accurately models the top homology of $\cX^{\lambda, \tau}_{\F}$, the critical feature being that the cycle class $[\cX^{\lambda, \tau}_{\F}]$ agrees with the limit cycle of some natural locus in the generic fiber of $\cX_\gamma^{\varepsilon = 1}$. This is an immediate consequence of the stronger ``mod $p$ local model theorem'' Theorem \ref{thm:mod p model}, which gives a complete description of $\cX^{\lambda, \tau}_{\F}$ itself in terms of a subscheme of $\rY_\gamma^{\varepsilon = 1}(\lambda)$ arising as a $p$-flat limit.
Results of this kind were first discovered in \cite{companion}, which is a substantial improvement of the theory in \cite{LLLM22} that is designed to work in the rather elaborate setup of non-generic $\tau$. In Appendix \ref{app: B}, we import these ideas in the much simpler generic setting relevant to this paper; thus \emph{this paper does not rely on \cite{companion}}.

\subsubsection{Deformation to affine Springer fibers} The $\Z_p$-scheme $\cX_\gamma^{\varepsilon  =1 }$ can be deformed to a family $\cX_\gamma^{\varepsilon}$ over $\Z_p[\varepsilon]$. The fiber over $\varepsilon = 0$ is a mixed characteristic degeneration of \emph{affine Springer fibers}, whose base change to $\Q_p$ we denote $\rX_\gamma^{\varepsilon = 0}$ and whose base change to $\F_p$ we denote $\rY_\gamma^{\varepsilon = 0}$. We remark that just as affine Springer fibers are local analogues of Hitchin spaces (moduli of Higgs bundles), the family $\cX_\gamma^{\varepsilon}$ is a local analogue of moduli of ``$\lambda$-connections'' (with $\varepsilon$ playing the role of $\lambda$; we reserve the notation $\lambda$ for other purposes). 

The entirety of Part 1 is devoted to analyzing the family $\cX^\varepsilon_\gamma$ along with the behavior of $\dZ(\cX^\varepsilon_\gamma)$ under degeneration in $\varepsilon$. In \S \ref{sec: deformed ASF} we define these families, construct the \emph{affine Springer action} on them, and tabulate their top-dimensional irreducible components at specific $\varepsilon$. In \S \ref{ssec: comparison of families} we analyze the specialization of irreducible components of $\cX_\gamma^\varepsilon$ from generic $\varepsilon$ to $\varepsilon = 0$ or $\varepsilon = 1$, which is needed to carry out the deformation of cycles from $\varepsilon = 0$ to $\varepsilon =1$ in \eqref{eq: intro BM cycle summary}.

\subsubsection{Microlocal support} The ``microlocal support'' map from \eqref{eq: intro BM cycle summary} is itself the composition of two steps, which are studied in Part 2. 
\begin{equation}\label{eq: intro microlocal support} 
\begin{tikzcd}[row sep = tiny, column sep = huge]
\substack{\text{Representation}\\\text{theory}} \ar[r, "\S \ref{sec: modular representation theory}"
] &  \substack{\text{Coherent} \\ \text{realization}} \ar[r, "\S \ref{sec: K-groups}"]  & \substack{\text{Affine} \\ \text{Springer fiber}}  \\
K_0(\Rep_0( \Ug^0, T)) \ar[r, "{\substack{\text{BMR}\\\text{localization}}}"]  &  K_0(\Coh_{\cB}^{T}(\wt{\cN}))  \ar[r, "\substack{\text{homological} \\ \text{mirror} \\ \text{symmetry}}"]  &  \dZ(\rY_\gamma^{\varepsilon = 0})  
\end{tikzcd}
\end{equation}

The first step, explained in \S \ref{sec: modular representation theory}, transforms graded Lie algebra representations to coherent sheaves on the Springer resolution $\wt{\cN}$ for $G$, with support conditions determined by central characters, via the localization functor of Bezrukavnikov-Mirkovic-Rumynin \cite{BMR08, BM13}. It should be thought of as the characteristic $p$ version of Beilinson-Bernstein localization, enhanced by observations related to $p$-curvature that allow one to describe $D$-modules on $X$ by coherent sheaves on the Frobenius twist of $T^*X$. 

The second step, explained in \S \ref{sec: K-groups}, transforms coherent sheaves into top-dimensional cycles on the affine Springer fiber $\rY_\gamma^{\varepsilon =0}$, using work of Bezrukavnikov--Boixeda Alvarez--McBreen--Yun \cite{BBMY2}. The key point here is that $\rY_\gamma^{\varepsilon = 0}$ has a natural realization as a Lagrangian inside a symplectic space $\cM_\psi$ of $\chG$-Higgs bundles (in fact $\rY_\gamma^{\varepsilon = 0}$ is homeomorphic to a Hitchin fiber in completely integrable Hitchin system for $\cM_\psi$), and the passage from coherent sheaves on $\wt{\cN}$ to Lagrangians in $\cM_\psi$ should be thought of as some incarnation of homological mirror symmetry. It is implemented by an instance of the geometric Langlands correspondence, relating coherent sheaves to constructible sheaves on a certain moduli space of $\chG$-bundles, which after forming singular support gives the ``microlocal support'' map alluded to above. We highlight that it is ultimately mirror symmetry which provides the passage from $G$ to its dual group $\chG$ in our story.

\subsubsection{Breuil--M\'ezard cycles} We consider the Lie algebra $\mf{g} = (\Lie G)_k$  in characteristic $p$. The category $\Rep_0(\Ug^0, T)$ is the category of finitely generated $X^*(T)$-graded $\mf{g}$-representations with trivial Harish-Chandra central character and nilpotent $p$-central character. It has well-known similarities to $\Rep_k(G(\F_p))$, under which the Serre weights $\sigma$ of $G$ correspond to certain irreducible representations $L(\sigma)$ in $\Rep_0(\Ug^0, T)$. We define $\cZ^{\EG}(\sigma)$ to be the image of $[L(\sigma)]$ under \eqref{eq: intro BM cycle summary}. The details are somewhat elaborate, and are explained in \S \ref{sec: existence of BM cycles}. 

\subsubsection{Verification of the Breuil--M\'ezard relations} We next turn to discuss the proof that the Breuil--M\'ezard cycles $\cZ^{\EG}(\sigma)$ satisfy the equations \eqref{eq: intro BM equations}, which is the content of \S \ref{sec: existence of BM cycles} (building on \S \ref{sec: Degeneration of affine Springer fibers} and Appendix \ref{app: A}). The equations concern the relation between $\Rep_k(G(\F_p))$ and the Emerton-Gee stack for $\chG$. By backtracking through \eqref{eq: intro BM cycle summary}, we can formulate equivalent equations relating $\Rep_0(\Ug^0, T)$ to the geometry of $\rY_\gamma^{\varepsilon = 0}$, which we ultimately prove by an analysis of the microlocal support map. 

We first focus on the case $\lambda = 0$, corresponding to minimal regular Hodge-Tate weights. Recall that $\rX_\gamma^{\varepsilon = 0}$ denotes the generic fiber of the scheme $\cX_\gamma^{\varepsilon=0} \rightarrow \Spec \Z_p$. It has an irreducible component $\rX_\gamma^{\varepsilon=0}(\rho)$, whose fundamental class we specialize to characteristic $p$, obtaining a cycle $\fsp_{p \rightarrow 0} [\rX_\gamma^{\varepsilon = 0}(\rho)]$ on the affine Springer fiber $\rY_\gamma^{\varepsilon = 0}$. This turns out to be the output of transferring ${[\cX_{\F}^{\rho, \tau}]}$ to the model $\rY_\gamma^{\varepsilon = 1}$ and then deforming from $\varepsilon = 1$ to $\varepsilon = 0$. 

On the other hand, \emph{Jantzen's generic decomposition formula} shows that the ``similarity'' between $\Rep_k(G(\F_p))$ and $\Rep_0(\Ug^0, T)$, which we used to make the correspondence of irreducible representations $\sigma \leftrightarrow L(\sigma)$, takes the inertial type $\sigma(\tau)$ to the \emph{baby Verma module} $\wh{Z}_1(p \rho)$. 
\[
\begin{tikzcd}
\text{Emerton-Gee stack} & {[\cX_{\F}^{\rho, \tau}]} \ar[r] \ar[d, dashed, "?"'] & {\fsp_{p \rightarrow 0} [\rX_\gamma^{\varepsilon = 0} (\rho)]} \ar[l] \ar[d, "?"]  & \text{Affine Springer fiber} \\
\Rep_k(G(\F_p)) &  {[\sigma(\tau)]} \ar[r]  \ar[u, dashed] & {[\wh{Z}_1(p \rho) ]} \ar[u, dashed] \ar[l]  & \Rep_0(\Ug^0, T) 
\end{tikzcd} 
\]
From these considerations we reduce \eqref{eq: intro BM condition} in the special case $\lambda = 0$ to the statement that 
\begin{equation}\label{eq: BM on ASF}
\textit{the microlocal support of $ {[\wh{Z}_1(p \rho) ]} $ is ${\fsp_{p \rightarrow 0} [\rX_\gamma^{\varepsilon = 0} (\rho)]}$.}
\end{equation}
We then prove this jointly with Roman Bezrukavnikov and Pablo Boixeda Alvarez in Appendix \ref{app: A}; it is highly non-obvious from the definition of the microlocal support map. Indeed, a major difficulty in the Breuil--M\'ezard Conjecture is that the geometry of a cycle obtained by degeneration (i.e., flat limit) is difficult to understand; the special fiber $\cX^{\rho, \tau}_{\F}$ is itself defined indirectly as a flat limit of an irreducible (and reduced) space $\cX^{\rho, \tau}_{E}$. From this perspective \eqref{eq: BM on ASF} appears at first glance to be as difficult as \eqref{eq: intro BM condition}. We will explain, however, that there is now key new traction provided by the \emph{presence of many symmetries} on $\dZ(\rY_\gamma^{\varepsilon = 0})$.

\subsubsection{Equivariant localization} A key structure of $\rY^{\varepsilon=1}_\gamma$ and $\rY^{\varepsilon = 0}_\gamma$, which is not present on the Emerton-Gee stack, is the existence of an action by a maximal torus $\chT \subset \chG$ that makes them equivariantly formal\footnote{For $\rY^{\varepsilon=1}_\gamma$ we need to restrict to certain bounded subschemes to guarantee equivariant formality.}. Then \emph{equivariant localization} allows us to compute their homology in terms of $\chT$-fixed points. This provides an alternative ``basis'' of homology in which it is easy to compute the effect of degeneration, because the degeneration of $\chT$-fixed points has a simple geometry. Thus, although we cannot compute $\fsp_{p \rightarrow 0} [\rX_\gamma^{\varepsilon = 0} (\rho)]$ in the geometrically natural basis of irreducible components, we are able to compute it in terms of equivariant localization.

\subsubsection{Two extended affine Weyl actions} It then remains to compute the left side of \eqref{eq: BM on ASF} in terms of equivariant localization. Because of the indirect way in which microlocal support is defined, we do not know how to do this directly. Instead, we leverage more symmetries on (the homology of) $\rY^{\varepsilon = 0}_\gamma$, which will enable us to \emph{recognize the microlocal support of $[\wh{Z}_1(p\rho)]$ by its symmetries.} 

More precisely, $\dZ(\rY_\gamma^{\varepsilon = 0})$ has two commuting actions of the extended affine Weyl group $\dWext$ for $\chG$: the \emph{centralizer-monodromy action} which we denote $(\dWext, \cdot)$ comes from the symmetries of the space $\rY_\gamma^{\varepsilon = 0}$, and the \emph{affine Springer action} $(\dWext, \bu)$ which is subtler and does not come from an action on the space. The microlocal support map is equivariant for these actions, where on $\Rep_0(\Ug^0, T) $ the action of $(\dWext, \cdot)$ comes from change of grading, and the action of $(\dWext, \bu)$ comes from wall-crossing functors. It turns out that the baby Verma has a very special ``eigenproperty'' with respect to these two actions, which allows us to characterize its microlocal support without first computing it (amusingly, this \emph{a posteriori} gives simple formulas for the equivariant localization of microlocal supports of the baby Vermas).

\subsubsection{Variation in Hodge-Tate weights} We have just sketched the proof of \eqref{eq: intro BM equations} for the special case $\lambda=0$. The extension to more general $\lambda$ involves an interesting new geometric observation. Using equivariant localization, we are able to express the cycle class $[\cX^{\lambda, \tau}_{\F}]$ as a linear combination of cycle classes of the form $[\cX^{\rho, \tau'}_{\F}]$ for certain $\tau'$, after transporting them to the model via the homological model theorem. Intriguingly, after equivariant localization the relation is just a geometric incarnation of \emph{Weyl's character formula}. This effectively reduces to the case $\lambda = 0$, which we already handled. 

The linear combination from the preceding paragraph involves the action of $(\dWext, \cdot)$ on the homology of the local model, which is ``invisible to the naked eye'' in the sense that it does not come from an action on the underlying space -- we can only ``see'' it through equivariant localization. Such a reduction does not seem to have been anticipated in previous approaches to the Breuil--M\'ezard Conjecture, perhaps because of the subtlety of the symmetries required to express it. We note a similarity to the main result of \cite{Bart23}, which proves an upper bound in an analogous (conjectural) equality where $\tau$ is trivial. In fact, one of the main results in that work, \cite[Theorem 1.1]{Bart23}, concerning cycle relations between special fibers of certain linear algebra moduli spaces, has a more conceptual interpretation in terms of equivariant localization. When combined with the appropriate homological model theorem, the inequality of cycles in \cite{Bart23} can be promoted to actual equalities, verifying the Breuil--M\'ezard equations in this setting. We will return to this in future work.

\subsection{Acknowledgments} The realization that Breuil--M\'{e}zard cycles must be related to microlocalization occurred to one of us (BLH) after listening to a talk by Roman Bezrukavnikov on the work \cite{BBMY} in the conference ``Modular Representation Theory'' (Clay Math Institute, 09/30-04/10/2019), and we would like to thank the organizers and the CMI for making that event happen.
We thank Roman Bezrukavnikov, Pablo Boixeda Alvarez, Gurbir Dhillon, Mark Goresky, Oscar Kivinen, Xin Jin, Simon Riche, Peng Shan, Geordie Williamson, Zhiwei Yun, and Xinwen Zhu for helpful conversations. We thank Christophe Breuil, Matt Emerton, Toby Gee, Mark Kisin, Zhongyipan Lin, and Jeff Manning for comments on a draft. TF acknowledges the hospitality of the Institute for Advanced Study in 2020-2021, where Bezrukavnikov and Boixeda Alvarez gave many patient explanations of their work. Parts of this paper were completed while the authors visited the Hausdorff Institute for Mathematics, funded by the Deutsche Forschungsgemeinschaft (DFG, German Research Foundation) under Germany's Excellence Strategy – EXC-2047/1 -- 390685813.

TF was supported by NSF Postdoctoral Fellowship DMS-1902927, NSF Grant DMS-2302520, and a Viterbi Postdoctoral Fellowship at the Simons-Laufer Mathematical Sciences Institute. B.LH. acknowledges support from the National Science Foundation under grant Nos.~DMS-1128155, DMS-1802037, DMS-2302619 and the Alfred P. Sloan Foundation.

\section{Notation and generalities}\label{ssec: notation}

\subsection{Notation on $p$-adic fields} Throughout this paper, $E$ denotes a finite extension of $\Q_p$, $\cO$ its ring of integers, and $\F$ its residue field. We fix a choice of algebraic closure $\ol{\Q}_p$ of $\Q_p$. When $K/\Q_p$ is an algebraic extension, we will denote the Galois group $\Gal(\ol{\Q}_p/K)$ by $G_{K}$, the Weil group $\Weil(\ol{\Q}_p/K)$ by $W_K$, and the inertia subgroup by $I_K$. 

If $X$ is a space (e.g., a scheme, stack, functor,...) over $\cO$, then we denote by $X_E$, $X_{\cO/p^k}$ and $X_\F$ its base change to $E$, $\cO/p^k$ and $\F$ respectively.

\subsection{Reductive groups}\label{sssec:reductive-groups}

We denote by $G$ an unramified reductive group over $\Q_p$ satisfying \cite[Hypothesis 9.1.1]{GHS18}: $G_{\mrm{der}}$ is simply connected, $Z(G)$ is connected, and $G$ admits a ``local twisting element'' (cf. Definition \ref{defn:twisting-elements} below). Since $G$ is unramified, it has a canonical model over $\Z_p$, which we will use tacitly throughout. 

Lie algebras of groups will be denoted with the corresponding lower-case fraktur letters (and will often be considered over $\F$, as indicated in the text). 

We denote by $\chG$ the Langlands dual group of $G$, regarded as a split reductive group over the ring of integers $\cO$. It is equipped with a canonical split maximal torus $\chT \subset \chG$, the dual torus to the ``abstract Cartan'' $A$ of $G$. Starting in \S \ref{sec: modular representation theory}, we also choose a split maximal torus $T \subset G$.

Furthermore, $(\chG, \chT)$ is equipped with a canonical pinning once we choose a direction of ``positivity''. Since different normalizations are found in the literature on this point, we will spell this out in excruciating detail. The reader may wish to ignore this on a first pass.

\subsubsection{Conventions on positivity}\label{sssec: positivity}
Given any Borel $B  = T \cdot N <G$ with unipotent radical $N$ and Cartan subgroup $T$, we obtain an isomorphism $T \xrightarrow{\sim} B/N \xrightarrow{\sim} A$. \emph{Our convention is that the roots of $A$ on $\mf{b}$ are negative}; this is perhaps a less standard convention but it is consistent with the references \cite{Jan03, BMR06, BMR08, BM13} that we will invoke. 

\begin{remark}
Let $\cB$ be the flag variety of $G$. Recall that $G$-equivariant line bundles on $\cB$ are identified with characters of $A$, according to the following construction. If $\cL$ is a $G$-equivariant line bundle on $\cB$, then any Borel subgroup $B < G$ corresponds to a point $[B] \in \cB$, and acts on the fiber $\cL|_{[B]}$ by a character, which is inflated from a character of $A$. For $\lambda \in X^*(A)$, we denote by $\cO(\lambda)$ the corresponding $G$-equivariant line bundle on $\cB$. Our convention on positivity is determined by the property that \emph{dominant weights of $A$ correspond to semi-ample line bundles on $\cB$} under this construction. 
\end{remark}

This choice equips $\chG$ with a canonical Borel subgroup $\chB$ containing $\chT$, such that the roots of $\chft$ on $\chfb$ are \emph{negative}. 
We denote by $\chPhi$ the roots of $\chT$ on $\chfg$, and by $\chPhi^+ \subset \chPhi$ the subset of positive roots, i.e., the roots of $\chT$ on $\chfg/\chfb$. We write $\chDelta$ for the simple positive roots. This induces a notion of standard parabolic, and for a coweight $\lambda \in X_*(\chT)$ we write $\chP_{\lambda}$ for the corresponding standard parabolic subgroup of $\chG$. 

\begin{defn}[Twisting elements]\label{defn:twisting-elements}If $\chG$ is semi-simple, then we denote by $\rho$ the half sum of its positive coroots. (Later we will want to also view $\rho$ as the half sum of positive roots for $G$, which explains this notation.) If $\chG$ is reductive, then we denote by $\rho$ any twisting element, i.e. a Galois invariant lift of $\rho$ from $\chG_{\ad}$ to $\chG$; this is ambiguous up to center, which will not affect the validity of our statements. 
\end{defn}

\begin{example} If $G = \GL_n$, then $\chB$ is the Borel subgroup of \emph{lower} triangular matrices in $\chG \cong \GL_n$. This is the \emph{opposite} Borel to the one chosen in \cite{LLLM22}; however, because we work with left cosets instead of right cosets in the affine flag variety of $\chG$, the upshot is that the combinatorial formulas here (e.g., for affine Weyl groups) will be \emph{consistent} with the formulas in \cite{LLLM22}.
\end{example}

\subsection{Root systems}
Starting in \S \ref{sec: modular representation theory}, we will have chosen a split maximal torus $T \subset G$ and an isomorphism $X^*(T) \cong X_*(\chT)$, where $X^*$ (resp. $X_*$) denote the \emph{geometric} character (resp. cocharacter) groups (i.e., formed after base change to a separable closure of a field.)

Then our choice of $\chB$ induces a choice of Borel $B \supset T$. We let $\Phi $ be the roots of $T$ on $\mf{g}$ and take the \emph{positive roots} $\Phi^+ \subset \Phi$ to be the roots of $T$ on $\mf{g}/\mf{b}$, so the roots $\Phi^- = - \Phi^+$ of $T$ on $\mf{b}$ are \emph{negative}. We let $\Delta \subset \Phi^+$ be the simple positive roots.

\subsubsection{Weyl groups}\label{sssec: weyl notation} Let $W$ be the finite Weyl group of $(\chG,\chT)$. We let $w_0 \in W$ be the longest element. 

Let $\chQ^\vee \subset X_*(\chT)$ be coroot lattice inside the coweight lattice of $\chG$. We let $\wt{W} \cong  X_*(\chT) \rtimes W $ be the \emph{extended affine Weyl group} (for $\chG$) 
and $W_{\aff}  \cong  \chQ^\vee \rtimes W$ be the \emph{affine Weyl group}. For $\lambda \in X_*(\chT)$, we write $t^\lambda$ for the corresponding element in $\wt{W}$.

For $(\alpha, k) \in \chPhi \times \Z$, the corresponding \emph{root hyperplane} is 
\[
H_{\alpha, k} := \{ \lambda \in X_*(\chT)_{\R} \co \tw{\lambda, \alpha}  = k \}.
\]
We abbreviate $H_\alpha := H_{\alpha, 0}$. The hyperplanes $H_{\alpha, k}$ divide $X_*(\chT)_{\R}$ into alcoves, which are acted upon simply transitively by $W_{\aff}$. Let $A_0$ be the \emph{dominant} base alcove anchored at 0. This choice determines the Bruhat order on $\wt{W}$, and a set of simple reflections for the Coxeter group $W_{\aff}$. Given $\wt{w}\in \wt{W}$, let $\wt{W}_{\leq \wt{w}}$ be the set of elements $\leq \wt{w}$ in the Bruhat order. For discussions related to modular representation theory (in characteristic $p$), we will also need $C_0=-\rho+pA_0$, the dominant $\rho$-shifted base $p$-alcove.

We denote by $\wt{W}^+$ the set of \emph{dominant} elements, i.e., those $\wt{w}$ such that $\wt{w}(A_0)$ is dominant. For $\wt{w}\in \wt{W}$, we denote by $\wt{w}_{\dom}$ the unique element in $W\wt{w}\cap \wt{W}^+$. It is the minimal length representative of $W\wt{w}$.

At times we will have chosen a split maximal torus $T \subset G$ and an identification of $T$ with $A$, i.e., an identification of $\chT$ with the Langlands dual torus of $T$. This entails an identification $X_*(\chT) \cong X^*(T) $, and we will sometimes use this to interpret subsets of one side as subsets of the other. For example, this identifies the coroots $\chPhi^\vee \subset X_*(\chT)$ with the roots $\Phi \subset X^*(T)$, hence the coroot lattice $\chQ^\vee \subset X_*(\chT)$ with the root lattice $Q \subset X^*(T)$. It also gives the alternate presentation $\wt{W} \cong  X^*(T) \rtimes W$ of the extended affine Weyl group.

Let $\Omega$ be the stabilizer of the $A_0$ for the action of $\wt{W}$. We have 
\[
\Omega \cong X_*(\chT) / \chQ^\vee\cong \wt{W}/W_{\aff},
\]
so that $\Omega$ also identifies with the group of central characters $X^*(Z)$ of $G$.

For $\lambda \in X_*(\chT)$, we define 
\[
\Adm(\lambda) :=  \{ \wt{w} \in \wt{W} \co \wt{w} \leq t^{w(\lambda)} \text{ for some } w \in W \}.
\]
 
Following the notation of \cite[\S 2.1.1]{LLLM22}, for $\alpha \in \chPhi$ we define the \emph{$m$th $\alpha$-strip}
\[
H_{\alpha}^{(m, m+1)} := \{ x \in X_*(\chT)_{\R} \co m < \tw{x, \alpha } < m+1\}.
\]
We say that an alcove $A \subset X^*(T)_{\R}$ is \emph{regular} if it does not lie in any strips  (for any $\alpha, m$) containing the base alcove. We say $\wt{w} \in \wt{W}$ is \emph{regular} if $\wt{w}(A_0)$ is regular. We define $\wt{W}^{\reg}$ to be the subset of regular elements in $\wt{W}$. We define $\Adm^{\reg}(\lambda) := \Adm(\lambda) \cap \wt{W}^{\reg}$.

The \emph{fundamental box} is the intersection of all $\alpha$-strips passing through $A_0$ where $\alpha$ is a simple root. Modulo the central directions, it is a fundamental domain for the action of $\chQ^\vee$, and we denote by $\wt{W}_1\subset \wt{W}$ the subset sending $A_0$ to the fundamental box. Elements of $\wt{W}_1$ can be enumerated as follows: for each $w\in W$, there is a unique (up to central translations) element $\wt{w}=t^{\rho_w}w$ such that $\wt{w}(A_0)$ is the unique translate of $w(A_0)$ contained in the fundamental box. One can choose $\rho_w=\sum_{w^{-1}\alpha<0} \omega_\alpha$, the sum of fundamental coweights of $\chG$ for simple roots $\alpha$ such that $w^{-1}(\alpha)<0$. 

\subsubsection{Actions} For $\wt{w} \in \wt{W}$ and $\lambda \in X_*(\chT)$, we write $\wt{w} \lambda$ or $\wt{w} \cdot \lambda$ for the natural action of $\wt{w}$ on $\lambda$. We write 
\[
\wt{w} \bu \lambda := \wt{w} (\lambda + \rho) - \rho
\]
for the dot action. 

Letting $\wt{w} = wt^{\nu}$, we write $\wt{w} \cdot_p \lambda := w t^{p \nu} \lambda$ for the natural action dilated by $p$ on the lattice part, and $\wt{w} \bup \lambda := w t^{p \nu} \bu \lambda$ for the dot action dilated by $p$ on the lattice part.

\subsubsection{Heights of weights}\label{sssec: height}
For $\lambda \in X_*(\chT) \cong X^*(T)$, we define its \emph{height}
\[
h_{\lambda} : =\max_{\alpha \in \chPhi} |\langle \lambda,\alpha \rangle|.
\] 
This can be generalized for $\wt{w} \in \wt{W}$: we define $h_{\wt{w}}$ to be the maximum over $\alpha \in \chPhi$ of the number of $\alpha$ root hyperplanes $H_{\alpha, k}$ separating $A_0$ and $\wt{w}A_0$.  

If $G$ is simple and simply connected, then the \emph{Coxeter number} of $G$ is $h_\rho+1$. In general, let $h$ be the maximum of the Coxeter numbers of simple factors of $G_{\ad}$. \textbf{We assume throughout that $p>h$.}

\subsubsection{Genericity} Let $m \in \N$. 
\begin{itemize}
\item We say that $\lambda \in X_*(\chT)$ is \emph{$m$-generic} if $|\tw{\lambda, \alpha} + p k| > m$ for all $\alpha \in \chPhi$ and $k \in \Z$. (This terminology is consistent with \cite[Definition 2.1.10]{LLLM22}.) 

Given a commutative ring $R$, we say that an element $\gamma \in X_*(\chT) \otimes_{\Z} R$ is \emph{$m$-generic} if $\tw{\gamma, \alpha^\vee}  + i + pk \in R^\times$ for all $k \in \Z$ and all $i  \in \{ 0, \pm 1, \ldots, \pm m\}$. Taking $R = \F_p$, note that $\lambda$ is $m$-generic (in the sense of the preceding paragraph) if and only if $\lambda \otimes 1 \in X_*(\chT) \otimes_{\Z} \F_p$ is $m$-generic.

\item For $m \geq 0$ and $\wt{w}  = wt^\nu \in \wt{W}$, we say that $\wt{w}$ is \emph{$m$-small} if $h_{\nu} \leq m$. 
\end{itemize}
\begin{example}
\begin{enumerate}
\item
We are interested in the examples $R_1  = \F[[t]]$ or $R_2 = \cO[\varepsilon]((t+p))$. In these cases $X_*(\chT) \otimes_{\Z} R_1 \cong \chft(\F)[[t]]$ or $X_*(\chT) \otimes_{\Z} R_2 \cong \chft(\cO)[\varepsilon]((t+p))$, and $\gamma$ is an element used to construct a ``deformation of an affine Springer fiber'' or a ``local model for a stack of potentially crystalline local Galois representations'', respectively. 
\item If $\wt{w}\in \Adm(\lambda)$ then $\wt{w}(0)$ lies in the convex hull of $W\lambda$, hence $h_{\wt{w}}\leq h_\lambda$. 
\end{enumerate}
\end{example}

\subsection{Specialization for Chow groups}\label{ssec: specialization for chow}
Throughout this paper, we denote $\Ch(\cX) := \Ch(\cX)_{\Q}$ for the \emph{rational} Chow group of an algebraic scheme or stack $\cX$. We will never consider integral structures on Chow groups. We denote by $Z_d(\cX)$ the group of $\Q$-linear combinations of $d$-equidimensional substacks of $\cX$. 

Let $S$ be a regular scheme, $i \co Z \inj S$ a closed regular embedding of codimension $d$, and $j \co U \inj S$ its open complement. Assume that $Z$ is regular and the normal bundle to $i$ has trivial top Chern class.

Let $f \co \cX \rightarrow S$ be a finitely presented map of schemes. Write $\cX_Z$ for the base change of $\cX$ to $Z$, and $\cX_U$ for the base change of $\cX$ to $U$. Recall the relative Chow group \cite[Chapter 20.1]{Ful98} $\Ch_*(\cX_U/U)$, etc. Under the assumption, \cite[Chapter 20.3]{Ful98} constructs a \emph{specialization map}
\[
\fsp \co \Ch_m(\cX_U/U) \rightarrow \Ch_m(\cX_Z/Z)
\]
as follows. There is a refined Gysin map $i^! \co \Ch_m(\cX_U/U) \rightarrow \Ch_m(\cX_Z/Z)$, and the Chern class assumption implies that $i^!$ factors over the restriction map $\Ch_m(\cX/S) \rightarrow \Ch_m(\cX_U/U) $, which is surjective. This factorization is by definition $\fsp$, as depicted in the diagram below. 
\[
\begin{tikzcd}
\Ch_m(\cX/S) \ar[r, twoheadrightarrow]  \ar[dr, "i^!"'] & \Ch_m(\cX_U/U) \ar[d, "\fsp ", dashed] \\
&  \Ch_m(\cX_Z/Z)
\end{tikzcd} 
\]

\subsubsection{Specialization over a DVR} Suppose $S$ is the spectrum of a discrete valuation ring, with generic point $\eta \in S$ and special point $s \in S$. Then $\fsp$ can be described as follows: for the cycle class $[W_\eta] \in \Ch_m(\cX_\eta)$ of some $W_\eta \subset \cX_\eta$, let $W \subset \cX$ be the Zariski closure of $W_\eta$. Then $W$ is flat over $S$, and 
\[
\fsp([W_\eta]) = [W_s] \in \Ch_m(\cX_s).
\]
This description shows the following effectivity of the specialization map. 

\begin{lemma}\label{lem: specialization preserves effectivity} Suppose $S$ is the spectrum of a DVR, with generic point $\eta \in S$ and special point $s \in S$. If $\alpha \in \Ch_m(\cX_\eta)$ is represented by an effective cycle, then $\fsp(\alpha) \in \Ch_m(\cX_s)$ is also represented by an effective cycle. 
\end{lemma}

\subsubsection{Iterated specialization}\label{sssec: iterated spc} 
 Suppose $S$ is regular scheme, and $f_1, f_2 \in \cO(S)$ are such that:
 \begin{itemize}
 \item For each $m \in \{1,2\}$, the closed subscheme $Z_m = V(f_m)$ is regular, realizing $i_m \co Z_m \inj S$ as a regular embedding of codimension $1$.
 \item The closed subscheme $Z := V(f_1, f_2)$ is regular, realizing $\ol{i}_m \co Z  \inj Z_m$ as a regular embedding of codimension $1$. 
 \end{itemize}
Note that our assumption implies that for each $m=1,2$ the normal bundle of $i_m$ is trivial and the normal bundle of $Z \inj Z_m$ is also trivial.

\begin{example}Let $k$ be a field. A prototypical example is $S = \Spec k[t_1, t_2]$ or $\Spec k[[t_1, t_2]]$ with $f_i = t_i$. 

We will be interested in a mixed-characteristic variant: $S = \Spec \cO[\varepsilon]$ or $S = \Spec \cO[\varepsilon]_{(\varepsilon)}$ with $f_1 = p$ and $f_2 = \varepsilon$. 
\end{example}

Now let $f \co \cX \rightarrow S$ be a finitely presented map. Under the assumptions above, we have specialization maps 
\begin{equation}\label{eq: iterated specialization diagram}
\begin{tikzcd}
\Ch(\cX[\frac{1}{f_1f_2}] / S[\frac{1}{f_1f_2}]) \ar[r, "\fsp_{f_1 \rightarrow 0}"] \ar[d, "\fsp_{f_2 \rightarrow 0}"]  &  \Ch(\cX_{Z_1}[\frac{1}{f_2}]/ Z_1[\frac{1}{f_2}])  \ar[d, "\fsp_{f_2 \rightarrow 0}"] \\
 \Ch(\cX_{Z_2}[\frac{1}{f_1}] / Z_2[\frac{1}{f_1}])  \ar[r, "\fsp_{f_1\rightarrow 0}"]  & \Ch(\cX_Z/Z)  \\
\end{tikzcd}
\end{equation}

\begin{lemma}\label{lem: iterated specialization}
Diagram \eqref{eq: iterated specialization diagram} commutes. 
\end{lemma}

\begin{proof} Consider the diagram 
\begin{equation}\label{eq: unfold iterated spc}
\begin{tikzcd}
\Ch(\cX/S) \ar[r, twoheadrightarrow] \ar[dr, "i_1^!"']  & \Ch(\cX[\frac{1}{f_2}]/S[\frac{1}{f_2}]) \ar[r, twoheadrightarrow] \ar[dr, "i_1^!"']  & 
\Ch(\cX[\frac{1}{f_1f_2}]/S[\frac{1}{f_1f_2}])   \ar[d, "\fsp_{f_1 \rightarrow 0}"] \\
&  \Ch(\cX_{Z_1}/Z_1) \ar[r, twoheadrightarrow]  \ar[r]  \ar[dr, "\ol{i}_2^!"'] & \Ch(\cX_{Z_1}[\frac{1}{f_2}]/Z_1[\frac{1}{f_2}]) \ar[d, "\fsp_{f_2 \rightarrow 0}"] \\
 & & \Ch(\cX_Z/Z)
\end{tikzcd}
\end{equation}
By definition of the specialization maps, the two triangles in the right column commute. By base change compatibility for the refined Gysin pullback $i_1^!$, the top left parallelogram commutes. Hence the whole diagram commutes. The right-then-down path in \eqref{eq: iterated specialization diagram} is the right column of \eqref{eq: unfold iterated spc}. The commutativity of \eqref{eq: unfold iterated spc} says that it can be computed by picking any lift in $\Ch_*(\cX/S)$ and applying $\ol{i}_2^! \circ  i_1^!$, which equals $i^!$ by the compositional property of the refined Gysin pullback. 

Now note that $i^! \co \Ch(\cX/S) \rightarrow \Ch(\cX_Z/Z)$ does not depend on the factorization of the embedding $i \co Z \inj S$. A symmetric argument shows that it computes the down-then-left path in \eqref{eq: iterated specialization diagram} in the same sense, hence verifying the commutativity of \eqref{eq: iterated specialization diagram}.
\end{proof}

 \subsection{Borel--Moore homology}\label{ssec: BM homology}
 Our convention is that all pullback/pushforward operations on $\ell$-adic sheaves are derived, so we write $\pi_* := R\pi_*$, etc. 

\subsubsection{Borel--Moore homology}\label{sssec: BM homology}Let $k$ be a field and $\pi \co X \rightarrow \Spec k$ be a finite type scheme. Suppose there is a square-root of the cyclotomic character $\Gal_k \rightarrow \Z_\ell^{\times}$; choose one to define the half Tate-twist $(1/2)$. (This is the case if $k$ is an algebraic extension of $\Q_p$ or $\F_p$, possibly after enlarging $\Z_\ell$, which covers all the cases we will consider.) Then we define the \emph{$m$th ($\ell$-adic) Borel--Moore homology group to be}
\[
H^{\BM}_m (X) := H^{-m}(X; \D_X(-m/2) )
\]
where $\D_X \cong \pi^!\Qll{\Spec k}$ is the dualizing sheaf on $X$. (We will only ever consider the case where $m$ is even, so we will never invoke the choice of square root of cyclotomic character.) 

We also define the \emph{geometric Borel--Moore homology groups} to be 
\[
\mBMg_m (X) := H^{\BM}_m(X_{\ol{k}}).
\]
\emph{We emphasize the difference in font used to distinguish the geometric and absolute Borel--Moore homology.} Note that the Tate twists do not affect the underlying group of $\mBMg_m (X) $. Pullback defines a map $H^{\BM}_m(X) \rightarrow \mBMg_m(X)$. We will mostly be interested in geometric Borel--Moore homology; the absolute Borel--Moore homology groups are only used as intermediate steps to make constructions with the geometric groups. 

Given a proper map $f \co X \rightarrow Y$ over $k$, the adjunction $f_* \D_X  = f_! \D_X \rightarrow \D_Y$ induces a map
\[
f_* \co H_*^{\BM}(X) \rightarrow H_*^{\BM}(Y).
\]

For an ind-scheme $X = \colim_i X_i$, we define 
\[
H^{\BM}_m(X) := \colim_i H^{\BM}_*(X_i)
\]
with the transition maps induced by the closed embeddings $X_i \inj X_j$ as above. Similarly we define the geometric Borel--Moore homology groups
\[
\mBM_m(X) := \colim_i \mBM_*(X_i).
\]

\subsubsection{Relative Borel--Moore homology} Given a map $\pi \co X \rightarrow S$, the \emph{relative Borel--Moore homology} of $X/S$ is 
\[
H^{\BM}_m(X/S) := H^{-m}(X; \pi^! \Qll{S}(-m/2)).
\]
When $S = \Spec k$, this recovers the previously defined (absolute) Borel--Moore homology groups.

\subsubsection{The cycle class map}\label{sssec: cycle class map}
For any scheme $X \rightarrow \Spec k$, there is a cycle class map 
\[
\Ch_m(X) \rightarrow H^{\BM}_{2m}(X).
\]
For a cycle $Z$ in a space $X$, we write $[Z] \in \Ch_{\dim Z}(X)$ for its Chow class or $[Z] \in \mBM_{2\dim Z}(X)$ for the $\ell$-adic realization of $[Z]$ in the \emph{geometric} Borel--Moore homology of $X$, depending on context to make it clear which version we refer to. 

\subsection{Specialization for Borel--Moore homology}\label{ssec: specialization for BM homology}

Let the setup be as in \S \ref{ssec: specialization for chow}. Then there is a specialization map (for example, apply \cite[\S 4.5.6]{DJK21} with coefficients being $\Q_\ell$)
\[
\spc \co H^{\BM}_*(\cX_U/U) \rightarrow H^{\BM}_*(\cX_{Z}/Z).
\]

\subsubsection{Functoriality}\label{sssec: functoriality for specialization maps} Specialization maps are functorial with respect to the following types of morphisms. 
\begin{itemize}
\item Let $f \co \cX \rightarrow \cY$ be proper. Then there are pushforward maps 
\[
f_* \co H^{\BM}_*(\cX_U/U)   \rightarrow H^{\BM}_*(\cY_U/U)  \hspace{1cm} \text{and} \hspace{1cm} 
f_* \co H^{\BM}_*(\cX_Z/Z)   \rightarrow H^{\BM}_*(\cY_Z/Z)  
\]
that fit into a commutative diagram 
\[
\begin{tikzcd}
 H^{\BM}_*(\cX_U/U) \ar[r, "\spc"] \ar[d, "f_*"] & H^{\BM}_*(\cX_{Z}/Z) \ar[d, "f_*"] \\
 H^{\BM}_*(\cY_U/U) \ar[r, "\spc"] & H^{\BM}_*(\cY_{Z}/Z) 
\end{tikzcd}
\]
\item Suppose $f \co \cX \rightarrow \cY$ promotes to a quasi-smooth map of derived schemes\footnote{See \cite{FH} for a primer on these notions, with references.} of virtual dimension $d(f)$. Then there are pullback maps 
\[
f^! \co H^{\BM}_*(\cY_U/U)   \rightarrow H^{\BM}_{*+2d(f)}(\cX_U/U)  \hspace{1cm} \text{and} \hspace{1cm} 
f^! \co H^{\BM}_*(\cY_Z/Z)   \rightarrow H^{\BM}_{*+2d(f)}(\cX_Z/Z)  
\]
that fit into a commutative diagram 
\[
\begin{tikzcd}
 H^{\BM}_{*+2d(f)}(\cX_U/U) \ar[r, "\spc"] & H^{\BM}_{*+2d(f)}(\cX_{Z}/Z)   \\
 H^{\BM}_*(\cY_U/U) \ar[r, "\spc"] \ar[u, "f^!"] & H^{\BM}_*(\cY_{Z}/Z) \ar[u, "f^!"]
\end{tikzcd}
\]
\end{itemize}

\subsubsection{Geometric specialization for DVRs}

Suppose $S$ is a discrete valuation ring, with generic point $\eta \in S$ and special point $s \in S$. Let $f \co \cX \rightarrow S$ be of finite presentation. Then there is a specialization map for the \emph{geometric} Borel--Moore homology groups, 
\[
\fsp \co \mBMg_*(\cX_\eta) \rightarrow \mBMg_*(\cX_s).
\]
It is obtained by taking the colimit of the specialization maps over finite extensions $S'/S$ (cf. \cite[Example 20.3.5]{Ful98}).

\subsubsection{Dimensions}

We set $d := \dim \check{\cB}$ to be the dimension of the flag variety of $\chG$. At various points we write $\topCh$ or $\topBM$ for the group of top-degree classes (in Chow groups or geometric Borel--Moore homology groups), which unless noted otherwise means $\topCh = \Ch_d$ and $\topBM = \mBM_{2d}$.

\part{Degeneration of local models}\label{part: I}

In this Part, we define and study families of spaces that interpolate between the \emph{local models} for potentially crystalline Emerton--Gee substacks, and \emph{affine Springer fibers}. These spaces we define will only depend on $\chG$, and we will not refer (in this Part) to their Frobenius endomorphisms. Therefore, everything we discuss in this Part is insensitive to the possible twisting of $G$ with respect to its split form.

\section{Deformations of affine Springer fibers} \label{sec: deformed ASF}

In this section, we define deformations of geometric objects called ``affine Springer fibers''. Let us give a roadmap to these constructions and their significance. We construct a family $\cX^{\varepsilon}_{\gamma}$ over the 2-dimensional base $\Spec \cO[\varepsilon]$, depending on $\chG$ and a parameter $\gamma \in \chfg(\cO)[\varepsilon] ((t+p))$. To give a feel for the family $\cX^{\varepsilon}_{\gamma} \rightarrow \Spec \cO[\varepsilon]$, we describe its fiber over various special loci in the base. 
\begin{itemize}
\item Over the divisor cut out by $\varepsilon=1$, $\cX^{\varepsilon}_\gamma$ is an $\cO$-scheme closely related to the local models of potentially crystalline substacks in the Emerton-Gee stack for $\chG$. More precisely, it is (for appropriate choices of $\gamma$) a ``naive local model'' of \cite[\S 3]{LLLM22}. 
\item Over the divisor cut out by $\varepsilon = 0$, $\cX^{\varepsilon}_\gamma$ is a mixed-characteristic degeneration of \emph{affine Springer fibers}, whose geometry will be seen to be closely connected to representation theory. 
\end{itemize}
The Breuil--M\'ezard Conjecture will be related to the divisor $\varepsilon = 1$, but we have better traction on the divisor $\varepsilon = 0$, thanks to geometric representation theory. The family $\cX^{\varepsilon}_\gamma$ interpolates between these, and will allow us to transfer information from one to the other. 


The construction of $\cX^{\varepsilon}_\gamma$ occupies \S \ref{ssec: affine flag families} and \S \ref{ssec: deformed ASF}. From \S \ref{ssec: translation} -- \S \ref{ssec: affine Springer action} we construct an \emph{affine Springer action} on its relative Borel--Moore homology over the locus $p=0$. Finally in \S \ref{ssec: top-irred-comp generic} and \S \ref{ssec: top-irred-comp special} we analyze the irreducible components in the fibers of $\cX^{\varepsilon}_\gamma$ over points of interest.

\subsection{Families of affine flag varieties}\label{ssec: affine flag families}

We regard $\chG$ as a reductive group over $\cO$. Recalling that we have fixed a Borel subgroup $\chB \subset \chG$, we let $\chcG$ be the Bruhat-Tits group scheme over $\A^1_{\cO} = \Spec \cO[t]$ obtained by dilatation of the Chevalley group scheme $\chG/ \A^1_{\cO}$ along $\chB_{\cO} \subset \chG_{\cO}$ in the fiber at the origin of $\A^1_{\cO}$. 

\begin{remark}[Comparison to conventions of \cite{LLLM22}]
When comparing to \cite{LLLM22}, our $\chB$ will be the \emph{opposite Borel} to the one of \emph{loc. cit.}. This comes from our convention to view the affine Grassmannian as a left coset space rather than a right coset space. For example, when $\chG = \GL_n$ our formulas will be based on the choice of $\chB$ as the \emph{lower-triangular} Borel subgroup. 
\end{remark}

Let $L \chcG$ be the functor on $\cO$-algebras $R$ sending $R \mapsto \chcG(R((t+p)))$ and $L^+ \chcG$ be the functor on $\cO$-algebras $R$ sending $R \mapsto \chcG(R[[t+p]])$. The fppf quotient $L\chcG / L^+ \chcG$ is represented by an ind-scheme $\Gr_{\chcG}$ over $\cO$, with the following properties: 
\begin{itemize}
\item The generic fiber $\Gr_{\chcG} \times_{\Spec {\cO}} \Spec E$ is isomorphic to $\Gr_{\chG, E}$, the affine Grassmannian for $\chG_{E}$.
\item The special fiber $\Gr_{\chcG} \otimes_{\Spec \cO} \Spec \F$ is isomorphic to $\Fl_{\chG, \F}$, the affine flag variety for $\chG_{\F}$. 
\end{itemize}
We therefore think of $\Gr_{\chcG}$ as a mixed-characteristic degeneration from the affine Grassmannian to the affine flag variety.

\subsection{Deformed affine Springer fibers} \label{ssec: deformed ASF}

Let $\gamma \in \chft(\cO)[\varepsilon] ((t+p)) \subset \chfg(\cO)[\varepsilon]((t+p))$. 

For $\A^1_{\cO} = \Spec \cO[\varepsilon]$, let $\cX_{\gamma}^{\varepsilon}$ be the sub-(ind-)scheme of $\Gr_{\chcG} \times_{\cO} \A^1_{\cO}$ defined as 
\[
\cX_{\gamma}^{\varepsilon} = \left\{ g \in \Gr_{\chcG} \times_{\cO} \A^1_{\cO}\co \Ad_{g^{-1}}(\gamma)  -  \varepsilon  t(t+p) \frac{dg^{-1}}{dt}  g \in  \Lie L^+ \chcG \right\},
\]
where the symbol $\frac{dg^{-1}}{dt} g =: d\log (g^{-1})$ is understood as in \cite{FZ10}, and is explained in \cite[\S 1.2.4]{F07}.\footnote{A small calculation is required to see that the defining equation is invariant for the right action of $L^+\chcG$.}
We note that
\[
- \frac{dg^{-1}}{dt} g = - d\log(g^{-1}) =  d\log(g) =   g^{-1} \frac{dg}{dt} 
\]
and occasionally we will use the latter form of the expression, for comparison to formulas in \cite{LLLM22}.

\begin{remark}The definition of $\cX_\gamma^{\varepsilon}$ makes sense when $\gamma$ lies more generally in $\chfg(\cO)[\varepsilon]((t+p))$, but we do not need this generality and it complicates the notion of translation by the ``stabilizer of $\gamma$'', so we restrict our attention to $\gamma$ of the stated form. 
\end{remark}

\begin{remark}
More generally, there are versions of $\cX_{\gamma}^{\varepsilon}$ with defining equation 
\[
\Ad_{g^{-1}} (\gamma)  -t^r(t+p) \frac{dg^{-1}}{dt}  g \in \Lie L^+ \chcG
\]
as soon as $r \geq 1$; the versions with $r>1$ will not be needed in this paper. 
\end{remark}

\subsubsection{Specializations} We introduce some notation for fibers of the family $\cX_{\gamma}^{\varepsilon} \rightarrow \Spec \cO$ over specific loci. 
\begin{itemize}
\item For a fixed $\varepsilon_0 \in \A^1_{\cO}$ we denote by $\cX_{\gamma}^{\varepsilon = \varepsilon_0}$ the fiber of $\cX_{\gamma}^{\varepsilon}$ over $\varepsilon_0$. 

\item We denote by $\rX_{\gamma}^{\varepsilon} := \cX_{\gamma}^{\varepsilon}|_{\A^1_{E}}$ the fiber of $\cX_{\gamma}^{\varepsilon}$ over $\Spec E$. For $\varepsilon_0 \in \A^1_{E}$ we denote by $\rX_\gamma^{\varepsilon = \varepsilon_0}$ the fiber of $\rX_{\gamma}^{\varepsilon} $ over $\varepsilon_0$. 

\item We denote by $\rY_\gamma^\varepsilon := \cX^{\varepsilon}_{\gamma}|_{\A^1_{\F}}$ the fiber of $\cX_{\gamma}^{\varepsilon}$ over $\Spec \F$. For $\varepsilon_0 \in \A^1_{\F}$ we denote by $\rY_\gamma^{\varepsilon = \varepsilon_0}$ the fiber of $\rY_{\gamma}^{\varepsilon} $ over $\varepsilon_0$.
\end{itemize}

\begin{example}[The specialization $\varepsilon= 0$]
Consider $\varepsilon_0=0 \in \A^1_{E}$. Then $\rX_{\gamma}^{\varepsilon = 0}$ is the (spherical) \emph{affine Springer fiber} associated to $\gamma$ (over $E$). In turn, $\rY_\gamma^{\varepsilon = 0}$ is the Iwahori affine Springer fiber associated to $\gamma$ (over $\F$). These notions of affine Springer fibers were originally introduced by Kazhdan-Lusztig in \cite{KL88}. Hence we may regard $\cX_{\gamma}^{\varepsilon}$ as a two-parameter deformation of affine Springer fibers, over the two-dimensional base $\A^1_{\cO} $.
\end{example}

\subsection{Translation action}\label{ssec: translation} For $\gamma \in \chft(\cO)[\varepsilon]((t+p))$, we write $\gamma_i \in  \chft(\cO)((t+p))$ for the coefficient of $\varepsilon^i$ in $\gamma$. The centralizer of $\gamma_0$ contains $L\chT$, and evidently acts on $\rX_{\gamma}^{\varepsilon = 0}  = \rX_{\gamma_0}^{\varepsilon = 0}$ and $\rY_{\gamma}^{\varepsilon = 0}= \rY_{\gamma_0}^{\varepsilon = 0}$ by (left) translation.

At the level of underlying reduced schemes we have $L\chT \cong T \rtimes X_*(\chT)$. The action of $L\chT$ on homology of $\rX^{\varepsilon = 0}_{\gamma}$ or $\rY^{\varepsilon = 0}_{\gamma}$ therefore factors through an action of $X_*(\chT)$, which we also refer to as the \emph{translation action}.

There is no translation action on deformed affine Springer fibers; instead translation elements take one deformed affine Springer fiber to another. Namely, a straightforward calculation shows that for $\gamma \in \chft(\cO)[\varepsilon]((t+p))$ and $h\in L\chT$, left multiplication by $h$ takes $\cX_{\gamma'}^{\varepsilon}$ isomorphically to $\cX_{\gamma}^{\varepsilon}$, where 
\begin{equation}\label{eq: translate DASF}
\gamma' := h^{-1}\gamma h - \varepsilon t(t+p)\frac{dh^{-1}}{dt} h.
\end{equation}

\begin{remark}The reason why we consider $\gamma$ having a dependence on $\varepsilon$ is to have a class of objects preserved by translations by $L\chT$. 
\end{remark}

\subsection{The Grothendieck alteration}\label{ssec: groth alteration}
We recall some facts from Springer theory; a reference is \cite[\S 1.2.2]{Yun17}. For a split reductive group $H$ over a field, the \emph{Grothendieck alteration} is the projection map $\pi \co \wt{\mf{h}} \rightarrow \mf{h}$, where $\mf{h}$ parametrizes pairs of $u \in \mf{h}$ and a Borel subgroup $B_H \subset H$ such that $u \in \Lie B_H$. The map $\pi$ is small, so $\pi_* \Qll{\wt{\mf{h}} }$ is an intermediate extension supported on all of $\mf{h}$. On the other hand, over the strongly regular semisimple locus of $\mf{h}$, $\pi$ is a torsor for the Weyl group $W_H$ of $H$. These two facts together equip $\pi_* \Qll{\wt{\mf{h}} }$ with a canonical action of $W_H$. Noting that $\wt{\mf{h}}$ is smooth, so the dualizing sheaf $\DD_{\wt{\mf{h}}}$ is isomorphic to a shift and twist of $\Qll{\wt{\mf{h}}}$, we equivalently get a $W_H$-action on $\DD_{\wt{\mf{h}}}$.

\begin{defn}\label{defn: cartesian up to nilpotents}
We say that a commutative diagram of schemes (or stacks)
\[
\begin{tikzcd}
A \ar[r] \ar[d] & B \ar[d] \\
C \ar[r] & D
\end{tikzcd}
\]
is \emph{Cartesian up to nilpotents} if the induced map on the underlying reduced subschemes from $A$ to the fibered product $B \times_D C$ is an isomorphism. 
\end{defn}

Suppose we have a diagram 
\[
\begin{tikzcd}
\wt{S} \ar[d, "\pi_S"] \ar[r, "\wt{f}"]  & \wt{\mf{h}} \ar[d, "\pi"] \\
S \ar[r, "f"] & \mf{h}
\end{tikzcd}
\]
which is Cartesian up to nilpotents. Then by proper base change, the $W_H$-action on $\pi_* \DD_{\wt{\mf{h}}}$ induces a $W_H$-action on 
\[
f^! \pi_* \DD_{\wt{\mf{h}}} \cong \pi_{S*} \wt{f}^! \DD_{\wt{\mf{h}}} \cong \pi_{S*} \DD_{\wt{S}}.
\]
In particular, after passing to cohomology we obtain a $W_H$-action on $\mBM_*(\wt{S})$. Actions constructed by this mechanism will generally be referred to as ``Springer actions''. 

Let us record the compatibility of Springer actions in a general situation. Given a commutative diagram 
\[
\begin{tikzcd}
\wt{S} \ar[d, "\pi_S"] \ar[r, "\wt{f}"] & \wt{T} \ar[r, "\wt{g}"] \ar[d, "\pi_T"] & \wt{\mf{h}} \ar[d, "\pi"] \\
S \ar[r, "f"] & T \ar[r, "g"] & \mf{h}
\end{tikzcd}
\]
in which all squares are Cartesian up to nilpotents, the $W_H$-action on $\pi_* \DD_{\wt{\mf{h}}}$ induces also a $W_H$-action on $\pi_* \DD_{\wt{T}}$ and $\pi_* \DD_{\wt{S}}$, hence on $\mBM_*(\wt{T})$ and $\mBM_*(\wt{S})$. 

\begin{lemma}\label{lem: Springer equivariant}
(1) If $f$ is proper, then the map $\mBM_*(\wt{S}) \xrightarrow{\wt f_*} \mBM_*(\wt{T})$ is equivariant for the $W_H$-action. 
 
 (2) If $f$ and $\wt{f}$ can be promoted to quasi-smooth maps of derived schemes with virtual dimension $d(f)$,\footnote{For example, this is the case if $f$ is LCI.} making the diagram commute up to nilpotents on classical truncations, then the map $\mBM_*(\wt{T}) \xrightarrow{f^!} \mBM_{*+2d(f)} (\wt{S})$ is equivariant for the $W_H$-action, where $d(f)$ is the \emph{virtual dimension} of $f$, i.e., the Euler characteristic of the cotangent complex of $f$. 
\end{lemma}

\begin{proof} (1) By base change, the map in question is obtained by taking global sections on $T$ of the map 
\[
f_* f^! g^! \pi_* \DD_{\wt{\mf{h}}}  =  f_! f^! g^! \pi_* \DD_{\wt{\mf{h}}}  \xrightarrow{\mrm{counit}} g^! \pi_* \DD_{\wt{\mf{h}}} ,
\]
and the $W_H$-action is induced by the Springer $W_H$-action on $\pi_*  \DD_{\wt{\mf{h}}}$. Naturality of the counit $f_! f^! \rightarrow \Id$ implies that the map is compatible with the $W_H$-action. 

(2) By base change, the map in question is obtained by taking global sections on $T$ of the composite map 
\[
g^! \pi_*  \DD_{\wt{\mf{h}}}  \xrightarrow{\mrm{unit}} f_* f^* g^! \pi_* \DD_{\wt{\mf{h}}}  \xrightarrow{[f]} f_* f^!g^!  \pi_* \DD_{\wt{\mf{h}}}  \tw{-d(f)}
\]
where $\tw{i} := [2i](2i)$ is a shift and Tate twist, and the natural transformation $[f] \co f^* \rightarrow f^!\tw{-d(f)}$ is induced by the relative fundamental class of $f$, as explained in \cite[\S 3.4]{FYZ3}. Naturality of $[f]$ and the unit $\Id \rightarrow f_* f^*$ imply that the map is compatible with the $W_H$-action. 
\end{proof}

\subsection{Affine Springer action}\label{ssec: affine Springer action} 
Let $\chI \subset L^+\chG_{\F}$ be the Iwahori subgroup corresponding to the fixed Borel subgroup $\chB_{\F} \subset \chG_{\F}$. 

There is an ``affine Springer action'' of $W_{\aff}$ on $\mBM_*(\rY_\gamma^{\varepsilon=0})$ and on $\mBM_*(\rY_\gamma^{\varepsilon = 1})$. For the affine Springer fibers, the action was constructed by Lusztig \cite{Lus96} and Sage \cite{Sa97} for $W_{\aff} \subset \dWext$. An exposition of the construction for $W_{\aff}$ can be found in \cite[\S 2.6.3]{Yun17}. For the deformed affine Springer fibers at $\varepsilon = 1$, an action of $W_{\aff}$ on $\mBM_*(\rY_\gamma^{\varepsilon = 1})$ was constructed by Frenkel-Zhu \cite[\S 6]{FZ10}, by imitating Lusztig's construction. 

\begin{remark}[Extended affine Springer action]
For $\varepsilon = 0$, the $W_{\aff}$-action on  $\mBM_*(\rY_\gamma^{\varepsilon=0})$ was extended to an action of $\dWext$ by Yun in \cite[Theorem 2.5]{Yun14}, which we also call the ``affine Springer action''. 
\end{remark}

We now construct an affine Springer action of $\dWaff$ on the relative Borel--Moore homology (cf. \S \ref{ssec: BM homology}) of $\rY_\gamma^{\varepsilon}$ over $\A^1_{\varepsilon}$, by a slight generalization of \cite{FZ10}. Our construction specializes to Lusztig's construction when $\varepsilon=0$, and Frenkel-Zhu's when $\varepsilon= 1$. 

\begin{notation}For $r \geq 2$, we introduce the notation 
\[
\Ad^{\varepsilon, r}_{g^{-1}}(\gamma) := \Ad_{g^{-1}}(\gamma)  - \varepsilon  t^r \frac{dg^{-1}}{dt}  g.
\]
We are primarily interested in $r=2$, in which case we abbreviate $\Ad^{\varepsilon} := \Ad^{\varepsilon, 2}$. 
\end{notation}

A straightforward calculation shows that 
\begin{equation}\label{eq: adh compose}
\Ad^{\varepsilon,r}_{g_1 g_2}(\gamma) = \Ad^{\varepsilon,r}_{g_1}(\Ad^{\varepsilon,r}_{g_2}(\gamma)).
\end{equation}
For each parahoric subgroup $\chI \subset \mbf{\chP} \subset L\chG_{\F}$, there is a corresponding affine Springer fiber 
\[
\rY_{\mbf{\chP}, \gamma}^\varepsilon := \{ g \in L\chG_{\F} \co \Ad^{\varepsilon,r}_{g^{-1}}(\gamma) \in \Lie \mbf{\chP}\} /\mbf{\chP}.
\] 
That $\rY_{\mbf{\chP}, \gamma}^{\varepsilon}$ is well-defined -- i.e., the condition $\Ad^{\varepsilon,r}_g(\gamma) \in \Lie \mbf{\chP}$ is preserved by right multiplication by $\mbf{\chP}$ (for $r \geq 2$) -- follows from \cite[Lemma 11]{FZ10}. 

Let $\mbf{\chP}^u$ be the pro-unipotent radical of $\mbf{\chP}$. 
Let $L_{\mbf{\chP}}$ be the Levi quotient of $\mbf{\chP}$ and $\mf{l}_{\mbf{\chP}} := \Lie L_{\mbf{\chP}}$. Then there is an \emph{evaluation map} sending $g \mbf{\chP}^u \in \rY_{\mbf{\chP}, \gamma}$ to the reduction of $\Ad^{\varepsilon,r}(g^{-1}) \gamma \in \Lie \mbf{\chP}$ in $\Lie \mbf{\chP} / \Lie \mbf{\chP}^u \cong \mf{l}_{\mbf{\chP}}$, which is well-defined up to the adjoint action of $L_{\mbf{\chP}}$, and thus defines a map
\begin{equation}\label{eq: parahoric evaluation hbar}
\ev \co \rY_{\mbf{\chP}, \gamma}^\varepsilon \rightarrow [\mf{l}_{\mbf{\chP}}/L_{\mbf{\chP}}].
\end{equation}

The following result is a variant of \cite[Proposition 12]{FZ10}. 

\begin{lemma}\label{lem: deformed cartesian square}
There is a natural (in $\mbf{\chP}$) Cartesian square
\begin{equation}\label{eq: affine springer cover hbar}
\begin{tikzcd}
\rY_{\gamma}^{\varepsilon} \ar[d, "\wt{\pi}"] \ar[r, "\ev"] &  {[\wt{\mf{l}}_{\mbf{\chP}}/L_{\mbf{\chP}}]} \ar[d, "\pi"] \\
\rY_{\mbf{\chP}, \gamma}^{\varepsilon} \ar[r, "\ev"] & {[\mf{l}_{\mbf{\chP}}/L_{\mbf{\chP}}]}
\end{tikzcd}
\end{equation}
where $\pi$ is the Grothendieck alteration for $\mf{l}_{\mbf{\chP}}$. 
\end{lemma}

\begin{proof}
Let $g\mbf{\chP}\in\rY_{\mbf{\chP}, \gamma}^{\varepsilon}$ and set
\[
\gamma':=\Ad_{g^{-1}}^{\varepsilon,r}(\gamma)\in \Lie \mbf{\chP}.	
\]
The fiber $\wt{\pi}^{-1}(g)$ consists of $gx\mbf{\chI}$ such that $x\in\mbf{\chP}$ and (using \eqref{eq: adh compose})
\[
\Ad_{(gx)^{-1}}^{\varepsilon,r}(\gamma)=\Ad_{x^{-1}}^{\varepsilon,r}(\gamma')\in \Lie \chI.
\]
But
\begin{equation}\label{eq:x-springer-condition}
\Ad_{x^{-1}}^{\varepsilon,r} (\gamma')= \Ad_{x^{-1}}(\gamma')  - \varepsilon t^r \frac{d (x^{-1} )}{dt} x
\end{equation}
and the second term on the RHS belongs to $\Lie \mbf{\chP}^u \subset \Lie \chI$ by \cite[Proposition 11]{FZ10}.
Thus the condition \eqref{eq:x-springer-condition} reduces to $\Ad_{x^{-1}}(\gamma')\in \Lie \chI$. Since $x \in \chP$ and $\gamma' \in \Lie \mbf{\chP}$, the element $\Ad_{x^{-1}}(\gamma')$ lies in $\mbf{\chP}$, and then it lies in $\Lie \chI$ if and only if it lies in $\mf{l}_{\mbf{\chP}}$ modulo $\Lie  \mbf{\chP}^u$. Since $\ev(g\mbf{\chP})$ is exactly the class of $\gamma'$ in $[\mf{l}_{\mbf{\chP}}/L_{\mbf{\chP}}]$, the result follows.

\end{proof}

As explained in \S \ref{ssec: groth alteration}, Lemma \ref{lem: deformed cartesian square} induces an action of the Weyl group $W_{\mbf{\chP}}$ associated to $\mbf{\chP}$ on $\mBM_*(\rY_{\gamma}^\varepsilon)$. Then just as in the case of affine Springer fibers \cite[\S 5.5]{Lus96}, these actions glue to an action of $W_{\aff}$ on $\mBM_*(\rY_{\gamma}^{\varepsilon})$. It will actually be more convenient for us to normalize the action as a \emph{right} action, using the antipode map on $W_{\aff}$. 

\begin{remark}[Compatibility with translation action]
It is immediate from the definitions that the resulting $W_{\aff}$-action commutes with the translation action on $\mBM_*(\rY_\gamma^{\varepsilon = 0})$. 
\end{remark}

Let us summarize the upshot of this construction. 

\begin{prop}[Affine Springer action]\label{prop: affine springer action}
For any $S \rightarrow \Spec \F[\varepsilon]$, letting $\rY_\gamma^\varepsilon|_S $ denote the base change of $\rY_\gamma^\varepsilon$ to $S$, there is a right action of $W_{\aff}$ on $H^{\BM}_*(\rY_\gamma^\varepsilon|_S /S)$, with the following properties. 
\begin{enumerate}
\item If $s$ is a geometric point over $\varepsilon = 0$, then it is the usual affine Springer action of Lusztig. 
\item If $s$ is a geometric point over $\varepsilon = 1$, then it is the action of Frenkel-Zhu \cite{FZ10}. 
\item It commutes with specialization in $\varepsilon$. 
\end{enumerate}
\end{prop}
 
\begin{proof}
The first two points are clear from the construction. The third follows from the construction of the specialization map, Lemma \ref{lem: Springer equivariant}, and the construction of the $W_{\aff}$-action. 
\end{proof}

\subsection{Parametrization of top-dimensional irreducible components: generic fiber}\label{ssec: top-irred-comp generic}  Here we study certain top-dimensional irreducible components of $\rX_{\gamma}^{\varepsilon}$ for $\gamma_0=(t+p)s$ where $s\in \chft$ is regular. (Recall that for $\gamma \in \chfg(\cO)[\varepsilon]((t+p))$, we write $\gamma_0 \in \chfg(\cO)((t))$ for the evaluation of $\gamma$ at $\varepsilon = 0$.)

For a dominant coweight $\lambda \in X_*(\chT)^+$, let $S^\circ(\lambda) \subset \Gr_{\chG, E}$ be the corresponding (open) Schubert cell and $S(\lambda)$ its closure.

For $\lambda \in X_*(\chT)$, let $\rX^{\varepsilon }_{\gamma}(\lambda) $ be the Zariski closure of $\rX_\gamma^{\varepsilon} \cap S^\circ(\lambda)$. For any $\varepsilon_0\in \A^1_{\varepsilon}$, we also define $\rX_\gamma^{\varepsilon=\varepsilon_0}(\lambda)$ as the Zariski closure of $\rX_\gamma^{\varepsilon=\varepsilon_0} \cap S^\circ(\lambda)$.

\begin{warning}\label{warning: closure base change} It is evident that $\rX_\gamma^{\varepsilon=\varepsilon_0}(\lambda)$ is a closed subscheme of the fiber at $\varepsilon_0$ of $\rX_\gamma^{\varepsilon}(\lambda)$, which we denote $\rX_\gamma^{\varepsilon}(\lambda)|_{\varepsilon_0}$, but it is not clear (at least to the authors) whether this closed embedding is an isomorphism. 
\end{warning}

\begin{lemma}\label{lem: X lambda irred e=0}
Assume $\gamma_0=(t+p)s$ with $s\in \chft$ regular. Then there is an open subscheme $V\subset \A^1_{E}$ containing $0$ and an open subscheme $U\subset S^\circ(\rho)$ such that 
\begin{itemize}
\item $(U\times V)\cap \rX^{\varepsilon}_{\gamma}\cong \A^d\times V$
where we recall that $d=\dim (\chG/\chB)$, and
\item Every fiber of $\rX^{\varepsilon}_{\gamma}(\rho)\setminus U\times V$ over $V$ has dimension strictly less than $d$.
\end{itemize}
In particular, if $\lambda < \rho$ then $\rX^{\varepsilon = 0}_{\gamma}(\lambda)$ has dimension strictly less than $d$, and $\rX^{\varepsilon = 0}_{\gamma}(\rho) $ is irreducible of dimension $d$. 
\end{lemma}

\begin{proof}
We explain how this essentially follows from the computations in the proof of \cite[Proposition 3.3.4]{LLLM22}, in particular how to convert notations in ~\emph{loc.cit.} to our situation. 

We have a standard affine open cover of $S^\circ(\lambda)$ by translates of Iwahori orbits, 
\[
\{ U_w(\lambda) :=w^{-1}I(\lambda)(t+p)^{\lambda} \}_{w\in W}.
\]
Here $I(\lambda)$ is the affine space of dimension $\langle \lambda, 2\rho^\vee\rangle$ given by a root group decomposition (after fixing any ordering of the roots)
\[
I(\lambda) \cong \prod_{\alpha \in \chPhi}\prod_{i=0}^{\langle \lambda, \alpha \rangle-1}U_{\alpha,i},
\]
where $U_{\alpha, i} \subset L\chG$ is the affine root group corresponding to the affine root $(\alpha, i) \in \chPhi \times \Z$. 
(In the notation of the proof of \cite[Proposition 3.3.4]{LLLM22}, $I(\lambda)$ is the transpose of $(v-t)^{-\lambda}\wt{N}_\lambda$, noting that $v$, $t$, $\mathrm{Diag}(\bf{a})$ in ~\emph{loc.cit.} correspond to $t,-p, s$ in our situation.)

Choosing coordinates $x_{\alpha,i}$ for the affine root groups appearing in $I(\lambda)$ identifies $(U_w(\lambda) \times \A^1_{\varepsilon})\cap \rX^{\varepsilon}_{\gamma}$ with the subspace of $I(\lambda)$ cut out by the equations
\[
 (\varepsilon i-\langle \mathrm{Ad}_w(s),\alpha \rangle )x_{\alpha,i}-p(i+1)\varepsilon x_{\alpha,i+1}=O(X_{\beta,j})\
 \]
for $0\leq i<\langle \lambda ,\alpha \rangle-1$. Here the right-hand side is an expression in the $x_{\beta,j}$ with $0<\beta<\alpha$. Over the locus $V\subset \A^1_{\varepsilon}$ where $\varepsilon i-\langle \Ad_w(s),\alpha \rangle$ is invertible, these equations cut out an affine subspace of dimension equal to $\dim (\chG/\chP_{\lambda})$, with coordinates in the $x_{\alpha,\langle \lambda, \alpha\rangle-1
}$. In particular, we see that $\rX^{\varepsilon}_{\gamma}(\lambda)|_V$ has fiberwise dimension $<d$ for $\lambda<\rho$.

Now, for each simple root $\alpha>0$ of $\chG$, the coordinate $x_{\alpha,0}$ is fiberwise over $V$ non-vanishing on $\rX^{\varepsilon}_{\gamma}\cap (U_w(\rho)\times V)$ since $\langle \rho, \alpha \rangle = 1 $, hence the intersection of
$\rX^{\varepsilon}_{\gamma}\cap (U_{w}(\rho)\times V)$ and $\rX^{\varepsilon}_{\gamma}\cap (U_{s_\alpha w}(\rho)\times V)$ is fiberwise open dense in either space. It follows that $\rX^{\varepsilon}_{\gamma}\cap (U_{w}(\rho)\times V)$ is fiberwise open dense in $\rX^{\varepsilon}_{\gamma}\cap (S^\circ(\rho)\times V)$, hence the complement has fiberwise dimension strictly less than $d$.

\begin{remark}\label{rmk:AFS:char 0} \begin{enumerate}

\item
For $\lambda\neq \rho$ regular, $\rX^{\varepsilon = 0}_\gamma(\lambda)$ is never irreducible, which is in stark contrast to the behavior of the deformed affine Springer fiber (see Lemma \ref{lem: X lambda irred e neq 0} below).
\item By \cite[Corollary 3.10.2]{Ngo10}, $\rX^{\varepsilon = 0}_{\gamma}$ is equidimensional, and the translation action by $L\chT_{E}$ is transitive on the set of irreducible components. In fact, one can check that $\rX^{\varepsilon = 0}_{\gamma}\cap S^\circ(\rho)$ is a fundamental domain for the $L\chT_{E}$-action on the regular locus of $\rX^{\varepsilon = 0}_\gamma$, and that the regular locus consists of exactly the $(L\chT/L^+\chT)_{E}=X_*(\chT)$ translates of $\rX^{\varepsilon = 0}_\gamma\cap S^\circ(\rho)$. This shows that for regular $\lambda\neq \rho$, $\rX^{\varepsilon = 0}(\lambda)$ is a (reducible) union of translates of $\rX^{\varepsilon = 0}(\rho)$.
\end{enumerate}
\end{remark}
\end{proof}

\begin{lemma}\label{lem: X lambda irred e neq 0}  Let $\rX_\gamma^{\varepsilon \neq 0}$ be the fiber of $\rX_\gamma^{\varepsilon}$ over $\G_m \subset \A^1_{\varepsilon}$. For $\lambda \in X_*(\chT)^+$, the Bialynicki-Birula map induces an isomorphism $S^\circ(\lambda) \cap \rX_\gamma^{\varepsilon \neq 0 } \xrightarrow{\sim} \chG/\chP_\lambda \times \G_{m, E}$ as a family over $\G_{m,E} = \Spec E[\varepsilon^{\pm 1}]$. 
\end{lemma}

\begin{proof} This is a consequence of the computation in the proof of Lemma \ref{lem: X lambda irred e=0}. Since $\varepsilon p$ is invertible under the hypotheses, in the charts defined in the proof of Lemma \ref{lem: X lambda irred e=0} we can solve all the $x_{\alpha,i}$ in terms of $x_{\alpha,0}$. But the Bialynicki-Birula map on these charts is exactly the map extracting the $x_{\alpha,0}$.

\end{proof}

In particular, the Bialynicki-Birula map induces an isomorphism $S^\circ(\lambda) \cap \rX_\gamma^{\varepsilon = 1} \xrightarrow{\sim} (\chG/\chP_\lambda)_{E}$. Therefore, the map 
\[
\lambda \mapsto S^\circ(\lambda) \cap \rX_\gamma^{\varepsilon = 1}
\]
induces a bijection from the set of regular dominant weights $X_*(\chT)^+$ to the top-dimensional irreducible components of $\rX_\gamma^{\varepsilon = 1}$. A similar discussion applies for $\varepsilon = \eta$ (the generic point of $\A^1_{\varepsilon}$).

For $\varepsilon_0 \in \G_{m, E}$, define
\[
\rX^{\varepsilon = \varepsilon_0}_{\gamma}(\leq \lambda) = \rX_\gamma^{\varepsilon = \varepsilon_0} \cap S(\lambda);
\]
which is a disjoint union of partial flag varieties $\rX^{\varepsilon=\varepsilon_0}_\gamma(\lambda')$ for $\lambda'\leq \lambda$. We also define $\cX^{\varepsilon=\varepsilon_0}_\gamma(\le \lambda)$, $\cX^{\varepsilon=\varepsilon_0}_\gamma(\lambda')$ to be the corresponding closures in $\cX^{\varepsilon=\varepsilon_0}_{\gamma}$.
We are particularly interested in $\varepsilon_0 \in \{ 1, \eta\}$.

\subsection{Parametrization of top-dimensional irreducible components: special fiber}\label{ssec: top-irred-comp special}
Here we establish a combinatorial parametrization of the top-dimensional irreducible components of $\rY_{\gamma}^{\varepsilon =  \varepsilon_0}$ for any $\varepsilon_0$.

 For each $\wt{w} \in \dWext$, let $S^\circ(\wt{w}) =  \chI \wt{w} \chI / \chI  \subset \Fl_{\chG, \F}$ be the corresponding (open) Schubert cell and $S(\wt{w}) \subset \Fl_{\chG, \F}$ its closure. 
 
\begin{defn}
For $\varepsilon_0\in \A^1_{\varepsilon}$ and $\wt{w} \in \dWext$, let 
\[
\rY^{\varepsilon = \varepsilon_0}_{\gamma}(\wt{w})^\circ := \rY_\gamma^{\varepsilon = \varepsilon_0} \cap S^\circ(\wt{w})
\]
and $\rY^{\varepsilon = \varepsilon_0}_{\gamma}(\wt{w}) $ be its closure in $\rY^{\varepsilon = \varepsilon_0}_{\gamma}$. (As in Warning \ref{warning: closure base change}, this construction may not commute with the base change in $\varepsilon$.)
\end{defn}

For the statements below, recall the notation $h_{\wt{w}}$ from \S \ref{sssec: weyl notation}. 

\begin{lemma}\label{lem: Y lambda irred e=0} Let $\gamma=t(s+\varepsilon r)$ with $r,s\in \chft$. Assume $s$ is regular semisimple.
\begin{enumerate}
\item $\rY^{\varepsilon = 0}_{\gamma}(\wt{w})^\circ$ is an affine space over $\F$ of dimension 
\[
\dim \rY^{\varepsilon = 0}_{\gamma}(\wt{w})^\circ = d - \# \{ \alpha \in \chPhi^{+} \mid \wt{w}(A_0) \subset H^{(0, 1)}_{\alpha} \}.
\]
In particular $\dim \rY^{\varepsilon = 0}_{\gamma}(\wt{w})^\circ=d$ if and only if $\wt{w}$ is regular.

\item $\rY^{\varepsilon = \eta}_{\gamma}(\wt{w})^\circ$ is an affine space over $\Spec \F(\varepsilon)$ of dimension 
\[
\dim \rY^{\varepsilon = \eta}_{\gamma}(\wt{w})^\circ =  d - \# \{ \alpha \in \chPhi^{+} \mid \wt{w}(A_0) \subset H^{(0, 1)}_{\alpha} \}.
\]

\item If $s+r$ is $h_{\wt{w}}$-generic, then $\rY_{\gamma}^{\varepsilon = 1}(\wt{w})^\circ$ is an affine space over $\F$ of dimension 
\[
\dim \rY_{\gamma}^{\varepsilon = 1}(\wt{w})^\circ =  d - \# \{ \alpha \in \chPhi^{+} \mid \wt{w}(A_0) \subset H^{(0, 1)}_{\alpha} \}.
\]
\end{enumerate}
\end{lemma}

\begin{proof}
As in the proof of Lemma \ref{lem: X lambda irred e=0}, we explain how this follows from modifying the proof of \cite[Theorem 4.2.4]{LLLM22}.

The main point is that in \cite[Equation (4.6)]{LLLM22}, the coefficient $(i + \delta_{\alpha > 0} + \langle \ol{\bf{a}},\alpha \rangle)$ becomes 
\[
\varepsilon(i+\delta_{\alpha>0})+\langle s+\varepsilon r,\alpha \rangle
\]
and we need this to be invertible for all $0\leq i <d_{\alpha,\wt{w}}$ and all roots $\alpha$. Our hypotheses are exactly arranged for this to be the case.

\end{proof}

\begin{remark}
The case $\varepsilon = 0$ of Lemma \ref{lem: Y lambda irred e=0} can also be found in \cite[proof of Lemma 2.6 (c)]{BBSV}.

\end{remark}

\begin{example}\label{ex: adm set h} 
If $\wt{w} \in \Adm(\lambda)$, then $h_{\wt{w}} \leq h_\lambda$. In this case, we find that $S^\circ(\wt{w}) \cap \rY_\gamma^{\varepsilon =1}$ has top dimension (equal to $d$) if and only if $\wt{w} \in \Adm^{\reg}(\lambda)$.
\end{example}

For $\varepsilon_0 \in \A^1_{\varepsilon}$, define the ``$\lambda$-admissible'' part of $\rY^{\varepsilon = \varepsilon_0}_{\gamma}$ to be 
\[
\rY^{\varepsilon = \varepsilon_0}_{\gamma}(\leq \lambda) :=  \rY_\gamma^{\varepsilon = \varepsilon_0} \cap  \left( \bigcup_{\wt{w} \in W \cdot \lambda} S(\wt{w}) \right). 
\]


\begin{cor}\label{cor: top cycles} Let $\gamma=t(s+\varepsilon r)$ with $r,s\in \mathfrak{t}$. 
\begin{enumerate}
\item Assume $s$ is regular. Then for $\varepsilon_0 \in \{0, \eta\}$, the map $\wt{w} \mapsto \rY_\gamma^{\varepsilon=\varepsilon_0}(\wt{w})$ induces a bijection between $\dWext^{\reg}$ and the top-dimensional irreducible components of $\rY_\gamma^{\varepsilon = \varepsilon_0}$. 

\item If $s+r$ is $h_\lambda$-generic, then the map $\wt{w} \mapsto \rY_\gamma^{\varepsilon=1}(\wt{w})$ induces a bijection between $\Adm^{\reg}(\lambda)$ and the top-dimensional irreducible components of $\rY_\gamma^{\varepsilon = 1}(\leq \lambda)$. 

\end{enumerate}
\end{cor}

\begin{proof}
Follows immediately from Lemma \ref{lem: Y lambda irred e=0}.
\end{proof}

We establish some basic properties of $\rY^{\varepsilon=0}_{\gamma}(\wt{w})$. The reader may need to review \S \ref{sssec: weyl notation}  for notation. 

\begin{lemma}\label{lem:transformation of components} Let $\gamma=t(s+\varepsilon r)$ with $s,r\in \chft$, $s$ regular. Below, equalities are as subschemes of $\Fl_{\chG}$.
\begin{enumerate}
\item Left multiplication by $t^{\nu}w \in \wt{W}$ induces
\[
t^{\nu}w \rY^{\varepsilon}_{\gamma}= \rY^{\varepsilon}_{w^{-1}\gamma+t\varepsilon \nu}.
\]
\item  Let $\wt{u}\in\wt{W}^{\reg}$ and factorize $\wt{u}$ uniquely as $\wt{u}=   \wt{v}^{-1}w_0 \wt{w}$ where $\wt{w}\in \wt{W}_1$ is restricted and $\wt{v}=v t^{\nu}\in \wt{W}^+$ is dominant, cf. (the proof of) \cite[Proposition 2.1.5]{LLLM22}. 
\begin{enumerate}
\item For $\varepsilon_0\in \{0,\eta\}$, we have
\[
\rY^{\varepsilon=\varepsilon_0}_{\gamma}(\wt{u})=t^{-\nu}v^{-1}\rY^{\varepsilon=\varepsilon_0}_{v\gamma-t\varepsilon v\nu}(w_0\wt{w}).
\]
\item If furthermore $\wt{u}\in \Adm^{\reg}(\rho)$, and both $s+r$, $s+r-\nu$ are $h_\rho$-generic, then we have
\[
\rY^{\varepsilon=1}_{\gamma}(\wt{u})=t^{-\nu}v^{-1} \rY^{\varepsilon=1}_{v\gamma-t v\nu}(w_0\wt{w}) .
\]
\end{enumerate}
\item Let $\wt{w}\in \wt{W}_1$ be restricted and $\sigma\in W$. For any $\varepsilon_0\in \{0,1,\eta\}$, we have
\[
\rY^{\varepsilon=\varepsilon_0}_{\gamma}(w_0\wt{w})=\sigma^{-1}\rY^{\varepsilon=\varepsilon_0}_{\sigma^{-1}\gamma}(w_0\wt{w}).
\]
\item For $\wt{w}=t^{\rho_w}w \in \wt{W}_1$ with $w\in W$, we have
\[
\rY^{\varepsilon=\varepsilon_0}_{\gamma}(w_0\wt{w})^\circ \subset (L\chN\cap \chI)t^{w_0\rho_w}w_0w\chI/\chI\subset t^{w_0\rho_w}(L\chN) w_0w\chI/\chI.
\] 
\end{enumerate}
\end{lemma}
\begin{proof}Assertion (1) follows from a direct computation, using that 
\[
\Ad_{t^{\nu}w}^{\varepsilon}(\gamma)=w^{-1}(\gamma)+t\varepsilon \nu.
\]

Assertions (2) and (3) follow from the proof of \cite[Proposition 4.3.5]{LLLM22} and \cite[Proposition 4.3.6]{LLLM22} respectively, noting that:
\begin{itemize}
\item The cited proofs did not use the running assumption that $\chG=\GL_n$ in \emph{loc.cit.}
\item The role of the genericity assumptions in the cited proofs were only used to guarantee that intersecting $\rY^{\varepsilon=\varepsilon_0}_\gamma$ with certain open affine Schubert cells are affine spaces of the correct dimension. This holds for $\varepsilon_0=0$ (and hence also for $\varepsilon_0=\eta$) due to the regularity of $s$, and holds for $\varepsilon_0=1$ by our assumptions, cf. Lemma \ref{lem: Y lambda irred e=0}.
\end{itemize}

Finally, (4) follows from writing standard representatives for $S^\circ(w_0\wt{w})$ and the fact that $w_0\wt{w}$ is anti-dominant, see for example \cite[Corollary 4.2.15]{LLLM22}.

\end{proof}

\section{Specialization of cycles}\label{ssec: comparison of families}

In the previous section we parametrized some irreducible components of $\rX^{\varepsilon}_{\gamma}$ and $\rY^{\varepsilon}_{\gamma}$ over $\varepsilon \in \{0,1,\eta\}$. In this section we study the behavior of these irreducible components under specialization in $\varepsilon$. 

\subsection{Bases for top homology}
We set up some degeneration problems. We have a family $\cX_{\gamma}^{\varepsilon}\rightarrow  \A^1_{\cO}$, 
\[
\begin{tikzcd}
\rX_{\gamma, E}^{\varepsilon} \ar[r, hook] \ar[d] & \cX_{\gamma}^{\varepsilon} \ar[d] & \rY_{\gamma, \F}^{\varepsilon} \ar[d]  \ar[l, hook'] \\
\Spec E  \ar[r, hook] & \Spec \cO & \Spec \F \ar[l, hook'] 
\end{tikzcd}
\]

 Below when we say a subscheme of $\rY^{\varepsilon = \varepsilon_0}_{\gamma}$ is ``top-dimensional'', we mean that its dimension is equal to that of the ambient space $\rY^{\varepsilon = \varepsilon_0}_{\gamma}$, which is $d = \dim \chG/\chB$. We write $\topCh(-)$ for the top-degree Chow group and $\topBM(-)$ for the top-degree \emph{geometric} Borel--Moore homology group; in the latter case the ``top'' degree is $2d$. Recall our notation for cycle classes from \S \ref{sssec: cycle class map}. 
 
\begin{lemma}\label{lem: bases}  Let $\gamma=t(s+\varepsilon r)$ with $r,s\in \mathfrak{t}$. 
\begin{enumerate}

\item Assume $s$ is regular. Then for $\varepsilon_0 \in \{0, \eta\}$, the cycle classes $[\rY_{\gamma}^{\varepsilon=\varepsilon_0}(\wt{w})]$ form a basis for $\topCh(\rY_{\gamma}^{\varepsilon=\varepsilon_0})$ as $\wt{w}$ ranges over $\dWext^{\reg}$. 

\item If $s+r$ is $h_\lambda$-generic. Then the cycle classes $[\rY_{\gamma}^{\varepsilon=1}(\wt{w})]$ form a basis for $\topCh(\rY_{\gamma}^{\varepsilon=1}(\leq \lambda))$ as $\wt{w}$ ranges over $\Adm^{\reg}(\lambda)$. 
\end{enumerate}
\end{lemma}

\begin{proof}
Follows immediately from Corollary \ref{cor: top cycles}. 
\end{proof}

\subsection{Specialization in $\varepsilon$} In this subsection we analyze the behavior of specialization in $\varepsilon$. 

\subsubsection{Generic fiber}

Below for a closed point $\varepsilon_0 \in \A^1_k$ over a field $k$, we write $\A^1_{(\varepsilon_0)}$ for the localization of $\A^1$ at $\varepsilon_0$, which is evidently a discrete valuation ring. For any family over $\A^1_{(\varepsilon_0)}$, there is then a specialization map from the Chow groups or Borel--Moore homology groups of the generic fiber to those of the special fiber (cf. \S \ref{ssec: specialization for chow} and \S \ref{ssec: specialization for BM homology}). 

Localizing the family $\rX^{\varepsilon}_{\gamma}$ over $\A^1_{(0)}$ (over $E$), we have a specialization map 
\begin{equation}\label{eq: gr hbar specialization 0}
\fsp_{\varepsilon \rightarrow 0} \co \topCh(\rX^{\varepsilon = \eta}_{\gamma}) \rightarrow 
\topCh(\rX^{\varepsilon=0}_{\gamma}).
\end{equation}
Localizing the family $\rX^{\varepsilon}_{\gamma}$ over $\A^1_{(1)}$ (over $E$), we have a specialization map 
\begin{equation}\label{eq: gr hbar specialization 1}
\fsp_{\varepsilon \rightarrow 1} \co \topCh(\rX^{\varepsilon = \eta}_{\gamma}) \rightarrow 
\topCh(\rX^{\varepsilon=1}_{\gamma}).
\end{equation}

\begin{prop}\label{prop: gr hbar spc}  Let $\gamma = (t+p)s $ with $s \in \chft$ regular. 

 (1) The map \eqref{eq: gr hbar specialization 1} is an isomorphism and sends $[\rX_{\gamma}^{\varepsilon =\eta}(\lambda)] \mapsto [\rX_{\gamma}^{\varepsilon=1}(\lambda)]$ for all $\lambda \in X_*(T)^+$. 

(2) The map \eqref{eq: gr hbar specialization 0} sends $[\rX_{\gamma}^{\varepsilon =\eta}(\rho)] \mapsto [\rX_{\gamma}^{\varepsilon=0}(\rho)]$. 
\end{prop}

\begin{proof}
(1) This follows from Lemma \ref{lem: X lambda irred e neq 0}. 

(2) By Lemma \ref{lem: X lambda irred e=0}, there are open subsets $V\subset \A^1_{E}  = \Spec E[\varepsilon]$, $U\subset S^\circ(\rho)$ such that $\cU^{\varepsilon} :=\rX_{\gamma}^{\varepsilon}\cap (V\times U)\cong V\times \A_{E}^d$, and such that its complement has fiberwise dimension less than $d$ over $V$. Consequently, we have a commutative diagram
\begin{equation}\label{eq: generic spc eq2}
\begin{tikzcd}
\Ch_{\mrm{top}}(\cU^{\varepsilon=\eta})   \ar[r, "\spc_{\varepsilon \rightarrow 0}"]  &\Ch_{\mrm{top}}(\cU^{\varepsilon=0})  \\
\Ch_{\mrm{top}}(\rX_{\gamma}^{\varepsilon=\eta}(\rho))  \ar[r, "\spc_{\varepsilon \rightarrow 0}"] \ar[u] &  \Ch_{\mrm{top}}(\rX_{\gamma}^{\varepsilon = 0}(\rho))  \ar[u]
\end{tikzcd}
\end{equation}
The bound on the fiber dimension of the complement of $\cU^{\varepsilon}$ implies that the vertical restriction maps are isomorphisms. Since $\cU^{\varepsilon} \cong V\times \A^d_{E}$ is the trivial family, we have $\spc_{\varepsilon \rightarrow 0}[\cU^{\varepsilon=\eta}]=[\cU^{\varepsilon=0}]$, and the result follows.

\end{proof}

\subsubsection{Special fiber} Let $\gamma = t(s+\varepsilon r)$ with $s,r \in \chft$. Assume $s$ is regular and let $\lambda \in X_*(\chT)^+$ be a dominant coweight such that $s+r$ is $h_\lambda$-generic. 

Localizing the family $\rY^{\varepsilon}_{\gamma}$ over $\A^1_{(0)}$ (over $\F)$, we have a specialization map
\begin{equation}\label{eq: fl hbar specialization 0}
\fsp_{\varepsilon \rightarrow 0} \co \topCh(\rY^{\varepsilon = \eta}_\gamma (\leq \lambda))  \rightarrow 
\topCh(\rY^{\varepsilon=0}_{\gamma}(\leq \lambda)).
\end{equation}
Localizing the family $\rY^{\varepsilon}_{\gamma}$ over $\A^1_{(1)}$ (over $\F)$, we have a specialization map
\begin{equation}\label{eq: fl hbar specialization 1}
\fsp_{\varepsilon \rightarrow 1} \co \topCh(\rY^{\varepsilon =\eta}_\gamma (\leq \lambda)) \rightarrow 
\topCh(\rY^{\varepsilon=1}_{\gamma}(\leq \lambda)).
\end{equation}
Lemma \ref{lem: bases} implies that for $\varepsilon_0 \in \{0, 1, \eta\}$, the classes $[\rY_\gamma^{\varepsilon = \varepsilon_0}(\wt{w}')] \in \topCh(\rY_\gamma^{\varepsilon = \varepsilon_0}(\leq \lambda))$ form a basis as $\wt{w}'$ varies over $\Adm^{\reg}(\lambda)$. Therefore, for each $\varepsilon_0 \in \{0,1\}$ there exists a unique matrix 
\[
M^{\varepsilon = \varepsilon_0} := (m_{\wt{w} \wt{w}'}^{\varepsilon = \varepsilon_0} \in \Z_{\geq 0})_{\wt{w},\wt{w}' \in \Adm^{\reg}(\lambda)}
\]
(the non-negativity by Lemma \ref{lem: specialization preserves effectivity}) such that 
\[
\fsp_{\varepsilon \rightarrow 0}  [\rY_\gamma^{\varepsilon = \eta}(\wt{w})] = \sum_{\wt{w}'  \in \Adm^{\reg}(\lambda)} m_{\wt{w} \wt{w}'}^{\varepsilon = \varepsilon_0} [\rY_\gamma^{\varepsilon = \varepsilon_0}(\wt{w}')] \quad \text{ for all } \wt{w} \in \Adm^{\reg}(\lambda).
\]

\begin{remark}[Independence of $\lambda$] The definition of $m_{\wt{w} \wt{w}'}^{\varepsilon= \varepsilon_0}$ appears to depend on $\lambda$. However, there is the following sense in which it is independent of $\lambda$ as long as it is defined: if $\lambda \leq \lambda'$ and $\wt{w}, \wt{w}' \in \Adm^{\reg}(\lambda) \subset \Adm^{\reg}(\lambda')$, then each $m_{\wt{w} \wt{w}'}^{\varepsilon = \varepsilon_0}$ is the same whether defined in terms of $\lambda$ or $\lambda'$. This is a consequence of pushforward functoriality (cf. \S \ref{sssec: functoriality for specialization maps}) for the closed embedding $\rY^{\varepsilon}_\gamma (\leq \lambda) \inj \rY^{\varepsilon}_\gamma (\leq \lambda')$. 
\end{remark}

\begin{lemma}\label{lem: kappa specialization upper triangular}  Maintain the running assumptions on $s,r$ and $\lambda$. Then the matrix $M^{\varepsilon = \varepsilon_0} :=  (m_{\wt{w} \wt{w}'}^{\varepsilon = \varepsilon_0})$ is unipotent and upper-triangular with respect to the Bruhat order on $\dWext$, with all entries non-negative and $m_{\wt{w} \wt{w}}^{\varepsilon = \varepsilon_0}=1$ for all $\wt{w} \in \Adm^{\reg}(\lambda)$. 


In particular, \eqref{eq: fl hbar specialization 0} and \eqref{eq: fl hbar specialization 1} are isomorphisms. 
\end{lemma}

\begin{proof} We prove the statement for $\varepsilon_0 = 0$, the argument for $\varepsilon_0=1$ being similar. The closure of $\rY_\gamma^{\varepsilon = \eta}(\wt{w})$ is a closed subset of $\A^1_{\F}\times S(\wt{w})$, hence its fiber over $\varepsilon=0$ is a subscheme of $S(\wt{w})$. Since $m_{\wt{w} \wt{w}'}^{\varepsilon = 0}\neq 0$ if and only if $\rY_\gamma^{\varepsilon = 0}(\wt{w})$ occurs in this closure, upper triangularity follows. The fact that $m_{\wt{w} \wt{w}}^{\varepsilon = 0}=1$ follows from the fact that $\rY_\gamma^{\varepsilon = \eta}(\wt{w})\cap (V\times S^\circ(\wt{w}))\cong V\times \A^d_{\F}$ is the trivial family for some open $V\subset \A^1_{\varepsilon}$ containing $0$ (as in the proof of Proposition \ref{prop: gr hbar spc}).
\end{proof}




\begin{defn}[Deforming from $\varepsilon = 1$ to $\varepsilon =0$]\label{def: specialization 1 to 0}
Suppose $\gamma$ is $h_\lambda$-generic. Then we abuse notation by defining 
\[
\fsp_{\varepsilon \rightarrow 0} \co \topCh(\rY^{\varepsilon = 1}_\gamma (\leq \lambda))  \rightarrow 
\topCh(\rY^{\varepsilon=0}_{\gamma}(\leq \lambda))
\]
to be the composition of \eqref{eq: fl hbar specialization 0} with the inverse of \eqref{eq: fl hbar specialization 1}. 
\end{defn}

\part{Microlocal analysis} \label{part: II}
Most of this Part focuses on affine Springer fibers (and not their deformations), so we will use the abbreviations $\cX_{\gamma} := \cX^{\varepsilon = 0}_{\gamma}$, $\rX_\gamma := \rX_\gamma^{\varepsilon = 0}$, and $\rY_\gamma := \rY_\gamma^{\varepsilon = 0}$. 

Henceforth we use $\rK(-) := K_0(-)$ for the Grothendieck group of an exact category, in order to improve the readability of the notation, because there will be many subscripts (including $0$) on the categories, and we never consider higher K-groups anyways. 

Recall that $h = h_\rho+1$ is the maximum of the Coxeter numbers of the simple factors of $\chG$, and that we assume throughout that $p>h$.  

\section{Equivariant homology of affine Springer fibers} 
In this section we recall some notions from equivariant (co)homology, which will be applied to the (deformed) affine Springer fibers. We begin with a review of generalities such as equivariant Borel--Moore homology, equivariant formality, the equivariant localization theorem, and equivariant Euler classes in \S \ref{ssec: equiv BM homology} -- \S \ref{ssec: equivariant Euler}. 

Then in \S \ref{ssec: GKM} we recall the ``GKM description'' (named after Goresky--Kottwitz--MacPherson), which gives a combinatorial description of the equivariant Borel--Moore homology of spaces satisfying the so-called ``GKM conditions''. We explicate the GKM description for the affine Springer fibers $\rY_\gamma$ and their variants over the complex numbers. This is used in particular to fix identifications $\topBM(\rY_{\gamma})$ for all $\gamma = ts$, where $s \in \chft_{\F}$ is regular, as well as to establish a bridge between the homology of affine Springer fibers in characteristic $p$ and the homology of their variants over $\C$. 

Finally, in \S \ref{ssec: equivariant actions} we discuss actions on the equivariant Borel--Moore homology of the $\rY_\gamma$, which in particular provide equivariant lifts of the translation action and affine Springer action defined in Part 1. 

\subsection{Equivariant Borel--Moore homology}\label{ssec: equiv BM homology}

Suppose $X$ is a finite type scheme over a field $k$, equipped with the action of an algebraic group $H$. Recall that the $H$-equivariant cohomology of $X$, denoted $H^*_H(X)$, is the cohomology of the quotient stack $[X/H]$. We define $\rH^*_H(X)$ to be the \emph{geometric} cohomology of $[X/H]$. 

Similarly, the \emph{$H$-equivariant ($\ell$-adic) Borel--Moore homology} $H^{\BM, H}_*(X)$ is defined as the relative Borel--Moore homology of the quotient stack $X/H$ relative to $[\Spec k / H]$,
\[
H^{\BM, H}_*(X) := H^{\BM}_*([X/H] / [\Spec k /H]). 
\] 
We denote by $\rH^{\BM, H}_*(X)$ the \emph{geometric} $H$-equivariant Borel--Moore homology, i.e., the same definition after base changing to a separable closure of $k$. 


For an ind-scheme $X = \colim_i X_i$ equipped with compatible $H$-actions on each $X_i$, we define 
\[
\rH^{\mrm{BM}, H}_*(X) := \colim_i \rH^{\mrm{BM}, H}_*(X_i)
\]
with transition maps induced by the closed embeddings $X_i \inj X_j$.

\begin{example}[Torus actions]\label{ex: torus} In this paper, we will only ever consider equivariant Borel--Moore homology with respect to the action of a split torus $T$. Equivariant Borel--Moore homology has a natural module structure over equivariant cohomology. Therefore $\rH^{\BM,T}_*(X)$ is naturally a module over 
\[
\rH_T^*(\pt ;\Q_\ell) = \Sym_{\Ql} (\mf{t}^*) = \cO(\mf{t}) =: \sph_T
\]	
where $\mf{t} := (\Lie T)_{\Ql}$. When the acting torus $T$ is understood, we will simply abbreviate $\sph := \sph_T$. 

From the definitions, we have $\rH^{\BM, T}_*(\pt)  \cong \rH^{-*}_T(\pt)$, so there is a natural pairing 
\[
\rH^*_T(\pt) \otimes \rH^{\BM, T}_*(\pt) \rightarrow \rH^0_T(\pt) \cong \Q_\ell
\]
which realizes $\rH^{\BM, T}_*(\pt)$ as the graded $\Q_{\ell}$-dual of $\rH^*_T(\pt) \cong  \cO(\mf{t})$. This naturally identifies $\rH^{\BM, T}_*(\pt)$ with $\Omega_{\mf{t}}^{\wedge \dim \mf{t}}$, the $\cO(\mf{t})$-module of top-dimensional differential forms on $\mf{t}$. Then any choice of generator of $\omega_{\mrm{top}} \in \Omega_{\mf{t}}^{\wedge \dim \mf{t}}$ over $\cO(\mf{t})$ is equivalent to a choice of trivialization of $\rH^{\BM, T}_*(\pt)$ as an $\rH^*_T(\pt)$-module. 
\end{example}

\subsection{Equivariant formality}
Let $H$ be a group acting on a variety $X$. The graded $H$-equivariant Borel--Moore homology group $\rH^{\BM, H}_*(X)$ has a graded action of $\rH^*_H(\pt)$. 

There is an augmentation homomorphism $\rH^*_H(\pt) \rightarrow \Q_\ell$. This induces a map 
\begin{equation}\label{eq: de-eq BM}
\rH^{\BM, H}_{*}(X) \otimes_{\rH^*_H(\pt)} \Q_\ell \rightarrow \mBM_*(X).
\end{equation}
Recall that $X$ is \emph{$H$-equivariantly formal} if this map is an isomorphism. 

\begin{remark}[Equivariant classes in top homological degree]\label{remark: de-equivariant BM homology}
If $X$ is $H$-equivariantly formal and equidimensional of dimension $d$, then the map 
\[
\rH^{\BM, H}_{2d}(X) \rightarrow \rH^{\BM}_{2d}(X)
\]
induced by \eqref{eq: de-eq BM} is an isomorphism; in other words,``top-dimensional cycles are equipped with a canonical $H$-equivariant structure''.
\end{remark}

\begin{example}[Equivalued affine Springer fibers]\label{ex: ASF equivariantly formal}
According to \cite[Theorem 0.2]{GKM06}, the \emph{equivalued} affine Springer fibers admit a paving by affine spaces, and are therefore pure. Hence they are equivariantly formal with respect to any group action. This applies in particular to the affine Springer fibers $\rY_{\gamma}$ and $\rX_{\gamma}$ where $\gamma = ts$ for $s \in \mf{t}^*$, under the translation action of $T$. 
\end{example}

\subsection{The equivariant localization theorem}

Let $T$ be a torus acting on a variety $X$. We recall the localization theorem for the $T$-equivariant Borel--Moore homology of $X$.  Recall that $\rH^{\BM, T}_*(X)$ has a natural module structure over $\sph = \sph_T = \cO(\mf{t})$. Hence we may regard $\rH^{\BM,T}_*(X)$ as a quasicoherent sheaf on $\mf{t}$. 

Let $\iota \co X^T \rightarrow X$ denote the inclusion of the $T$-fixed points. The following is the famous \emph{Equivariant Localization Theorem} for torus actions, which goes back to work of Atiyah-Bott; a modern reference is \cite[Theorem A]{AKLPR}.

\begin{thm}[Equivariant Localization Theorem]
The kernel and cokernel of the map 
\[
\iota_* \co \rH^{\BM,T}_*(X^T) \rightarrow \rH^{\BM,T}_*(X)
\]
are supported (as quasicoherent sheaves) on the union of $\Lie K \subset \mf{t}$ as $K$ runs over proper stabilizer subtori $K \subset T$. In particular, $\iota_*$ is an isomorphism after tensoring over $\sph$ with $\Frac(\sph)$. 
\end{thm}

The theorem extends to ind-varieties in the obvious way. Note that the structure of $\rH^{\BM,T}_*(X^T) $ is simple: since $T$ acts trivially on $X^T$, we simply have $\rH^{\BM,T}_*(X^T) = \sph \otimes_{\Q_\ell} \mBM_*(X^T)$, which is free over $\sph$. 

We give some examples below, in which we maintain the notation of Part \ref{part: I}. 

\begin{example}\label{ex: fixed points X}
Let $\varepsilon_0 \in \A^1_{E}$ and let $X = \rX_{\gamma}^{\varepsilon = \varepsilon_0}$. Since $\gamma \in \chft[[t]]$, there is an action of $\chT$ via left translation on $X$ (cf. \S \ref{ssec: translation}), and we may identify $X^{\chT}$ with the discrete scheme $X_*(\chT)$, with $\lambda \in X_*(\chT)$ corresponding to $t^\lambda  L^+\chG \in \Gr_{\chG ,E}$. 
\end{example}


\begin{example}\label{ex: fixed points Y}
Let $\varepsilon_0 \in \A^1_{\F}$ and $X = \rY_\gamma^{\varepsilon = \varepsilon_0}$. Since $\gamma \in \chft[[t]]$, there is an action of $\chT$ via left translation on $X$. Then we may identify $X^{\chT}$ with the discrete scheme $\dWext$, with $w t^\lambda \in \dWext$ corresponding to $ wt^\lambda \chI \in \rY_{\gamma}^{\varepsilon = \varepsilon_0}$. 
\end{example}


\begin{defn}\label{eq:equivariant-localization}
Suppose $X$ is equivariantly formal for the action of a split torus $T$. Then $\rH^{\BM, T}_*(X) \cong \rH^{\BM}_*(X) \otimes_{\Ql} \sph_T$ is free over $\sph_T$, so we have an inclusion 
\begin{equation}\label{eq: GKM eBM 1} 
\Loc^T \co  \rH^{\BM, T}_*(X) \inj \rH^{\BM, T}_*(X^T) \otimes_{\sph_T} \Frac(\sph_T) 
\end{equation}
which is the dashed arrow in the diagram  
\[
\begin{tikzcd}
\rH^{\BM, T}_*(X^T) \ar[r, "i_*"] \ar[d, hook] & \rH^{\BM, T}_*(X) \ar[d, hook] \ar[dl, hook, dashed, "\Loc^T"']  \\
\rH^{\BM, T}_*(X^T) \otimes_{\sph_T} \Frac(\sph_T)  \ar[r, "i_* \otimes \Id", "\sim"'] & \rH^{\BM, T}_*(X) \otimes_{\sph_T} \Frac(\sph_T)
\end{tikzcd}
\]
Recall that we set $\mf{t} := (\Lie T)_{\Ql}$. Fix a generator $\omega_{\mrm{top}}$ of $\Omega_{\mf{t}}^{\wedge \dim \mf{t}}$, which induces a $\sph_T$-module trivialization $\rH^{\BM, T}_*(\pt) \cong \sph_T$ (cf. Example \ref{ex: torus}). Using this, we have
\begin{equation}\label{eq:GKM-eBM-1.5}
\rH^{\BM, T}_*(X^T) \otimes_{\sph_T} \Frac(\sph_T) \cong  \bigoplus_{x \in X^T} \Frac(\sph_T) [x].
\end{equation}
\end{defn}

\begin{defn}
For $\alpha \in \rH^{\BM, T}_*(X)$, the \emph{equivariant support} of $\alpha$ is the subset of $X^T$ at which $\Loc^T(\alpha)$ has non-zero component in the direct sum decomposition \eqref{eq:GKM-eBM-1.5}. 
\end{defn}

\subsection{Equivariant Euler classes}\label{ssec: equivariant Euler} Suppose $X$ is an ind-variety with an action of a split torus $T$, such that $X^T$ consists of isolated points in $X$. Let $x \in X^T$ be a \emph{smooth} point of $X$. Decompose the tangent space $\rT_xX$ as a representation of the torus $T$ into a sum of characters of $T$,
\[
\rT_x X \cong \bigoplus_i \lambda_i, \quad \lambda_i \in X^*(T).
\]
We may view $d\lambda_i$ as (linear) elements of $\cO(\mf{t})$. Then the \emph{equivariant Euler class} of $X$ at $x$ is defined to be $\prod_i d\lambda_i \in \cO(\mf{t})$, and we denote its inverse by 
\[
e_T(x,X) := \frac{1}{\prod_i d\lambda_i} \in \Frac(\sph_T).
\]

If $X$ is $T$-equivariantly formal of pure dimension $d$, then by Remark \ref{remark: de-equivariant BM homology} the fundamental class of $X$ admits a unique $T$-equivariant lift, which we denote $[X]_T \in \rH^{\BM,T}_{2d}(X)$. 

\begin{lemma}
Assume that we are in the situation of Definition \ref{eq:equivariant-localization}. Suppose $X$ has isolated $T$-fixed points. If $x \in X^T$ is a smooth point of $X$, then the image of $[X]_T$ in 
\[
\rH^{\BM,T}_*(X^T) \otimes_{\sph} \Frac(\sph_T) \cong \bigoplus_{x\in X^T} \Frac(\sph_T)[x]
\]
under the map $\Loc^T$ from \eqref{eq: GKM eBM 1} has coefficient of $[x]$ equal to $e_T(x,X)  \in \Frac(\sph)$. 
\end{lemma}

\begin{proof}This follows immediately from the Atiyah-Bott localization formula, for which a modern reference (in much more generality) is \cite[Theorem D]{AKLPR}. 
\end{proof}

\begin{example}\label{ex: flag variety}
Consider $X = \chG/\chB$ with the left translation action of $\chT$. Then the $\chT$-fixed points of $X$ may be identified with the discrete scheme $W$, where $w \in W$ corresponds to the fixed point $w\chB := \dot{w} \chB$ for any lift $\dot{w} \in N(\chT)$ of $w$ (note that the coset $w\chB$ does not depend on the choice of lift). The tangent space $\rT_{w\chB} X$ is $\chT$-equivariantly identified as (under our conventions on positive roots from \S \ref{sssec: positivity})
\[
\rT_{w\chB} X \cong \bigoplus_{\alpha \in \chPhi^+} \mf{g}_{w \alpha}.
\]
Letting $\beta := \prod_{\alpha \in \chPhi^+} d \alpha$, the equivariant fundamental class $[\chG/\chB]_{\chT}$ satisfies
\[
\Loc^{\chT} ([\chG/\chB]_{\chT}) = 
\sum_{w \in W}  \frac{\sgn(w)}{\beta} [w] \in \bigoplus_{w \in W} \Frac(\sph_{\chT})[w] = \rH^{\BM,T}_*(X^T) \otimes_{\sph} \Frac(\sph_T).
\]
\end{example}



\subsection{GKM description of equivariant Borel--Moore homology}\label{ssec: GKM} Consider the situation of Definition \ref{eq:equivariant-localization}, and suppose that $X$ furthermore satisfies the \emph{GKM conditions}\footnote{Named after Goresky, Kottwitz, and MacPherson who initiated these ideas in \cite[Theorem 7.5]{GKM04}.}: on any quasicompact subset of $X$, there are only finitely many $T$-fixed points and finitely many one-dimensional $T$-orbits. We review the so-called \emph{GKM description} of $\rH^{\BM, T}_*(X)$, which applies under the GKM conditions. Following the formulation of \cite[Proposition 4.3]{BL21}, it says that the image of \eqref{eq: GKM eBM 1} is the space of 
\[
\sum_{x \in X^T} f_x   [x]  \in  \bigoplus_{x\in X^T} \Frac(\sph_T) [x]
\]
satisfying the following conditions:
\begin{itemize}
\item The poles of $f_x$ are of order $\leq 1$ and contained in the union of the \emph{singular hyperplanes}, meaning hyperplanes of the form $\ker(d \chi)$ for a character $\chi \co T \rightarrow \G_m$ such that $X^{\ker \chi} \neq X^T$. 
\item For every \emph{singular} character $\chi$ (meaning that $X^{\ker \chi} \neq X^T$) and every connected component $Y \subset X^{\ker(\chi)}$, we have 
\[
\Res_{\ker(d \chi)} \left( \sum_{x \in X^T \cap Y} f_x  \ \topom  \right) = 0
\]
where we recall that $\topom$ is the fixed generator of $\Omega_{\mf{t}}^{\wedge \dim \mf{t}}$. Note that the collection of such $Y$ form the $1$-dimensional orbits of $T$ on $X$. 
\end{itemize}

\begin{example}[Borel--Moore homology of equivalued unramified affine Springer fibers]\label{ex: GKM description for Y}Let $\gamma = ts$ with $s \in \chft_{\C}$ regular. One can write down an analogous version of $\rY_\gamma$ over the complex numbers, which we denote $\rY_{\gamma, \C}$. Then $\rY_{\gamma, \C}$ satisfies the GKM conditions for the translation action of $\chT$, has the same fixed points as the positive characteristic version analyzed in Example \ref{ex: fixed points Y}, and the 1-dimensional $\chT$-orbits are calculated in \cite[\S 14]{GKM04}. 

The GKM description of $\rH^{\BM, \chT}_*(\rY_\gamma)$ is described explicitly in \cite[Corollary 4.8]{BL21}, and is manifestly independent of $\gamma$ (satisfying the hypotheses), and is used to identify $\rH^{\BM, \chT}_*(\rY_{\gamma, \C})$ for all $\gamma$ satisfying the hypotheses. By equivariant formality of the $\chT$-action, this also identifies $\rH^{\BM}_*(\rY_{\gamma,\C})$ for all $\gamma$ satisfying the hypotheses. 

To bootstrap these results to positive characteristic, observe that the same analysis shows that for $\gamma = ts$ such that $s \in \chft_{\F}$ is regular, then $\rY_{\gamma}$ also satisfies the GKM conditions, and with the same combinatorics of fixed points, singular characters, and 1-dimensional orbits, hence $\rH^{\BM, \chT}_*(\rY_\gamma)$ has the same GKM description. In particular, we use this description to identify $\rH^{\BM, \chT}_*(\rY_\gamma)$ for all $\gamma$ satisfying the hypotheses. By equivariant formality of the $\chT$-action on $\rY_\gamma$, this also identifies $\mBM_*(\rY_\gamma)$ for all $\gamma$ satisfying the hypotheses.

%

\end{example}

\begin{remark}[Deformed affine Springer fibers in the admissible region] 
Let $\gamma = ts$ with $s \in \chft_{\F}$ regular and let $\lambda \in X_*(\chT)^+$ be a dominant coweight such that $\gamma$ is $h_{\lambda+\rho}$-generic. Then it follows from \cite[Proposition 3.3.4]{LLLM22} that:
\begin{itemize}
\item $\rY_\gamma^{\varepsilon} (\leq \lambda)$ satisfies the GKM description for $\chT$ (and in particular is $\chT$-equivariantly formal). 
\item The $\chT$-fixed points of the family $\rY_\gamma^{\varepsilon} (\leq \lambda) \rightarrow \Spec \F[\varepsilon]$ are identified with the constant family $\Adm(\lambda) \times \Spec \F[\varepsilon]$, with $wt^\mu \in \Adm(\lambda)$ corresponding to $wt^\mu \chI \in \Fl_{\chG, \F[\varepsilon]}$. 
\item For each $\varepsilon_0 \in \Spec \F[\varepsilon]$, the singular characters are independent of $\varepsilon_0$ and the residue conditions they induce are independent of $\varepsilon_0$. 
\end{itemize}
Hence the GKM descriptions of $\mBMT_*(\rY_\gamma^{\varepsilon = \varepsilon_0} (\leq \lambda) )$ are independent of $\varepsilon_0 \in \Spec \F[\varepsilon]$. This gives compatible identifications of $\mBMT_*(\rY_\gamma^{\varepsilon = \varepsilon_0} (\leq \lambda))$ for all $\varepsilon_0 \in \Spec \F[\varepsilon]$, with respect to which the specialization maps in $\varepsilon$ are the identity map. In particular this gives \emph{some} bases of $\mBMT_*(\rY_\gamma^{\varepsilon = \varepsilon_0} (\leq \lambda) )$ with respect to which the specializations in $\varepsilon$ are the identity maps (but it is completely unclear whether this is the case for the basis comprised by the cycle classes of top-dimensional irreducible components.)
\end{remark}

\subsection{Actions on equivariant Borel--Moore homology}\label{ssec: equivariant actions} We define equivariant lifts of the affine Springer action and translation action on $\mBM_*(\rY_\gamma)$. 

\subsubsection{Monodromy-centralizer action}\label{sssec: equivariant monodromy-centralizer action}

 We will define an action of $(\dWext, \cdot)$ on $\mBMT_*(\rY_\gamma)$ using the GKM description from Example \ref{ex: GKM description for Y}, which we call the \emph{monodromy-centralizer action}. 

Fix a generator $\topom \in \Omega^{\wedge \dim \chft}_{\chft}$. Then $\Loc^{\chT}$ embeds $\mBMT_*(\rY_\gamma) \inj \bigoplus_{x \in \dWext} \Frac(\sph_{\chT})[x]$. There is a left action of $\wt{w} \in \dWext$ on $\bigoplus_{x \in \dWext} \Frac(\sph_{\chT})[x]$ via 
\[
\wt{w} \cdot \sum_{x \in X^{\chT}}  f_x  [x]   =  \sum_{x \in X^{\chT}} (\wt{w} f_x) [\wt{w} x]
\]
where $\wt{w} f_x$ refers to the natural $\wt{W}$-action on $\Frac(\sph_{\chT})$. 
One checks from the GKM description in Example \ref{ex: GKM description for Y} that this action preserves the subspace $\rH^{\BM, \chT}_*(\rY_\gamma)$. 

Finally, by Remark \ref{remark: de-equivariant BM homology}, it induces an action $(\dWext, \cdot)$ on the non-equivariant Borel--Moore homology $\topBM(\rY_\gamma)$, for which the action of $X_*(\chT) \subset \wt{W}$ agrees with the translation action of \S \ref{ssec: translation}. We refer to (all these variants of) the action $(\dWext, \cdot)$ as the \emph{monodromy-centralizer action}.

\begin{remark}[Explanation of terminology]
This action is also defined (for $Y_{\gamma, \C}$) in \cite[\S 5.2]{BL21}, up to converting between left and right actions. As explained in \cite[Remark 5.4]{BL21}, the action of the $X_*(\chT) \subset \dWext$ on $\mBMT_*(\rY_\gamma)$ is induced by the translation action by the centralizer of $\gamma$ (cf. \S \ref{ssec: translation}).

The action of $W \subset \dWext$ is subtler in that it is \emph{not} induced by an action of $W$ on $\rY_\gamma$, but could be thought of informally as the the monodromy action coming from a local system on $\chft^{\reg}$ whose fiber over $\gamma$ is $\rH^{\BM, \chT}_*(\rY_\gamma)$; indeed, Example \ref{ex: GKM description for Y} suggests that the $\rY_\gamma$ are ``$\chT$-equivariantly homotopy equivalent''. The cleanest way we know to make this precise is to just use the GKM description, as above.
\end{remark}

\subsubsection{Affine Springer action}\label{sssec: affine Springer equivariant} 

 We will define a \emph{right} action of $(\dWext, \bu)$ on $\mBMT_*(\rY_\gamma)$ that lifts the affine Springer action from \S \ref{ssec: affine Springer action}. The action of $\wt{w} \in \dWext$ on $\bigoplus_{x \in \dWext} \Frac(\sph_{\chT})[x]$ is given by  
\[
\wt{w} \bu \sum_{x \in X^{\chT}} f_x  [x]  = \sum_{x \in X^{\chT}} f_x [x \wt{w}].
\]
One checks from the GKM description in Example \ref{ex: GKM description for Y} that this action preserves $\mBMT_*(\rY_\gamma)$ embedded as a subspace of $\bigoplus_{x \in \dWext} \Frac(\sph_{\chT})[x]$ via $\Loc^{\chT}$. By equivariant formality, it induces an action on $\mBM_*(\rY_\gamma)$, which agrees with the affine Springer action constructed in \S \ref{ssec: affine Springer action}, according to \cite[\S 14.4]{GKM04}\footnote{For the version over $\C$, but this translates to the version over $\F$ as in Example \ref{ex: GKM description for Y}.}. We refer to (all these variants of) the action $(\dWext, \bu)$ as the \emph{affine Springer action}.

\begin{remark}\label{rem: actions commute}
It is immediate from the definitions that the actions $(\dWext, \cdot)$ and $(\dWext, \bu)$ commute with each other. 
\end{remark}

\section{Modular representation theory}\label{sec: modular representation theory}

Let $G$ be a reductive group over $\F_p$ with simply connected derived subgroup. In this section, which has no original results due to us, we recall some facts about the representation theory of $\mf{g} := \Lie G$ as well as of the \emph{Frobenius kernel} (of the relative Frobenius $F$)
\[
G_1 := \ker(F \co G_{\ol \F_p} \rightarrow G_{\ol \F_p}),
\]
and their graded variants. Roughly speaking, these will be used in Part 3 to ``approximate'' the representation theory of $G(\F_p)$, which is more directly related to Serre weights. We also review in \S \ref{ssec: BMR localization} the crucial theory of \emph{Bezrukavnikov-Mirkovic-Rumynin localization}, which provides the bridge between the categories of such representations and the categories of coherent sheaves that feature into the instance of mirror symmetry which is relevant for us. 

Finally, in \S \ref{ssec: rep m-c action} and \S \ref{ssec: braid action} we collect some natural symmetries of these categories. 

\subsection{Choice of torus}
Since reductive groups over finite fields are quasi-split, $G$ has a Borel subgroup. Let $B_0 \subset G$ be a Borel subgroup conforming to the conventions of \S \ref{sssec: positivity}, and $T_0 \subset B_0$ its Levi subgroup. Let $k := \ol \F_p$ and $B := (B_0)_k$, $T := (T_0)_{ k}$. Then we have an identification $X^*(T) \cong X_*(\chT)$, realizing $\chT$ as the Langlands dual group of $T$, for which $B \leftrightarrow \chB$ are corresponding standard Borel subgroups. 

Let $B^+$ be the opposite Borel to $B$. Recall that conventions from \S \ref{sssec: positivity} are that the positive roots $\Phi^+$ are the roots of $T$ on $\mf{b}^+$, equivalently on $\mf{g}/\mf{b}$, and that $\rho$ is the half-sum of the positive roots. 

The choice of $B$ induces an isomorphism of the flag variety $\cB$ with $G/B$. Recall that our positivity conventions are normalized so that the $G_k$-equivariant line bundles $\cO(\lambda)$ on $G_k/B$ associated to dominant weights $\lambda \in X^*(T)^+$ are semi-ample.

\subsection{Center of the universal enveloping algebra}\label{ssec: center of Ug} Fix an algebraic closure $k = \ol{\F}_p$. For a variety $X$ over $k$, we write $X^{(1)} := X \times_{k, \Frob_p}k$ for the \emph{(first) Frobenius twist} of $X$. We regard Lie algebras $\mf{g}$, etc. over $k$, and write $\Ug$ for the universal enveloping algebra of $\mf{g}$ over $k$. By $\Rep(\mf{g}) = \Rep(\Ug)$, etc. we mean the category of \emph{finitely generated} representations of $\Ug$. 

\subsubsection{The Harish-Chandra center} We recall some facts about the center of the universal enveloping algebra $\Ug$. Let $\mf{h} = \Lie A_{k}$ be the abstract Cartan Lie algebra of $\mf{g}$. An argument of Harish-Chandra produces a map $\Sym_k(\mf h)^{(W, \bu)} \inj Z(\Ug)$. Its image is called the \emph{Harish-Chandra center} $\zhc$. 

\subsubsection{The Frobenius center} In the analogous characteristic zero story, the Harish-Chandra center comprises the entirety of the center of $\Ug$. But in characteristic $p$, the center of $\Ug$ is much larger: there is also the so-called ``Frobenius center''
\[
\Sym_k \mf{g}^{(1)} \inj Z(\Ug)
\]
induced by the map sending $X \in  \mf{g}$ to $X^p-X^{[p]}$, whose image we denote $\zp$. Here $X \mapsto X^{[p]}$ is the $p$-operation on a Lie algebra in characteristic $p$, e.g., for $\mf{g} = \mf{gl}_n$ it sends a matrix to its $p$th power.

\subsubsection{The full center $Z(\Ug)$} Under our assumption that $p>h$, the center $Z(\Ug)$ is generated by the Harish-Chandra center and the Frobenius center, and has the more precise geometric description (cf. \cite{MR99}) 
\begin{equation}\label{eq: center Ug}
\Spec Z(\Ug) \cong \mf{h}^* \sslash W \times_{\mf{h}^{*(1)}\sslash W} \mf{g}^{*(1)}.
\end{equation}
Here:
\begin{itemize}
\item The map $\mf{g}^{*(1)} \rightarrow \mf{h}^{*(1)}\sslash W$ is the composition 
\[
\mf{g}^{*(1)} \rightarrow \mf{g}^{*(1)}\sslash G^{(1)} \xrightarrow{\mrm{Chevalley}} \mf{h}^{*(1)}\sslash W.
\]
\item The map $\mf{h}^*\sslash W \rightarrow \mf{h}^{*(1)}\sslash W$ is induced by the ``Artin-Schreier map'' $t \mapsto t^p - t^{[p]}$, where $t \mapsto t^{[p]}$ is the $p$-operation on $\mf{h}$.
\end{itemize}

\subsubsection{Representations with central conditions}
By the preceding discussion, a character of $Z(\Ug)$ is given by a compatible pair $(\lambda \in \mf{h}^*, \chi \in \mf{g}^{*(1)})$. For such a compatible pair $(\lambda, \chi)$, we define:
\[
\Ug^{\lambda} := (\Ug) \otimes_{\zhc} \lambda, \quad \Ug_{\chi} := (\Ug) \otimes_{\zp} \chi, \quad \Ug^{\lambda}_{\chi} := (\Ug) \otimes_{Z(\Ug)} (\lambda, \chi).
\]
We also make the following definitions.
\begin{itemize}
\item Define $\Rep^{\lambda}(\Ug)$ to be the full subcategory of $\Rep(\Ug)$ where $\zhc$ acts with \emph{generalized} eigenvalue $\lambda$. 
\item Define $\Rep_{\chi}(\Ug)$ to be the full subcategory of $\Rep(\Ug)$ where $\zp$ acts with \emph{generalized} eigenvalue $\chi$. 
\item Define $\Rep^\lambda(\Ug_\chi) := \Rep^\lambda(\Ug) \cap \Rep(\Ug_{\chi})$, and $\Rep^\lambda_{\chi}(\Ug):= \Rep^\lambda(\Ug) \cap \Rep_{\chi}(\Ug)$, etc.  
\end{itemize}

\subsection{The Frobenius kernel}\label{sssec: Frob kernel} Recall that the \emph{Frobenius kernel} $G_1$ is the kernel of the relative Frobenius morphism $F \co G_k \rightarrow G_k^{(1)}$. (Note that since $G$ is defined over $\F_p$, there is an isomorphism $G_k \cong G_k^{(1)}$.) Then $\cO(G_1)$ is a finite-dimensional commutative Hopf algebra over $k$, which is $k$-dual to $\Ug_0$ as a Hopf algebra. This induces a symmetric monoidal equivalence of categories \cite[I.9.6]{Jan03}
\begin{equation}\label{eq: frob kernel}
\Rep(G_1) \cong \Rep(\Ug_0),
\end{equation}
where we remind that ``$\Rep$'' means finitely generated representations in all contexts. We will freely use this equivalence of categories to transport statements between $\Rep(G_1)$ and $\Rep(\Ug_0)$.

 Recall that the simple representations of $G_k$ are in bijection with $X^*(T)^+$: for each $\lambda \in X^*(T)^+$ there is a unique simple representation of $G_k$ with highest weight $\lambda$, which we denote $L(\lambda)$.

Recall from \S \ref{sssec: mod rep theory} that a dominant weight $\lambda$ is \emph{$p$-restricted} if $0 \leq \tw{\lambda, \alpha^\vee}<p$ for all simple roots $\alpha$; the set of $p$-restricted weights is denoted $X_1^*(T)$. The simple representations of $G_1$ are exactly the restrictions of simple representations $L(\lambda)$ of $G_k$ with highest weight $\lambda \in X_1^*(T)$. 
 

\subsection{Graded representations}\label{ssec: graded representations}  Now we invoke the chosen torus $T \subset G_k$ to define $T$-graded representations. 

\subsubsection{Graded Lie algebra representations} If the adjoint action of $T$ on $\mf{g}$ fixes a central character $\chi$ of $\zp$, then we define $\Rep(\Ug_{\chi}, T)$ to be the Harish-Chandra category of representations $V$ of $\Ug_\chi$ together with a lift of $V|_{\mf{t}}$ to a representation of $T$. We define in an analogous way $\Rep_{\chi}(\Ug,T)$, $\Rep_{\chi}(\Ug^\lambda,T)$, etc.  

\begin{example}[$\chi=0$]\label{ex: graded simple} Take $\chi=0$. Then $\Rep(\Ug_0, T)$ is the category of graded $\Ug_0$-representations in the sense of \cite[\S D.5]{Jan04}. Concretely, $\Ug_{0}$ has a natural $X^*(T)$-grading where $X_\alpha \in \mf{g}$ has weight $\alpha$, for which $\Rep(\Ug_0, T)$ is the category of $X^*(T)$-graded representations of $\Ug_0$. 

The simple representations of $\Rep(\Ug_0,T)$ (which are the same as the simple representations in $\Rep_0(\Ug, T)$) are in bijection with $X^*(T)$, indexed by their highest weights, and we denote by $\wh{L}(\lambda)$ the simple representation of $\Rep(\Ug_0,T)$ with highest weight $\lambda \in X^*(T)$. 
\end{example}

\begin{example}[Graded baby Verma representations]\label{ex: graded baby Verma} Take $\chi = 0$. Let $\mf{b} =  \mf{n} \oplus \mf{t} \subset \mf{g}$ be the Lie algebra of $B$. It induces an isomorphism\footnote{This is (tautologically) compatible with our identifications $X^*(T) \cong X^*(A) \cong X_*(\chT)$.} $\mf{t} \xrightarrow{\sim} \mf{h}$, so that for $\lambda \in X^*(T)$ we may regard its derivative $d\lambda \in \mf{t}^*$ as an element of $\mf{h}^*$. Then $(d\lambda, \chi)$ forms a character of $Z(\Ug)$. We may also regard $d\lambda$ as a character of $\mf{b}$ by inflation. The \emph{graded baby Verma module} $\wh{Z}_{\mf{b}}(\lambda) \in \Rep(\Ug_0,T)$ is defined as
\[
\wh{Z}_{\mf{b}}(\lambda) := \Ug_0 \otimes_{\cU \mf{b}} d\lambda
\]
where $d\lambda$ has graded weight $\lambda$ and the universal enveloping algebras are equipped with their natural gradings. 
\end{example}


\subsubsection{Graded Frobenius kernel representations}\label{sssec: graded frob kernel}

Since the Frobenius kernel $G_1 \triangleleft G$ is normal, it generates along with $T$ a subgroup scheme $G_1 T < G_k$, isomorphic to the pushout of $G_1$ and $T$ along $T_1$. Its representation theory is studied (for example) in \cite[II.9]{Jan03}. The equivalence \eqref{eq: frob kernel} and its version for $T$ combine to give a monoidal equivalence of categories
\begin{equation}\label{eq: graded frob kernel}
\Rep(\Ug_0 , T) \cong \Rep(G_1 T).
\end{equation}

We denote by $\wh{L}_1(\lambda) \in \Rep(G_1 T)$ the simple representation of highest weight $\lambda$. To define baby Vermas of $G_1T$, we must make a choice of Borel subgroup. For compatibility with the literature on $G_1T$ representations that we will cite later, \textbf{we normalize the definition in the following way}, which is ``opposite'' to Example \ref{ex: graded baby Verma}: we denote by $\wh{Z}_1(\lambda)  \in \Rep(G_1 T)$ the graded baby Verma module of highest weight $\lambda$ for the Borel $B^+ = w_0B \subset G$, i.e. $\wh{Z}_1(\lambda)$ corresponds to 
\[
\wh{Z}_{w_0 \mf{b}}(\lambda)  := \Ug_0 \otimes_{\cU \mf{b}^+} d\lambda \in \Rep(\Ug_0 , T)
\]
under \eqref{eq: graded frob kernel}. This definition is made for compatibility with \cite{Jan03, GHS18}: our notation agrees with \cite[II.9.1 equation (2)]{Jan03}.

Recall that the \emph{linkage class} of $\lambda \in X^*(T)$ is $W_{\aff} \bup \lambda$, We denote by $\Rep^{\lambda}(G_1T)$ the Serre subcategory generated by simples $\wh{L}_1(\lambda')$ for $\lambda' \in W_{\aff} \bup \lambda$. The \emph{linkage principle} says that if $\lambda' \notin W_{\aff} \bup \lambda$, then $\Rep^{\lambda}(G_1 T)$ and $\Rep^{\lambda'}(G_1T)$ lie in different blocks of $\Rep(G_1 T)$.

\begin{example}[The extended principal block] 
The equivalence \eqref{eq: frob kernel} intertwines 
\[
\Rep^\lambda(\Ug_0 , T) \cong \bigoplus_{\lambda' \in \wt{W} \bup \lambda / W_{\aff} \bup \lambda} \Rep^{\lambda'}(G_1 T),
\]
which is the \emph{extended block of $\lambda$} in $\Rep(G_1 T)$. When $\lambda$ is regular for the action of $(\wt{W}, \bup)$, the indexing set for the sum is naturally a torsor for $\Omega = X^*(T)/Q$. This applies in particular when $\lambda = 0$, which is the case of most interest to us; then the RHS is called the \emph{extended principal block} of $\Rep(G_1T)$ and abbreviated $\Rep^{\emptyset}(G_1 T)$. 
\end{example}

\subsection{BMR Localization}\label{ssec: BMR localization}
Let $\cB$ be the flag variety of $G_k$. For nilpotent $\chi \in \mf{g}^*$, let $\cB_{\chi}$ be the inverse image of $\chi$ under the Grothendieck alteration\footnote{In Part I we viewed the Grothendieck alteration as a map $\wt{\mf{g}} \rightarrow \mf{g}$. It is related to this one via a $G$-equivariant isomorphism $\mf{g} \cong \mf{g}^*$, which exists under our characteristic hypotheses (cf. \cite[\S 3.1.2]{BMR08}).} $\wt{\mf{g}} \rightarrow \mf{g}^*$ (i.e., the \emph{Springer fiber} associated to $\chi$). Below we state in a special case a localization theorem for Lie algebras in positive characteristic due to Bezrukavnikov--Mirkovic. Recall that the Springer resolution $\wt{\cN} := T^* \cB$ is the cotangent bundle of the flag variety of $G_k$, and the superscript $(-)^{(1)}$ denotes Frobenius twist over $k$. 

\begin{thm}[{\cite[Theorem 1.6.7]{BM13}}]\label{thm: BM} Let $\lambda \in X^*(T) \cap p A_0$ and $\chi \in \mf{g}^*$ be nilpotent and fixed by $T$ under the coadjoint action. Then there is an equivalence 
\[
\gamma^{\lambda}_{\chi}\co D^b (\Rep_{\chi}(\Ug^\lambda,T)) \xrightarrow{\sim} D^b (\Coh_{\cB_{\chi}^{(1)}}^{T^{(1)}} (\wt{\cN}^{(1)})),
\]
where $\Coh_{\cB_{\chi}^{(1)}}^{T^{(1)}}( \wt{\cN}^{(1)})$ denotes the category of ${T^{(1)}}$-equivariant coherent sheaves on $\wt{\cN}^{(1)}$ with set-theoretic support on $\cB_{\chi}^{(1)}$ (with $T^{(1)}$-action on $\wt{\cN}^{(1)}  = (T^*\cB)^{(1)}$ induced by the $T$-translation action on $\cB$). 

Moreover, the equivalence $\gamma^{\lambda}_{\chi}$ is t-exact for the usual t-structure on $D^b (\Rep_{\chi}(\Ug^\lambda,T))$ and the exotic t-structure on $D^b (\Coh_{\cB_{\chi}^{(1)}}^{T^{(1)}} (\wt{\cN}^{(1)}))$. 
\end{thm}

Theorem \ref{thm: BM} is an equivariant enhancement of the localization theorem from {\cite[Theorem 5.3.1]{BMR08}}, which gives an equivalence 
\begin{equation}\label{eq: BMR}
 D^b (\Rep_\chi(\Ug^\lambda)) \xrightarrow{\sim} D^b (\Coh_{\cB_{\chi}^{(1)}} (\wt{\cN}^{(1)})).
\end{equation}

\begin{example}\label{ex: principal block BMR}
We will only apply Theorem \ref{thm: BM} for $\chi = 0$ and $\lambda = 0$. In this case $\cB_{\chi} = \cB$ is the full flag variety, and Theorem \ref{thm: BM} supplies an equivalence
\begin{equation}\label{eq: BM 0 0}
\gamma^0_0 \co D^b (\Rep_0(\Ug^0,T))  \xrightarrow{\sim} D^b (\Coh_{\cB^{(1)}}^{T^{(1)}}(\wt{\cN}^{(1)})).
\end{equation}
In this case, the non-equivariant version \eqref{eq: BMR} goes as follows. Since $\lambda = 0$, the characteristic $p$ analogue of Beilinson-Bernstein localization identifies $D^b(\Rep(\Ug^0))$ with the derived category of coherent $\cD$-modules on $\cB$. The condition that $\chi = 0$ translates into the condition that the $p$-curvature (of the $\cD$-modules obtained via Beilinson-Bernstein localization) is nilpotent. 
The ring of differential operators on $\cB$ pushes forward to an Azumaya algebra $\Fr_* \cD_{\cB}$ on $T^*\cB^{(1)}$, and $\cD$-modules on $\cB$ with nilpotent $p$-curvature push forward exactly to $\Fr_* \cD_{\cB}$-modules on $(T^*\cB)^{(1)} \cong \cN^{(1)}$ with set-theoretic support on $\cB^{(1)}$. Finally, the Azumaya algebra $\Fr_* \cD_{\cB}$ splits canonically on the formal neighborhood of $\cB^{(1)}$ by Cartier descent, giving a Morita equivalence between coherent $\Fr_* \cD_{\cB}$-modules supported on $\cB^{(1)}$ and  $\Coh_{\cB^{(1)}}(\wt{\cN}^{(1)})$. 

The graded version is bootstrapped from the non-graded one by relating $D^b (\Rep_0(\Ug^0,T))$ to $T$-equivariant $D$-modules on $\cB$, and then tracking the equivariant structure through the Morita equivalence. 
\end{example}




\begin{example}[The trivial representation]\label{ex: graded trivial BMR} Take $\chi = 0$ and $\lambda = 0$. The equivalence \eqref{eq: BM 0 0} sends the trivial representation $\wh{L}(0) \in \Rep_0(\Ug^0,T)$ to $i_* \cO_{\cB^{(1)}}$ where $\cO_{\cB^{(1)}}$ is equipped with its natural $T^{(1)}$-equivariant structure induced by the translation action of $T^{(1)}$ on $\cB^{(1)}$.
\end{example}

\begin{example}[Graded baby Verma modules]\label{ex: graded baby verma BMR}
The choice of Borel $B \subset G$ induces an isomorphism of the flag variety $\cB$ with $G/B$. Let $\mf{b} =  \mf{n} \oplus \mf{t} \subset \mf{g}$ be the Lie algebra of $B$. Then from \cite[\S 3.1.4]{BMR08} we see that \eqref{eq: BM 0 0} sends $\wh{Z}_{\mf{b}}(2\rho)$ to $\delta_{\mf{b}}$, the $T$-equivariant skyscraper sheaf on $\cB^{(1)}$ supported at $\mf{b}$ regarded as a point of $\cB^{(1)}$, with its natural $T$-equivariant structure induced by the translation action of $T^{(1)}$ on $\cB^{(1)}$.

\end{example}

\subsection{Monodromy-centralizer action}\label{ssec: rep m-c action} 
We will define an action of $\dWext$ on $\rK(\Rep_0(\Ug^0, T))$ which is parallel to the monodromy-centralizer action of \S \ref{ssec: equivariant actions} in a direct sense. 

\subsubsection{Graded Lie algebra and Frobenius kernel} Let $N(T) < G$ be the normalizer of $T$. Then the adjoint action of $N(T)$ on $\mf{g}$ induces an action of $N(T)$ on $\Rep_0(\Ug^0, T)$. Under the equivalence \eqref{eq: frob kernel}, this corresponds to the action of $N(T)$ on $\Rep^\emptyset(G_1 T)$ induced by conjugation action of $N(T)$ on $G_1 T <G$. 

Now note that at the level of Grothendieck groups, the subgroup $T \triangleleft N(T)$ acts trivially on $\rK(\Rep^\emptyset(G_1 T))$ since its action is inner. Hence we obtain an action of $N(T)/T \cong W$ on $\rK(\Rep^\emptyset(G_1 T)) \cong \rK(\Rep_0(\Ug^0, T))$. 

In addition, via the obvious quotient map $G_1T \surj T^{(1)}$ we have an action of $\Rep(T^{(1)})$ on $\Rep^\emptyset(G_1 T)$ by inflation and tensoring. On $\Rep_0(\Ug^0, T)$ this corresponds to the action of changing the grading: there is an identification $T^{(1)} \cong T$ since $T$ is defined over $\F_p$, inducing $X^*(T^{(1)} ) \cong X^*(T)$. With respect to this identification, the relative Frobenius $\Fr_T \co T \rightarrow T^{(1)}$ induces the second map in the sequence
\[
X^*(T) \cong X^*(T^{(1)}) \xrightarrow{\Fr_T^*} X^*(T) 
\]
whose composite is multiplication by $p$. In particular, $\lambda \in X^*(T) = X^*(T^{(1)})$ acts on $\Rep_0(\Ug^0, T)$ by translating the grading by $p\lambda$. 

At the level of Grothendieck groups, this gives an action of $X^*(T)$ on $\rK(\Rep^\emptyset(G_1 T)) \cong \rK(\Rep_0(\Ug^0, T))$. Together with the earlier $W$-action, these define an action of $\dWext \cong X^*(T) \rtimes W$ on $\rK(\Rep_0(\Ug^0, T))$, which we call the \emph{monodromy-centralizer action}. (This terminology is still unexplained, but is made to be parallel to \S \ref{sssec: equivariant monodromy-centralizer action}.) We denote the action of $\wt{w} \in \dWext$ by $\wt{w} \cdot_p (-)$. 

\begin{example}[Graded simple modules]\label{ex: monodromy action on simples}
Recall that for $\lambda \in X_1^*(T)$ (cf. \S \ref{sssec: Frob kernel} for the notation), $\wh{L}_1(\lambda)$ is the simple representation in $\Rep^{\emptyset}(G_1T)$ with highest weight $\lambda$. For $\mu \in X^*(T)$, we have 
\begin{equation}\label{eq: translate act on simple}
t^\mu \cdot_p [\wh{L}_1(\lambda)] = [\wh{L}_1(\lambda + p \mu)] \in \rK (\Rep^\emptyset(G_1T)).
\end{equation}
For $w \in W$, we have 
\[
w\cdot_p [\wh{L}_1(\lambda)]  = [\wh{L}_1(\lambda)]  \in \rK (\Rep^\emptyset(G_1T))
\]
because for such $\lambda$, the simple representation $\wh{L}_1(\lambda)$ extends to a simple representation of $G$, where the conjugation action of $N(T)$ is inner, hence trivial on $\rK(\Rep(G))$. Together with \eqref{eq: translate act on simple}, this determines the $(\dWext, \cdot)$-action on all simples, and shows that it permutes the classes of simples.
\end{example}

\begin{example}[Graded baby Vermas]
Recall that $\wh{Z}_{\mf{b}}(\lambda)$ is the baby Verma representation in $\Rep_0(\Ug^0, T)$ with highest weight $\lambda$ and Borel $\mf{b}$. Note that the central character conditions force $\lambda \in pX^*(T)$. For any $\mu \in X^*(T)$, we have 
\begin{equation}\label{eq: translate act on verma}
\mu \cdot_p [\wh{Z}_{\mf{b}}(\lambda)] = [\wh{Z}_{\mf{b}}(\lambda + p \mu)] \in \rK(\Rep_0(\Ug^0, T)).
\end{equation}
For any $w \in W$, we have according to \cite[\S 9.3]{Jan03}
\begin{equation}\label{eq: monodromy act on verma}
w\cdot_p  [\wh{Z}_{\mf{b}}(\lambda)]  = [\wh{Z}_{w \mf{b}}(w \lambda)]  \in \rK(\Rep_0(\Ug^0, T))
\end{equation}
where $w \mf{b}$ is the translate of $\mf{b}$ by any $\dot{w} \in N(T)$ lifting $w$. Together with \eqref{eq: translate act on verma}, this determines the $(\dWext, \cdot)$-action on all baby Vermas, and shows that it permutes the classes of baby Vermas. 
\end{example}

\subsubsection{Coherent sheaves}

There is an obvious action of $\Rep(T^{(1)})$ on $\Coh^{T^{(1)}}_{\cB^{(1)}}(\wt{\cN}^{(1)})$ by tensoring with equivariant representations. At the level of Grothendieck groups, this induces an action of $X^*(T^{(1)}) \cong X^*(T)$ on $\rK (\Coh^{T^{(1)}}_{\cB^{(1)}}(\wt{\cN}^{(1)})) \cong  \rK (\Coh^{T^{(1)}}(\cB^{(1)}))$.

 Also, writing $\Coh^{T^{(1)}}_{\cB^{(1)}}(\wt{\cN}^{(1)}) = \Coh_{\cB^{(1)}}(T^{(1)} \bs \wt{\cN}^{(1)})$, we see that there is an action of $N(T^{(1)})$ by left multiplication on the quotient $T^{(1)} \bs \wt{\cN}^{(1)}$. At the level of Grothendieck groups, it factors over $T^{(1)}$, inducing an action of $W \cong N(T^{(1)})/T^{(1)}$ on $\rK (\Coh(T^{(1)} \bs \cB^{(1)}))$ by left translation. 

Together, these combine into an action of $\dWext \cong X^*(T^{(1)}) \rtimes W$ on $\rK (\Coh^{T^{(1)}}_{\cB^{(1)}}(\wt{\cN}^{(1)}))$, which we denote $(\dWext, \cdot)$. 

 \begin{lemma}\label{lem: BMR monodromy-centralizer}
At the level of Grothendieck groups, the equivalence \eqref{eq: BM 0 0} intertwines the action of $(\dWext, \cdot_p)$ on $\rK(\Rep_0(\Ug^0, T))$ with the action of $(\dWext, \cdot)$ on $\rK (\Coh^{T^{(1)}}_{\cB^{(1)}}(\wt{\cN}^{(1)}))$.
 \end{lemma}

\begin{proof}The action of $X^*(T^{(1)})  \cong X^*(T) \subset \dWext$ on both sides can be described as tensoring with representations of $T^{(1)}$, and it is clear from the construction that the equivalence \eqref{eq: BM 0 0} intertwines these operations. 
 
 It follows from chasing through the construction that \eqref{eq: BM 0 0} also intertwines the actions of $W$. Let us sketch why: the point is that the map $\Ug^0 \rightarrow \Gamma(\cB, \cD)$ is $G_k$-equivariant for the adjoint action of $G_k$ on $\Ug^0$ and the left translation action of $G_k$ on $\cB$. This in turn follows from the fact that the action map 
\[
G_k \times_k \cB \rightarrow \cB
\]
is equivariant for the conjugation action of $G_k$ on itself and the left translation action on $\cB$. Hence a fortiori the map $\Ug^0 \rightarrow \Gamma(\cB, \cD)$ is $H$-equivariant for any subgroup $H<G_k$. Apply this to $H = N(T)$ and the result follows, using that the adjoint action of $N(T)$ induced the $W$-action on $\rK(\Rep_0(\Ug^0, T))$, while the translation action of $N(T)$ induced the $W$-action on $\rK (\Coh^{T^{(1)}}_{\cB^{(1)}}(\wt{\cN}^{(1)}))$.

\end{proof}

\subsection{Braid action}\label{ssec: braid action}

Recall that the \emph{extended affine Braid group} $\wt{\BB}$ has generators $T_{\wt{w}}$ for $\wt{w} \in \dWext$, and relations $T_{\wt{w}} T_{\wt{w}'} = T_{\wt{w}\wt{w}'}$ if $\len(\wt{w}\wt{w}')  = \len(\wt{w}) + \len(\wt{w}')$. We have $\wt{\BB} \cong \BB_{\aff} \rtimes \Omega$, where $\BB_{\aff}$ is the affine Braid group associated to the Coxeter group $W_{\aff}$.

We will define an action of $\dWext$ on $\rK(\Rep_0(\Ug^0, T))$ which is parallel to the affine Springer action of \S \ref{ssec: equivariant actions} in a direct sense.

\subsubsection{Graded Lie algebra}\label{sssec: braid action on lie}

Let $(\chi, \lambda)$ be as in Theorem \ref{thm: BM}. In \cite[\S 2]{BMR06}, Bezrukavnikov-Mirkovic-Rumynin constructed an action of $\wt{\BB}$ on $D^b (\Rep^\lambda(\Ug,T))$\footnote{Strictly speaking, they did not consider the graded version, but the same construction goes through essentially verbatim.}. We will describe the action of $T_{\wt{w}}$ for $\wt{w} \in \dWext$. First, recall that for any $\mu, \nu \in X^*(T)$ there is a \emph{translation functor}
\[
T^\nu_{\mu} \co \Rep^\mu(\Ug,T) \rightarrow \Rep^{\nu}(\Ug,T).
\]
For $\omega \in \Omega$, the action of $T_\omega \in \wt{\BB}$ on $D^b(\Rep^\mu(\Ug, T))$ is via $T_\mu^{\omega \mu}$; in particular, it is exact (since it does not pass through a wall).

If $\mu$ lies in the interior of the alcove and $\nu$ lies on a codimension-1 face, then we define the \emph{reflection functor} 
\[
R_{\mu|\nu} := T_{\nu}^\mu \circ  T_{\mu}^\nu \co \Rep^\mu(\Ug,T) \rightarrow \Rep^{\mu}(\Ug,T).
\]
The functors $R_{\mu | \nu}$ are naturally isomorphic for different choices of $\nu$ in the interior of the codimension-1 face, so we fix a choice and denote the reflection functor (also known as \emph{wall-crossing functor}) by $R_s$ where $s \in \wt{W}$ is the reflection through the codimension-1 face. There is a distinguished triangle of functors on $D^b(\Rep^\mu(\Ug, T))$,
\[
\Id \rightarrow R_{s} \rightarrow \II^*_{s}.
\]
Then the action of $T_{s} \in \wt{\BB}$ on $D^b(\Rep^{\mu}(\Ug,T))$ is via $\II^*_{s}$. 

These functors are compatible with the Frobenius-center, hence the same formulas induce a $\wt{\BB}$-action on $D^b(\Rep_{\chi}(\Ug, T))$.

\subsubsection{Geometric braid action on coherent sheaves} There is also a $\wt{\BB}$-action on $D^b (\Coh(\wt{\mf{g}}^{(1)}))$ constructed in \cite[Theorem 1.4.1]{Ric08} which preserves each $D^b (\Coh_{\cB_{\chi}^{(1)}}(\wt{\mf{g}}^{(1)}))$. Then Riche shows in \cite[Theorem 5.4.1]{Ric08} that the equivalence of Theorem \ref{thm: BM} intertwines the two $\wt{\BB}$-actions. For the graded case, it is also true that the equivalence Theorem \ref{thm: BM} intertwines the two $\wt{\BB}$-actions, by the same argument as for Riche's result in the non-graded case. 

\subsubsection{Converting to right actions} Henceforth we convert the $\wt{\BB}$-action to a \emph{right} action by the anti-involution $\wt{\BB} \xrightarrow{\sim} \wt{\BB}^{\mrm{opp}}$ given by the inverse map. This is for compatibility with how we normalized the affine Springer action on $\mBM_*(\rY_\gamma)$ to be a right action. 

\subsubsection{Steinberg action}

Our reductive group $G$ extends canonically to an unramified reductive group over $\Z_p$. Let $\St_G := \wt{\cN} \times_{\cN} \wt{\cN}$. This is defined in characteristic $p$ by our convention, but there is an analogous construction over $\C$ that we denote $\St_{G, \C}$. Recall that the Kazhdan-Lusztig isomorphism (conjectured by Deligne-Langlands) identifies $\rK (\Coh^{G \times \G_m}(\St_{G,\C}))$, as an algebra under convolution, with the affine Hecke algebra of $\chG$. This induces an isomorphism $\rK (\Coh^{G}(\St_{G,\C})) \cong  \Z[\dWext]$. A small argument, which we presently give, shows that the same holds in characteristic $p$.

Recall that there is a specialization map in $K$-theory for a Noetherian scheme $X$ flat over a DVR (cf. \cite[\S 7.1.3]{BMR08}). When the DVR is $\Z_p$, we denote it as 
\begin{equation}\label{eq: spc Coh}
\fsp_{p \rightarrow 0} \co  \rK(\Coh(X_{\Q_p}))  \rightarrow \rK(\Coh(X_{\F_p})).
\end{equation}
It is defined on a coherent sheaf $\cF/X_{\Q_p}$ by choosing a $\Z_p$-lattice and then taking the (derived) tensor product with $\F_p$. If a split reductive group scheme $H/\Z_p$ acts on $X/\Z_p$, then there is also an equivariant version. 
 Let $\St_{G, \Z_p} = \wt{\cN}_{\Z_p} \times_{\cN_{\Z_p}} \wt{\cN}_{\Z_p}$ and define $\St_{G, \Q_p}$, etc. analogously. 

\begin{lemma}\label{lem: steinberg fsp} 
The map 
\begin{equation}\label{eq: steinberg fsp}
\fsp_{p \rightarrow 0} \co  \rK(\Coh^{G \times \G_m}(\St_{G, \Q_p}))  \rightarrow	 \rK(\Coh^{G \times \G_m}(\St_{G, \F_p}))
\end{equation}
induced by $\St_{G, \Z_p}$ is an algebra isomorphism (both sides being equipped with the convolution product). 
\end{lemma}

\begin{proof}
Over $\Z_p$, the identification 
\[
\wt{\cN}_{\Z_p} \times_{\cN_{\Z_p}} \wt{\cN}_{\Z_p}  \cong T^* (\cB_{\Z_p}) \times_{\cN_{\Z_p}} T^*(\cB_{\Z_p})
\]
realizes the Steinberg variety $\St_G$ as the conormal space (relative to $\Z_p$) to the Bruhat stratification of $\cB_{\Z_p} \times \cB_{\Z_p}$ by diagonal $G_{\Z_p}$-orbits. In particular, $\St_{G, \Z_p}$ admits a stratification into affine spaces over $\Z_p$. This equips both sides of \eqref{eq: steinberg fsp} with filtrations such that each graded is the specialization map for an affine space, which is an isomorphism  (cf. the Cellular Fibration Lemma \cite[\S 5.5]{CG10}). Therefore, $\fsp_{p \rightarrow 0}$ is an isomorphism of groups. 

The algebra structure on either side is given by convolution, which is defined because over any field $\wt{\cN}$ is smooth and the projection $\wt{\cN} \rightarrow \cN$ is proper. Since $\wt{\cN}_{\Z_p}$ is smooth over $\Z_p$ and the projection $\wt{\cN}_{\Z_p} \rightarrow \cN_{\Z_p}$ is proper, the operations constituting convolution are compatible with $\fsp_{p \rightarrow 0}$. Therefore the map \eqref{eq: steinberg fsp} is compatible with convolution. 
\end{proof}

Note that the flat base change from $\Q_p$ to $\C$ induces an isomorphism 
\[
\rK(\Coh^{G \times \G_m}(\St_{G, \Q_p})) \xrightarrow{\sim} \rK(\Coh^{G \times \G_m}(\St_{G, \C})),
\]
again because of the stratification into affine spaces. We deduce that $\rK(\Coh^{G \times \G_m}(\St_{G, \F_p}))$ is also isomorphic to the affine Hecke algebra for $\chG$. 

We resume working over a field $k$ of characteristic $p$, and omit the subscripts indicating the base field. For any nilpotent $\chi \in \mf{g}^*$, there is a natural (right) \emph{Steinberg action} of $\rK (\Coh^{G \times \G_m}(\St_G))$ on $\rK (\Coh^T(\cB_{\chi}))$, the $K$-theory of the corresponding Springer fiber, by convolution on the right. On the other hand, the $T$-equivariant version of the construction of \cite[Theorem 1.4.1]{Ric08} induces a right action of $\wt{\BB}$ on $D^b (\Coh^T_{\cB_\chi}(\wt{\cN}))$. At the level of Grothendieck groups, according to \cite[Theorem 1.3.2]{BM13}, if $\chi \in \mf{g}^*$ is nilpotent then the action of $\wt{\BB}$ on $\rK (\Coh^T_{\cB_\chi}(\wt{\cN})) \cong \rK (\Coh^T(\cB_\chi))$ factors through the Steinberg action of $\dWext$. 

Repeating this discussion with appropriate Frobenius twists, we obtain a Steinberg action of $\dWext$ on $\rK (\Coh^{T^{(1)}}_{\cB^{(1)}}(\wt{\cN}^{(1)}))$, which we denote $(\dWext, \bu)$. It follows from (the graded version of) Riche's Theorem that: 

\begin{lemma}\label{lem: K-groups BMR} The action of $\wt{\BB}$ on $\rK(\Rep_0(\Ug^0, T))$ induced by \S \ref{sssec: braid action on lie} factors through an action of $\dWext$ that we denote $(\dWext, \bup)$.

Furthermore, at the level of Grothendieck groups the equivalence \eqref{eq: BM 0 0} intertwines the action of $(\dWext, \bup)$ on $\rK(\Rep_0(\Ug^0, T))$ with the action of $(\dWext, \bu)$ on $\rK (\Coh^{T^{(1)}}_{\cB^{(1)}}(\wt{\cN}^{(1)}))$. 
\end{lemma}

\subsubsection{Comparison to characteristic zero} We want to compare the $K$-group and actions in Lemma \ref{lem: K-groups BMR} to the analogous situation in characteristic zero. The flat family $\cB_{\Z_p}/\Z_p$ induces a specialization map in equivariant $K$-theory, 
\begin{equation}\label{eq: spc flag}
\fsp_{p \rightarrow 0} \co  \rK (\Coh^{T_{\Q_p}}(\cB_{\Q_p})) \rightarrow \rK(\Coh^{T_{\F_p}}(\cB_{\F_p})).
\end{equation}

\begin{lemma}\label{lem: char p to 0}
The map \eqref{eq: spc flag} is an isomorphism. It has the following properties: 
\begin{enumerate}
\item It is equivariant for the action of $(\wt{W}, \cdot)$ and also for the Steinberg action of $(\wt{W}, \bu)$. 
\item It takes $[\cO_{\cB_{\Q_p}}]$ to $[\cO_{\cB_{\F_p}}]$. 
\item For each Borel subgroup $B<G$ defined over $\Z_p$, it takes the skyscraper class $[\delta_{B_{\Q_p}}] \in  \rK(\Coh^{T_{\Q_p}}(\cB_{\Q_p}))$ to the skyscraper class $[\delta_{B_{\F_p}}]  \in \rK(\Coh^{T_{\F_p}}(\cB_{\F_p}))$.
\end{enumerate}
\end{lemma}

\begin{proof}
The Bruhat stratification decomposes $\cB_{\Z_p}$ into affine spaces over $\Z_p$. Therefore \eqref{eq: spc flag} is an isomorphism by the same cellular fibration argument as in the proof of Lemma \ref{lem: steinberg fsp}. 

The compatibility with the action $(\wt{W}, \cdot)$ is evident from the definitions. The compatibility with the action $(\wt{W}, \bu)$ follows from that the fact that the action maps come from tensor products on smooth ambient spaces or pushforward along proper morphisms defined over $\Z_p$, like in the proof of Lemma \ref{lem: steinberg fsp}, which are therefore compatible with specialization. 

The computation of $\fsp$ on the structure sheaf and skyscrapers is evident from the definition.

\end{proof}

We note again that the flat base change from $\cB_{\Q_p}$ to $\cB_{\C}$ induces an isomorphism of K-groups, by the Bruhat stratification into affine spaces.

\begin{remark}
The specialization map $\fsp_{p \rightarrow 0}$ and base change from $\Q_p$ to $\C$ are analyzed for the $K$-theory of more general Springer fibers (without the $T$-equivariance) in \cite[Proposition 7.1.7]{BMR08}. 
\end{remark}

\subsection{Upshot} Summarizing, we have:

\begin{thm}\label{thm:mod-rep-theory-summary}
There is an isomorphism 
\[
\rK(\Rep_0(\Ug^0, T)) \rightarrow \rK(\Coh^{T_{\C}}(\cB_{\C}))
\]
which has the following properties:
\begin{enumerate}
\item It intertwines the left action $(\wt{W}, \cdot_p)$ on the LHS with the left action $(\wt{W}, \cdot)$ on the RHS. 
\item It intertwines the right action $(\wt{W}, \bup)$ on the LHS with the right action $(\wt{W}, \bu)$ on the RHS. 
\item It sends $[\wh{L}(0)] \in 
\rK(\Rep_0(\Ug^0, T))$ to $[\cO_{\cB_{\C}}] \in \rK(\Coh^{T_{\C}}(\cB_{\C}))$. 
\item It sends $[\wh{Z}_{\mf{b}}(2\rho)] \in 
\rK(\Rep_0(\Ug^0, T))$ to $[\delta_{B_{\C}}]\in \rK(\Coh^{T_{\C}}(\cB_{\C}))$. 
\end{enumerate}
\end{thm}

\begin{proof}
Combine Lemma \ref{lem: char p to 0}, Lemma \ref{lem: K-groups BMR}, and Lemma \ref{lem: BMR monodromy-centralizer}.
\end{proof}


\section{Shadows of mirror symmetry}\label{sec: K-groups}

Our goal is to compare the representation-theoretic information of $G$ measured in $\rK(\Rep_0(\Ug^0, T))$ with the geometric information of $\chG$ measured in $\topCh(\rY_\gamma)$. In the preceding section we explained a more geometrical incarnation of $\rK(\Rep_0(\Ug^0, T))$ in terms of coherent sheaves on the flag variety for $G$. In this section, we will connect this with $\topCh(\rY_\gamma)$. This connection may be viewed as some manifestation of homological mirror symmetry, which relates Lagrangians on a symplectic manifold and coherent sheaves on a mirror variety. \emph{This provides in particular the passage from $G$ to its dual group $\chG$.}

\begin{prop}[Bezrukavnikov--Boixeda Alvarez--McBreen--Yun \cite{BBMY2}]\label{prop: HMS}Let $s \in \chft$ be regular semisimple and $\gamma = ts \in \chft[[t]]$.  There is a map 
\begin{equation}\label{eq: HMS}
 \rK (\Coh^{T_{\C}}_{\cB_{\C}} (\wt{\cN}_{\C})) \rightarrow \topCh(\rY_{\gamma}),
\end{equation}
with the following properties: 
\begin{enumerate}
\item It intertwines the actions $(\wt{W}, \cdot)$  defined in \S \ref{ssec: braid action} for the LHS and \S \ref{ssec: equivariant actions} for the RHS.
\item It intertwines the actions $(\dWext, \bu)$ defined in \S \ref{ssec: braid action} for the LHS and in \S \ref{ssec: equivariant actions} for the RHS. 
\item It sends $[\cO_{\cB_{\C}}] \in \rK(\Coh^{T_{\C}}_{\cB_{\C}} (\wt{\cN}_{\C}))$ to the fundamental class of the unique (top-dimensional) irreducible component of $\rY_{\gamma}$ which is the pre-image of $[t^0] \in \Gr_{\chG, \F}$ under the projection map $\Fl_{\chG, \F} \rightarrow \Gr_{\chG, \F}$. 
\end{enumerate}
\end{prop}


\begin{remark}The left side of \eqref{eq: HMS} does not depend on $\gamma$ while the right side seems to depend on it. Recall however that in Example \ref{ex: GKM description for Y} we saw that $\topCh(\rY_{\gamma})$ is ``independent of $\gamma$'' in a suitable sense. With this said, we may choose $\gamma$ to come from a regular semisimple element of $\mf{t}_{\Z_p}[[t]]$, and by Example \ref{ex: GKM description for Y} again it is equivalent to prove the analogous statement with $\rY_{\gamma, \C}$, the complex version of the affine Springer fiber, in place of $\rY_{\gamma}$. 
\end{remark}

Proposition \ref{prop: HMS} is a consequence of constructions in \cite{BBMY2}. For the sake of being self-contained, we sketch the relevant constructions below, while emphasizing that they are entirely due to \cite{BBMY2}; \emph{we do not claim any original results in this section}.

\begin{remark}\label{rem: BBMY}
To explain the title of this section: the map \eqref{eq: HMS} is the shadow of an instance of \emph{homological mirror symmetry} which predicts in some form that given a symplectic manifold $M$ with a Lagrangian skeleton $L$, there should be an equivalence between a ``Fukaya category''\footnote{This incarnation of the Fukaya category is in the spirit of \cite{NS20} rather than Kontsevich's original formulation.} of microlocal sheaves $\muSh_L(M)$ with supports on $L$, and the derived category of coherent sheaves on a mirror algebraic variety. We do not attempt to be precise, because the technicalities are complicated and prevent currently existing general conjectures from covering the case at hand. The paper \cite{BBMY2}  establishes an equivalence of categories which should be interpreted as ``homological mirror symmetry'' in this instance. In this case the coherent category is the obvious one; on the other side $\rY_\gamma$ will be seen to be Lagrangian in a certain ``de Rham moduli space'' $\cM_{\psi}$ studied in \cite{BBMY}, and the relevant ``Fukaya category'' is a certain subcategory of microlocal sheaves on $\cM_{\psi}$ supported on $\rY_\gamma$. In particular, this will imply that the maps in \eqref{eq: HMS} are isomorphisms after tensoring with $\Q$. 

We emphasize that \emph{our applications do not require} the deep categorical equivalences forthcoming in \cite{BBMY2}; we require only the construction of a functor (explained below), which is relatively easy (given existing technology). Interestingly, although Proposition \ref{prop: HMS} is a completely decategorified statement, we do not know how to produce the map except by categorical considerations. 
\end{remark}

The next two subsections, \S \ref{ssec: constructible realization} and \S \ref{ssec: microlocalization}, will not be logically used in the rest of this paper. They consist solely of summarizing certain excerpts from \cite{BBMY2}, in order to elucidate the nature of Proposition \ref{prop: HMS}.

\subsection{Constructible realization}\label{ssec: constructible realization} We work over the complex numbers $\C$. We will define a ``constructible realization'' of $\Coh_{\cB}^T(\wt{\cN})$. 

Abbreviate $\chK = L^+\chG$ for the arc group of $\chG$ and let $\chI \subset \chK$ be the Iwahori subgroup corresponding to $\chB \subset \chG$.

\subsubsection{Equivariant sheaves} Let $L^-_1 \chG$ be the first congruence subgroup of the negative loop group, whose value on a $\C$-algebra $R$ is 
\[
L^-_1\chG(R) := \ker (\chG(R[t^{-1}]) \rightarrow \chG(R)).
\]
Fix a regular semi-simple $s \in \chft$, and write $\gamma := ts \in \chfg[[t]]$. Using the Killing form, we may also view $s$ as an element of $\chfg^*$. There is a filtration of $L^-_1\chG$ by congruence subgroups, and the quotient by the second congruence subgroup gives a surjection $L^-_1\chG \surj \chfg$. We write $\psi$ for the additive character of $L^-_1\chG$ induced by inflating $s \in \chfg \cong \chfg^*$ along this map. We write $\Psi$ for the exponential $D$-module on $\A^1$ pulled back to a character sheaf on $L^-_1\chG$.

For a space with an action of $L^-_1\chG$, there is a derived category of sheaves equivariant with respect to $(L^-_1\chG, \Psi)$. More precisely, we write $\cD(\ldots)$ for presentable stable $\infty$-categories and $D(\ldots)$ for the corresponding homotopy category, and $\cD_{(L^-_1 \chG, \Psi)}(\ldots)$ or $D_{(L^-_1 \chG, \Psi)}(\ldots)$ for the respective equivariant derived categories. We apply these considerations to $\Gr_{\chG} := L\chG/\chK$ and $\Fl_{\chG} := L\chG / \chI$. More generally, for any parahoric subgroup $\chI \subset \mbf{\chP} \subset \chK$ we consider the equivariant derived categories
\[
\cD_{\psi, \mbf{\chP}} := \cD_{(L^-_1\chG, \Psi)}(L\chG / \mbf{\chP}). 
\]
We will construct a functor
\begin{equation}\label{eq: GL}
\cD \Coh^{T}_{\cB}(\wt{\cN}) \rightarrow \cD_{\psi, \chI} 
\end{equation}
following \cite[\S 6]{BBMY2}.

\subsubsection{Ingredients} To begin, we tabulate some categorical equivalences. 
\begin{enumerate}
\item Bezrukavnikov-Finkelberg constructed in \cite{BF08} a monoidal equivalence
\begin{equation}\label{eq: BF equivalence}
\cD_{\chK}(\Gr_{\chG}) \cong \cD\Coh^G(0 \dtimes_{\mf{g}} 0)
\end{equation}
where the monoidal structure is given by convolution on both sides. Here and throughout, $\dtimes$ means the derived fibered product. 

\item Bezrukavnikov constructed in \cite{Bez16} a monoidal equivalence
\begin{equation}\label{eq: B16 equivalence}
\cD_{\chI}(\Fl_{\chG}) \cong \cD\Coh^G(\wt{\cN} \dtimes_{\mf{g}} \wt{\cN})
\end{equation}
where the monoidal structure is given by convolution on both sides. 

\item Arkhipov-Bezrukavnikov-Ginzburg constructed in \cite{ABG04} an equivalence 
\begin{equation}\label{eq: ABG}
 \cD_{\chK}(\Fl_{\chG}) \cong \cD\Coh^G(0 \dtimes_{\mf{g}} \wt{\cN}).
 \end{equation}
 Moreover, the LHS carries a left convolution action $\cD_{\chK}(\Gr_{\chG}) $ and a right convolution action by $\cD_{\chI}(\Fl_{\chG})$, while the RHS carries a left convolution action by $\cD\Coh^G(0 \dtimes_{\mf{g}} 0)$ and a right convolution action by $\cD\Coh^G(\wt{\cN} \dtimes_{\mf{g}} \wt{\cN})$. The work \cite{Bez16} shows that the equivalence \eqref{eq: ABG} respects these actions via \eqref{eq: BF equivalence} and \eqref{eq: B16 equivalence}.
\item Bezrukavnikov--Boixeda Alvarez--McBreen--Yun construct in \cite[Theorem 5.5.2]{BBMY2} an equivalence
\begin{equation}\label{eq: BBMY eq 1}
\cD_{\psi, \chK} \xrightarrow{\sim} \cD(\Rep(T)).
\end{equation}
 Moreover, the LHS carries a lattice translation action on the left by $X_*(\chT)$ and a right convolution action by $\cD_{\chK}(\Gr_{\chG})$, while the RHS carries a tensoring action by $\Rep(T)$ on the left and a right convolution action by $\cD\Coh^G(0 \dtimes_{\mf{g}} 0)$. In \cite[\S 6.6]{BBMY2} it is proved that the equivalence \eqref{eq: BBMY eq 1} respects the left actions under the identification between $\Rep(T)$ and $X^*(T)\cong X_*(\chT)$-graded vector spaces, and in \cite[\S 5.4]{BBMY2} it is proved that the equivalence \eqref{eq: BBMY eq 1} respects the right actions under \eqref{eq: BF equivalence}. 
\end{enumerate}

\subsubsection{Convolution} Convolution induces a functor
\begin{equation}\label{eq: conv}
\cD_{\psi, \chK} \otimes_{\cD_{\chK}(\Gr_{\chG})} \cD_{\chK}(\Fl_{\chG}) \inj \cD_{\psi, \chI}.
\end{equation}
(In fact, this is fully faithful because $\Gr_{\chG}$ and $\Fl_{\chG}$ are ind-proper, but we will not need this.)

Convolution induces a fully faithful functor
\[
\cD(\Rep(T)) \otimes_{\cD\Coh^G(0 \dtimes_{\mf{g}} 0)}\cD\Coh^G(0 \dtimes_{\mf{g}} \wt{\cN})  \rightarrow \cD \Coh^{T}(\wt{\cN}),
\]
whose essential image is precisely $\cD \Coh^{T}_{\cB}(\wt{\cN})$. Combining this with the equivalences in (1), (2), (3) and using \eqref{eq: conv} gives the desired functor \eqref{eq: GL}. Moreover, by construction the functor is equivariant for the left action of $\Rep(T)$, and right convolution action of \eqref{eq: B16 equivalence}.

\begin{remark}[Relation to Geometric Langlands] The category $\cD_{\psi, \chI} $ is denoted $\wt{\cD}_{\psi}$ in \cite{BBMY2}. There is a pullback functor $\cD_{\psi, \chK} \rightarrow \cD_{\psi, \chI}$. Define $\cD_{\psi} \subset \cD_{\psi, \chI}$ to be the full subcategory generated by $\cD_{\psi, \chK}$ under the right convolution action of $\cD_{\chI}(\Fl_{\chG})$. Then \cite[Theorem 6.2.1]{BBMY2} shows that the functor \eqref{eq: GL} is an equivalence onto $\cD_{\psi}$. The category $\cD_{\psi}$ may be interpreted as a certain category of sheaves on the moduli stack of $\chG$-bundles on $\PP^1$ with Iwahori level structure at $0$ and certain wild ramification at $\infty$, while the category $\cD \Coh^{T}_{\cB}(\wt{\cN})$ may be interpreted as a category of coherent sheaves on a corresponding moduli space of (Betti) $G$-local systems. As such, the equivalence \eqref{eq: GL} may be viewed as proving a certain wildly ramified instance of global \emph{Geometric Langlands} -- see \cite[\S 5]{BBMY}. 
\end{remark}

\subsection{Microlocalization}\label{ssec: microlocalization} We continue to work over the complex numbers $\C$. Let $\gamma = ts$ for regular semisimple $s \in \chft$ be as in the previous subsection, and $\Fl_{\gamma} := \rY_{\gamma, \C}$, the $\C$-version of the affine Springer fiber associated to $\gamma$ (cf. Example \ref{ex: GKM description for Y}). The upshot of the constructible realization is a map 
\begin{equation}\label{eq: constructible K-map}
\rK (\Coh^{T_{\C}}_{\cB_{\C}} (\wt{\cN}_{\C})) \rightarrow \rK(D^b_{(L_1^-\chG, \Psi)}( \Fl_{\chG})).
\end{equation}
To obtain the map of Proposition \ref{prop: HMS}, we compose \eqref{eq: constructible K-map} with the \emph{singular support} map. In this situation the singular support will be a Lagrangian in the ``twisted cotangent bundle'' $T^*_\psi(L^-_1 \chG \bs \Fl_{\chG})$, and it turns out that the affine Springer fiber $\Fl_{\gamma}$ is essentially such a Lagrangian. 

\subsubsection{Moduli stack of bundles} We will now be more precise, following \cite[\S 3]{BBMY2}. We fix the curve $\PP^1$ over $\C$. Let $\chI_0 = \chI$ be the Iwahori group at $0$ and $\chK_\infty^2$ be the second congruence subgroup at $\infty$. We consider $\bupn_{\chG}(\chI_0, \chK_\infty^2)$ with $S$-points $(\cE, \tau_{\infty}, \tau_{0, \chB})$ where: 
\begin{itemize}
\item $\cE$ is a $\chG$-bundle over $\PP^1_S$. 
\item $\tau_{\infty}$ is a trivialization of $\cE$ along the divisor $2\infty_S \inj \PP^1_S$. 
\item $\tau_{0,\chB}$ is a $\chB$-reduction of $\cE$ along the divisor $0_S \inj \PP^1_S$.
\end{itemize}
In particular, the group $\chK_{\infty}^1/\chK_{\infty}^2 \cong \chfg$ acts on $\bupn_{\chG}(\chI_0, \chK_\infty^2)$ through changing the level structure $\tau_{\infty}$, inducing a moment map 
\[
\mu \co T^* \bupn_{\chG}(\chI_0, \chK_{\infty}^2) \rightarrow \chfg^*.
\]
We write $\psi \in \chfg^*$ for the image of $s \in \chfg \cong \chfg^*$ given by the Killing form. We write $\Psi$ for the additive character sheaf on $\chfg$ pulled back from the exponential $D$-module via $\psi$. 

\subsubsection{Hitchin stack} Define the Hitchin stack $\cM(\chI_0, \chK_\infty^2) := T^* \Bun_{\chG}(\chI_0, \chK_\infty^2)$, viewed with the canonical symplectic structure coming from its nature as a cotangent bundle. Then define $\cM_{\psi}$ as the symplectic reduction
\[
\cM_{\psi} :=  \left( \cM(\chI_0, \chK_\infty^2)\sslash_{\psi} \chfg \right) := [\mu^{-1}(\psi)/\chfg].
\]
In preparation for describing its functor of points, we trivialize $\omega_{\PP^1}([\infty]+[0])$ with the differential form $dt/t$, thus identifying $\omega_{\PP^1}(2[\infty]+[0]) \cong \omega_{\PP^1}([\infty])$. Then $\cM_{\psi}$ has $S$-points the groupoid of $(\cE, \tau_{\infty}, \tau_{0, \chB}, \varphi)$ where $(\cE, \tau_{\infty}, \tau_{0,\chB}) \in \bupn_{\chG}(\chI_0, \chK_\infty^2)(S)$, and $\varphi \in H^0(\PP^1_S, \Ad^* \cE \otimes \omega_{\PP^1}([\infty]))$ such that 
\begin{itemize}
\item Around $\infty_S$, $\varphi = \psi dt $ plus higher order terms under the trivialization $\tau_\infty$. 
\item $\Res_{0_S}(\varphi) \in \chfn$ under the trivialization $\tau_{0, \chB}$. 
\end{itemize}

\subsubsection{Hitchin fibration} By its nature as a symplectic reduction, $\cM_{\psi}$ carries a natural symplectic structure. There is a Hitchin fibration $f_{\psi}  \co \cM_{\psi} \rightarrow \cA_{\psi}$, which is a completely integrable system \cite[Lemma 3.2.3]{BBMY2} and a $\G_m$-action on $A_{\psi}$ contracting it to a central point $a_{\psi} \in \cA_\psi$. The central \emph{Hitchin fiber} $f_\psi^{-1}(a_\psi)$ is Lagrangian in $\cM_\psi$. 

\subsubsection{Microlocalization}
By the uniformization of $\Bun_{\chG}$ for $\PP^1$, there is a canonical equivalence 
\[
\cD_{\psi, \chI}  \xrightarrow{\sim} \cD_{(\chfg, \Psi)}(\bupn_{\chG}(\chI_0, \chK_\infty^2)).
\]
Then the formation of singular support induces a map 
\begin{equation}\label{eq: SS D}
\rK(D^b_{(\chfg, \Psi)}( \bupn_{\chG}(\chI_0, \chK_\infty^2 )) \rightarrow \topCh(f_\psi^{-1}(a_\psi)). 
\end{equation}

\begin{remark}
The map \eqref{eq: SS D} is categorified by the \emph{microlocalization functor} 
\[
\muloc \co \cD_{ (\chfg, \Psi)} (\bupn_{\chG}(\chI_0, \chK_\infty^2 )) \rightarrow \muSh_{f_\psi^{-1}(a_\psi)} (\cM_\psi),
\]
where the right hand side is the category of \emph{microlocal sheaves with support in $f_\psi^{-1}(a_\psi)$}. The singular support of $\cF \in D^b_{(\chfg, \Psi)}( \bupn_{\chG}(\chI_0, \chK_\infty^2 ))$ is the naive support of $\muloc(\cF)$. 
\end{remark}

\subsubsection{Relation to affine Springer fiber} By \cite[Proposition 3.4.2]{BBMY2}, there is a canonical homeomorphism $\Fl_\gamma \rightarrow f_\psi^{-1}(a_\psi)$, which in particular induces an isomorphism 
\begin{equation}\label{eq: homeomorphism}
\topCh(f_\psi^{-1}(a_\psi)) \rightarrow \topCh(\Fl_\gamma)
\end{equation}
compatible with the two actions of $\wt{W}$. Finally, composing \eqref{eq: SS D} with \eqref{eq: homeomorphism} and \eqref{eq: constructible K-map} gives the desired map of Proposition \ref{prop: HMS}.

\subsection{The microlocal support map}\label{ssec: musupp} We may now define a map that casts representations of $G_1T$ onto cycles in $\topCh(\rY_\gamma)$. 

\begin{defn}\label{defn: musupp}
We define the \emph{microlocal support}\footnote{The notation ``SS'' is an abbreviation for ``singular support''. This is synonymous with ``microlocal support'' which is often denoted $\mu\mrm{supp}$; we do not use the latter notation because it would render certain equations too long.} map 
\[
\musupp \co \rK(\Rep_0( \Ug^0, T)) \rightarrow \topCh(\rY_\gamma)
\] 
to be the composition of the maps from Theorem \ref{thm:mod-rep-theory-summary}, the identification $\rK(\Coh^{T_{\C}}_{\cB_{\C}}(\wt \cN_{\C})) \cong \rK(\Coh^{T_{\C}}(\cB_{\C}))$, and Proposition \ref{prop: HMS}. 

Note that we have 
\[
\rK(\Rep_0( \Ug^0, T)) \xleftarrow{\sim} \rK(\Rep( \Ug^0_0, T)) \xrightarrow{\sim} \rK(\Rep^0(\Ug_0, T)) \cong \rK(\Rep^\emptyset (G_1 T)),
\]
where the first two isomorphisms come from the fact that the centers act by scalars on simple representations, and the last isomorphism comes from \eqref{eq: frob kernel}. Therefore we will also view $\musupp$ as being defined on any of these other groups. In Part 3, we will mostly view it as a map
\[
\musupp \co \rK(\Rep^\emptyset (G_1 T)) \rightarrow \topCh(\rY_\gamma).
\]
\end{defn}

Since $\topCh(\rY_\gamma) \subset \topBM(\rY_\gamma)$, we will also regard $\musupp$ as having target in $\topBM(\rY_\gamma)$ at times. By Theorem \ref{thm: BM} and Proposition \ref{prop: HMS}, we know that $\musupp$ has the following properties. 
\begin{enumerate}
\item It intertwines the action of $(\wt{W}, \cdot_p)$ on the LHS (defined in \S \ref{ssec: rep m-c action}) with the action of $(\wt{W}, \cdot)$ on the RHS (defined in \S \ref{ssec: equivariant actions}). 
\item It intertwines the action of $(\wt{W}, \bup)$ on the LHS (defined in \S \ref{ssec: braid action}) with the action of $(\wt{W}, \bu)$ on the RHS (defined in \S \ref{ssec: equivariant actions}). 
\item It sends $[\wh{L}(0)] \in \rK(\Rep_0( \Ug^0, T)) $ to the fundamental class of the unique (top-dimensional) irreducible component of $\rY_{\gamma}$ which is the pre-image of $[t^0] \in \Gr_{\chG, \F}$ under the projection map $\Fl_{\chG, \F} \rightarrow \Gr_{\chG, \F}$.
\end{enumerate}

Since the map $\musupp$ is essential, and its definition meandered through a rather serpentine construction, we recapitulate it in the diagram below, where the left (resp. right) column pertains to $G$ (resp. $\chG$). 
\[
\begin{tikzcd}[column sep = huge]
\rK(\Rep^\emptyset (G_1 T)) \ar[drrr, "\musupp", dashed] \ar[d, equals] \\ 
\rK(\Rep_0( \Ug^0, T)) \ar[rrr, "\musupp"', dashed] \ar[d, "\text{BMR localization}"']  &  & &    \topCh(\rY_\gamma) \\
\rK(\Coh_{\cB^{(1)}}^{T^{(1)}}(\wt{\cN}^{(1)})) \ar[d, "(\fsp_{p \rightarrow 0})^{-1}"'] & & & \topCh(\rY_{\gamma,\C})  \ar[u, "\fsp_{p \rightarrow 0}"']  \\
 \rK (\Coh^{T_{\C}}_{\cB_{\C}} (\wt{\cN}_{\C})) \ar[rrr, "\text{geometric Langlands}"'] \ar[urrr, "\text{mirror symmetry}"] &  & &  \rK(\cD_\psi) \ar[u, "\text{microlocalization}"'] 
\end{tikzcd}
\]

\begin{remark}
The identification $\topCh(\rY_{\gamma,\C}) \xrightarrow{\sim} \topCh(\rY_\gamma)$ from Example \ref{ex: GKM description for Y} was not defined via $\fsp_{p \rightarrow 0}$. However, from the observation that the degeneration from characteristic zero to characteristic $p$ is constant on $\chT$-fixed points, it is clear that the identification using the GKM description coincides with $\fsp_{p \rightarrow 0}$. 
\end{remark}

\section{Degeneration of affine Springer fibers}\label{sec: Degeneration of affine Springer fibers}

For the rest of the paper, we view $G$ as an unramified reductive group over $\Z_p$, satisfying the hypotheses of \S \ref{sssec:reductive-groups}:
 \begin{hypothesis}\label{hyp:part3-hypotheses} $G_{\mrm{der}}$ is simply connected, $Z(G)$ is connected, and $G$ admits a ``local twisting element'', which we denote $\rho$. 
 \end{hypothesis}

In this section, we assemble the ingredients from the preceding sections in order to finally do geometric calculations related to the Breuil--M\'ezard Conjecture (although the precise connection will not be explained until Part 3). Our goal here is to ``understand'' the limit cycle $\fsp_{p \rightarrow 0} [\rX_\gamma(\lambda)] \in \topCh(\rY_\gamma)$. We will express this cycle in terms of representation theory via the microlocal support map from Definition \ref{defn: musupp}. 

The first difficulty is that the specialization process is a priori mysterious, at least in terms of the basis of irreducible components, since specialization does not (in general) interact well with the properties of being reduced or irreducible. A key tool for us is equivariant localization, which allows us to calculate the specialization in terms of torus-fixed points instead. We illustrate this in \S \ref{ssec: equivariant class of limit cycle}, where we calculate explicitly the equivariant fundamental class of $\fsp_{p \rightarrow 0} [\rX_\gamma(\rho)]$. We then go on in \S \ref{ssec: general Schubert cells} to express $\fsp_{p \rightarrow 0} [\rX_\gamma(\lambda)]$, for general $\lambda$, in terms of the case $\lambda = \rho$. In Part 3, this calculation will be used to reduce the Breuil--M\'ezard Conjecture for all sufficiently small $\lambda$ to the fundamental case $\lambda = \rho$. 

At the next stage, we want to identify $\fsp_{p \rightarrow 0} [\rX_\gamma(\rho)]$ with the microlocal support of a particular baby Verma module in $\Rep^0(G_1T)$. Although we know an explicit formula for the equivariant fundamental class of the former object, it does not seem to be easy to compute the latter object in these terms. Hence we have to take a more indirect approach. In \S \ref{ssec: recognition principle} we prove a ``Recognition Principle'' that characterizes the class $\fsp_{p \rightarrow 0} [\rX_\gamma(\rho)]$ in somewhat more conceptual terms. Then in \S \ref{ssec: microlocal support of baby verma} and Appendix \ref{app: A}, which is joint with Bezrukavnikov and Boixeda Alvarez, we check that this Recognition Principle applies to the desired baby Verma module. 


\subsection{Equivariant class of limit cycles}\label{ssec: equivariant class of limit cycle}

We fix, once and for all, a generator $\topom$ of $\Omega_{\chft}^{\wedge \dim \chft}$. We also abbreviate $\sph := \sph_{\chT}$. For each root $\alpha \in \chPhi$, we have its derivative $d\alpha \in \Sym_{\Ql}( \chft^*) \cong \sph$. Define 
\begin{equation}\label{eq: Delta}
\beta := \prod_{\alpha \in \chPhi^+} d\alpha  \in \Sym_{\Ql}( \chft^*) \cong \sph. 
\end{equation}
Recall from Example \ref{ex: fixed points X} that there is an isomorphism $X_*(\chT) \xrightarrow{\sim} (\rX_\gamma^{\varepsilon = \varepsilon_0})^{\chT}$ sending $\lambda \mapsto t^\lambda$, which identifies (implicitly using $\topom$)
\begin{equation}\label{eq: gr fixed H_*}
\mBMT_*((\rX_\gamma^{\varepsilon = \varepsilon_0})^{\chT}) \cong \bigoplus_{\lambda \in 
X_*(\chT)} \sph[t^\lambda].
\end{equation}

Recall that if $X$ is $\chT$-equivariantly formal, then for $\alpha \in \topBM(X)$ we write $\alpha_{\chT} \in \topBMT(X)$ for its equivariant lift, which exists and is unique by Remark \ref{remark: de-equivariant BM homology}. This applies to $X = \rX_\gamma^{\varepsilon = \varepsilon_0}$, whose $\chT$-fixed points may be identified with $X_*(\chT)$ as in Example \ref{ex: fixed points X}. 

\begin{lemma}\label{lem: equivariant class} Let $\varepsilon_0  \in \A^1_{E}$ be non-zero and let $\lambda \in X_*(\chT)^+$ be regular. Then, with respect to \eqref{eq: gr fixed H_*}, we have
\begin{equation}\label{eq: equivariant class}
\Loc^{\chT}([\rX^{\varepsilon=\varepsilon_0}_\gamma(\lambda)]_{\chT}) = \frac{1}{\beta} \sum_{w \in W} \sgn(w) [t^{w \lambda}]  \in 
\mBMT_*((\rX_\gamma^{\varepsilon = \varepsilon_0})^{\chT}) \otimes_{\sph} \Frac(\sph).
\end{equation}
\end{lemma}

\begin{proof}
Since $\lambda$ is regular, by Lemma \ref{lem: X lambda irred e neq 0} the Bialynicki-Birula map takes $\rX^{\varepsilon=\varepsilon_0}_\gamma(\lambda)$ isomorphically to $\chG/\chB$. Therefore, the $\chT$-fixed points of $\rX^{\varepsilon=\varepsilon_0}_\gamma(\lambda)$ are identified with $\{[t^{w\lambda}]\}_{w \in W} \subset X_*(\chT) = (\Gr_{\chG, E})^{\chT}$. 

As explained in Example \ref{ex: flag variety}, the component at $wB \in (\chG/\chB)^{\chT}$ of $\Loc^{\chT}([\chG/\chB]_T) \in \mBMT_{2d}(\chG/\chB)$ is $(\frac{\sgn(w)}{\beta})$. Then we conclude using that $wB \in \chG/\chB$ corresponds to the fixed point $[t^{w \lambda }] \in \rX^{\varepsilon=\varepsilon_0}_\gamma(\lambda)$ under the Bialynicki-Birula isomorphism $\rX^{\varepsilon=\varepsilon_0}_\gamma(\lambda) \xrightarrow{\sim} \chG/\chB$. 
\end{proof}

\begin{cor}\label{cor: basic cpt e=0 equivariant class} With $\beta$ defined as \eqref{eq: Delta}, we have 
\[
\Loc^{\chT}([\rX^{\varepsilon=0}_\gamma(\rho)]_{\chT}) = \frac{1}{\beta} \sum_{w \in W} \sgn(w) [t^{w  \rho}] \in \mBMT_*((\rX_\gamma^{\varepsilon = 0})^{\chT}) \otimes_{\sph} \Frac(\sph).
\]
\end{cor}
\begin{proof} We use Lemma \ref{lem: equivariant class} to obtain the analogous description of $[\rX^{\varepsilon=\eta}_\gamma(\rho)]_{\chT}$. Then the identity follows from (the equivariant version of) Proposition \ref{prop: gr hbar spc}, and the fact that $\fsp_{\varepsilon\rightarrow 0}$ is the identity map in the common GKM description of $\rH^{\BM, T}_*(\rX^{\varepsilon}_\gamma)$.
\end{proof}

\subsection{General Schubert cells}\label{ssec: general Schubert cells}

In this part we consider $\rX^{\varepsilon = \eta}_{\gamma}(\lambda+\rho) \inj \rX^{\varepsilon = \eta}_\gamma$, the closure of $\rX^{\varepsilon \neq 0}_{\gamma}\cap S^\circ (\lambda+\rho)$ indexed by general $\lambda \in X_*(T)^+$. We will calculate the fundamental class $[\rX^{\varepsilon = \eta}_{\gamma}(\lambda+\rho) ]$ in terms of the ``basic'' case $\lambda=0$, which was analyzed in Lemma \ref{lem: equivariant class}. This will have significance for the Breuil--M\'ezard Conjecture with Hodge-Tate weights $\lambda + \rho$.

By the Equivariant Localization Theorem and Example \ref{ex: fixed points X}, we have 
\begin{equation}\label{eq: eq loc X}
\mBMT_*(\rX_\gamma^{\varepsilon = \varepsilon_0}) \otimes_{\sph} \Frac(\sph) \cong \bigoplus_{\mu \in X_*(\chT)} \Frac(\sph)[t^\mu].
\end{equation}
This has an obvious (left) action of $X_*(\chT)$, through left translation on the indexing set.

\begin{thm}\label{thm: generic fiber cycle decomposition}
Let $\varepsilon_0 \in \A^1_{E}$ be non-zero and $\lambda \in X_*(\chT)^+ \cong X^*(T)^+$. Then we have
\[
\Loc^{\chT} [\rX_{\gamma}^{\varepsilon = \varepsilon_0 }(\lambda+\rho)]_{\chT}  = \sum_{\mu \leq \lambda}  m_{\mu}(\lambda)  t^{\mu} \cdot  \Loc^{\chT}[\rX_{\gamma}^{\varepsilon = \varepsilon_0}(\rho)]_{\chT}  \in \mBMT_*(\rX_\gamma^{\varepsilon = \varepsilon_0}) \otimes_{\sph} \Frac(\sph),
\]
where $m_{\mu}(\lambda)$ is the multiplicity of the weight $\mu$ in the Weyl module $W(\lambda)$ of $G_{\ol\Q_p}$. Here the action $t^{\mu} \cdot $ is the one defined just above, through \eqref{eq: eq loc X}. 
\end{thm}

\begin{proof}
According to Lemma \ref{lem: equivariant class}, we have
\begin{equation}\label{eq: weights 1}
\Loc^{\chT} ([\rX^{\varepsilon=\varepsilon_0}_\gamma(\lambda+\rho)]_{\chT}) = \frac{1}{\beta} \sum_{w \in W} \sgn(w) [t^{w (\lambda + \rho)}] \in \bigoplus_{\mu \in X^*(T)} \Frac(\sph)[t^\mu]
\end{equation}
and
\begin{equation}\label{eq: weights 1.5}
\Loc^{\chT}([\rX^{\varepsilon=\varepsilon_0}_\gamma(\rho)]_{\chT}) = \frac{1}{\beta} \sum_{w \in W} \sgn(w) [t^{w \rho}] \in \bigoplus_{\mu \in X^*(T)} \Frac(\sph)[t^\mu].
\end{equation}
From \eqref{eq: weights 1.5} we find that 
\begin{equation}\label{eq: weights 2}
\sum_{\mu \leq \lambda} m_{\mu}(\lambda) t^\mu \cdot [\rX^{\varepsilon=\varepsilon_0}_\gamma(\rho)]_{\chT} = \frac{1}{\beta} \sum_{ \mu \leq \lambda} m_{\mu}({\lambda}) t^\mu \cdot \left( \sum_{w \in W} \sgn(w) [t^{w \rho}] \right) \in \mBMT_*(\rX_\gamma^{\varepsilon = \varepsilon_0})  \otimes_{\sph} \Frac(\sph).
\end{equation}



Below we recall the Weyl character formula. To set notation, we regard the characters of representations of $G_{\ol \Q_p}$ as elements of the group ring $\Ql[X^*(T)]  = \Ql[X_*(\chT)]$; recall that $X^*(T)$ is the \emph{geometric} character group. When writing characters, we use $e^\lambda \in \Ql[X^*(T)]$ to represent the group element $\lambda \in X^*(T) \cong X_*(\chT)$. In these terms, the Weyl character formula says that
\[
\sum_{\mu \leq \lambda} m_{\mu}(\lambda) e^{\mu} = \frac{\sum_{w \in W} \sgn(w) e^{w(\lambda + \rho)}}{\prod_{\alpha \in \Phi^+} (e^{\alpha/2} - e^{-\alpha/2})} \in \Ql[X^*(T)]
\]
and the Weyl denominator formula says that
\[
\prod_{\alpha \in \Phi^+} (e^{\alpha/2} - e^{-\alpha/2}) = \sum_{w \in W} \sgn(w) e^{w \rho}  \in \Ql[X^*(T)].
\]
Combining them, we find that 
\[
\sum_{w \in W} \sgn(w) e^{w(\lambda + \rho)} = \sum_{\mu \leq \lambda} m_{\mu}(\lambda) e^\mu \left(\sum_{w \in W} \sgn(w) e^{w \rho} \right) \in  \Ql[X^*(T)].
\]
Hence we have an identity
\begin{equation}\label{eq: weights 3}
\sum_{ w \in W} \sgn(w) [t^{w(\lambda+ \rho)}]  = \sum_{ \mu \leq \lambda} m_{\mu}(\lambda) [t^{\mu}] \left(\sum_{w \in W} \sgn(w) [t^{w \rho}] \right) \in \bigoplus_{\mu \in X_*(\chT)} \Frac(\sph)[t^\mu].
\end{equation}
Now the desired equality follows from comparing \eqref{eq: weights 1}, \eqref{eq: weights 2}, and \eqref{eq: weights 3}.
\end{proof}

Below, for ease of notation we abbreviate $C_{\wt{w}}$ for the irreducible component $\rY_{\gamma}(\wt{w})$ of $\rY_\gamma$ from Corollary \ref{cor: top cycles}, for all $\wt{w} \in \wt{W}^{\reg}$. As usual, $[C_{\wt{w}}]$ is its cycle class and $[C_{\wt{w}}]_{\chT}$ its $\chT$-equivariant lift. 

\begin{lemma}\label{lem:Euler class of component} The irreducible component $C_{t^{w\rho}} \subset \rY_\gamma(\leq \rho)$ has $t^{w\rho}$ as a smooth point, and the equivariant Euler class of $[C_{t^{w\rho}}]_{\chT}$ at $t^{w\rho}$ is 
\[
e_T(t^{w\rho},C_{t^{w\rho}})=\frac{\sgn(w)}{\beta}.
\]
\end{lemma}
\begin{proof} This follows from the proof of Lemma \ref{lem: Y lambda irred e=0} and the explicit description of the affine space chart around $t^{w\rho}$.
\end{proof}

\subsection{Recognition principle}\label{ssec: recognition principle} Let $\sph := \sph_{\chT}$. We will refer to the GKM description of $\mBMT_*(\rY_{\gamma}) $ inside
\[
\mBMT_*(\rY_{\gamma}^{\chT})  \otimes_{\sph} \Frac(\sph) \cong \bigoplus_{w \in \dWext} \Frac(\sph)  [\wt{w}].
\]

\begin{prop}[Recognition principle]\label{prop: uniqueness of cycle}
There is a unique class $[Z] \in 
\topBMT(\rY_{\gamma}) $ which in the GKM description has the following properties: 
\begin{enumerate}
\item (Eigenclass) $[Z]$ has equivariant support contained in $X_*(\chT) \subset \dWext$. 
\item (Support bound) $[Z]$ has equivariant support contained in the regular admissible set 
\[
\mrm{Adm}(\rho) := \{ \wt{w} \in \dWext \co \wt{w} \leq t^{w \rho} \quad \text{ for some } w \in W \}.
\]
\item (Normalization) The component of $[Z]$ at $t^{\rho}$ is $1/\beta  \in \Frac(\sph)$.
\end{enumerate}
\end{prop}

\begin{remark}
Condition (1) can be formulated alternatively as an eigenproperty (cf. Proposition \ref{prop: eigenproperty verma} below) for the lattice parts of the two actions $(\dWext, \cdot)$ and $(\dWext, \bu)$, which can be thought of as ``Hecke actions'' (because they are literally given by convolving with Hecke operators in the constructible realization of \S \ref{ssec: constructible realization}). This explains why we call (1) an ``eigenclass'' property.\footnote{In an earlier draft of this paper, the main result of this section was proved using a different Recognition Principle formulated in terms of a more abstract eigenproperty, before we realized that we could instead work with a more explicit characterization in terms of translation elements.} 
\end{remark}

\begin{proof}
Let $[Z]$ and $[Z']$ be two elements of $\topBMT(\rY_{\gamma}) $ satisfying all of these conditions. Then $[Z]-[Z']$ is a class $\delta =\sum a_{\wt{w}}[\wt{w}]\in \mBMT_*(\rY_\gamma)$ with equivariant support contained in $X_*(\chT) \cap \mrm{Adm}(\rho) $, but not supported at $t^\rho$. 

Note that as $C_{\wt{w}}$ is the closure of $C_{\wt{w}}\cap S^\circ(\wt{w})$, $[C_{\wt{w}}]_{\chT}$ has equivariant support only on elements $\leq \wt{w}$. We have
\[
\delta=\sum_{\wt{w} \in \wt{W}^{\reg}} m_{\wt{w}}[C_{\wt{w}}]_{\chT}, \hspace{1cm} m_{\wt{w}} \in \Z.
\]
 If $\wt{w}$ is a maximal element such that $m_{\wt{w}}\neq 0$, then $\delta$ has a non-trivial coefficient at $[\wt{w}]$, so $\wt{w}\in \Adm(\rho)$. We conclude that only $\wt{w}\in \Adm^{\reg}(\rho)$ contributes in the above sum.

Let $w\in W$, and $\alpha$ be a simple root of $\chT$ such that $ws_{\alpha}>w$, and set $\wt{u}=wt^{\rho}s_{\alpha}w^{-1}$. Then:
\begin{itemize}
\item (cf. \cite[Proposition 4.6]{BL21}\footnote{This is a reference to \cite[\S 5.11]{GKM04}, which focuses on the version of the affine Springer fiber over $\C$, but the same analysis applies in characteristic $p$, as explained in Example \ref{ex: GKM description for Y}.}) The $1$-dimensional $\chT$-orbits of $\Fl_{\chG, \F}$ joining $t^{w\rho}$ and $\wt{u}$, resp. $t^{ws_\alpha\rho}$ and $\wt{u}$ belong to $\rY_\gamma(\leq\rho)$. Furthermore, both these orbits are associated to the character $w\alpha\in X^*(\chT)$. 
\item By \cite[Proposition 2.2.6]{LLLMextremal}, the connected component of $\rY_{\gamma}^{\ker w\alpha}$ passing through $t^{w\rho}$ intersects $\rY_{\gamma}(\leq \rho)^{\chT}=\Adm(\rho)$ at exactly $t^{w\rho},\wt{u}$ and $t^{ws_\alpha\rho}$.
\end{itemize}
The GKM description thus gives
\[
\Res_{d(w\alpha)}(a_{t^{w\rho}})+\Res_{d(w\alpha)}(a_{\wt{u}})+\Res_{d(w\alpha)}(a_{t^{ws_\alpha\rho}})=0.
\]
But our hypothesis implies that the middle term vanishes (since $\wt{u}$ is not a translation), hence the two outer terms must either both vanish or both not vanish. Since $t^{w\rho}$ is maximal in $\Adm(\rho)$, $[C_{t^{w\rho}}]_{\chT}$ is the only irreducible component that contributes to $a_{t^{w\rho}}$. Hence by Lemma \ref{lem:Euler class of component} we have 
\[a_{t^{w\rho}}=m_{t^{w\rho}}\frac{\sgn{w}}{\beta},\]
and we have an analogous formula for $a_{t^{ws_\alpha\rho}}$. 
Combining this with the previous observation, we learn that $m_{t^{w\rho}}$ and $m_{t^{ws_\alpha\rho}}$ are either both zero or both non-zero. Since we also have $m_{t^{\rho}}=0$, this gives $m_{t^{w\rho}}=0$ for all $w\in W$, so $m_{t^{w\rho}} = m_{t^{ws_\alpha\rho}} = 0$. 

Finally, if $\delta\neq 0$, then there must be a maximal element $\wt{w}$ such that $m_{\wt{w}}\neq 0$. Then $\delta$ has non-trivial coefficient at $[\wt{w}]$, so $\wt{w}$ is a translation in $\Adm^{\reg}(\rho)$. However, Lemma \ref{lem:regular admissible translation} below shows that such a translation must be of the form $t^{w\rho}$, contradicting what we showed in the previous paragraph that $m_{t^{w\rho}}=0$ for all $w\in W$.

\end{proof}

\begin{lemma}\label{lem:regular admissible translation} We have 
\[
\Adm^{\reg}(\rho)\cap X_*(\chT)=\{t^{w\rho}  \co  w\in W\}.
\]
\end{lemma}
\begin{proof} Suppose $t^{\mu}\in \Adm^{\reg}(\rho)$. Then by \cite[Corollary 2.1.7]{LLLM22} we have
\[t^{\mu}=\wt{w}_2^{-1}w_0\wt{w}_1\]
where 
\begin{itemize}
\item $\wt{w}_1\in \dWext_1$ is restricted, i.e., $\wt{w}_1(A_0)$ belongs to the fundamental box (cf. \S \ref{sssec: weyl notation}),
\item $\wt{w}_2\in \dWext^+$ is dominant, i.e., $\wt{w}_2(A_0)$ belongs to the dominant cone, 
\item $\wt{w}_2 \uparrow w_0t^{-\rho}\wt{w}_1$. (See \cite[II.6]{Jan03} for the definition of the $\uparrow$ order.) 
\end{itemize}
The fact that $\wt{w}_2^{-1}w_0\wt{w}_1$ is a translation in $X_*(\chT)$ shows that 
\[
t^{\nu}w_0t^{-\rho}\wt{w}_1= \wt{w}_2 \quad \text{for some $\nu\in X_*(\chT)$}.
\]
Since any dominant alcove is uniquely a dominant translation of a restricted alcove, $\nu$ is dominant. But then $ w_0t^{-\rho}\wt{w}_1\uparrow \wt{w}_2$, so the third bullet point above forces $\nu=0$, and 
\[\wt{w}_2^{-1}w_0\wt{w}_1=( w_0t^{-\rho}\wt{w}_1)^{-1}w_0\wt{w}_1=t^{w_1^{-1}\rho},\]
where $w_1$ is the image of $\wt{w}_1$ in $W$.
\end{proof}


\begin{remark}
The proof of Proposition \ref{prop: uniqueness of cycle} applies also to deformed affine Springer fibers as long as they satisfy the GKM conditions. By Lemma \ref{lem: Y lambda irred e=0}, this holds for $\gamma = ts$ when $s$ is regular semi-simple and $\varepsilon = \eta$, and if $s$ is furthermore $h_\rho$-generic then it holds for $\rY_\gamma^{\varepsilon = 1}(\leq \rho)$. The only step in the proof that requires additional commentary is the calculation of 1-dimensional $\chT$-orbits, which in the deformed case is done as in the proof of \cite[Proposition 4.3.]{LLLM22}. 
\end{remark}

\begin{lemma}\label{lem: limit cycle recognition}
The class
\[
\fsp_{p \rightarrow 0} [\rX^{\varepsilon=0}_\gamma(\rho)]_{\chT} \in \topBMT(\rY_\gamma)
\]
satisfies the conditions of Proposition \ref{prop: uniqueness of cycle}. 
\end{lemma}

\begin{proof}
By \S \ref{sssec: functoriality for specialization maps} we have a commutative diagram 
\[
\begin{tikzcd}
\mBMT_*( \rX_{\gamma}^{\chT}) \ar[d] \ar[r, "\fsp_{p \rightarrow 0}"] & \mBMT_*(\rY_\gamma^{\chT}) \ar[d] \\
\mBMT_*(\rX_{\gamma}) \ar[r, "\fsp_{p \rightarrow 0}"] & \mBMT_*(\rY_\gamma)
\end{tikzcd}
\]
From this and Corollary \ref{cor: basic cpt e=0 equivariant class} we conclude 
\begin{equation}\label{eq: loc limit cycle}
\Loc^{\chT}\left( \fsp_{p \rightarrow 0} [\rX^{\varepsilon=0}_\gamma(\rho)]_{\chT}  \right) = \frac{1}{\beta} \sum_{w \in W} \sgn(w) [t^{w\rho}],
\end{equation}
which visibly satisfies the conditions of Proposition \ref{prop: uniqueness of cycle}.
\end{proof}

\subsection{The microlocal support of baby Verma modules}\label{ssec: microlocal support of baby verma} Recall the microlocal support map $\musupp$ from Definition \ref{defn: musupp}. Henceforth we mostly focus on $\rK(\Rep^\emptyset(G_1T))$ instead of $\rK(\Rep_0(\Ug^0, T))$; although these are canonically isomorphic we remind that they carry different normalizations for their baby Verma modules. The following result, obtained jointly with Bezrukavnikov and Boixeda Alvarez, is established in Appendix \ref{app: A}. 

\begin{thm}[Joint with Bezrukavnikov--Boixeda Alvarez]\label{thm: baby verma equals limit cycle} We have
\[
\musupp[\wh{Z}_1(p\rho)]_{\chT}=\fsp_{p\rightarrow 0}[\rX^{\varepsilon=0}_\gamma(\rho)]_{\chT}  \in \topBMT(\rY_\gamma)
\]
and under $\Loc^{\chT}$ they are given by $\frac{1}{\beta} \sum_{w \in W} \sgn(w) [t^{w \rho}]$. 
\end{thm}

For the proof of Theorem \ref{thm: baby verma equals limit cycle}, and also for other purposes later, we need the following technical lemma.

\begin{lemma}\label{lem:supp of simple} Let $\wt{w}\in \dWext_1$ be restricted (cf. \S \ref{sssec: weyl notation} for notation). Then  $\musupp [\wh{L}_1(\wt{w} \bup 0 )]_{\chT}$ has equivariant support in $\dWext_{\leq w_0\wt{w}}=\{\wt{u} | \wt{u}\leq w_0\wt{w}\}$. Furthermore, $W\wt{w}$ occurs in the equivariant support.
\end{lemma}
\begin{proof} We induct on $\ell(\wt{w})$, the case $\ell(\wt{w})=0$ being true by the normalization condition (3) of \S \ref{ssec: musupp}. 

If $\ell(\wt{w})>0$, we can find a simple affine reflection $s$ such that $\wt{w}s\in \dWext_1$ and $\wt{w}s<\wt{w}$. By \cite[II.7.15--21, II.9.22]{Jan03} (for description of the wall-crossing functors for $G$ and then their variants for $G_r T$, respectively), the wall-crossing functor $R_s$ ($=\Theta_s$ in the notation of \emph{loc. cit.}) satisfies
\[
[\wh{L}_1(\wt{w} \bup 0  )] = [R_s (\wh{L}_1( \wt{w}s \bup 0 ))] +\sum_{\substack{\wt{u}\in \dWext^+ \co \\ \wt{u}\leq \wt{w}s}} m_{\wt{u}} [\wh{L}_1(\wt{u} \bup 0)] \in \rK(\Rep^0(G_1T))
\]
for some $m_{\wt{u}} \in \Z$ (note that $\wt{u}$ need not be restricted). We will show that the equivariant support of $\musupp$ applied to each of the terms on the right-hand side belongs to $W_{\leq w_0\wt{w}}$:
\begin{enumerate}
\item We have $\musupp [R_s(\wh{L}_1( \wt{w}s \bup 0 ))]_{\chT}= \musupp[\wh{L}_1(\wt{w}s \bup 0 )]_{\chT} \bup s +\musupp[\wh{L}_1(\wt{w}s \bup 0 )]_{\chT}$. By induction, the equivariant support of $\musupp[\wh{L}_1(\wt{w}s \bup 0 )]_{\chT}$ belongs to $\wt{W}_{\leq w_0\wt{w}s}$, and the equivariant support of $\musupp[\wh{L}_1(\wt{w}s \bup 0 )]_{\chT} \bup s$ belongs to $\wt{W}_{\leq w_0\wt{w}s}s\subset \wt{W}_{\leq w_0\wt{w}}$, since $w_0\wt{w}=(w_0\wt{w}s)s$ is a reduced factorization and $s$ is simple. Furthermore, by induction $w_0\wt{w}$ does occur in the equivariant support of this term. 
\item For $\musupp[\wh{L}_1(\wt{u} \bup  0)]_{\chT}$: Decompose $\wt{u}=t^{\nu}\wt{v}$ with $\wt{v}\in \dWext_1$ and $\nu$ dominant. Note that $\wt{v}(A_0)\uparrow t^{\nu}\wt{v}(A_0)\uparrow \wt{w}s(A_0)$, so $\ell(\wt{v})\leq \ell(\wt{w}s)<\ell(\wt{w})$ (note however that $\wt{v}$ and $\wt{w}$ may be incomparable since $\nu$ may fail to be in $Q$). By induction, the equivariant support belongs to
\[
t^\nu\dWext_{\leq w_0 \wt{v}}\subset \dWext_{\leq w_0\wt{w}s} \subset \dWext_{\leq w_0\wt{w}}
\]
where the first inclusion is Lemma \ref{lem: support bound simple} below. We also note that $w_0\wt{w}$ cannot occur in the equivariant support of this term, since $\ell(\wt{v})<\ell(\wt{w})$.
\end{enumerate}
Finally, note that we have shown that $w_0\wt{w}$ does occur in the equivariant support. Example \ref{ex: monodromy action on simples} then shows that all of $W\wt{w}$ occurs in the equivariant support.
\end{proof}

\begin{lemma}\label{lem: support bound simple} Let $\wt{w}\in \dWext_1$, $\nu \in X^*(T)^+$, and $\wt{v}\in \dWext_1$ such that
$t^{\nu}\wt{v}\leq \wt{w}$. Then  we have 
\[
t^{\nu}\dWext_{\leq w_0 \wt{v}}\subset \dWext_{\leq w_0\wt{w}}.
\]
\end{lemma}
\begin{proof} Let $\wt{u}\leq w_0\wt{v}$, so $\wt{u}_{\dom}\leq \wt{v}$. Let $\sigma\in W$ be the unique element such that $\sigma t^{\nu}\wt{u}$ is dominant. It suffices to show $\sigma t^{\nu}\wt{u}\leq \wt{w}$.

Since $\nu$ is dominant, we have $\nu-\sigma\nu \geq 0$, and 
\[
\sigma t^{\nu} \wt{u} \uparrow t^{-\sigma\nu+\nu}\sigma t^{\nu}\wt{u}=t^{\nu}\sigma \wt{u}.
\]
We also have $\sigma \wt{u} \uparrow \wt{u}_{\dom}$, hence
\[
\sigma t^{\nu}\wt{u}\uparrow t^{\nu}\wt{u}_{\dom}.
\]
But for dominant elements of $\wt{W}$, the $\uparrow$ order and the Bruhat order coincide, hence
\[\sigma t^{\nu}\wt{u}\leq t^{\nu}\wt{u}_{\dom}\leq t^{\nu}\wt{v}\leq \wt{w},\]
as desired.
\end{proof}

 \part{Applications to the Breuil-M\'{e}zard Conjecture}
 
 Finally, we assemble the preceding ingredients for application to the Breuil--M\'ezard Conjecture. 
 
\section{The Breuil-M\'{e}zard Conjecture}
In this section we set up the formulation of the geometric, refined version of the Breuil--M\'ezard Conjecture. In \S \ref{ssec: Serre weights} we discuss the parametrization of Serre weights. In \S \ref{ssec: tame inertial types} we define tame inertial types and their parametrization by Deligne-Lusztig representations. In \S \ref{ssec: EG stack} we recall the generalization of Emerton-Gee stack for tame groups due to Lin \cite{Lin23b} and state Emerton-Gee's formulation of the Breuil--M\'ezard Conjecture and its (partial) generalization for tame groups. 

\subsection{Frobenius automorphism}\label{ssec:Frobenius}
Let $F$ be the relative Frobenius of $\ulG_{\ol{\F}_p}$, so $\ul{G}^F \cong \ul{G}(\F_p)$. Let $\pi$ be the finite order automorphism of $(\ul{G}_{\ol{\F}_p}, \ul{B}_{\ol{\F}_p}, \ul{T}_{\ol{\F}_p})$ from \cite[Proof of Lemma 9.2.4]{GHS18}. It is characterized by the property that $F \circ \pi^{-1}$ is the relative Frobenius for the $\F_p$-structure on $\ul{G}_{\ol{\F}_p}$ induced by the split form of $G_{\F_p}$.

 The automorphism $\pi$ induces an automorphism of $X^*(\ulT) \cong X_*(\ul{\chT})$ which we denote by the same name. It gives rise to an endomorphism
$\varphi=p\pi^{-1}$ of $\wt{W}$, characterized by the formula
\begin{equation}\label{eq:varphi}
\varphi(wt^{\mu})=\pi^{-1}(w)t^{p\pi^{-1}(\mu)}.
\end{equation}

 \subsection{Serre weights}\label{ssec: Serre weights}
 
 Fix an algebraic closure $k = \ol{\F}_p$. The simple representations of $G(\F_p)$ over $k$ are called the \emph{Serre weights} of $\ulG$.  
 
 For $\lambda \in X^*(T)$, write $L(\lambda)$ for the simple representation of $\ul{G}_k$ with highest weight $\lambda$.

 \subsubsection{Parametrization of Serre weights}
Recall from \S \ref{sssec: Frob kernel} that the \emph{$p$-restricted weights} $X_1^*(\ulT) \subset X^*(\ulT)^+$ are defined as 
\[
 X_1^*(\ulT) := \{ \lambda \in X^*(\ulT)^+ \co 0 \leq \tw{\lambda, \alpha^\vee}  < p \text{ for all } \alpha \in \ul{\Delta}\}.
\]
Let
 \[
 X^0(\ulT) := \{ \lambda \in X^*(\ulT) \co \tw{\lambda, \alpha^\vee} = 0 \text{ for all } \alpha \in \ul{\Delta}\}.
 \]
Then $X_1^*(\ulT)$ is a finite union of $X^0(\ulT)$-cosets. For $\lambda \in X_1^*(\ulT)$, we write $F(\lambda) := L(\lambda)|_{G(\F_p)}$. Then $F(\lambda)$ is a simple representation of $\ulG(\F_p)$, and the map $\lambda \mapsto F(\lambda)$ induces a bijection
\[ 
\frac{X_1^*(\ulT)}{(F-1) X^0(\ulT)} \leftrightarrow \{\text{Serre weights of }\ulG\}.
\]

\subsubsection{The Frobenius kernel} As in \S \ref{sssec: Frob kernel}, we write $\ulG_1 := \ker(\ulG_k \xrightarrow{F} \ulG_k^{(1)})$ for the Frobenius kernel of $\ulG$ and $\ulG_1\ulT < \ulG_k$ for the subgroup scheme generated by $\ulG_1$ and $\ulT$. For $\lambda \in X^*(T)$, let $\wh{L}_1(\lambda)$ be the simple representation of $\ulG_1\ulT$ with highest weight $\lambda$. (As a representation of $\ulG_1$, we have $\wh{L}_1(\lambda):= L(\lambda_0)|_{\ulG_1}$ where $\lambda_0$ is any $pX^*(\ulT)$-translate of $\lambda$ into the fundamental box.) Then $\wh{L}_1(\lambda)$ is simple and the map $\lambda \mapsto \wh{L}_1(\lambda)$ induces a bijection 
 \[
 X^*(\ulT) \leftrightarrow  \left\{ \text{simple $\ulG_1\ulT$-representations} \right\}.
 \]

 \subsubsection{Reparametrization of Serre weights}\label{sssec: reparametrization of serre weights} Recall the following notations:
 \begin{itemize}
 \item $\ul{A}_0 \subset X^*(T)$ is the \emph{dominant base alcove} anchored at $0$, and $\ul{C}_0 := - \rho  + p \ul{A}_0$ is the \emph{dominant $\rho$-shifted base $p$-alcove}. 
 \item The \emph{dominant} affine Weyl group elements are 
 \[
\wt{\ulW}^+ := \{ \wt{w} \in \wt{\ulW} \co \wt{w} \bup \ulC_0 \text{ is dominant}\}
\]
and the \emph{dominant $p$-restricted} affine Weyl group elements are 
\[
\wt{\ulW}_1 := \{ \wt{w} \in \wt{\ulW}^+ \co \wt{w} \bup \ulC_0 \text{ is $p$-restricted (and dominant)}\}.
\]
\end{itemize}

We will also use the reparametrization from \cite[(2.5)]{LLLM22} of Serre weights by pairs $\wt{w}_1 \in \wt W_1^+ :=  \wt{W}_1 \cap \wt{W}^+$ and $\omega \in X^*(\ulT)$ such that $\omega-\rho \in \ulC_0$:
\[
F_{(\wt{w}_1, \omega)} := F(\pi^{-1}(\wt{w}_1) \bup (\omega - \rho)).
\]
As explained in \cite[Lemma 2.2.4]{LLLM22}, there are many such parametrizations: the map 
\begin{equation}\label{eq: compatible}
(\wt{w}_1,\omega) \mapsto  \wt{w}_1 t^{\omega} \in  \wt{\ulW}/\ulW_{\aff}
\end{equation}
gives an injection from the set of such parametrizations to the set of central characters of $\ul G$, which is naturally identified with $X^*(Z)=\Omega=\wt{\ulW}/\ulW_{\aff}$ \footnote{This is called an \emph{algebraic central character} in \cite{LLLM22}.}. For the endomorphism $\varphi$ defined in \eqref{eq:varphi}, the image of \eqref{eq: compatible} is a fixed coset under $(\varphi-1)\Omega$. 

\begin{defn}\label{def: compatible pair}
We will say that such a pair $(\wt{w}_1, \omega)$ is \emph{compatible} with an element $\zeta \in \wt{\ulW}/\ulW_{\aff}$ if it maps to $\zeta$ under \eqref{eq: compatible}. 
\end{defn}

 \subsection{Tame inertial $L$-parameters}\label{ssec: tame inertial types} We establish some aspects of the ``tame part'' inertial Local Langlands correspondence for unramified groups. 
 
\begin{defn}\label{defn:tame-inertial-param} For $R\in \{\ol \Q_p,\ol \Z_p,\ol \F_p\}$, an \emph{inertial $L$-parameter over $R$} is a $\bG(R)$-conjugacy class of homomorphisms
\[\tau: I_{\Q_p}\to \bG(R)\]
which can be extended to an $L$-homomorphism $\wt{\tau}:W_{\Q_p}\to \LG$. Such a $\tau$ is called \emph{tame} if $\tau$ factors through the tame quotient of $I_{\Q_p}$.
\end{defn}
\begin{remark} \begin{enumerate}
\item
Our inertial types are actual homomorphism as opposed to $1$-cocycles since we assumed $G$ is unramified. 
\item There is a natural bijection between {\it tame} inertial $L$-parameters over $\ol \Q_p,\ol \Z_p$, and $\ol \F_p$.
\end{enumerate}
\end{remark}

\subsubsection{Combinatorial parametrization} By \cite[Proposition 9.2.1]{GHS18} and \cite[Propositions 5.7, 5.26]{DL76}, after choosing a generator $\psi \in\varprojlim \F^\times_{p^i}$,
inertial $L$-parameters are in natural bijection with the set
\[
((X_*(\bT)\otimes \Q/\Z)/W)^{\varphi=1},
\]
which in turn naturally identifies with 
\[\big(\bigcup_{\wt{w}\in\wt{W}} (X_*(\bT)\otimes \Q)^{\varphi=\wt{w}}\big)/\wt{W}.\]
In the proof of
Lemma \ref{lemma: existence of Kisin lattice},
we make this bijection explicit
by matching $x\in X_*(\bT)\otimes \Q$
with the inertial parameter ``$C_\gamma:\gamma\mapsto v^{-x}\gamma(v^x)$''; 
as this is only used in that Lemma, we leave the clarification to its proof.

\begin{const}\label{const:combinatorial tame type} Let $\wt{w}=t^{\mu}w\in\wt{W}$, where $w\in W$ and $\mu\in X_*(\chT)$. Then the element 
\begin{equation}\label{eq:x_w}
x_{\wt w}:=(\varphi-w)^{-1}(\mu)\in X_*(\chT)\otimes \Q
\end{equation}
is the unique solution to the equation
\[\varphi(x)=\wt{w}x\]
and hence induces an element of $((X_*(\bT)\otimes \Q/\Z)/W)^{\varphi=1}$, which we denote by $\tau(\wt{w}) =\tau(w,\mu)$.
This construction gives rise to a surjection
\[\wt{W}/_\varphi \wt{W} \surj ((X_*(\bT)\otimes \Q/\Z)/W)^{\varphi=1}\]
where on the left side, the quotient is via the $\varphi$-conjugation action $\wt{z}*_{\varphi}\wt{w}=\varphi(\wt{z})\wt{w}\wt{z}^{-1}$.

\end{const}
\begin{defn}\label{defn:presentation} Let $\tau$ be an inertial $L$-parameter. 
\begin{enumerate}
\item We say that $\wt{w}=t^{\mu}w\in\wt{W}$ (or equivalently, $(w,\mu)$) is a \emph{presentation} of $\tau$ if $\tau\ = \tau(\wt{w}) :=\tau(w,\mu)$, i.e., $\tau$ arises from $\wt{w}$ via Construction \ref{const:combinatorial tame type}. 
\item A presentation $\wt{w}=t^{\mu}w$, of $\tau$ induces a point $x_{\wt{w}}=(\varphi-w)^{-1}(\mu)\in X_*(\chT)\otimes \Q $, whose orbit under $\wt{W}$ corresponds to $\tau$.
\item A presentation $\wt{w}$ of $\tau$ is \emph{lowest alcove} if the induced point $x_{\wt{w}}$ is in the lowest alcove.
\item We say an element $\wt{w}=t^{\mu} w\in\wt{W}$ is \emph{$m$-generic} if for each root $\alpha \in \Phi$
\[pn_\alpha+m<\langle \mu,\alpha^{\vee} \rangle<p(n_\alpha +1)-m\]
for some integer $n_\alpha$.
\end{enumerate}
\end{defn}
\begin{remark} We warn the reader that when $\wt{w}=t^{\mu}w$ is a lowest alcove presentation of $\tau$, then we can only conclude that $\mu$ belongs to the closure $\bar C_0$ of the base $p$-alcove, i.e., it is only $(-1)$-generic with respect to the alcove $C_0$. 
\end{remark}

\begin{prop} Let $\wt{w}=t^\mu w \in \wt{W}$.
\begin{enumerate}
\item If $\mathrm{Stab}_{\wt{W}}(x_{\wt{w}})=\{1\}$, then the set of presentations of $\tau(\wt{w})$ coincides with the set $\wt{W}*_\varphi \wt{w}$ of $\varphi$-conjugates of $\wt{w}$.  
\item If $\mu$ is $1$-generic in the base $p$-alcove $C_0$, then $x_{\wt{w}}$ is in the base alcove $A_0$, so that $\wt{w}$ is a lowest alcove presentation of $\tau=\tau(\wt{w})$. If this is the case, then the set of lowest alcove presentations of $\tau$ is
\[\{\varphi(\wt{z})\wt{w}\wt{z}^{-1}\, \co  \wt{z} \in \wt{W} \text{ such that } \wt{z}(A_0)=A_0\}.\]
\end{enumerate}
\end{prop}
\begin{proof}(1) Suppose $\wt{z}$ is another presentation of $\tau(\wt{w})$, meaning there is $\wt{u}\in\wt{W}$ so that $x_{\wt{z}}=\wt{u}x_{\wt{w}}$ (recall the definition of $x_{\wt w}$ from \eqref{eq:x_w}). But then
\[\wt{w}x_{\wt{w}}=\varphi(x_{\wt{w}})=\varphi(\wt{u}^{-1})\varphi(x_{\wt{z}})=\varphi(\wt{u}^{-1})\wt{z}x_{\wt{z}}=\varphi(\wt{u}^{-1})\wt{z}\wt{u}x_{\wt{w}}\]
hence $\wt{z}=\wt{u}*_\varphi \wt{w}$.

(2) We have $\mu+wx_{\wt{w}}=p\pi^{-1} x_{\wt{w}}$. This gives
\[ph_{x_{\wt{w}}}\leq h_\mu+h_{x_{\wt{w}}}<p-1+h_{x_{\wt{w}}}\]
Hence $h_{x_{\wt{w}}}<1$. 
In turn this implies that for any $\alpha>0$ 
\[p\langle x_{\wt{w}},\alpha^{\vee}\rangle=\langle \mu+wx_{\wt{w}},\alpha^{\vee}\rangle>0\]
so that indeed $x_{\wt{w}}\in A_0$.

Now, any element $\wt{z}\in\wt{W}$ which stabilizes $x_{\wt{w}}\in A_0$ must stablize $A_0$, and hence $\wt{z}\in \Omega$. However, projecting $x_{\wt{w}}$ to $ X_*(\chT) / \chQ^\vee=X^*(Z)\cong \Omega$, we see that $\wt{z}$ must be trivial. The last conclusion now follows from the computation in part (1).
\end{proof}

In particular, similar to Definition \ref{def: compatible pair}, we get a notion of compatibility for presentations of $1$-generic inertial parameters:
\begin{defn} \label{def: compatible presentation} Let $\tau$ be a tame inertial $L$-parameter that admits at least one $1$-generic lowest alcove presentation. Then the map
\[\wt{w}\mapsto  \wt{w}\ulW_{\aff} \in  \wt{\ulW}/\ulW_{\aff}=\Omega\]
gives a bijection between the set of lowest alcove presentations of $\tau$ and a coset of $(\varphi-1)\Omega$.

We say a lowest alcove presentation $\wt{w}$ of $\tau$ is \emph{compatible} with $\zeta\in \Omega$ if the image of $\wt{w}$ in the above injection is $\zeta$. 
\end{defn}

\subsubsection{Tame inertial types} We assume for initial discussion that $G = \GL_n$. Then to a Weil-Deligne inertial parameter $\tau$ the inertial Local Langlands correspondence associates a smooth irreducible $G(\Z_p)$-representation $\sigma(\tau)$ called the \emph{inertial type} of $\tau$, with properties explained in \cite[Theorem 2.5.4]{LLLM22}. 

We focus on the case where $\tau$ is tame, and $\sigma(\tau)$ is then called a \emph{tame inertial type}. Let $(w, \mu) \in \ulW \times X^*(\ulT)$ be a \emph{good pair} in the sense of \cite[\S 2.2]{LLL19}. Then we have a Deligne-Lusztig representation $R(w,\mu)$ \cite[Proposition 9.2.1 and 9.2.2]{GHS18}. For $1$-generic $\mu$, we can take $R(w, \mu) = \sigma(\tau(w, \mu))$,  according to \cite[Proposition 2.5.5]{LLLM22}. 

Now we return to the general case where $G$ is merely assumed to be unramified (and satisfy Hypothesis \ref{hyp:part3-hypotheses}). Then we do not have the inertial Local Langlands correspondence in general, but we may still define the map $\tau \mapsto \sigma(\tau)$ for mildly generic tame $\tau$, because the ``tame part'' of the Local Langlands correspondence is understood for such $G$. If $\wt{w} = t^\mu w$ is a presentation for $\tau$ in the sense of Definition \ref{defn:presentation} (so that $\tau = \tau(w, \mu)$), and the Deligne--Lusztig representation $R(w,\mu)$ is irreducible (for example, when $\wt{w}$ is $1$-generic), then we define $\sigma(\tau)$ to be $R(w, \mu)$.
\begin{remark} One also has a notion of Weil--Deligne inertial parameter which is an enhancement of the notion of tame inertial parameter with a nilpotent operator as in \cite[Definition 2.5.1]{LLLM22}. Then for inertial Weil--Deligne parameter $\wt{\tau}$ whose inertial parameter is $\tau(w,\mu)$, one can define $\sigma(\wt{\tau})$ to be certain irreducible factors of $R(w,\mu)$ when the latter is not irreducible. This generalizes the above recipe to situations where $\tau$ is tame but not sufficiently generic. Since our main results do not concern this situation, we do not elaborate further on the precise recipe.
\end{remark}

\subsection{The Emerton-Gee stack}\label{ssec: EG stack}
Emerton-Gee have constructed a moduli stack over $\Spf \cO$ of $(\varphi, \Gamma)$-modules, whose groupoid of points in any finite $\cO$-algebra $A$ is naturally equivalent to the groupoid of continuous $A$-valued representations of $G_{\Q_p}=\Gal(\ol{\Q}_p/\Q_p)$. For the $L$-group $\LG$ of a tamely ramified $G/\Q_p$, \cite{Lin23b} generalizes this construction to produce a corresponding moduli stack $\cX^{\EG}_{\LG}$ of representations valued in $\LG$. 

\subsubsection{Potentially crystalline substacks} For a pair $(\lambda, \tau)$ where $\lambda \in X_*(\chT)$ and $\tau$ is an inertial $L$-parameter, there is a substack $\cX^{\lambda, \tau}_{\LG} \subset \cX^{\EG}_{\LG}$ which is finite type and flat over $\Spf \Z_p$, and then determined by the property that for any finite flat $\cO$-algebra $A^\circ$, $\cX^{\lambda, \tau}_{\LG}(A^\circ)$ is the groupoid of $A^\circ$-lattices in potentially crystalline $G_{\Q_p}$-representations with Hodge-Tate weights $\lambda$ and inertial parameter $\tau$. This construction is due to Emerton--Gee for general linear groups, and is carried out in \cite[\S 2]{Lin23b} for more general $L$-groups (a brief discussion of the cases relevant to this paper occurs in Appendix \ref{app: EG stacks}). When $\lambda$ is regular, the special fiber $\cX^{\lambda, \tau}_{\LG, \F}$ is an algebraic stack of dimension $\dim \chG/\chB$; when $\lambda$ is irregular, the dimension of its special fiber is strictly smaller.

\subsubsection{Irreducible components}\label{sssec:irred-cpnts} 
When $G = \Res_{K/\Q_p} \GL_n$, recall \cite[Theorem 1.2.1, Theorem 6.5.1]{EG23} that $\cX_{\LG}^{\EG})_{\red}$ is equidimensional of dimension $[K:\Q_p] n(n-1)/2 = \dim \check\cB$ (the flag variety of $\chG$), and its irreducible components are in bijection with Serre weights $\sigma$.

\begin{example} When $\sigma=F(\kappa-\rho)$ with $2$-generic $\kappa\in X_1^*(\ulT)$, the irreducible component $\cC_\sigma$ corresponding to $\sigma$ is characterized as the scheme theoretic image of a monomorphism
\[\cX_{\lsup LB,\F}^{\kappa,\mathrm{mns}}\to \cX_{\lsup LG,\F}^{\EG}\]
where $\cX_{\lsup LB,\F}^{\kappa,\mathrm{mns}}$ is the stack parametrizing parameters $\ol \rho: G_{\Q_p}\to \lsup LB(\ol \F)$ such that
\begin{itemize}
\item The composite $\ol \rho^{\semis}:I_{\Q_p}\to \lsup LB(\ol \F) \to \lsup LT(\ol \F)$ corresponds to the tame inertial $L$-parameter $\tau(1,\kappa)$, or in other words, the semisimplification $\ol \rho^{\semis}$ has inertial weight $\kappa$.
\item $\ol \rho$ is ``maximally non-split'' in the sense that it does not factor through $(\lsup LB\cap \lsup LB')(\ol \F)$ for another Borel $B'\neq B$.
\end{itemize}
\end{example}

When $G$ is an unramified group, \cite{Lin23b} studies the structure of the irreducible components of $(\cX_{\LG}^{\EG})_{\red}$. It is expected, but only established when $\chG$ has no simple factors of type E or F (in a forthcoming work of Zhongyipan Lin), that $(\cX_{\LG}^{\EG})_{\red}$ is equidimensional of dimension $d=\dim \check{\cB}$ (the flag variety of $\chG$) and admits a natural decomposition into a union of $d$-dimensional closed substacks $\cC_\sigma$ parameterized by Serre weights $\sigma$ of $G(\F_p)$. 
In general, for $\sigma=F(\kappa-\rho)$ with $2$-generic $\kappa\in X_1^*(\ulT)$, \cite{Lin23b} constructs the expected closed substack $\cC_\sigma$ as the scheme theoretic images of the $\lsup LB$-valued family $\cX_{\lsup LB,\F}^{\kappa}\to  \cX_{\lsup LG,\F}^{\EG}$ whose semisimplification has inertial weight $\kappa$ in the above sense.


\subsubsection{The Breuil--M\'ezard Conjecture}
For a Serre weight $\sigma$, let $W(\lambda)$ be the Weyl module with highest weight $\lambda$ and 
\begin{equation}\label{eq:n-sigma-lambda-tau}
n_{\sigma}(\lambda, \tau) := [\ol{W(\lambda) \otimes \sigma(\tau) }: \sigma]
\end{equation}
be the Jordan-H\"{o}lder multiplicity of $\sigma$ in a mod $p$ reduction of $W(\lambda) \otimes \sigma(\tau)$. Recall that $d=\dim \check{\cB}$ (the flag variety of $\chG$). The \emph{geometric Breuil--M\'ezard Conjecture} due to Emerton-Gee, predicts: 

\begin{conj}[{cf. \cite[Conjecture 8.2.2]{EG23}}]
Let $G = \GL_n$. There is a collection of effective cycles $\cZ(\sigma) \in \dZ((\cX_{\LG}^{\EG})_{\red})$ indexed by Serre weights $\sigma$ of $G(\F_p)$, such that for all $\lambda \in X^*(\ulT)^+$ and inertial parameters $\tau$, we have 
\[
[\cX^{\lambda + \rho, \tau}_{\LG, \F}] = \sum_\sigma n_{\sigma}(\lambda, \tau) \cZ(\sigma) \in \dZ((\cX_{\LG}^{\EG})_{\red}).
\]
\end{conj}

For general unramified groups $G/\Q_p$ satisfying Hypothesis \ref{hyp:part3-hypotheses}, there should ideally be an analogous conjecture, but its formulation is currently problematic since we do not have the (inertial) Local Langlands correspondence in such generality. Nevertheless, we can formulate the ``tame'' part of the Breuil--M\'ezard Conjecture, cf. \cite[Theorem D, Conjecture 1]{Lin23c}.

\begin{conj}[Tame Breuil-M\'{e}zard Conjecture]\label{conj:tame-BM-Conj}
Let $G$ be an unramified group over $\Q_p$ satisfying Hypothesis \ref{hyp:part3-hypotheses}. Then there is a collection of effective cycles $\cZ(\sigma) \in \dZ((\cX_{\LG}^{\EG})_{\red})$ indexed by Serre weights $\sigma$ of $G(\F_p)$, such that for every $\lambda \in X^*(T)^+$ and every tame inertial parameter $\tau$, we have
\[
[\cX^{\lambda + \rho, \tau}_{\LG, \F}] = \sum_\sigma n_{\sigma}(\lambda, \tau) \cZ(\sigma) \in \dZ((\cX_{\LG}^{\EG})_{\red}).
\]
where $[\cX^{\lambda+\rho, \tau}_{\LG, \F}]$ is the cycle class of $\cX_{\LG, \F}^{\lambda+\rho, \tau}$. 
\end{conj}

\begin{remark}
The effectivity is not part of the original formulations in \cite{BM02, EG23}, but it has become part of the conjectural picture relating the Breuil--M\'ezard Conjecture to patching, and is conjectured explicitly in \cite[\S 9.4]{EGSurvey}. 
\end{remark}

\section{Existence of Breuil-M\'{e}zard cycles}\label{sec: existence of BM cycles}

Fix an unramified group $G$ satisfying hypothesis \ref{hyp:part3-hypotheses}, and abbreviate $\cX^{\EG} := \cX^{\EG}_{\LG}$, $\cX^{\lambda, \tau} := \cX^{\lambda, \tau}_{\LG}$. In this section we will prove the following Theorem, which gives partial evidence to Conjecture \ref{conj:tame-BM-Conj} 

\begin{thm}\label{thm: BM 13}There exists a collection $\{ \cZ^{\EG}(\sigma) \in \dZ(\cX^{\EG}_{\red}) \}$, indexed by Serre weights $\sigma$ for $G(\F_p)$, such that for all $\lambda \in X^*(T)^+$ and $\tau = \tau(w, \mu)$ with $\mu \in \ulC_0$ being $2h_{\lambda+\rho}$-generic, we have 
\begin{equation}\label{eq: BM relation}
[\cX^{\lambda+\rho, \tau}_{\F}] = \sum_{\sigma} n_{\sigma}(\lambda, \tau)  \cZ^{\EG}(\sigma),
\end{equation}
where $n_\sigma(\lambda, \tau)$ is defined as in Conjecture \ref{conj:tame-BM-Conj}. 
\end{thm}

\begin{remark}
It will not be immediate from the construction whether the cycles $\cZ^{\EG}(\sigma)$ are effective, or whether they are uniquely characterized by equations \eqref{eq: BM relation}. However, they will satisfy a weaker ``bounded support'' property, and we show in \S \ref{sec: complements} that this property together with equations \eqref{eq: BM relation} does characterize them uniquely.
\end{remark}

Here is the outline of this section. In \S \ref{ssec: X_R} we discuss the relationship between the deformed affine Springer fibers $\ul{\cX}_\gamma^{\varepsilon=1}$ and the potentially crystalline substacks $\cX^{\lambda, \tau}$ of $\cX^{\EG}$. Roughly speaking, we cook up $\gamma$ from $\tau$ so that $\ul{\cX}_\gamma^{\varepsilon=1}(\leq \lambda)$ is a mod $p$ local model for $\bigcup_{\lambda' \leq \lambda} \cX^{\lambda', \tau}$. In \S \ref{ssec: gluing} we discuss the comparison of the special fibers $\ul{\rY}_\gamma^{\varepsilon=1}$ as $\gamma$ varies. In \S \ref{ssec: BM cycles on local model} we construct the incarnations $\cZ_\gamma^{\varepsilon=1}(\sigma)$ of Breuil--M\'ezard cycles on the model spaces $\ul{\rY}_\gamma^{\varepsilon=1}$, and in \S \ref{ssec: BM relations on local model} we prove that they satisfy relations which correspond to \eqref{eq: BM relation}. Finally, in \S \ref{ssec: BM relations} we construct the cycles $\cZ^{\EG}(\sigma)$ of Theorem \ref{thm: BM 13}, and deduce that they satisfy the relations \eqref{eq: BM relation}. 

For proving Theorem \ref{thm: BM 13}, we may and do fix once and for all a central character $\zeta\in \wt{\ulW}/\ulW_{\aff}$ of $\ul G$. The significance of this choice is two-fold: First, it determines a central character of $\ulG(\F_p)$, which geometrically means we work on an open-closed substack of $\cX^{\EG}$. (We also note that the choice of $\zeta$ also corresponds to a choice of connected component for the affine flag variety for $\ul\chG$.)  Second, it pins down consistently the choice of lowest alcove presentations for the tame inertial types $\tau$ and parametrizations of Serre weights that we work with, which feeds into the local model theorems relating parts of $\cX^{\EG}$ with deformed affine Springer fibers. 

\subsection{Modeling potentially crystalline loci}\label{ssec: X_R} Recall the endomorphism $\varphi$ of $\wt W$ defined in \eqref{eq:varphi}. 

\begin{defn}
Suppose $\tau(w, \mu)$ is a lowest alcove presentation of a tame inertial parameter $\tau$. Then we define 
\[
\gamma(w, \mu) = -(t+p)x,
\]
where $x=(\varphi-w)^{-1}\mu\in \chft$ is the unique solution to $\varphi(x)=\mu+w(x)$. (In the notation of Construction \ref{const:combinatorial tame type}, this $x$ is $x_{\wt{w}}$ for $\wt w = t^\mu w$.) 
\end{defn}

\begin{remark}
It follows from the definitions that $\gamma(w,\mu)\equiv tw^{-1}(\mu)$ mod $p$, and in particular $x\in\chft $ is regular semisimple if $\mu$ is $0$-generic.
\end{remark}

\subsubsection{The homological model theorem} The following theorem is the main $p$-adic Hodge theoretic input that we need. It shows that one can model the homology of certain unions of potentially crystalline loci with regular Hodge-Tate weights in terms of suitable deformed affine Springer fibers: 

\begin{thm}[Homological model theorem]\label{thm: companion} Let $\lambda\in X^*(\ulT)^+$ and $\tau = \tau(w, \mu)$ be a tame inertial parameter with lowest alcove presentation $\wt{w} := t^\mu w$. 
Suppose $\mu$ is $2h_{\lambda+\rho}$-generic, and set $\gamma :=\gamma(w,\mu)$. Then for $d := \dim \check{\cB}$ (the flag variety of $\chG$), there is an injection 
\[
\transfer_\gamma \co  \dZ(\ulY_\gamma^{\varepsilon=1}(\leq \lambda+\rho)) \inj \dZ(\cX^{\EG}_{\red})
\]
with the following properties. 
\begin{enumerate}
\item For each dominant $\lambda'\leq \lambda+\rho$,
\[
\transfer_\gamma  \left( \fsp_{p \rightarrow 0} [\ulX_\gamma^{\varepsilon =1}(\lambda')] \right) = [\cX^{\lambda', \tau}_{\F}] \in  \dZ(\cX^{\EG}_{\red}).
\]
\item Let $\sigma\in \JH( \ol{W(\lambda)\otimes R(w, \mu)})$. There is a unique pair $(\wt{u},\wt{v})\in \ul{\wt{W}}_1\times \ul{\wt{W}}^+$ with $\wt{w}=\wt{v}^{-1}w_0\wt{u}\in \Adm^{\reg}(\lambda+\rho)$ such that $\sigma=F(\pi^{-1}(\wt{u})\bup (t^\mu w \wt{v}^{-1}(0) - \rho))$. Then
\[
\transfer_\gamma \left( [\ulY^{\varepsilon=1}_{\gamma}(\wt{w})] \right)  =  [\cC_{\sigma}] \in  \dZ(\cX^{\EG}_{\red}),
\]
where we recall that $\cC_\sigma$ is the irreducible component of $\cX^{\EG}_{\red}$ corresponding to $\sigma$ (\S \ref{sssec:irred-cpnts}).
\end{enumerate}
\end{thm}

\begin{remark}
The proof appears in Appendix \ref{app: B}, and shows moreover that under the assumptions of the theorem, 
\begin{equation}\label{eq:cycle-contain}
\transfer_\gamma \Big(\dZ \Big(\ulY_\gamma^{\varepsilon=1}(\leq \lambda') \Big)\Big) \supset Z_d(\cX^{\lambda', \tau}_{\F}) \quad \text{as subspaces of $\dZ(\cX^{\EG}_{\red})$}.
\end{equation}
\end{remark}

\begin{remark} 
Suppose $G=\Res_{K/\Q_p} \GL_n$ for some unramified $K$. Then:
\begin{itemize} \item When $\tau$ is \emph{very generic} relative to $\lambda$ in the sense that $\mu$ avoids some universal closed subvariety in $X^*(\chT)\otimes_{\Z} \,\F_p$ depending on $\lambda$, the Theorem \ref{thm: companion}(1) can be deduced from \cite[Theorem 7.3.2]{LLLM22}. However, the inexplicit nature of this very generic condition prevents us from making this deduction when $\lambda$ has a dependency on $p$.
\item
When $\mu$ is around $\max\{h_{\lambda+\rho}+4h_\rho,2h_{\lambda+\rho}\}$-generic, \cite[Theorem 7.4.2]{LLLM22} shows that the inclusion \eqref{eq:cycle-contain} is an equality. Using this as input, one can show this also holds under our weaker genericity hypothesis by a degeneration argument.
\end{itemize}
We remark that the arguments of \cite{LLLM22} make essential use of Taylor--Wiles patching, and in particular depend on global inputs which do not extend to general groups.


\end{remark}

\subsubsection{Strategy for constructing Breuil--M\'ezard cycles}\label{sssec:bm-cycle-strat} We want to define candidate Breuil--M\'ezard cycles $\cZ^{\EG}(\sigma)$ on $\cX^{\EG}$, as $\sigma$ ranges over Serre weights. We will do this by choosing some auxiliary tame type $\tau$ that contains $\sigma$ as a Jordan-H\"older factor, constructing a cycle $\cZ_{\gamma}^{\varepsilon = 1}(\sigma)$ on the model $\ulY_\gamma^{\varepsilon = 1}$, and then taking 
\[
\cZ^{\EG}(\sigma) = \transfer_\gamma \Big(\cZ_{\gamma}^{\varepsilon = 1}(\sigma) \Big).
\]

\begin{reminder}We recall that for each $\varepsilon_0 \in \A^1_{\varepsilon}$, $d$ is the (equi)dimension of $\ulY_\gamma^{\varepsilon=1}$ and $\ulY_\gamma^{\varepsilon=0}$, so that $\topCh(-) = \dZ(-)$ for these spaces. This is also expected to be true for $\cX_{\red}^{\EG}$, but not known in the exceptional types E and F. 
\end{reminder}

In turn, the cycle $\cZ_{\gamma}^{\varepsilon = 1}(\sigma)$ will be constructed by applying the inverse of $\fsp_{\varepsilon \rightarrow 0} \co (\ulY_\gamma^{\varepsilon=1}) \rightarrow \topCh(\ulY_\gamma^{\varepsilon=0})$ to the microlocal support of simple representations of $\ulG_1\ulT$.

\subsubsection{Change of types} For this strategy to work, we need to know that the recipe produces cycles that are independent of the choice of the auxiliary tame type $\tau$ (or equivalently, $\gamma$). This requires understanding how the intersection of $\cX^{\lambda, \tau}_{\F}$ and $\cX^{\lambda', \tau'}_{\F}$ behaves under the two transfer maps. We refer to the explicit answer as a ``change-of-type'' formula. 

Write $\delta := \mu - \mu'$ and let $\gamma := \gamma(w, \mu)$ and $\gamma' := \gamma(w',\mu')$. 
By Lemma \ref{lem:transformation of components}(1), left translation by $(w')^{-1} t^{\delta}w$ maps $\ulY^{\varepsilon = 1}_{\gamma}$ isomorphically to $\ulY_{\gamma'}^{\varepsilon =1}$ as subschemes of $\Fl_{\ul{\chG}}$, 
\begin{equation}\label{eq: translate type}
\tr_{(w')^{-1} t^{\delta}w}  \co \ulY_\gamma^{\varepsilon =1} \xrightarrow{\sim} \ulY_{\gamma'}^{\varepsilon =1}.
\end{equation}
This induces an isomorphism 
\begin{equation}\label{eq: translate DASF homology}
\tr_{(w')^{-1} t^{\delta}w} \co \topCh(\ulY_{\gamma}^{\varepsilon=1}) \xrightarrow{\sim} \topCh(\ulY_{\gamma'}^{\varepsilon=1})
\end{equation}
and similarly on Borel--Moore homology. Given $\lambda, \lambda' \in X^*(\ulT)^+$, \eqref{eq: translate DASF homology} induces a partially defined map 
\begin{equation}\label{eq: partial transfer map}
\begin{tikzcd}[column sep = huge]
\topCh(\ulY_{\gamma}^{\varepsilon=1}(\leq \lambda)) \ar[r, "{\tr_{(w')^{-1} t^{\delta}w}}", dashed] & \topCh(\ulY_{\gamma'}^{\varepsilon=1}(\leq \lambda'))
\end{tikzcd}
\end{equation}
where the domain of definition is generated by the classes of top-dimensional irreducible components contained in $\ulY_\gamma^{\varepsilon =1}(\leq \lambda)$ which are mapped by $\tr_{(w')^{-1} t^{\delta}w}$ to top-dimensional irreducible components lying in $\ulY_{\gamma'}^{\varepsilon=1}(\leq \lambda')$.

\begin{prop}[Change-of-type formula]\label{prop: companion compatibility}
Let $\lambda \in X^*(\ulT)^+$,  $\lambda' \in X^*(\ulT)^+$, and $\delta := \mu -  \mu'$. Let $\tau = \tau(w, \mu)$, and $\tau' = \tau(w', \mu')$ be two lowest alcove presentations compatible with $\zeta-\lambda$ and $\zeta-\lambda'$, respectively. Assume that $\mu$ is $2h_{\lambda+\rho}$-generic and $\mu'$ is $2h_{\lambda'+\rho}$-generic. Then the diagram 
\[
\begin{tikzcd}
\dZ(\ulY_\gamma^{\varepsilon=1}(\leq \lambda+\rho)) \ar[rr, dashed, "\tr_{(w')^{-1} t^\delta w}"] \ar[dr, "\transfer_\gamma"'] & & 
\dZ(\ulY_{\gamma'}^{\varepsilon=1}(\leq \lambda'+\rho))  \ar[dl, "\transfer_{\gamma'}"]  \\
& \dZ(\cX^{\EG}_{\red}) 
\end{tikzcd}
\]
commutes.
\end{prop}
\begin{proof} From Theorem \ref{thm: companion} we have 
\[
\transfer_\gamma \left( [\ulY^{\varepsilon=1}_{\gamma}(\wt{w})] \right)  =  [\cC_{\sigma}]
\]
where $\wt{w}=\wt{v}^{-1}w_0\wt{u}\in \Adm^{\reg}(\lambda+\rho)$ and $\sigma=F(\pi^{-1}(\wt{u})\bup (t^\mu w \wt{v}^{-1}(0) - \rho))$.
The condition that $\tr_{(w')^{-1} t^\delta w}[\ulY^{\varepsilon=1}_{\gamma}(\wt{w})]$ is defined guarantees that $\sigma\in \JH( \ol{W(\lambda')\otimes R(w, \mu')})$. By Theorem \ref{thm: companion}(2), there is $\wt{w}'=(\wt{v}')^{-1}w_0\wt{u}'\in \Adm^{\reg}(\lambda'+\rho)$ such that 
\[
\transfer_\gamma \left( [\ulY^{\varepsilon=1}_{\gamma'}(\wt{w}')] \right)  =  [\cC_{\sigma}],
\]
hence $F(\pi^{-1}(\wt{u})\bup (t^\mu w \wt{v}^{-1}(0) - \rho))=F(\pi^{-1}(\wt{u}')\bup (t^{\mu'} w' (\wt{v}')^{-1}(0) - \rho))$. By our compatibility assumptions on $(w,\mu)$ and $(w',\mu')$, the two parametrizations of the Serre weight coincide. In particular, we have $\wt{u}=\wt{u}'$, hence 
\[t^\mu w \wt{v}^{-1}(0)=t^{\mu'} w' (\wt{v}')^{-1}(0)\]
so that 
\[t^\mu w \wt{v}^{-1}=t^{\mu'} w' (\wt{v}')^{-1}s\]
for some $s\in \ul{W}$.
But by Lemma \ref{lem:transformation of components}, for an appropriate $\gamma_1$ we have 
\[t^{\mu}w \ulY^{\varepsilon=1}_{\gamma}(\wt{w})=t^{\mu}w\wt{v}^{-1}\ulY^{\varepsilon=1}_{\gamma_1}(w_0\wt{u})=t^{\mu'} w' (\wt{v}')^{-1}s\ulY^{\varepsilon=1}_{\gamma_1}(w_0\wt{u})=t^{\mu'} w' (\wt{v}')^{-1}\ulY^{\varepsilon=1}_{\gamma_1}(w_0\wt{u})=t^{\mu'} w' \ulY^{\varepsilon=1}_{\gamma'}(\wt{w}')\]
so that 
\[\tr_{(w')^{-1} t^\delta w}[\ulY^{\varepsilon=1}_{\gamma}(\wt{w})]=[\ulY^{\varepsilon=1}_{\gamma'}(\wt{w}')]\]
as required.
\end{proof}

\subsection{Gluing homology}\label{ssec: gluing} As has been discussed, the spaces $\ulY_\gamma^{\varepsilon=1}$ contain models for the special fibers $\bigcup_\lambda \cX^{\lambda, \tau}_{\F} \subset \cX^{\EG}_{\F}$. As $\tau$ varies, the special fibers of potentially crystalline substacks overlap. In this section, we will see how to implement the corresponding overlapping for the irreducible components of $\ulY_\gamma^{\varepsilon=1}$.

Recall that in Example \ref{ex: GKM description for Y}, we used the GKM description to fix (in particular) an identification
\begin{equation}\label{eq: identify equivariant types e = 0}
\begin{tikzcd}[column sep = huge]
\topBMTg(\ulY_{\gamma}^{\varepsilon = 0})  \ar[r, "{\sim}"', "\mrm{GKM}"]  & \topBMTg(\ulY_{\gamma'}^{\varepsilon = 0}) 
\end{tikzcd}
\end{equation}
characterized by the commutativity of the diagram (where $\sph := \sph_{\ul{\chT}}$)
\[
\begin{tikzcd}
\topBMulT(\ulY_{\gamma}^{\varepsilon = 0})  \ar[r, "{\sim}"', "\mrm{GKM}"] \ar[d, hook] &   \topBMulT(\ulY_{\gamma'}^{\varepsilon = 0})  \ar[d, hook]   \\
\Frac(\sph) \otimes_{\sph} \topBMulT(\ulY_{\gamma}^{\varepsilon = 0})  \ar[r, "{\sim}"]   &  \Frac(\sph) \otimes_{\sph}  \topBMulT(\ulY_{\gamma'}^{\varepsilon = 0})  \\
\Frac(\sph) \otimes_{\sph} \topBMulT((\ulY_{\gamma}^{\varepsilon = 0})^{\chT})  \ar[r, equals]  \ar[u, "{\sim}"] \ar[d, equals]   & \Frac(\sph) \otimes_{\sph}  \topBMulT((\ulY_{\gamma'}^{\varepsilon = 0})^{\chT})   \ar[u, "{\sim}"]  \ar[d, equals] \\
\bigoplus_{\wt{w} \in \wt{\ulW}} \Frac(\sph)[\wt{w}] \ar[r, equals] & \bigoplus_{\wt{w} \in \wt{\ulW}} \Frac( \sph)[\wt{w}]
\end{tikzcd}
\]
By the equivariant formality of $\ulY_\gamma^{\varepsilon = 0}$ and $\ulY_{\gamma'}^{\varepsilon = 0}$ and Remark \ref{remark: de-equivariant BM homology}, \eqref{eq: identify equivariant types e = 0} induces in turn an identification of the \emph{non-equivariant} top Borel--Moore homology groups, 
\begin{equation}\label{eq: identify types e = 0}
\topBM(\ulY_{\gamma}^{\varepsilon = 0}) =   \topBM(\ulY_{\gamma'}^{\varepsilon = 0}) .
\end{equation}

\begin{defn}[Renormalized specialization]\label{defn: renormalized specialization} Suppose $\gamma=\gamma(w,\mu)$ where $\mu$ is $h_\lambda$-generic and $t^\mu w$ is compatible with $\zeta-\lambda$. Then we define the \emph{renormalized specialization} map
\[
\fsp_{\gamma}^{\mrm{ren}} \co \topCh(\ulY^{\varepsilon = 1}_{\gamma} (\leq \lambda)) \inj \topCh(\ulY^{\varepsilon = 0}_{\gamma})
\]
to be the composition of $\fsp_{\varepsilon \rightarrow 0}$ in the sense of Definition \ref{def: specialization 1 to 0} followed by the monodromy-translation action of $t^{\mu} w$. 
\end{defn}

\begin{figure}
  \centering
\begin{tikzpicture}
  \draw (0,4) -- (4,4) -- (2,0) -- cycle;
  
  \draw[dashed] (2.5,2.7) circle (0.7) node[font=\small]  {$\gamma$};
  \draw[dashed] (1.7,2.2) circle (0.6)  node[font=\small] {$\gamma'$};
  \draw[dashed] (2.2,1.7) circle (0.5)  node[font=\small] {$\gamma''$};

  \draw (7,3) circle (0.7) node[font=\small]  {$\rY_\gamma^{\varepsilon=1}$};
  \draw (5,2) circle (0.6)  node[font=\small] {$\rY_{\gamma'}^{\varepsilon=1}$};
  \draw (6,1) circle (0.5) node[font=\small]  {$\rY_{\gamma''}^{\varepsilon=1}$};

  \draw[<-, >=stealth] (3,3) -- node[above, font=\small] {$\transfer_\gamma$} (6.3,3);
  \draw[<-, >=stealth] (2.3,2.2) -- (4.4,2);
  \draw[<-, >=stealth] (2.5,1.5) --node[below, font=\small] {$\transfer_{\gamma''}$} (5.5,1);

  \draw[dashed] (10.5,2.7) circle (0.7) node[font=\small]  {$\gamma$};
  \draw[dashed] (9.7,2.2) circle (0.6)  node[font=\small] {$\gamma'$};
  \draw[dashed] (10.2,1.7) circle (0.5)  node[font=\small] {$\gamma''$};

  \draw[<-, >=stealth] (10,3) --node[above, font=\small] {$\fsp_\gamma^{\mrm{ren}}$} (7.7,3);
  \draw[<-, >=stealth] (9.2,2.2) -- (5.6,2);
  \draw[<-, >=stealth] (9.9,1.5) --node[below, font=\small] {$\fsp_{\gamma''}^{\mrm{ren}}$} (6.5,1);

  \node[font=\small, above] at (2,4.5) {Emerton-Gee stack};
  
  \def\numColumns{2}
  \def\numRows{3}
  
  \def\columnSpacing{2.0} 
  \def\rowSpacing{1.0}    
  
  \foreach \col in {1,2}{
    \coordinate (Label\col) at ($(2*\col*\columnSpacing + 2, 4.5)$);
    \node[font=\small, above] at (Label\col) {%
      \ifnum\col=1
        Local models
      \else
        Affine Springer fiber
      \fi
    };
  }
\end{tikzpicture}

  \caption{{\small This cartoon (produced with the aid of ChatGPT after much coaxing) depicts the construction of Breuil--M\'ezard cycles. The triangle represents the $p$-dilated fundamental alcove $\ulC_0$, whose weights are the $\mu$ in the lowest-alcove presentation $R(w, \mu)$ for tame types; the choice of $w$ and $\mu$ together is captured by $\gamma = \gamma(w, \mu)$. The models $\rY_\gamma^{\varepsilon=1}(\leq \lambda)$ provide charts for corresponding portions of $\cX^{\EG}_{\red}$; the sizes of the charts, depicted by the dashed circles, depend on the genericity of $\gamma$, which is measured by the distance from the walls. Cycles are initially produced on the affine Springer fiber via microlocal support, deformed to the models at $\varepsilon = 1$, and then transferred to the Emerton-Gee stack. Well-definedness of the cycles comes from the fact that $\transfer_\gamma$ and $\fsp_{\gamma}^{\mrm{ren}}$ effect the ``same amount'' of gluing. }}
  \label{fig: bm cycle}
\end{figure}
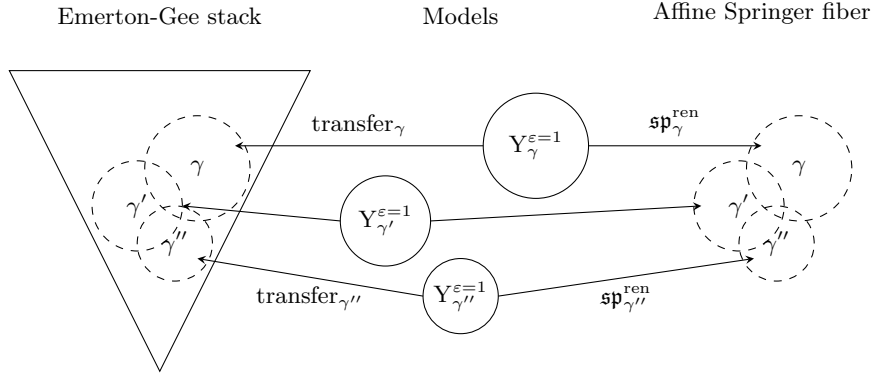

The significance of the renormalized specialization map is seen in the Proposition below, whose statement should be compared to Proposition \ref{prop: companion compatibility}, and interpreted as saying that 
\begin{quote}
``the maps $\fsp_{\gamma}^{\mrm{ren}}$ implement the same combinatorial gluing relations among $\{\topBM(\rY_\gamma^{\varepsilon=1})\}_{\gamma}$ as the maps $\transfer_\gamma$''.
\end{quote}
See Figure \ref{fig: bm cycle} for a visual depiction of this slogan. 

\begin{lemma}\label{lem: e=1 equivariant formality}
If $\gamma$ is $h_\lambda$-generic, then the spaces $\ulY_\gamma^{\varepsilon = \varepsilon_0}(\leq \lambda)$ are equivariantly formal for $\varepsilon_0   = 1$ or $\eta$. 
\end{lemma}

\begin{proof}
According to Lemma \ref{lem: Y lambda irred e=0}, the hypotheses imply that $\ulY_\gamma^{\varepsilon = 1}(\leq \lambda)$ has an affine paving, and therefore that its cohomology is pure. This implies the degeneration of the Grothendieck spectral sequence for $H$-fixed points, which implies equivariant formality. 
\end{proof}

\begin{prop}\label{prop: type change of cycles} Let $\gamma = \gamma(w, \mu)$ and $\gamma' = \gamma(w', \mu')$. Assume that $\gamma$ is $h_\lambda$-generic and $\gamma'$ is $h_{\lambda'}$-generic. Then the diagram 
\begin{equation}\label{eq: top change of cycles}
\begin{tikzcd}[column sep = huge]
\topBMg(\ulY_{\gamma'}^{\varepsilon = 1}(\leq \lambda')) \ar[r, " \fsp_{\gamma'}^{\mrm{ren}}"] & \topBMg(\ulY_{\gamma'}^{\varepsilon = 0}) \ar[d, equals, "\eqref{eq: identify types e = 0}"] \\ 
\topBMg(\ulY_{\gamma}^{\varepsilon = 1}(\leq \lambda))  \ar[u, dashed, "\tr_{(w')^{-1} t^{\delta} w}"] \ar[r, "
\fsp_{\gamma}^{\mrm{ren}}"] & \topBMg(\ulY_{\gamma}^{\varepsilon = 0})
\end{tikzcd}
\end{equation}
commutes, where $\delta := \mu - \mu'$ and the left vertical arrow is the partially defined map induced by \eqref{eq: partial transfer map}.  
\end{prop}

\begin{proof}
By Lemma \ref{lem: e=1 equivariant formality} and Remark \ref{remark: de-equivariant BM homology}, if $\gamma$ is $h_\lambda$-generic then we have 
\begin{equation}\label{eq: pass to equivariant class}
\topBMg(\ulY_{\gamma}^{\varepsilon = \varepsilon_0}(\leq \lambda)) \xleftarrow{\sim} \topBMulT(\ulY_\gamma^{\varepsilon = \varepsilon_0}(\leq \lambda)) \inj \Frac(\sph) \otimes_{\sph} \topBMulT(\ulY_\gamma^{\varepsilon = \varepsilon_0}(\leq \lambda))
\end{equation}
for any $\varepsilon_0 \in \{0,1, \eta\}$. Therefore it suffices to check the commutativity of the diagram in rationalized $\ulchT$-equivariant Borel--Moore homology. Then we may apply the equivariant localization theorem to compare the rationalized $\ulchT$-equivariant homology version of \eqref{eq: top change of cycles} with the corresponding diagram for the rationalized equivariant homology version for the $\ulchT$-fixed points. This embeds the corresponding diagram for $\chT$-equivariant Borel--Moore homology 
 \begin{equation}\label{eq: t-equiv compatibility}
\begin{tikzcd}[column sep = huge]
\topBMulT(\ulY_{\gamma'}^{\varepsilon = 1}(\leq \lambda')) \ar[r, "\fsp_{\gamma'}^{\mrm{ren}}"] & \topBMulT(\ulY_{\gamma'}^{\varepsilon = 0}) \ar[d, equals, "\eqref{eq: identify types e = 0}"] \\ 
\topBMulT(\ulY_{\gamma}^{\varepsilon = 1}(\leq \lambda))  \ar[u, dashed,  "\tr_{(w')^{-1} t^{\delta} w}"] \ar[r, "\fsp_\gamma^{\mrm{ren}}"] & \topBMulT(\ulY_{\gamma}^{\varepsilon = 0})
\end{tikzcd}
\end{equation}
into the similar diagram for $\ulchT$-equivariant Borel--Moore homology of the $\ulchT$-fixed points, tensored with $\Frac(\sph)$. Since the family $\ulY_\gamma^{\varepsilon}$ restricts to the constant family $\wt{\ulW} \times \Spec \F_p[\varepsilon]$ on $\ulchT$-fixed points, the similar diagram on fixed points reads 
\[
\begin{tikzcd}
\bigoplus_{\wt{w} \in \wt{W}} \sph [ \wt{w}]  \ar[r, "t^{\mu'}w'\cdot" , hook] & \bigoplus_{\wt{w} \in \wt{W}} \sph [ \wt{w}]\ar[d, equals] \\
\bigoplus_{\wt{w} \in \wt{W}} \sph [ \wt{w}] \ar[u, "\tr_{(w')^{-1} t^{\delta} w}", dashed] \ar[r,  "t^{\mu}w \cdot ", hook] & \bigoplus_{\wt{w} \in \wt{W}} \sph [ \wt{w}]
\end{tikzcd}
\]
for which the commutativity is evident. 


\end{proof}

\subsection{Breuil-M\'{e}zard cycles on the local model}\label{ssec: BM cycles on local model} 
We now undertake the construction of Breuil--M\'ezard cycles. As explained in \S \ref{sssec:bm-cycle-strat}, the strategy is to define cycles on the affine Springer fibers $\ul{\rY}_\gamma^{\varepsilon=0}$ using the microlocal support map (cf. Definition \ref{defn: musupp}), then to deform them to the $\ul{\rY}_\gamma^{\varepsilon=1}$, and finally to transfer them to the Emerton-Gee stack using the maps $\transfer_\gamma$ from Theorem \ref{thm: companion}, as in the diagram below. 
\[
\begin{tikzcd}
\rK(\Rep^{\emptyset} (G_1 T)) \ar[r, "\musupp"]  &  \topCh(\ul{\rY}_\gamma^{\varepsilon = 0}) \ar[r, "(\fsp_{\varepsilon \rightarrow 0})^{-1}"]  & \topCh(\ul{\rY}_\gamma^{\varepsilon = 1}) \ar[r, "\transfer_\gamma"] & \topCh(\cX_{\F}^{\lambda, \tau}) \ar[r, hook] & \dZ(\cX^{\EG}_{\red})
 \end{tikzcd}
\]

\subsubsection{Cycles on affine Springer fibers} We begin by constructing the cycles on $\rY_\gamma^{\varepsilon = 0}$. Recall the definition of $\wt{\ulW}_1$ from \S \ref{sssec: reparametrization of serre weights}. 

\begin{defn}\label{def: BM cycles on local model e=0}Let $\gamma  = ts $ with $s \in \ul{\chft}(\ol{\F}_p)$ regular. Any ($1$-generic) Serre weight $\sigma = F(\lambda)$ has a unique lowest alcove presentation of its highest weight $\lambda=\pi^{-1}(u)\bup (\xi-\rho)$, for $u \in \wt{\ulW}_1$ and $\xi \in \ulC_0$ so that $(u,\xi)$ is compatible with our fixed choice of central character $\zeta$ in the sense of Definition \ref{def: compatible pair}. We define
\begin{equation}\label{eq: translated version BM Z}
\cZ_{\gamma}^{\varepsilon = 0}(\sigma)= \cZ_{\gamma}^{\varepsilon = 0}(u, \xi) :=  t^{\xi} \cdot  \musupp[\wh{L}_1(u \bup 0)] \in \topCh(\ulY_\gamma^{\varepsilon = 0})
\end{equation}
where we remind that ``$t^{\xi} \cdot$'' refers to the monodromy-centralizer action $t^{\xi} \in \wt W$. 

In what follows, whenever the notation $\cZ_{\gamma}^{\varepsilon=0}(u, \xi)$ appears, the implicit assumption is always that $(u, \xi)$ is compatible with $\zeta$. 
\end{defn}

For later use, we record the effect of the monodromy action on the cycles just defined. 

\begin{lemma}[Monodromy action on the cycles]\label{lem: w-action on cycle} Let $\mu \in  \rho +  \ulC_0$ and $u \in \wt{\ulW}_1$. If $\nu \in X^*(\ulT)$ is sufficiently small so that $\mu - \rho + \ulW \nu \subset \ulC_0$, then 
\begin{equation}\label{eq: monodromy w action on BM cycle}
w  t^{-\mu } \cdot \cZ_\gamma^{\varepsilon = 0}(u , \mu +\nu) = t^{-\mu } \cdot \cZ_{\gamma}^{\varepsilon = 0}(u,  \mu + w \nu) \in \topCh(\ulY^{\varepsilon = 0}_{\gamma}).
\end{equation}
\end{lemma}

\begin{proof}Plugging in the definition \eqref{eq: translated version BM Z} to both sides, we find that the LHS of \eqref{eq: monodromy w action on BM cycle} is 
\[
w t^{\nu} \cdot \musupp [\wh{L}_1(u \bup 0 )]  = t^{w \nu} w\cdot \musupp [\wh{L}_1(u \bup 0 )]
\]
and the RHS of \eqref{eq: monodromy w action on BM cycle} is $t^{w \nu} \cdot \musupp [\wh{L}_1(u \bup 0 )]$. The two agree because $w\cdot  [\wh{L}_1(u \bup 0 ) ]= [\wh{L}_1(u \bup 0 )] $, as explained in Example \ref{ex: monodromy action on simples}.
\end{proof}

\subsubsection{Cycles on the local models} Our next task is to transport the cycles to the model $\ulY^{\varepsilon = 1}_{\gamma}$. First we analyze when the cycle $\cZ_{\gamma}^{\varepsilon = 0}(u, \xi )$ lies in the image of the renormalized specialization map $\fsp_\gamma^{\mrm{ren}}$ on $\topCh(\ulY^{\varepsilon = 1}_{\gamma}(\leq \lambda + \rho))$, and can therefore be deformed to $\varepsilon = 1$.

\begin{defn}[Admissible tuples] \label{defn:admissible tuples} Let $u \in \wt{\ulW}_1, w \in \ulW$ and $\nu \in X^*(\ulT)$. We say that $(u, w, \nu)$ is \emph{admissible} for $\lambda \in X^*(\ulT)^+$ if 
\[
(t^{-w^{-1}\nu})_{\dom} \leq w_0t^{-\lambda-\rho}u.
\]
\end{defn}

\begin{remark}\label{rem: admissible remark}
The representation-theoretic significance of this definition is as follows. Suppose $\mu$ is $h_{\lambda +2\rho}$-generic. Then by \cite[Proposition 2.3.7]{LLLM22}, the Serre weight $F(\pi^{-1}(u) \bup (\mu -\rho+\nu))$ appears as a Jordan-H\"{o}lder factor of $W(\lambda) \otimes R(w, \mu)$ if and only if $(u, w, \nu)$ is admissible for $\lambda$. 
\end{remark}

\begin{lemma}\label{lem: admissible support}
When defined, we have $t^{-\mu+\rho}  \cdot \cZ_{\gamma}^{\varepsilon = 0}(u, \mu + w^{-1} \nu) \in \topCh(\ulY^{\varepsilon = 0}_{\gamma}(\leq \lambda + \rho))$ if and only if $(u, w, \nu)$ is admissible for $\lambda$. 
\end{lemma}

\begin{proof}
This follows from Lemma \ref{lem:supp of simple}.
\end{proof}

If $\mu$ is $h_{\lambda+\rho}$-generic and $\gamma = \gamma(w, \mu)$, then 
\begin{equation}\label{eq: fsp gamma lambda+rho}
\fsp_{\gamma}^{\mrm{ren}} \co \topCh(\ulY^{\varepsilon = 1}_{\gamma}(\leq \lambda + \rho)) \inj \topCh(\ulY^{\varepsilon = 0}_{\gamma})
\end{equation}
is defined. By Lemmas \ref{lem: w-action on cycle} and \ref{lem: admissible support}, $w^{-1} t^{-\mu+\rho} \cdot \cZ_\gamma^{\varepsilon = 0}(u, \mu+\nu)$ lies in the image of \eqref{eq: fsp gamma lambda+rho} precisely when $(u, w, \nu)$ is admissible for $\lambda$. 

We can now define the incarnation of Breuil--M\'ezard cycles on the models.

\begin{defn}\label{def: BM cycles on local model} Assume $\mu \in \rho  + \ulC_0$ is $h_{\lambda +\rho}$-generic and let $\gamma = \gamma(w, \mu)$. If $(u, w, \nu)$ is admissible for $\lambda$ and $(u, \mu + \nu)$ is compatible with $\zeta$, we define $\cZ_\gamma^{\varepsilon = 1}(u , \mu +\nu)_{\lambda} \in \topCh(\ulY^{\varepsilon = 1}_\gamma(\leq \lambda+\rho))$ to be the unique cycle satisfying
\[
\fsp_\gamma^{\mrm{ren}}\left(\cZ_\gamma^{\varepsilon = 1}(u, \mu+\nu)_{\lambda} \right)
 =   \cZ_\gamma^{\varepsilon = 0}(u, \mu+\nu)   \in \topCh(\ulY^{\varepsilon = 0}_\gamma).
\]
(This notation reflects that in practice $\nu$ is ``small'' relative to $\mu$, so we visualize the $\{\xi\}$ for which $\cZ_{\gamma}^{\varepsilon =1}(u, \xi)_\lambda$ is defined as forming a constellation of weights orbiting $\mu$.) 

For $\sigma = F( \pi^{-1}u \bup ({\mu-\rho+\nu}))$, we denote 
\[
\cZ_\gamma^{\varepsilon = 1}(\sigma)_{\lambda} := \cZ_\gamma^{\varepsilon=1}( u,\mu+\nu)_{\lambda}.
\]
\end{defn}

\begin{remark}[Independence of $\lambda$]\label{rem: ind of lambda}
If $\lambda \leq \lambda'$ and $\gamma$ is $h_{\lambda' + \rho}$-generic, then $\gamma$ is also $h_{\lambda+\rho}$-generic and the inclusion $\ulY^{\varepsilon = 1}_\gamma(\leq \lambda+\rho) \inj \ulY^{\varepsilon = 1}_\gamma(\leq \lambda'+\rho)$ sends $\cZ_\gamma^{\varepsilon = 1}(u, \mu+\nu)_{\lambda} \mapsto \cZ_\gamma^{\varepsilon = 1}(u, \mu+\nu)_{\lambda'}$. We therefore have a well-defined class 
\[
\cZ_\gamma^{\varepsilon = 1}(u, \mu+\nu)\in \topCh(\ulY^{\varepsilon = 1}_\gamma)
\]
as long as there exists $\lambda$ such that $\mu$ is $h_{\lambda+\rho}$-generic and $(u, w, \nu)$ are admissible for $\lambda$. In this case we say that the Serre weight $\sigma = F( \pi^{-1}u \bup ({\mu-\rho+\nu}))$ is \emph{admissible} for $(\gamma, \lambda)$, and we abbreviate $\cZ_\gamma^{\varepsilon = 1}(\sigma) := \cZ_\gamma^{\varepsilon = 1}(u, \mu+\nu)$. With these definitions, the identity 
\[
\fsp_{\gamma}^{\mrm{ren}} \left( \cZ_\gamma^{\varepsilon = 1}(u, \mu+\nu) \right) =  
\cZ^{\varepsilon = 0}_\gamma(u, \mu+\nu)  \in \topCh(\ulY^{\varepsilon = 0}_\gamma)
\]
holds as long as all of its terms are defined. \end{remark}

We will next establish a certain ``independence of $\gamma$'' for the cycles $\cZ_{\gamma}^{\varepsilon = 1}(\sigma)$, in preparation for transferring these cycles to the Emerton-Gee stack.

\begin{lemma}[Independence of $\gamma$]\label{lem: translate BM cycle} Let $\mu, w, \mu', w', \delta, \gamma, \gamma'$ be as in Proposition \ref{prop: type change of cycles}. Let $u \in \wt{\ulW}_1$ and $\xi \in \rho + \ulC_0$ such that $\cZ_{\gamma}^{\varepsilon = 1}(u, \xi)$ and $\cZ_{\gamma'}^{\varepsilon = 1}(u, \xi) $ are both defined. Then $\tr_{(w')^{-1} t^{\delta} w} (\cZ_{\gamma}^{\varepsilon = 1}(u, \xi))$ is defined and 
\begin{equation}\label{eq: translate BM cycle}
\tr_{(w')^{-1} t^{\delta} w} (\cZ_{\gamma}^{\varepsilon = 1}(u , \xi)) = \cZ_{\gamma'}^{\varepsilon = 1}(u , \xi)  \in \topCh(\ulY_{\gamma'}^{\varepsilon = 1}).
\end{equation}
\end{lemma}

\begin{proof}Note that it is enough to check the equality in $\topBMg$, since we have $\topCh \inj \topBMg$. 

By definition, if $\cZ_{\gamma}^{\varepsilon = 1}(u, \xi)$ is defined then $\fsp_{\gamma}^{\mrm{ren}}(\cZ_{\gamma}^{\varepsilon = 1}(u , \xi)) \in \topBMg(\ulY_\gamma^{\varepsilon=0})$ is defined and equals $t^{\xi} \cdot \musupp[\wh{L}_1(u \bup 0)]$.

Similarly, if $\cZ_{\gamma'}^{\varepsilon = 1}(u, \xi)$ is defined then  $\fsp_{\gamma'}^{\mrm{ren}}(\cZ_{\gamma'}^{\varepsilon = 1}(u, \xi))\in \topBMg(\ulY_{\gamma'}^{\varepsilon=0})$ is defined and equals $ t^{\xi} \cdot \musupp[\wh{L}_1(u \bup 0)]$.

Taking $\lambda$ and $\lambda'$ so that $\sigma = F( \pi^{-1}u \bup ({\mu-\rho+\nu}))$ is admissible for $(\gamma, \lambda)$ and $(\gamma', \lambda')$, the partially defined map 
\[
\begin{tikzcd}[column sep = huge]
\topBMg(\ulY_{\gamma}^{\varepsilon = 1}(\leq \lambda + \rho)) \ar[r, dashed, "\tr_{(w')^{-1} t^{\delta} w}"] & \topBMg(\ulY_{\gamma'}^{\varepsilon = 1}(\leq \lambda' + \rho))
\end{tikzcd}
\]
includes $\cZ_{\gamma}^{\varepsilon = 1}(u, \xi)$ in its domain, so $\tr_{(w')^{-1} t^{\delta} w}(\cZ_{\gamma}^{\varepsilon = 1}(u, \xi))$ is defined in $\topBMg(\ulY_{\gamma'}^{\varepsilon = 1}(\leq \lambda'+\rho))$, and the equality \eqref{eq: translate BM cycle} then follows from the previous paragraphs plus the commutativity of the diagram \eqref{eq: top change of cycles}.

\end{proof}

\subsection{Breuil--M\'ezard relations on the model}\label{ssec: BM relations on local model} We will next verify a collection of relations on the models, which will later be seen to correspond to \eqref{eq: BM relation} under $\transfer_\gamma$. In terms of Figure ~\ref{fig: bm cycle}, we will show that the cycles $\cZ_\gamma^{\varepsilon=1}(\sigma)$ verify all the Breuil--M\'ezard relations which concern cycles contained in a single ``chart'' labeled by a single $\gamma$. 

\subsubsection{The case $\lambda = 0 $} The fundamental case is where $\lambda = 0$, corresponding to potentially crystalline $L$-parameters with minimal regular Hodge-Tate weights, which is handled by the Theorem below.

\begin{thm}\label{thm: BM equation}
Suppose that $\mu \in \rho + \ulC_0$ is $2h_\rho$-generic. Let $\gamma = \gamma(w, \mu)$ for a lowest alcove presentation $(w, \mu)$ compatible with $\zeta$. Then we have 
\[
\fsp_{p \rightarrow 0} [\ulX_\gamma^{\varepsilon = 1}(\rho)] = \sum_{\sigma}  [\ol{R(w, \mu)} \co \sigma] \cZ_\gamma^{\varepsilon=1}(\sigma) \in \topCh(\ulY_{\gamma}^{\varepsilon = 1}).
\]
In this expression, we understand the summand to be $0$ whenever $[\ol{R(w, \mu)} \co \sigma] = 0$ (even when $\cZ^{\varepsilon=1}_\gamma(\sigma)$ is undefined -- we are implicitly claiming that $\cZ^{\varepsilon=1}_\gamma(\sigma)$ is defined whenever $[R(w, \mu) \co \sigma] \neq 0$). 
\end{thm}

Before giving the proof of Theorem \ref{thm: BM equation} we record some representation-theoretic preliminaries. We will relate the multiplicities $[\ol{R(w, \mu)} \co \sigma]$ to decomposition multiplicities for $\ulG_1\ulT$ using Jantzen's generic decomposition pattern (found in this generality in \cite[Proposition 10.1.2]{GHS18})\footnote{There is a discrepancy with the formula appearing in \cite[Proposition 10.1.2]{GHS18} because they use the dot action of the $p$-dilated extended affine Weyl group, while we use the $p$-dilated dot action of the extended affine Weyl group. Also note that their ``$R(w, \mu+\rho)$'' is our $R(w, \mu)$.}: if $\mu$ is $2h_\rho$-generic, then 
\begin{equation}\label{eq: GHS}
\ol{R(w, \mu)} = \sum_{u \in \wt{\ulW}_1/X^0(\ulT)}  \sum_{\nu \in X^*(\ulT)} [\wh{Z}_1(\mu - \rho +p\rho) \co \wh{L}_1(u \bup (\mu - \rho)+ p \nu)] F(u \bup (\mu - \rho +w \pi \nu)).
\end{equation} 

Let $m_{u, \nu}^{\mu} := [\wh{Z}_1(\mu - \rho +p\rho) \co \wh{L}_1( u \bup (\mu - \rho) + p\nu)]$. Note that by the translation principle, this is independent of $\mu$ as long as $\mu - \rho $ is regular in $\ulC_0$. In particular, we have 
\begin{align*}
m_{u, \nu}^{\mu} &= [\wh{Z}_1(p \rho): \wh{L}_1(  u \bup 0 + p \nu )].
\end{align*}

We also note for future use that 
\begin{equation}\label{eq: same multiplicity}
m_{\pi u,\pi \nu}^{\mu}=m_{u,\nu}^{\mu},
\end{equation}
since we have  
\begin{align*}
m_{u,\nu}^{\mu} &=  [\wh{Z}_1(\mu - \rho +p\rho):\wh{L}_1( u \bup (\mu-\rho) + p\nu)] \\
&  = [\wh{Z}_1(\pi(\mu - \rho)+p\rho):\wh{L}_1( \pi u \bup \pi(\mu-\rho) + p\pi\nu)] \\
&  = [\wh{Z}_1(\mu - \rho +p\rho):\wh{L}_1( \pi u \bup (\mu-\rho) + p\pi\nu)]  =m_{\pi u, \pi \nu}^{\mu},
\end{align*}
where on the third line we use the translation principle to replace $\pi( \mu-\rho)$ by $\mu-\rho$. 

\begin{lemma}\label{lem: genericity} 
(1) If $\mu$ is $m$-generic for some $m \geq 2h_\rho$, then every $\sigma$ such that $[\ol{R(w, \mu)} \co \sigma] \neq 0$ is of the form $\sigma=F(\lambda)$ where $\lambda + \rho$ is $(m-h_\rho)$-generic. 

(2) If $\lambda + \rho$ is $m$-generic and $[\sigma(\tau) \co F(\lambda)] \neq 0$ then $\sigma(\tau)=R(w,\mu)$ where $\mu$ is $(m-h_\rho)$-generic. 
\end{lemma}

\begin{proof} 
Part (1) follows from \cite[Proposition 2.3.7]{LLLM22}, and part (2) follows from \cite[Lemma 2.3.4]{LLLM22}.
\end{proof}

\begin{proof}[Proof of Theorem \ref{thm: BM equation}]
Since
\begin{equation}\label{eq: covering weights}
\wt{\ulW} \bup 0 = \bigcup_{u \in \wt{\ulW}_1/X^0(\ulT)} \bigcup_{\nu \in X^*(\ulT)} \{u \bup 0  + p \nu\}
\end{equation}
we have
\begin{equation}\label{eq: BM equation 1}
 \musupp[\wh{Z}_1(p \rho)] = \sum_{u \in \wt{\ulW}_1/X^0(\ulT)} \sum_{\nu \in X^*(\ulT)} m_{u, \nu}^{\mu} \, \musupp[\wh{L}_1(   u \bup 0 + p\nu)] \in \topCh(\ulY_{\gamma}^{\varepsilon = 0}).
\end{equation}
By translation equivariance of $\musupp$, we have $\musupp[\wh{L}_1(u \bup 0 + p\nu )] = t^{\nu} \cdot \musupp[\wh{L}_1(u \bup 0)]$, so we can rewrite \eqref{eq: BM equation 1} as 
\begin{equation}\label{eq: BM equation 2}
\musupp[\wh{Z}_1(p \rho)] = \sum_{u, \nu} m_{u, \nu}^{\mu}  t^{\nu } \cdot \musupp[\wh{L}_1(u \bup 0)] \in \topCh(\ulY_{\gamma}^{\varepsilon = 0}).
\end{equation}

For later comparison, we reparametrize $u \mapsto \pi u$ and $\nu \mapsto \pi \nu$ so that inserting Theorem \ref{thm: baby verma equals limit cycle} into \eqref{eq: BM equation 2} gives
\begin{equation}\label{eq: BM equation 3.1}
\fsp_{p \rightarrow 0} [\ulX_{\gamma}^{\varepsilon = 0}(\rho)] = \sum_{u \in \wt{\ulW}_1/X^0(\ulT)}  \sum_{\nu \in X^*(\ulT)} m_{\pi u, \pi \nu}^{\mu} t^{\pi   \nu}  \cdot  \musupp[\wh{L}_1(\pi u  \bup 0)] \in \topCh(\ulY_{\gamma}^{\varepsilon = 0}).
\end{equation}
Substituting the observation \eqref{eq: same multiplicity} that $m_{\pi u, \pi \nu}^{\mu} = m_{u,\nu}^{\mu}$ into \eqref{eq: BM equation 3.1} gives
\begin{equation}\label{eq: BM equation 3.2}
\fsp_{p \rightarrow 0} [\ulX_{\gamma}^{\varepsilon = 0}(\rho)] = \sum_{u \in \wt{\ulW}_1/X^0(\ulT)}  \sum_{\nu \in X^*(\ulT)} m_{ u,  \nu}^{\mu} t^{\pi   \nu}  \cdot  \musupp[\wh{L}_1(\pi u \bup 0)] \in \topCh(\ulY_{\gamma}^{\varepsilon = 0}).
\end{equation}
The assumptions on genericity and compatibility with $\zeta$ are such that Definition \ref{def: BM cycles on local model} says, using also \eqref{eq: monodromy w action on BM cycle}, 
\begin{equation}\label{eq: translation property}
\fsp_{\varepsilon \rightarrow 0} \cZ^{\varepsilon = 1}_\gamma( \pi u,   \mu  +  w \pi  \nu) =  t^{\pi \nu}\cdot  \musupp[\wh{L}_1(\pi u \bup 0)] \in \topCh(\ulY_{\gamma}^{\varepsilon = 0}).
\end{equation}
Using this, we may rewrite \eqref{eq: BM equation 3.2} as 
\begin{equation}\label{eq: BM equation 3.3}
\fsp_{p \rightarrow 0} [\ulX_{\gamma}^{\varepsilon = 0}(\rho)] = \sum_{u \in \wt{\ulW}_1/X^0(\ulT)}  \sum_{\nu \in X^*(\ulT)} m_{ u,  \nu}^{\mu} \fsp_{\varepsilon \rightarrow 0} \cZ^{\varepsilon = 1}_\gamma( \pi u ,  \mu  +  w \pi \nu) \in \topCh(\ulY_{\gamma}^{\varepsilon = 0}).
\end{equation}
The assumptions imply that 
\[
\fsp_{\varepsilon \rightarrow 0} \co \topCh(\ulY_\gamma^{\varepsilon=1} (\leq  2\rho)) \rightarrow \topCh(\ulY_\gamma^{\varepsilon = 0}(\leq 2 \rho))
\]
is an isomorphism (cf. Lemma \ref{lem: kappa specialization upper triangular}). Therefore each $\cZ^{\varepsilon = 0}_{\gamma}(\pi u,  \mu  + w \pi \nu )$ for which $m_{u, \nu}^{\mu} \neq 0$ lies in the domain where $\fsp_{\varepsilon \rightarrow 0}$ is an isomorphism, so we may apply $\fsp_{\varepsilon \rightarrow 0}^{-1}$ to \eqref{eq: BM equation 3.3} along with Proposition \ref{prop: gr hbar spc} (and also use that specializations in $\varepsilon$ and $p$ commute by Lemma \ref{lem: iterated specialization}) to find that\footnote{Here again, if $m_{u, \nu}^{\mu} = 0$ then the corresponding summand is interpreted as $0$, and we are implicitly asserting that $m_{u, \nu}^{\mu} = 0$ if $\pi u \cdot (\mu - \rho + w \pi   \nu) \notin X_1^*(T)$ or if $\cZ^{\varepsilon = \eta}_\gamma(u ,  \mu + w \pi  \nu)$ is undefined.}
\begin{equation}\label{eq: BM equation 5}
\fsp_{p \rightarrow 0} [\ulX_{\gamma}^{\varepsilon = 1}(\rho)] = \sum_{u \in \wt{\ulW}_1/X^0(\ulT)}  \sum_{\nu \in X^*(\ulT)} m_{u, \nu}^{\mu} \cZ^{\varepsilon = 1}_{\gamma}(\pi u , \mu + w \pi   \nu).
\end{equation}
Note that $\cZ^{\varepsilon = 1}_\gamma(\pi u , \mu + w  \pi  \nu)  = \cZ^{\varepsilon = 1}_{\gamma}(\sigma)$ for $\sigma = F(u \bup (\mu - \rho+  w \pi \nu))$. Putting this into \eqref{eq: BM equation 5} completes the proof. 
\end{proof}

\subsubsection{Higher Hodge--Tate weights} We will next show that our Breuil-M\'{e}zard cycles $\cZ_{\gamma}^{\varepsilon = 1}(\sigma)$ satisfy the further relations expected from higher Hodge-Tate weights (the precise sense in which this is related to higher Hodge-Tate weights will be explained later in \S \ref{ssec: BM relations}). We begin with a purely representation-theoretic lemma. 

\begin{lemma}\label{lem: translation R} Assume $p \geq 2 h_\rho$. Let $\lambda \in X^*(\ulT)^+$ be such that $\mu+\kappa \in \rho + \ulC_0$ is $2h_\rho$-generic for all weights $\kappa$ of $W(\lambda)$. For $\nu \in X^*(\ulT)$, write $m_{\kappa}(\lambda)$ for the multiplicity of $\kappa$ as weight of the Weyl module $W(\lambda)$. Then we have
\[
[\ol{W(\lambda) \otimes R(w,\mu)}] =  \sum_{\kappa \in X^*(\ulT)} m_{\kappa}(\lambda)  [\ol{R(w, \mu+\kappa)}] \in \rK(\Rep_k(\ul{G}(\F_p))).
\]
\end{lemma}

\begin{proof}
Given $\xi \in X_1^*(\ulT)$ such that $\xi + \kappa$ lies in the same alcove as $\xi$ for all weights $\kappa$ of $W(\lambda)$, \cite[Lemma 4.2.4]{LLLM20} implies that  we have 
\begin{equation}\label{eq: small Weyl tensor}
[W(\lambda) \otimes_k L(\xi)] = \sum_{\kappa \in X^*(\ulT)} m_{\kappa}(\lambda)  [L(\xi+ \kappa)] \in \rK(\Rep_k(\ul{G}(\F_p))).
\end{equation}
By \eqref{eq: GHS} we have
\[
[\ol{R(w,\mu)}] = \sum_{u \in \wt{\ulW}_1/X^0(T)} \sum_{\nu \in X^*(\ulT)} m_{u,\nu} [F(u \bup (\mu - \rho+ w \pi \nu))] \in \rK(\Rep_k(\ul{G}(\F_p))),
\]
where $m_{u, \nu}  = [\wh{Z}_1(\mu - \rho + p \rho): \wh{L}_1(u \bup (\mu - \rho)+ p \nu)]$ as before (recall that we saw it was independent of $\mu$ by the translation principle). Under our assumptions, $m_{u, \nu} \neq 0 $ implies that $\xi := u \bup (\mu - \rho + w \pi  \nu)$ satisfies the condition needed to apply \eqref{eq: small Weyl tensor}. We therefore have 
\begin{equation}\label{eq: weyl tensor 1}
[\ol{W(\lambda) \otimes R(w,\mu)}] = \sum_{\kappa \in X^*(\ulT)} m_{\kappa}(\lambda)  \sum_{u \in \wt{\ulW}_1/X^0(T)} \sum_{\nu \in X^*(\ulT)}  m_{u,\nu} [F(u \bup (\mu - \rho + w \pi  \nu) + \kappa)] \in \rK(\Rep_k(\ul{G}(\F_p))). 
\end{equation}
Since the character of $W(\lambda)$ is invariant under the action of $W$, we may rewrite
\begin{equation}\label{eq: weyl tensor 2}
 \eqref{eq: weyl tensor 1} = \sum_{\kappa \in X^*(\ulT)} m_{\kappa}(\lambda)  \sum_{u \in \wt{\ulW}_1/X^0(T)} \sum_{\nu \in X^*(\ulT)}  m_{u,\nu} [F(u \bup (\mu - \rho + \kappa + w \pi  \nu) )] \in \rK(\Rep_k(\ulG(\F_p))). 
\end{equation}
Then we obtain the desired equality upon rearranging terms and applying \eqref{eq: GHS}.

\end{proof}

The next theorem is the generalization of Theorem \ref{thm: BM equation} that handles higher Hodge-Tate weights, although we will see that the proof is a reduction to Theorem \ref{thm: BM equation}. 

\begin{thm}[Breuil-M\'{e}zard relations on the models]\label{thm: BM equation with HT weight} Let $\lambda \in X^*(\ulT)^+$. Suppose that $\mu$ is $(h_\lambda + 2h_\rho)$-generic. Let $\gamma = \gamma(w, \mu)$ for a lowest alcove presentation $(w, \mu)$ compatible with $\zeta-\lambda$. Then we have
\begin{equation}\label{eq: BM equation with HT weight}
\fsp_{p \rightarrow 0} [\ulX^{\varepsilon=1}_{\gamma}(\lambda+\rho)]= \sum_{\sigma}  [\ol{W( \lambda) \otimes R(w, \mu)} \co \sigma] \cZ^{\varepsilon=1}_\gamma(\sigma) \in \topCh(\ulY_{\gamma}^{\varepsilon=1}).
\end{equation}
\end{thm}

\begin{proof} 

By Lemma \ref{lem: kappa specialization upper triangular}, the assumptions imply that 
\[
\fsp_{\varepsilon \rightarrow 0} \co  \topCh(\ulY_\gamma^{\varepsilon= 1} (\leq \lambda  + \rho)) \rightarrow \topCh(\ulY_\gamma^{\varepsilon = 0}(\leq \lambda + \rho))
\]
is defined and is an isomorphism, so the identity can be checked after applying $\fsp_{\varepsilon \rightarrow 0}$. We first analyze what happens upon doing this to the LHS of \eqref{eq: BM equation with HT weight}. Applying Proposition \ref{prop: gr hbar spc} and then Theorem \ref{thm: generic fiber cycle decomposition}, we deduce that 
\begin{equation}\label{eq: BM relation ht weight 1}
\fsp_{\varepsilon \rightarrow 0}  [\ulX_{\gamma}^{\varepsilon=1}(\lambda + \rho)]  = \sum_{\kappa \in X^*(\ulT)} m_{\kappa}(\lambda) \fsp_{p \rightarrow 0} t^{\kappa}  \cdot [\ulX_\gamma^{\varepsilon = 0}(\rho)] \in \topCh(\ulY_\gamma^{\varepsilon = 0}).
\end{equation}

Putting Theorem \ref{thm: BM equation} into \eqref{eq: BM relation ht weight 1} yields 
\begin{equation}\label{eq: BM relation ht weight 1.75}
\fsp_{\varepsilon \rightarrow 0}  [\ulX_{\gamma}^{\varepsilon=1}(\lambda + \rho)]   = \sum_{\kappa \in X^*(\ulT)} m_{\kappa}(\lambda)\sum_{\sigma} [\ol{R(w, \mu)}: \sigma]   t^{\kappa}  \cdot \cZ^{\varepsilon = 0}_\gamma(\sigma)\in \topCh(\ulY_\gamma^{\varepsilon = 0}).
\end{equation} 
For $u \in \wt{\ulW}_1$ and $\nu \in X^*(\ulT)$, consider the contribution of $\sigma := F(u \bup (\mu - \rho  + w \pi \nu))$ on the RHS of \eqref{eq: BM relation ht weight 1.75}. Set $\sigma' := F(u \bup (\mu - \rho + \kappa   + w \pi \nu))$. By construction we have $ t^{\kappa} \cdot \cZ_\gamma^{\varepsilon = 0} (\sigma) = \cZ_\gamma^{\varepsilon = 0} (\sigma')$. Also we saw from the translation principle that $ [\ol{R(w, \mu)}: \sigma] =  [\ol{R(w, \mu+\kappa)}: \sigma']$, so we may rewrite \eqref{eq: BM relation ht weight 1.75} as 
\begin{equation}\label{eq: BM relation ht weight 2}
\fsp_{\varepsilon \rightarrow 0}  [\ulX_{\gamma}^{\varepsilon=1}(\lambda + \rho)]  = \sum_{\kappa \in X^*(\ulT)} m_{\kappa}(\lambda) \sum_\sigma [R(w, \mu + \kappa ): \sigma'] \cZ^{\varepsilon = 0}_\gamma(\sigma').
\end{equation}
By Lemma \ref{lem: translation R}, the coefficient of $\cZ_\gamma^{\varepsilon = 0}(\sigma')$ in \eqref{eq: BM relation ht weight 2} is 
\[
\sum_{\kappa \in X^*(\ulT)} m_{\kappa}(\lambda)  [\ol{R(w, \mu + \kappa )}: \sigma'] = [\ol{W(\lambda) \otimes R(w, \mu)}:\sigma'],
\]
so that \eqref{eq: BM relation ht weight 2} agrees with 
\[
\sum_{\sigma}  [\ol{W( \lambda) \otimes R(w, \mu  )} \co \sigma'] \cZ^{\varepsilon=0}_\gamma(\sigma')\in \topCh(\ulY_{\gamma}^{\varepsilon=0}),
\]
as desired.

\end{proof}

\subsection{Proof of Theorem \ref{thm: BM 13}}\label{ssec: BM relations}

First, we (finally!) construct Breuil-M\'{e}zard cycles on the Emerton-Gee stack $\cX^{\EG}$.

\begin{lemma}[Independence of $\gamma$]\label{lem: coordinate change}
Consider the map $\transfer_\gamma$ from Theorem \ref{thm: companion}. For all $\gamma$ such that $\cZ_{\gamma}^{\varepsilon = 1}(u,  \xi)$ is defined, the classes 
\[
\transfer_{\gamma}(\cZ_{\gamma}^{\varepsilon = 1}(u, \xi)) \in \dZ(\cX^{\EG}_{\red})
\]
coincide.
\end{lemma}

\begin{proof} This follows from Proposition \ref{prop: companion compatibility} and Lemma \ref{lem: translate BM cycle}.
\end{proof}

\begin{defn}[Construction of Breuil-M\'{e}zard cycles]\label{def: bm cycles} Suppose that $\sigma = F(\pi^{-1}u \bup (\xi-\rho))$ occurs in a tame type $\tau = R(w, \mu)$ where 
\begin{itemize}
\item $\mu$ is $2h_\rho$-generic\footnote{This is guaranteed whenever $\xi + \rho$ is $(2h_\rho+1)$-generic: adjusting $u$ by an element in $\ul{\Omega}$, we can make $u\in t^{\rho}\ulW_{\aff}$ and so we can then take $\nu=0$ in Remark \ref{rem: admissible remark}. Then any lowest alcove presentation for the resulting inertial parameter is $2h_\rho$-generic.}, and 
\item $(u, \xi)$ and $(w, \mu)$ are compatible with $\zeta$. 
\end{itemize}
Let $\gamma = \gamma(w,\mu)$. Then the cycle $\cZ_{\gamma}^{\varepsilon = 1}(\sigma) \in \topCh(\ulY_\gamma^{\varepsilon = 1})$ is defined in Definition \ref{def: BM cycles on local model}, and we define
\[
\cZ^{\EG}(\sigma) := \transfer_\gamma (\cZ_{\gamma}^{\varepsilon = 1}(\sigma) ) \in \dZ(\cX^{\EG}_{\red}).
\]
A priori this definition seems to depend on the choice of $(w, \mu)$, but the independence of this choice was established in Lemma \ref{lem: coordinate change}. 
\end{defn}

We may now complete the proof of Theorem \ref{thm: BM 13}. 

\begin{proof}[Proof of Theorem \ref{thm: BM 13}] Let $\lambda \in X^*(\ulT)^+$ and $\tau = \tau(w, \mu)$ be a lowest alcove presentation of a tame inertial parameter such that $\mu$ is $(2h_\rho+h_\lambda)$-generic. 

Let $\gamma := \gamma(w, \mu)$. Then from Theorem \ref{thm: BM equation with HT weight} and the definition \eqref{eq:n-sigma-lambda-tau} of $n_\sigma(\lambda, \tau)$, we have that
\begin{equation}\label{eq: thm bm 1}
\fsp_{p \rightarrow 0} [\ulX^{\varepsilon=1}_{\gamma}(\lambda+\rho)]= \sum_{\sigma}  n_\sigma(\lambda, \tau) \cZ^{\varepsilon=1}_\gamma(\sigma) \in \topCh(\ulY_{\gamma}^{\varepsilon=1}).
\end{equation}
Consider applying $\transfer_\gamma$ to this identity. Since $\transfer_{\gamma}(\fsp_{p \rightarrow 0} [\ulX_\gamma^{\varepsilon = 1}(\lambda + \rho)])  = 
[\cX^{\lambda+\rho,\tau}_{\F}] $ by Theorem \ref{thm: companion}(1) and $\cZ^{\EG}(\sigma) = \transfer_\gamma(\cZ_{\gamma}^{\varepsilon = 1}(\sigma))$ by definition, \eqref{eq: thm bm 1} becomes
\[
[\cX^{\lambda+\rho,\tau}_{\F}] =\sum_{\sigma}  n_\sigma(\lambda, \tau) \cZ^{\EG}(\sigma) \in \dZ(\cX^{\EG}_{\red}),
\]
which is exactly what we wanted to show. 
\end{proof}

\section{Uniqueness of Breuil-M\'{e}zard cycles}\label{sec: complements} Let $G$ be an unramified group over $\Q_p$ satisfying Hypothesis \ref{hyp:part3-hypotheses}, and continue to abbreviate $\cX^{\EG} := \cX^{\EG}_{\LG}$ and $\cX^{\lambda, \tau} := \cX^{\lambda, \tau}_{\LG}$.  Here we establish the second part of Theorem \ref{thm: intro main}, asserting that if Conjecture \ref{conj:tame-BM-Conj} is true, then the ``true'' Breuil--M\'ezard cycles must agree with the cycles $\cZ^{\EG}(\sigma)$ as soon as $\sigma$ is sufficiently generic (which can be quantified effectively). 

The main result of this subsection is: 
\begin{thm}\label{thm: uniqueness} 
If there are effective cycles $\cZ(\sigma) \in \dZ(\cX^{\EG}_{\red})$ satisfying Conjecture \ref{conj:tame-BM-Conj}, then $\cZ(\sigma)$ agrees with the $\cZ^{\EG}(\sigma)$ from Definition \ref{def: bm cycles} whenever $\sigma = F(\lambda)$ such that $\lambda$ is $6h_\rho$-generic. 
\end{thm}

\begin{remark}\label{remark:bounded-support-property}The proof of Theorem \ref{thm: uniqueness} only uses the following ``bounded support'' property of $\cZ(\sigma)$: if $[\ol{R(w, \mu)} \co \sigma ] \neq 0$ for some $2h_\rho$-generic $\mu$, then the support of $\cZ(\sigma)$ is contained in $\cX^{\rho, \tau}_{\F}$. Effectivity of  $\cZ(\sigma)$ clearly implies this property. Note that the $\cZ^{\EG}(\sigma)$ constructed by Theorem \ref{thm: BM 13} have this bounded support property, but are not obviously effective. Thus even in the case $G=\Res_{K/\Q_p}\GL_n$ for unramified $K$, Theorem \ref{thm: uniqueness} is stronger than the similar uniqueness result in \cite[Proposition 8.6.5]{LLLM22}.
\end{remark}

\subsection{Reformulation in equivariant homology} \label{subsect:uniqueness reduction}
Fix an algebraic central character $\zeta$, as in \S \ref{sec: existence of BM cycles}. By the containment \eqref{eq:cycle-contain} and the injectivity of $\transfer_\gamma$ from
Theorem \ref{thm: companion}: for each $\tau = \tau(w, \mu)$ and $\gamma = \gamma(w, \mu)$ such that $\mu$ is $2h_\rho$-generic and $(w, \mu)$ is compatible with $\zeta$, 
we have an injection
\[
\transfer_\gamma^{-1} \co \dZ(\cX^{\leq \rho,\tau}_{\F}) \inj \dZ(\ulY_\gamma^{\varepsilon =1}(\leq \rho)),
\]

As explained in Remark \ref{remark:bounded-support-property}, the cycles $\cZ(\sigma)$ from the hypothesis of Theorem \ref{thm: uniqueness} must have support in $\cX^{\leq \rho,\tau}_{\F}$. Then applying the map $\transfer_\gamma^{-1}$ to the cycles $\cZ(\sigma)$ gives a collection of cycles $\cZ_{\gamma}^{\varepsilon =1}(\sigma)^\dagger \in \topCh(\ulY_\gamma^{\varepsilon =1}(\leq \rho))$ for each $\sigma \in \JH(R(w, \mu))$. Recall the renormalized specialization map from Definition \ref{defn: renormalized specialization},
 \[
 \fsp_{\gamma}^{\mrm{ren}} \co \topCh(\ulY^{\varepsilon = 1}_{\gamma} (\leq \rho))\cong   \topCh(\ulY^{\varepsilon = 0}_{\gamma} (\leq \rho)) \inj\topCh(\ulY^{\varepsilon = 0}_{\gamma}).
 \]
Applying this, we obtain cycles $\cZ_{\gamma}^{\varepsilon =0}(\sigma)^\dagger :=\fsp_{\gamma}^{\mrm{ren}} (\cZ_{\gamma}^{\varepsilon =1}(\sigma)^\dagger) \in \topCh(\ulY^{\varepsilon = 0}_{\gamma})$. On the other hand, the Breuil--M\'ezard cycles $\cZ_{\gamma}^{\varepsilon =1}(\sigma)$ we constructed in Definition \ref{def: BM cycles on local model} were characterized by the property (cf. also Definition \ref{def: BM cycles on local model e=0})
\[
\fsp_{\gamma}^{\mrm{ren}}\left(\cZ_{\gamma}^{\varepsilon =1}(F(\pi^{-1}u\bup \xi)) \right) =t^{\xi+\rho} \cdot \musupp[\wh{L}_1(u \bup 0)] \in \topCh(\ulY_\gamma^{\varepsilon = 0}).
\]
Therefore Theorem \ref{thm: uniqueness} is equivalent to the statement 
\begin{equation}\label{eq: equivariant reformulation of uniqueness}
\cZ_{\gamma}^{\varepsilon =0}(F(\pi^{-1}u\bup \xi))^\dagger=t^{\xi+\rho}  \cdot \musupp[\wh{L}_1(u \bup 0)] \in \topCh(\ulY_\gamma^{\varepsilon = 0}).
\end{equation}

Now recall that the spaces $\topCh(\ulY^{\varepsilon = 0}_{\gamma})$ (resp. $\topBM(\ulY^{\varepsilon = 0}_{\gamma})$, $\topBMulT(\ulY^{\varepsilon = 0}_{\gamma})$) are canonically identified via \eqref{eq: identify equivariant types e = 0} as we vary $\gamma=\gamma(w,\mu)$ over possible choices of $(w,\mu)$. It follows from Proposition \ref{prop: companion compatibility} and Proposition \ref{prop: type change of cycles} that under these identifications $\cZ_{\gamma}^{\varepsilon =0}(\sigma)^\dagger$ is independent of the choice of $\gamma$.

\subsection{Equivariant support bounds}  We establish some technical statements on the equivariant support of $
\cZ_{\gamma}^{\varepsilon =0}(F(\pi^{-1}u\bup \xi))^\dagger$ for later use.

\begin{lemma}\label{lem:effective support bound} Assume all types containing $\sigma=F(\pi^{-1}u\bup \xi)$ are $2h_\rho$-generic. Then $\Loc^{\ulchT} (\cZ_{\gamma}^{\varepsilon=0}(\sigma)^{\dagger}_{\ulchT})$ has equivariant support in $t^{\xi+\rho}(\ul{\wt{W}}_{\leq w_0u})$.
\end{lemma}

In particular, the lemma applies whenever $\xi+\rho$ is $3h_\rho$-generic.

\begin{proof} The tame types $\sigma(\tau)=R(w,\mu)$ containing $\sigma$ are exactly those such that for $\kappa:=\mu-\rho-\xi$, the triple $(u,w,-\kappa)$ is admissible for $0$ in the sense of Definition \ref{defn:admissible tuples}, i.e.,
\[
(t^{w^{-1}\kappa})_{\dom} = (t^{\kappa}w)_{\dom}\leq w_0t^{-\rho}u.
\]
It follows from the definitions that any class in $\fsp_{\gamma(w,\mu)}\circ \transfer^{-1}_{\gamma(w,\mu)}\big(\topCh(\cX^{\leq \rho,\tau}_{\F})\big)$ has equivariant support in $t^{\mu}w \cdot \Adm(\rho)$. Hence (using the effectivity assumption) we learn that the equivariant support of $\cZ_{\gamma}^{\varepsilon=0}(\sigma)^{\dagger}_{\ulchT}$ belongs to
\[
\bigcap_{\substack{\kappa \in X^*(\ulT), \ w \in \ulW: \\(t^{\kappa}w)_{\dom}\leq w_0t^{-\rho}u}} t^{\xi+\rho}t^{\kappa}w \cdot \Adm(\rho).
\]
But Lemma \ref{lem:intersection of admissible} below shows this intersection is exactly $t^{\xi+\rho}\wt{\ulW}_{\leq w_0u}$.
\end{proof}

\begin{lemma}\label{lem:intersection of admissible}
Let $u\in\ul{\wt{W}}_1$. Then we have 
\[
\bigcap_{\substack{\kappa \in X^*(\ulT), \ w \in \ulW: \\(t^{\kappa}w)_{\dom}\leq w_0t^{-\rho}u}} t^{\kappa}w\Adm(\rho)= \ul{\wt{W}}_{\leq w_0u}.
\]
\end{lemma}
\begin{proof} 
Suppose $\wt{x}$ belong to the LHS. Since $(\ul{W}t^{\kappa}w)_{\dom}=(t^{\kappa}w)_{\dom}$, we see that 
\[(w_0t^{-\rho}u)^{-1}\sigma\wt{x}\in \Adm(\rho)\]
for all $\sigma \in \ulW$. In particular, for an appropriate choice of $\sigma$, we get
\[(w_0t^{-\rho}u)^{-1}\sigma\wt{x}=(w_0t^{-\rho}u)^{-1}w_0\wt{y}\]
with $\wt{y}\in\ul{\wt{W}}^+$. But now \cite[Proposition 2.1.6 (4)]{LLLM22} implies that 
$\wt{y}\leq u$, hence $\wt{x}=\sigma^{-1}w_0\wt{y}\leq w_0u$.

Conversely, if $\wt{x}\in \ul{\wt{W}}_{\leq w_0u}$, then for any $\kappa,w$ indexing the intersection in the RHS, there is a suitable $\sigma\in \ul{W}$ such that
\[
(t^{\kappa}w)^{-1}\wt{x}=((t^{\kappa}w)_{\dom})^{-1}\sigma \wt{x}\leq (w_0t^{-\rho}u)^{-1}(w_0u)=u^{-1}t^{\rho}u
\]
since $(t^{\kappa}w)_{\dom}\leq w_0t^{-\rho}u$, $\sigma\wt{x}\leq w_0u$, and the second-to-last product is a reduced factorization by \cite[Lemma 2.1.4]{LLLM22}. This shows that $\wt{x}\in t^{\kappa}w\Adm(\rho)$.
\end{proof}

\subsection{A reconstruction algorithm} The idea to prove \eqref{eq: equivariant reformulation of uniqueness} is that we should be able to reconstruct each side from the Breuil-M\'{e}zard relations. Here we explain how to carry out this reconstruction process.
\begin{situation}\label{situation:algorithm}
Suppose we have a collection of cycles $Z(u,\xi)\in \mBMulT(\ulY^{\varepsilon = 0}_{\gamma})$ parametrized by a subset of $\ul{\wt{W}}_1^+\times \ulC_0$ and an integer $h$ such that: 
\begin{enumerate}
\item For any $(w,\mu)$ such that $\mu$ is $2h_\rho$-generic, we have for some $m_{u, \nu}^{\mu, w} \in \sph$,
\begin{equation} \label{eq:algorithm}
t^{\mu}w\cdot \musupp[\wh{Z}_1(p\rho)]_{\ulchT}=\sum_{u,\nu} m_{u,\nu}^{\mu,w} Z(u,\mu-\rho+w\nu)
\end{equation}
where
\begin{itemize}
\item The sum runs over $(u,\nu)\in \ul{\wt{W}}_1\times X^*(\ulT)$ such that $(t^{-\nu})_{\dom}\leq w_0t^{-\rho}u$, and each term in the sum is defined. 
\item When $(t^{-\nu})_{\dom}=t^\rho w_0u$, we have $m_{u,\nu}^{\mu,w}=1$. 
\end{itemize}
\item \label{it:situation 2}For every $u$ and every $h$-generic $\xi$, $Z(u,\xi)$ is defined and has equivariant support in $t^{\xi+\rho} \cdot \ul{\wt{W}}_{\leq w_0u}$.
\end{enumerate}
\end{situation}

\begin{prop} \label{prop:uniqueness from support bound} In Situation \ref{situation:algorithm}, for any $\xi$ which is $(h+3h_\rho)$-generic, $Z(u,\xi)$ is uniquely determined.
\end{prop}

Proposition \ref{prop:uniqueness from support bound} essentially follows from the recursive algorithm to compute generic Breuil--M\'ezard cycles in terms of Emerton-Gee stacks in \cite[\S 8.6.1]{LLLM22}. For the convenience of the reader, we will adapt this algorithm to the more combinatorial setting of equivariant homology. This gets rid of the inputs from patching in \emph{loc.cit.}, and our reformulation should be more practical for computer implementation.

By item (\ref{it:situation 2}), it is enough to compute the component of $\Loc^{\ulchT}(Z(u,\xi))$ at each $t^{\xi}\wt{x}\in t^{\xi}\ul{\wt{W}}_{\leq w_0u}$. If $C$ is an equivariant cycle class and $\wt{z}\in \ul{\wt{W}}$, we will use the short hand $m_{\wt{z}}(C)$ for the component of $\Loc^{\ulchT}(C)$ at $\wt{z}$.

Our algorithm will be based on recursion for the following notion of defect. 

\begin{defn} Let $\wt{z}=t^{\xi}\wt{x}\in \ul{\wt{W}}$ where $\wt{x}\leq w_0t^{-\rho}u $. The \emph{defect} of $\wt{z}$ with respect to $Z(u,\xi)$ is defined to be $\delta_{\wt{z}}(u,\xi) :=\ell(u)-\ell(\wt{x}_{\dom})\geq 0$.
\end{defn}

We have the following key recursion relation, which expresses $m_{\wt{z}}(Z(u,\xi))$ in terms of the $m_{\wt{z}}$ for lower defect situations:

\begin{lemma} \label{lem:key recursion} Suppose we are given $\wt{z}\in \ul{\wt{W}}$, $(u,\xi)\in \ul{\wt{W}}_1\times \ulC_0$ such that $\wt{z}\in t^{\xi+\rho} \ul{\wt{W}}_{\leq w_0u}$. Write
\[
\wt{z}=t^{\xi+\rho}\sigma w_0\wt{x} \hspace{.5cm} \text{with} \quad \sigma\in \ulW, \quad \wt{x}\in \ul{\wt{W}}^+.
\]
Assume $\wt{z}(0)$ is $(\max\{h,h_\rho\}+2h_\rho)$-generic.

Let $\kappa \in X^*(\ulT)$ and $w \in \ulW$ be such that 
\[ \sigma w_0t^{-\rho}u=t^\kappa w \]
and set $\mu:=\xi+\kappa - \rho$.
Then
\[
m_{\wt{z}}(t^{\mu}w\cdot \musupp[\wh{Z}_1(p\rho)]_{\ulchT})=m_{\wt{z}}(Z(u,\xi))+\sum_{u',\xi'}  m^{\mu,w}_{u',w^{-1}(\xi'-\mu+\rho)}m_{\wt{z}}(Z(u',\xi'))
\]
where the sum runs over $(u', \xi') \in \wt{\ul{W}}_1 \times \ulC_0$ such that $\wt{z}\in t^{\xi'+\rho}\ul{\wt{W}}_{\leq w_0u'}$ and $\delta_{\wt{z}}(u',\xi')<\delta_{\wt{z}}(u,\xi)$.
\end{lemma}
\begin{proof}
Note that 
\[
w^{-1}(-\mu+\wt{z}(0))=(t^{\xi+\kappa}w)^{-1}\wt{z}(0)=
(w_0t^{-\rho}u)^{-1}(t^{\xi + \rho }\sigma)^{-1}t^{\xi + \rho}\sigma w_0\wt{x}(0)=
(w_0t^{-\rho}u)^{-1}w_0\wt{x}(0).
\]
By \cite[Proposition 2.1.6]{LLLM22}, we have $(w_0t^{-\rho}u)^{-1}w_0\wt{x}\in \Adm(\rho)$, hence $\mu-\wt{z}(0)$ is $h_\rho$-small. It follows that $\mu$ is $2h_\rho$-generic. In particular, we have equation (\ref{eq:algorithm}) for our choice of $(w,\mu)$.

We now observe:
\begin{enumerate}
\item $Z(u,\xi)$ contributes to the right-hand side of (\ref{eq:algorithm}) with coefficient $m_{u,-w^{-1}\kappa}^{\mu,w}=1$. This is because $\xi=\mu-\rho-\kappa$ and
\[
(t^{\kappa}w)_{\dom}=(t^{w^{-1}\kappa})_{\dom}=w_0t^{-\rho}u.
\]
\item Any pair $(u',\nu')$ contributing to (\ref{eq:algorithm}) for our choice $(w,\mu)$ has the property that $\xi'=\mu-\rho+w \nu'$ is $h$-generic. This is because $(t^{-\nu'})_{\dom}\leq w_0t^{-\rho}u'$ belongs to $\ul{\wt{W}}_1$, hence $w\nu'\in -\ul{W}(t^{-\nu'}_{\dom})(0)$ is $h_\rho$-small.
\end{enumerate}
The second item in particular implies that the equivariant support of any contributing $Z(u',\xi')$ belongs to $t^{\xi' + \rho}\wt{\ulW}_{\leq w_0 u'}$.

It remains to check $\delta_{\wt{z}}(u',\xi')<\delta_{\wt{z}}(u,\xi)$ when $(u',\xi')\neq (u,\xi)$. Assume $\wt{z}$ belongs to the equivariant support of some $Z(u',\xi')$ (thus $\xi'=\xi+\kappa+w\nu'$). Then we can write
\[
\wt{z}=t^{\xi'+\rho}\wt{y}
\]
with $\wt{y}\leq w_0u'$. In particular, $\wt{y}_{\dom}\leq u'$, and $\delta_{\wt{z}}(u',\xi')=\ell(u')-\ell(\wt{y}_{\dom})$.

We have
\[
t^\xi \sigma w_0 \wt{x}=t^{\xi+\kappa+w\nu'}\wt{y}=t^{\xi}t^{\kappa}wt^{\nu'}w^{-1}\wt{y}=t^{\xi}\sigma w_0t^{-\rho}ut^{\nu'}w^{-1}(\wt{y})\]
and
\begin{equation} \label{eq:key comparison}
(w_0t^{-\rho}u)^{-1}w_0\wt{x}=t^{\nu'}w^{-1}\wt{y}=((t^{-\nu'})_{\dom})^{-1}\sigma' \wt{y}_{\dom}
\end{equation}
for some $\sigma'\in \ulW$.

Now by \cite[Lemma 2.1.4]{LLLM22}, since $\wt{x}\in \ul{\wt{W}}^+$ we have
\[
\ell((w_0t^{-\rho}u)^{-1}w_0\wt{x})=\ell((w_0t^{-\rho}u)^{-1}w_0)+\ell(\wt{x})
\]
and 
\[
\ell(t^{\rho})=\ell(u^{-1}t^{\rho}u)=\ell((w_0t^{-\rho}u)^{-1}w_0u)=\ell((w_0t^{-\rho}u)^{-1})+\ell(w_0)+\ell(u)=\ell((w_0t^{-\rho}u)^{-1}w_0)+\ell(u)
\]
so that
\[
\ell((w_0t^{-\rho}u)^{-1}w_0\wt{x})=\ell(t^{\rho})-\delta_{\wt{z}}(u,\xi).
\]
But this quantity also equals 
\begin{align*}
\ell((t^{-\nu'}_{\dom})^{-1}\sigma' \wt{y}_{\dom}) & \leq \ell((t^{-\nu'})_{\dom})+\ell(\sigma')+\ell(\wt{y}_{\dom}) \\
&\leq \ell(w_0t^{-\rho}u')+\ell(w_0)+\ell(u')-\delta_{\wt{z}}(u',\xi')=\ell(t^{\rho})-\delta_{\wt{z}}(u',\xi').
\end{align*}
We conclude that $\delta_{\wt{z}}(u',\xi')\leq \delta_{\wt{z}}(u,\xi)$. If equality occurs, then we must have 
\[
(t^{-\nu'})_{\dom}=t^{\rho}w_0 u', \quad  \sigma'=w_0.
\]
But applying \cite[Proposition 2.1.5]{LLLM22}, to the leftmost and rightmost factorization in equation \eqref{eq:key comparison} forces $u'=u$ and hence $(t^{-\nu'})_{\dom}=(t^{w^{-1}\kappa})_{\dom}=w_0t^{-\rho}u$. But this implies
\[\nu'=(w_0t^{-\rho}u)^{-1}(0)=-w^{-1}\kappa\]
so $\xi'=\xi+\kappa+w\nu'=\xi$.

\end{proof}
\begin{cor} \label{cor:recursion} In Situation \ref{situation:algorithm}, the quantity $m_{\wt{z}}(Z(u,\xi))$ is uniquely determined whenever
\begin{itemize}
\item $\wt{z}\in \ul{\wt{W}}$ satisfies $\wt{z}(0)$ is $(\max\{h,h_\rho\}+2h_\rho)$-generic, and
\item $(u,\xi)\in \ul{\wt{W}}_1\times \ulC_0$ satisfies $\wt{z}\in t^{\xi + \rho}\ul{\wt{W}}_{\leq w_0 u}$.
\end{itemize}
\end{cor}
\begin{proof} Lemma \ref{lem:key recursion} gives a recursive formula for $m_{\wt{z}}(Z(u,\xi))$ with respect to the defect $\delta_{\wt{z}}(u,\xi)\geq 0$.
\end{proof}

\begin{proof}[Proof of Proposition \ref{prop:uniqueness from support bound}]
For any $(u,\xi)$ such that $\xi+\rho$ is $(\max\{h,h_\rho\}+3h_\rho)$-generic, we know the equivariant support of $Z(u,\xi)$ is bounded by $t^{\xi+\rho}\wt{\ulW}_{\leq w_0u}$. But if $\wt{z}\in t^{\xi+\rho}\wt{\ulW}_{\leq w_0u}$, then $\xi+\rho-\wt{z}(0)$ is $h_\rho$-small, so $\wt{z}(0)$ is $(\max\{h,h_\rho\}+2h_\rho)$-generic. Thus Corollary \ref{cor:recursion} shows that each $m_{\wt{z}}(Z(u,\xi))$ are uniquely determined, hence so is $Z(u,\xi)$.
\end{proof}

\subsection{Proof of Theorem \ref{thm: uniqueness}}
Finally, we conclude the proof of Theorem \ref{thm: uniqueness}. After the reformulation in \S \ref{subsect:uniqueness reduction}, we need to check the equation
\[
\cZ_{\gamma'}^{\varepsilon =0}(F(\pi^{-1}u\bup \xi))^\dagger_{\ulchT}= t^{\xi+\rho} \cdot \musupp[\wh{L}_1(u \bup 0)]_{\ulchT}
\]
whenever $\xi$ is $6h_\rho$-generic.

We are in Situation \ref{situation:algorithm} for the cycles $Z(u,\xi) :=\cZ_{\gamma}^{\varepsilon =0}(F(\pi^{-1}u\bup \xi))^\dagger_{\ulchT}$ with the choice of $h=3h_\rho$ by Lemma \ref{lem:effective support bound}. On the other hand, we are also in Situation \ref{situation:algorithm} for the cycles $Z'(u,\xi) :=t^{\xi+\rho} \cdot \musupp[\wh{L}_1(u \bup 0)]_{\ulchT}$ and the choice $h=0$, with the same structure constants $m_{u,\nu}^{\mu,w}\in \Z\subset \sph$. Thus Proposition \ref{prop:uniqueness from support bound} applies and gives the desired equality.

\appendix

\part{Appendices}



\section{The microlocal support of baby Vermas \\ by R. Bezrukavnikov, P. Boixeda Alvarez, T. Feng, B. Le Hung}\label{app: A}
We will prove Theorem \ref{thm: baby verma equals limit cycle} by showing that $\musupp [\wh{Z}_1(p \rho)]$ satisfies the conditions of the Recognition Principle \ref{prop: uniqueness of cycle}. Conditions (2) and (3) are not difficult to verify by direct computation. One difficulty is that the ``eigenclass'' condition is less accessible; indeed, at first glance it seems to require calculating the equivariant fundamental class, which was the difficulty in the first place. However, it admits a more conceptual interpretation in terms of the interplay between the two actions of $\wt{W}$ restricted to $X^*(T)$. The key point is then that under BMR localization, baby Vermas are localized to \emph{skyscraper sheaves} in the coherent realization. This suggests that they are ``Hecke eigensheaves'' in some sense, and making this precise leads to the ``eigenclass'' condition. 

\subsection{Eigenclass condition}
We use the notation of \S \ref{sec: existence of BM cycles}; in particular, we regard $T<B<G_k$, etc. as being over $\F_p$ in this part. The $G_k$-action on our fixed Borel $B < G_k$ induces an isomorphism $G_k/B \xrightarrow{\sim} \cB$. The fixed points for the $T$-action on $\cB$ are identified with $W$, with $w \in W$ corresponding to $wB \in \cB$. We let $\cO_{wB}$ be the skyscraper sheaf at $wB$ viewed as a point in the zero section of $\wt{\cN} := T^*(\cB)$. Since $w B$ is $T$-fixed, $\cO_{wB}$ carries a native $T$-equivariant structure, which we use to view it as an object of $\Coh_{\cB}^T(\wt{\cN})$. 

\begin{example}[BMR localization of baby Verma]\label{ex: bmr of Z(0)}
As explained in Example \ref{ex: graded baby verma BMR}, $\wh{Z}_{\mf{b}}(2 \rho)$ localizes to the class of the skyscraper sheaf $[\cO_{B}] \in \Coh_{\cB}^T(\wt{\cN})$. We will calculate the BMR localization of $[\wh{Z}_{w_0 \mf{b}}(p \rho)]$. To this end, observe by comparing characters that we have 
\[
[\wh{Z}_{\mf{b}}(2\rho)] =  [\wh{Z}_{w_0 \mf{b}}(2p \rho)] \in \rK (\Coh_{\cB}^T(\wt{\cN}))
\]
where $w_0 \in W$ is the longest Weyl element. Hence we may write $\wh{Z}_{\mf{b}}(p \rho) = t^{-\rho} \cdot_p \wh{Z}_{\mf{b}}(2p\rho)$. Then using the equivariance of $\musupp$ for the $(\wt{W}, \cdot_p)$-action, we find that 
\begin{equation}\label{eq: ss of z}
\musupp [\wh{Z}_{w_0\mf{b}}(p \rho)]  = t^{-\rho} \cdot [\cO_{\mf{b}}] \in \rK(\Coh_{\cB}^T(\wt{\cN})).
\end{equation}
Recall from \S \ref{sssec: graded frob kernel} that for $\lambda \in X^*(T)$ we defined the representation $\wh{Z}_1(\lambda) \in \Rep^0(G_1 T)$ to correspond to $\wh{Z}_{w_0 \mf{b}}(\lambda)$. Hence we may rewrite \eqref{eq: ss of z} as 
\begin{equation}\label{eq: ss of z_1}
\musupp [\wh{Z}_1(p\rho)]  = t^{-\rho} \cdot [\cO_{\mf{b}}] \in \rK(\Coh_{\cB}^T(\wt{\cN})).
\end{equation}
\end{example}

Recall that $X^*(T)$ acts on $\rK(\Coh^{T}_{\cB}(\wt{\cN})) \cong \rK(\Coh^T(\cB))$ in two ways through the embedding $X^*(T) \inj \wt{W}$. 
\begin{enumerate}
\item For the $\cdot$-action, $\lambda \in X^*(T)$ acts by tensoring with the $T$-equivariant line bundle $\cO\tw{\lambda}$, which is the pullback of the $T$-equivariant line bundle on a point corresponding to the character $\lambda$ of $T$.  
\item For the $\bu$-action, $\lambda \in X^*(T)$ acts by tensoring with the $G$-equivariant line bundle $\cO(\lambda)   = G \times^{B} \A^1$ on $\cB$, where $B$ acts on $\A^1$ through the character $\lambda$. 

\end{enumerate}
Moreover, it is evident that the two actions of $X^*(T)$ commute with each other.

\begin{lemma}\label{lem: two actions on skyscraper} For all $w \in W$ and $\mu \in X^*(T)$, on the object $t^{\mu} \cdot \cO_{wB} \in \Coh^{T}(\cB)$, we have 
\[
t^{w\lambda} \cdot (t^{\mu} \cdot [\cO_{wB}]) = (t^{\mu} \cdot [\cO_{wB}]) \bu t^{\lambda}  \in \rK (\Coh^T(\cB)).
\]
\end{lemma}

\begin{proof} We immediately reduce to the case $\mu = 0$ since the two actions of $X^*(T)$ commute. By definition, $t^\lambda \cdot \cO_{wB}  = \cO\tw{\lambda} \otimes \cO_{wB}$ is $\cO_{wB}$ equipped with $T$-equivariant structure placing it in graded degree $\lambda$.

For $g \in T$, we have $g w B = w (w^{-1} g w B ) \in \cB$. Therefore, the left translation of $g \in T$ on $\cO_{wB} \bu t^{\lambda} = \cO_{wB} \otimes \cO(-\lambda)$ acts as multiplication by $\lambda (w^{-1} g w) = (w \lambda)(g)$. Hence $\cO_{wB} \otimes \cO(-\lambda)$ is $\cO_{wB}$ in graded degree $w\lambda$ for the left translation action of $t$. 
\end{proof}

\begin{prop} \label{prop: eigenproperty verma} The translation action of $X_*(\chT)$ on 
$\musupp[\wh{Z}_1(p \rho)]$ agrees with the affine Springer action of $X_*(\chT)$ on $\musupp[\wh{Z}_1(p \rho)]$. In particular, the equivariant support of $\musupp[\wh{Z}_1(p \rho)]_{\chT}$ is concentrated on the translation elements $X_*(\chT) \subset \dWext$.
\end{prop}
\begin{proof}
By \eqref{eq: ss of z_1} and Lemma \ref{lem: two actions on skyscraper} we see that 
\[
t^\lambda \cdot  \musupp[\wh{Z}_1(p \rho)] = \musupp[\wh{Z}_1(p \rho)] \bu t^{\lambda}  \in \rK (\Coh^T(\cB)) 
\]
for all $\lambda \in X^*(T)$. This gives the first statement of the Proposition. This implies that the two actions of $X_*(\chT)$ -- one restricted from $(\dWext, \cdot)$ and the other restricted from $(\dWext, \bu)$ -- also agree after lifting to $\chT$-equivariant Borel--Moore homology, so that
\[
t^{-\nu} \cdot  \musupp[\wh{Z}_1(p \rho)]_{\chT} \bu t^{\nu} = \musupp[\wh{Z}_1(p \rho)]_{\chT}.
\]
By the description of the two $\dWext$-actions in $\chT$-equivariant Borel--Moore homology (cf. \S \ref{ssec: equivariant actions}), this implies that if $\wt{w}$ belongs to the equivariant support of $\musupp[\wh{Z}_1(p \rho)]_{\chT}$ then so does $t^{\nu}\wt{w}t^{-\nu}=t^{\nu-w\nu}\wt{w}$ for any $\nu\in X_*(\chT)$, where $w$ is the projection of $\wt{w}$ to $W$. But if $w\neq 1$, the set
 $t^{\nu-w\nu}\wt{w}$ is unbounded, a contradiction; this concludes the proof. 
\end{proof}


\subsection{Support bound and normalization conditions}

\begin{prop}\label{prop: supp of baby verma}
Assume that $p > 2 h_\rho$. Then $\musupp[\wh{Z}_1(p\rho)]_{\chT}$ has equivariant support in $\Adm(\rho) \subset \dWext$. 
\end{prop}

\begin{proof}By \cite[Lemma 10.1.5]{GHS18}, the simple constituents of $\wh{Z}_1(p\rho)$ consists of $\wh{L}_1(p \nu +  \wt{w} \bup 0  )$ where
\[\sigma t^{-\nu}\uparrow w_0t^{-\rho}\wt{w}\]
for all $\sigma\in W$. This is equivalent to $\sigma t^{-\nu} \uparrow w_0t^{-\rho}\wt{w}$ for the choice of $\sigma$ such that $\sigma t^{-\nu}$ is dominant. In turn, this is equivalent to $\sigma t^{-\nu}\leq w_0t^{-\rho} \wt{w}$ for this choice. 
By Lemma \ref{lem: support bound simple}, it suffices to show that $t^{\nu}\dWext_{\leq w_0\wt{w}}\subset \Adm(\rho)$ for such $\wt{w},\nu,\sigma$.

Let $w$ be the projection of $\wt{w}$ to $W$. Then we have the reduced factorization
\[
t^{w^{-1}\rho}=  (w_0t^{-\rho}\wt{w})^{-1} w_0\wt{w}.
\]
We have 
\[
(\sigma t^{-\nu})^{-1}\leq (w_0t^{-\rho}\wt{w})^{-1}, \quad \sigma\leq w_0
\]
hence for any $u \leq w_0 \wt{w}$, 
\[
t^{\nu}u=(\sigma t^{-\nu})^{-1}\sigma u \leq (w_0t^{-\rho}\wt{w})^{-1} w_0\wt{w}=t^{w^{-1}\rho}
\]
belongs to $\Adm(\rho)$. This shows the desired inclusion $t^{\nu}\dWext_{\leq w_0\wt{w}}\subset \Adm(\rho)$.

\end{proof}

\begin{prop}\label{prop: normalization of baby verma} The coefficient of $\Loc^{\chT}\left(\musupp[\wh{Z}_1(p\rho)]_{\chT}\right)$ at $[t^{\rho}]$ is $1/\beta \in \Frac(\sph)$.
\end{prop}
\begin{proof} From the proof of Proposition \ref{prop: supp of baby verma}, the only composition factor $\wh{L}_1(p \nu +  \wt{w} \bup 0  )$ of $\wh{Z}_1(p\rho)$ whose microlocal support can contribute to the coefficient of $[t^{\rho}]$ is the one with $\wt{w}=1$ and $w_0t^{-\nu}=w_0t^{-\rho}$, i.e., $\wh{L}_1(p \nu + \wt{w} \bup 0 )=\wh{L}_1(p\rho)$.
Since $\wh{L}_1(p\rho)$ occurs in $\wh{Z}_1(p\rho)$ with multiplicity one, and $\musupp[\wh{L}_1(p\rho)]_{\chT}=t^\rho[\chG/\chB]_{\chT}$ has coefficient $1/\beta$ by Example \ref{ex: flag variety}, we are done. 
\end{proof}

\begin{proof} [Proof of Theorem \ref{thm: baby verma equals limit cycle}]
The equality 
\[
\Loc^{\chT}\left( \fsp_{p\rightarrow 0}[\rX^{\varepsilon=0}_\gamma(\rho)]_{\chT} \right)  = \frac{1}{\beta} \sum_{w \in W} \sgn(w) [t^{w \rho}]
\]
comes from \eqref{eq: loc limit cycle}. As observed in Lemma \ref{lem: limit cycle recognition}, this shows that $\fsp_{p\rightarrow 0}[\rX^{\varepsilon=0}_\gamma(\rho)]_{\chT}$ satisfies the properties of Proposition \ref{prop: uniqueness of cycle}. 

Propositions \ref{prop: eigenproperty verma}, \ref{prop: supp of baby verma} and \ref{prop: normalization of baby verma} show that $\musupp[\wh{Z}_1(p\rho)]_{\chT}$ also satisfy the characterizing properties in Proposition \ref{prop: uniqueness of cycle}. Hence they coincide. 
\end{proof}

\section{Generic tame potential crystalline stacks for unramified groups \\ by B. Le Hung and Z. Lin}
\label{app: EG stacks}
Let $(G, B, T, \{X_{\alpha}\})$ be an unramified pinned group over $\Q_p$ whose Langlands dual pinned group is $(\bG, \bB, \bT, \{Y_{\alpha^\vee}\})$, which we regard over $\Z_p$. We assume that $G$ has connected center, which is equivalent to $\bG$ having simply connected derived subgroup. Since $G$ is unramified, the Galois action on $(\bG, \bB, \bT, \{Y_{\alpha^\vee}\})$ is equivalent to the data of a finite order automorphism, which we denote by $\pi^{-1}$. (The notation is compatible with that in \S \ref{ssec:Frobenius}.) 
Write $W=N_{\bG}(\bT)/\bT$
for the Weyl group. 
Fix a coefficient field $E$ that splits $G$
with ring of integers $\cO$ and residue field $\F=\cO/\varpi$.
Frobenius and Galois groups always act trivially on $E$.

Denote by $\cG$ the group scheme $\bG\otimes_{\Z_p}\Z_p[v]$.
Set
\begin{align*}
\bG_{\F_p(\!(v)\!)} & := \cG\otimes_{\Z_p[v], p \mapsto 0} \F_p((v))\text{, and}\\
\bG_{\Q_p} &:= \cG\otimes_{\Z_p[v], v\mapsto -p} \Q_p,
\end{align*}
which are split groups over $\F_p((v))$ and $\Q_p$ respectively.
Our choice of pinning determines
(through the Chevalley valuation)
an apartment
\[
A(\bG_{\F_p(\!(v)\!)})\xrightarrow{\cong}A(\bG_{\Q_p})\cong
X_*(\bT)\otimes \R
\]
of the enlarged Bruhat-Tits building, together with a hyperspecial vertex $\ol o\mapsto o$ therein.

The extended affine Weyl group $\wt{W}$ of $\bG_{\F_p(\!(v)\!)}$ can be identified with
$$N_{\bG(\F_p((v)))}(\bT(\F_p((v))))/\bT(\F_p[[v]]),$$
and thus acts on the affine apartment $A(\bG_{\F_p(\!(v)\!)})$.
The automorphism $\pi^{-1}$ of the pinned root datum gives rise to an endomorphism
$\varphi:=p\pi^{-1}$ of $\wt{W}$ characterized by the formula
\begin{equation}
\varphi(wv^{\mu})=\pi^{-1}(w)v^{p\pi^{-1}\mu}
\end{equation}
We also denote by $\varphi$ the endomorphism $p\pi^{-1}$ of $A(\bG_{\F_p(\!(v)\!)})$. Then the natural action of $\wt{W}$ on $A(\bG_{\F_p(\!(v)\!)})$ is compatible with $\varphi$.

\subsection{Breuil-Kisin modules with tame descent data}
From now on, we work in the following situation. 

\begin{situation} \label{running assumption} Fix a generator $\psi\in \varprojlim \F^\times_{p^i-1}$.
We fix a regular tame inertial $L$-parameter $\tau$ (cf. Definition \ref{defn:tame-inertial-param}), and assume it has a lowest alcove presentation. We note that this assumption is equivalent to asking that any lift of $\tau$ to $X_*(\chT)$ has trivial stabilizer in $\wt{W}$. We will fix a lowest alcove presentation $\wt{w}(\tau)=v^{\mu}w$, which fixes a lift $x=x_{\wt{w}(\tau)}\in A(\bG_{\F_p(\!(v)\!)})$ of $\tau\in((X_*(\bT)\otimes \Q/\Z)/W)^{\varphi=1}$. Under our assumption, the $x$ and $\wt{w}(\tau)$ uniquely determine each other, and we abusively say $x$ is $m$-generic for some integer $m$ if $\wt{w}(\tau)$ is $m$-generic in the sense of Definition \ref{defn:presentation}.

We have $x\in X_*(\bT)\otimes \frac{1}{p^f-1}\Z$
for a sufficiently large integer $f$, and we set $e=p^f-1$.
We assume $E$ is large enough to contain
an embedding of the Galois closure of
$L_x=\Q_{p^{f}}( (-p)^{\frac{1}{e}})$.
Note that our choice of $\psi$ gives a preferred embedding $\iota_{\psi}:\F_x\inj \F$, where $\F_x$ is the residue field of $L_x$.
\end{situation}

\subsubsection{Pappas-Zhu group schemes}\label{subsect: PZ group schemes}
By Pappas-Zhu \cite[Theorem 3.1]{PZ13}, $x$ gives rise to a group scheme $\cG_x$ over $\Z_p[v]$
such that
\begin{itemize}
\item the fiber $\cG_x$ at $p=0$ (resp. at $v=-p$) is the connected stabilizer of $x$ in $\bG_{\F_p(\!(v)\!)}$ (resp. of $\bG_{\Q_p}$), 
\item $\cG_x[v^{-1}]=\chG\otimes_{\Z_p} \Z_p[v,v^{-1}]$.
\end{itemize}
The group scheme $\cG_x$ gives rise to the positive loop groups $L^+\cG_x$ and $L\cG_x$ by completing $\Z_p[v]$ at $v+p$,
followed by base change along $\Spec\cO\to\Spec\Z_p$.
Their functor of points have the following simple description 
\begin{itemize}
\item 
$L^+\cG_x(R) = \{A\in \bG(R[\![v+p]\!])| A \mod v \in \bB (R)\}$,
\item
$L\cG_x(R) = \bG(R(\!(v+p)\!))$.
\end{itemize}
Note that on $p$-adically complete $R$, $L^+\cG_x(R)$ is the subgroup of $L\cG_x(R)$
consisting of elements $A$ such that $v^{-x} A v^{x}\in \bG([\![v^{\frac{1}{e}}]\!])$. 
Their Lie algebras have functor of points
\[\Lie L\cG_x(R)=\chfg(R(\!(v+p)\!))\]
\[\Lie L^+\cG_x( \F)=\chft[\![v+p]\!]\oplus (\bigoplus_{\alpha\in \Phi^+} v\chfg_\alpha[\![v+p]\!])\oplus (\bigoplus _{\alpha\in \Phi^-} \chfg_\alpha[\![v+p]\!])\]
Here we use the convention that for any $\cO$-algebra $A$, $\chfg(A):=\chfg\otimes_\cO A$ (note that as in section \ref{sssec: positivity}, our convention is that the root groups of $\bB$
correspond to negative roots $\Phi^-$, rather than $\Phi^+$).

When the test ring $R$ is $p$-adically complete, $R[\![v+p]\!]=R[\![v]\!]$, so that the endomorphism $\varphi_0$ that fixes $R$ and sends $v\mapsto v^p$ is defined. This induces endomorphisms of $L^+\cG_x(R), L\cG_x(R), \Lie L\cG_x(R),  L^+\cG_x(R)$ which we continue to denote by $\varphi_0$. We will also let $\varphi=\pi^{-1}\varphi_0=\varphi_0\pi^{-1}$ on all these objects. We note that this notation is compatible with the previously defined $\varphi$ on $A(\bG_{\F_p(\!(v)\!)})$.

By \cite[Proposition 5.3]{PZ13}, the twisted affine Grassmannian $\Gr_{\cG_x}=L\cG_x/L^+\cG_x$ is an ind-proper ind-scheme over $\cO$.
Given a dominant $\lambda\in X_*(\chT)$, we have the global affine Schubert variety $\Gr_{\cG_x}^{\leq \lambda}$, and define $L\cG_{x}^{\le \lambda}$ to be its pullback to $L\cG_x$, cf \cite{Ei}.
By the description of the special fiber $\Gr_{\cG_x,\ol \F}$ (c.f. \cite[Theorem 9.3]{PZ13}), we have an affine Bruhat decomposition
\[L\cG_x^{\le \lambda}(\ol \F)=
\bigcup_{\widetilde{z}\in \Adm(\lambda)}
L^+\cG_x(\ol \F)\widetilde{z}L^+\cG_x(\ol \F)\]
where we recall from section \ref{sssec: weyl notation} the $\lambda$-admissible set
\[
\Adm(\lambda) :=  \{ \wt{w} \in \wt{W} \co \wt{w} \leq t^{w(\lambda)} \text{ for some } w \in W \}.
\]

\begin{lemma}
\label{lem:amplitude bound} Let $A\in L\cG_x^{\le \lambda}(\ol \F)$. 
\[\Ad(A)(\Lie L^+\cG_{x}(\ol \F))\subset \frac{1}{v^{h_\lambda}}\Lie L^+\cG_x(\ol \F)\]
\end{lemma}
\begin{proof} The claim reduces to the case $A=\wt{w}\in \Adm(\lambda)$ where it follows by a direct computation on affine root spaces.
\end{proof}

\begin{lemma}
\label{lem:tangent space of Schubert} Let $A\in L\cG_x^{\le \lambda}(\ol \F)$. Then the tangent space of $L\cG_{x,\F}^{\le \lambda}A^{-1}$ at the identity belongs to the subspace
\[\frac{1}{v^{h_\lambda}}\Lie L^+\cG_x(\ol \F) \subset \Lie L\cG(\ol \F)\]
\end{lemma}
\begin{proof} Without loss of generality, we can assume $A=\wt{w}\in \Adm(\lambda)$. Then the tangent space in question is stable under the action of the extended torus $\chT\times \G_{m,\mathrm{rot}}$, hence is the direct sum of its intersection with $\chft(\!(v)\!)$ and $\chfg_\alpha(\!(v)\!)$. The argument is now similar to \cite[Proposition 5.1.1]{BL23}. Any tangent vector $X$ in our tangent space satisfies
\[[X, \Lie L^+\cG_x(\ol \F)]\subset \frac{1}{v^{h_\lambda}}\Lie L^+\cG(\ol \F)\]
If $X\in \chfg_\alpha(\!(v)\!)$, bracketing with elements of $\chft$ gives $X\in \frac{1}{v^{h_\lambda}}\Lie L^+\cG_x(\ol \F)$, while if $X\in \chft(\!(v)\!)$, bracketing with a suitable element in $\chfg_\alpha$ yields the same.
\end{proof}

\begin{defn}\label{defn:twisted conjugation}
Assume Situation \ref{running assumption}. Define $Y^{\le \lambda,\tau}$ to be the 
\[[(L\cG^{\le \lambda})^{\wedge_p}_{x}/_{\varphi, x}(L^+\cG_{x})^{\wedge_p}]\]
where the action of the $p$-adic completion $(L^+\cG_x)^{\wedge_p}$ of $L^+\cG_x$ is the $x$-twisted $\varphi$-conjugation action given by
\[g \cdot_{\varphi, x} h := \varphi_x(g) hg^{-1}.\]
where $\varphi_x=\Ad(\wt{w}(\tau)^{-1})\circ \varphi$.

\end{defn}

\subsubsection{The Emerton-Gee stacks}\label{subsect: EG stacks}
Write $\lsup LG=\bG\rtimes \{1,\sigma,\dots,\sigma^{f-1}\}$.
Here we choose $\sigma$ so that
under the embedding $\bG\subset \lsup LG$,
we have $\pi^{-1} = \Ad(1\rtimes \sigma^{-1})$.

Recall that we have $\varphi_0:R((v))\to R((v))$, the Frobenius endomorphism which acts trivially on $R$
and send $v$ to $v^p$. This gives rise to the moduli space $\cR_{\lsup LG}$ of \'etale $\varphi_0$-modules, such that $\cR_{\lsup LG}(\F_p)$ is the equivalent to the groupoid
of Langlands parameters $G_{\Q_{p,\infty}}=\Gal(\ol{\Q}_p/\Q_{p, \infty})\to\lsup LG(\F_p)$.
When $G$ is an unramified group, $\cR_{\lsup LG}$ has a simple description:
for any $p$-adically complete $\Z_p$-algebra $R$, $\cR_{\lsup LG}(R)$ is the groupoid of $(M, M^\circ, \varphi_M)$
where
\begin{itemize}
\item[i.] $M$ is a {\it left} $\lsup LG$-torsor over $R(\!(v)\!)^{\wedge_p}$,
\item[ii.] $M^\circ$ is a {\it left} $\bG$-torsor over $\big((R\otimes_{\Z_p}W(\F_{p^f}))(\!(v)\!)\big)^{\wedge_p}$ together with an identification $M^\circ\times ^{\bG}\lsup LG=M\otimes_{\Z_p}W(\F_{p^f})$,
\item[iii.]
$\varphi_M:\varphi_0^*M\xrightarrow{\cong}M$ is a $\lsup LG$-torsor isomorphism such that $\varphi_M(M^\circ)\subset M^\circ (1\rtimes \sigma)$.
\end{itemize}
We call objects of $\cR_{\lsup LG}(R)$ {\it \'etale $\varphi_0$-modules with rigidified $\lsup LG$-structure}.
The Emerton-Gee stack $\cX^{\EG}$ for $\lsup LG$ is obtained by a variant of the above construction: it is the stack of $\varphi_0$-module where $\Q_{p,\infty}$ is replaced by the cyclotomic extension $\Q_p(\zeta_{p^{\infty}})$, together with additional commuting action of $\Gamma=\Gal(\Q_p(\zeta_{p^{\infty}})/\Q_p)$.
As in \cite{EG23}, there is a natural map $\cX^{\EG}\to \cR_{\lsup LG}$ which correspond to ``restricting parameters from $G_{\Q_p}$ to $G_{\Q_{p,\infty}}$''.


We suppress the phrase ``with (rigidified) $\lsup LG$-structure'' when it is clear from the context.

\begin{remark}
Since we make frequent use of Tannakian formalism, the notion of \( \varphi_0 \)-modules is more convenient than that of \( \varphi \)-modules---especially because, for split \( \GL_n \), the geometric Frobenius \( \pi^{-1} \) is trivial.

To give the reader a sense of the conditions:
without condition (iii), $\cR_{\lsup LG}(\F_p)$
becomes the groupoid of all group
homomorphisms $G_{\Q_{p, \infty}}\to\lsup LG(\F_p)$,
up to $\bG(\F_p)$-conjugacy;
without the rigidification $M^\circ$ and condition (ii, iii),
$\cR_{\lsup LG}(\F_p)$
becomes the groupoid of all group
homomorphisms $G_{\Q_{p, \infty}}\to\lsup LG(\F_p)$,
up to $\lsup LG(\F_p)$-conjugacy.
\end{remark}

Let $R$ be a $p$-adically complete $\cO$-algebra and
let $A\in L^{\le \lambda}\cG_x(R)$.
We define an element $\iota_x(A)=(M_A, M_A^\circ, \varphi_{F_A})\in \cR_{\lsup LG}(R)$ as follows:
\begin{itemize}
\item $M_A=\lsup LG\otimes_{\Z_p}R(\!(v)\!)$ is the trivial $\lsup LG$-torsor over $\Spec R(\!(v)\!)$;
\item $M_A^\circ = \bG\otimes_{\Z_p}(R\otimes_{\Z_p}W(\F_{p^{f}}))(\!(v)\!)$;
\item $\varphi_{F_A}: h\mapsto \varphi_0(h) (1\rtimes \sigma)\wt{w}(\tau)A$.
\end{itemize}
By abuse of notation, we write $\iota_x(A)=(1\rtimes \sigma)\wt{w}(\tau)A$.

Note that if $A, B\in L^{\le \lambda}\cG_x(R)$ and $g\in L\cG(R)$,
then $B=g\cdot_{\varphi, x}A$ if and only if
$\iota_x(B) =g \cdot_{\varphi_0} \iota_x(A):=\varphi_0(g) \iota_x(A) g^{-1}$. In this situation, $g$ induces a natural arrow between $\iota_x(A)$ and $\iota_x(B)$ in $ \cR_{\lsup LG}(R)$.
Hence the morphism $\iota_x:(L^{\le \lambda}\cG_x)^{\wedge_p}\to \cR_{\lsup LG}$ factors through $[(L\cG^{\le \lambda}_x)^{\wedge_p}/_{\varphi, x}(L^+\cG_{x})^{\wedge_p}]$.

\begin{prop} Let $\lambda\in X_*(\bT)^+$, and assume that $x$ is $(h_\lambda+1)$-generic. 
Then the morphism
$$
[(L\cG^{\le \lambda}_x)^{\wedge_p}/_{\varphi, x}(L^+\cG_{x})^{\wedge_p}]
\to \cR_{\lsup LG}$$
is a monomorphism.
\label{prop: monomorphism}
\end{prop}
This is the generalization of \cite[Proposition 5.4.3]{LLLM22} in our setting.

As in \emph{loc.cit.}, we begin with some preparatory lemmas:
\begin{lemma}\label{lem: mod-p-monomorphism} Assume that $x$ is $(h_\lambda+1)$-generic.
 If $A_1,A_2\in L^{\leq \lambda}\cG_x(\ol \F)$, and $X\in L\cG_x(\ol \F)$ such that
\[
A_1X=\varphi_x(X)A_2
\]
Then $X\in L^+\cG(\ol \F)$.
\end{lemma}
\begin{proof} Write $X=I_1\wt{w}I_2$ for $I_1,I_2\in L^+\cG_x(\ol \F)$ and $\wt{w}\in \wt{W}$.
Write $\wt{w}(\tau)=v^{\mu}w$, so that $\mu$ is dominant.
Since $\mu$ is in the interior of the lowest dominant alcove, 
for each $Y\in v^{-\mu}\varphi(L^+\cG_x(\ol \F))v^{\mu}$,
we have $Y\mod v\in \bT(\F_p)$.
Thus
we have
$$
\varphi_x(L^+\cG_x(\ol \F))
=w^{-1} v^{-\mu}\varphi(L^+\cG_x(\ol \F))v^{\mu}w
\subset L^+\cG_x(\ol \F).
$$
Our hypothesis implies 
\[\varphi_x(\wt{w})\in L\cG^{\leq \lambda}(\ol \F)\wt{w}\big(L\cG^{\leq \lambda}(\ol \F)\big)^{-1}\].
Write $\wt{w}=v^{\kappa}\sigma$ for $\sigma\in W$, so that $\varphi_x(\wt{w})=w^{-1}v^{-\mu+p\pi^{-1}(\kappa)+\pi^{-1}(\sigma)(\mu)}\pi^{-1}(\sigma)w$. 

The above equation implies that $v^{p\pi^{-1}(\kappa)-\mu+\pi^{-1}(\sigma)(\mu)}\in  W\Adm(\lambda)v^{\kappa}W\Adm(-\lambda)W$, and hence
\begin{equation}\label{eq:convex hull}
p\kappa-\pi(\mu)+\sigma(\pi(\mu))\in \Conv(\pi(\lambda-w_0\lambda+\kappa_{dom}))
\end{equation}
First, this together with the fact that $\mu$ is $(h_\lambda+1)$-generic implies
\[ph_{\kappa}\leq  2h_{\mu}+2h_{\lambda}+h_{\kappa}< 2p-2+h_{\kappa}\]
so that $h_{\kappa}\leq 1$.
Using this and the fact that $\mu$ is $(h_\lambda+1)$-generic, equation \eqref{eq:convex hull} shows that for any positive root $\alpha$
\[p\langle \kappa,-\alpha^\vee \rangle+h_\lambda < \langle p\kappa-\pi(\mu),-\alpha^\vee\rangle\leq \langle \sigma(\pi(\mu)),\alpha^\vee\rangle +2h_\lambda+h_\kappa\leq p+h_\lambda\]
so that $0\leq \langle \kappa, \alpha^\vee\rangle\leq 1$.

If $\langle \kappa, \alpha^\vee\rangle=1$ for $\alpha>0$ then applying equation \eqref{eq:convex hull} again we get
\[ \langle \sigma(\pi(\mu)), \alpha^\vee\rangle\leq -p+h_\mu+2h_{\lambda}+1< h_\lambda,\]
so $\sigma^{-1}(\alpha)<0$ since $\mu$ is $(h_\lambda+1)$-generic.
Conversely, if $\alpha>0$ and $\sigma^{-1}(\alpha)<0$ then 
\[\langle \kappa,\alpha^\vee \rangle\geq \langle\pi(\mu),\alpha^\vee\rangle -\langle \pi(\mu),\sigma^{-1}(\alpha)\rangle-2h_\lambda-1>0 \]
so that $\langle \kappa,\alpha^\vee \rangle =1$.

Thus $\Ad(\wt{w})=\Ad(v^{\kappa}\sigma)$ sends $\Lie L^+\cG_x$ into $\Lie L^+\cG_x$, so that $\wt{w}$ stabilizes the base alcove. But since $G$ has connected center, the image of $\wt{w}$ in $\wt{W}/\wt{W}_a=X_*(\chT)/\chQ^\vee$ is fixed by $p\pi^{-1}$, hence must be trivial since the former group is torsion-free. Thus $\wt{w}\in \wt{W}_a$ and hence $\wt{w}=1$. We conclude that $X\in L^+\cG_x(\F)$.
\end{proof}
\begin{lemma}\label{lem: mod-p-monomorphism Lie}
Let $A\in L^{\le \lambda}\cG_x(\ol \F)$
and $X\in \Lie L\cG_x(\ol \F)$.
Assume that $x$ is $(h_\lambda+1)$-generic and
\[\Ad(A)(X) - \varphi_x(X)
\in \frac{1}{v^{h_\lambda}}\Lie L^+\cG_x(\ol \F).\]
Then
\[X\in \Lie L^+\cG_x(\ol \F).\]
\end{lemma}
\begin{proof} Let $k$ be the maximal integer such that $X\in \frac{1}{v^k}\Lie L^+\cG_x(\ol \F)$. Assume for contradiction that $k>0$. Since $\Ad(A)(\Lie L^+\cG_x(\ol \F))\subset \frac{1}{v^{h_\lambda}}\Lie L^+\cG_x(\ol \F)$, looking at the component with minimal $v$-valuation, the hypothesis implies
\[pk\leq \{\max_{\alpha>0}\{\langle \mu,\alpha^\vee\rangle, p-\langle \mu,\alpha^\vee\rangle\}+k+h_{\lambda}<p-1+k\]
by the $(h_\lambda+1)$-genericity of $\mu$. But this implies $k<1$, a contradiction.
\end{proof}
\begin{proof}[Proof of Proposition \ref{prop: monomorphism}]  We need to check that for any finite type $\cO/\varpi^a$-algebra $R$, $A_1, A_2\in L^{\leq \lambda}\cG_x(R)$ and $X\in L\cG_x(R)$ such that
\begin{equation}\label{eq:equivalent phi}
A_1X=\varphi_x(X)A_2
\end{equation}
then $X\in L^+\cG_x(R)$.

Since $R$ embeds into the product of its local Artinian quotients, it suffices to treat the case $R$ is Artinian local with residue field $\ol \F$. We induct on the length of such $R$.
The base case follows from Lemma \ref{lem: mod-p-monomorphism}.

Suppose now that $R\surj R/J$ is a small extension, i.e. $J\cong \ol \F$ as $R$-modules. By the induction hypothesis, we may assume that $A_1=A_2$ mod $J$ and $X=1$ mod $J$, so that $X$ corresponds to an element $Y\in \Lie L\cG_x(\ol \F)$. Thus \eqref{eq:equivalent phi} becomes
\[\varphi_x(Y)=\Ad(\cl{A})(Y)+\delta\]
where $\cl{A}\in L^{\leq \lambda}(\ol \F)$ the common image of $A_1$, $A_2$ mod $\mf{m}_R$ and $\delta$ is the tangent vector corresponding to $A_1A_2^{-1}$. Note that $\delta$ is the tangent vector in $L\cG^{\le \lambda}\cl{A}^{-1}$ at the identity, hence by Lemma \ref{lem:tangent space of Schubert} $\delta\in \frac{1}{v^{h_\lambda}}\Lie L^+\cG(\ol \F)$. But then Lemma \ref{lem: mod-p-monomorphism Lie} shows that $Y\in \Lie L^+\cG(J)$, so that $X\in L^+\cG(R)$.

\end{proof}

\subsection{Tame potentially crystalline stacks}
For each tame inertial parameter $\tau$ and Hodge--Tate weights $\lambda$,
there exists a unique $\cO$-flat closed substack
$\cX^{\lambda,\tau}$
of $\cX^{\EG}$ such that
for each finite flat $\cO$-algebra $\Lambda$,
$\cX^{\lambda, \tau}(\Lambda)$
consists of parameters
$G_{\Q_p}\to \lsup LG(\Lambda)$
which becomes de Rham 
of Weil-Deligne
type $\tau$ and Hodge--Tate weight $\lambda$
after inverting $p$ (cf. \cite[Theorem 2.6.3]{Lin23b}). Similarly we have the stack $\cX^{\le \lambda,\tau}$ which identifies with the scheme theoretic union of $\cX^{ \lambda',\tau}$ for all dominant $\lambda'\le \lambda$.
The goal of this subsection is to give a more group theoretic description of $\cX^{\le \lambda,\tau}$ when $\tau$ is tame and generic relative to $\lambda$.

\subsubsection{Breuil-Kisin-Fargues $G_{\Q_p}$-modules and Weil-Deligne types}
\label{subsect:BKF modules}
Recall that $L_x$ is a splitting field of $\tau$.
Let $\varpi_x$ be a uniformizer of $L_x$
such that $\varpi_x^e=-p$
and set $\F_x:=\cO_{L_x}/\varpi_x$. We fix a compatible family
$\varpi_{x}^\flat\in \cO_{\C^\flat}$ of $p$-power roots of $\varpi_{x}$, which gives rise to a compatible family $p^\flat=(\varpi_x^\flat)^e$ of $p$-power roots of $-p$ in the tilt $\C^{\flat}$ of $\C=\widehat{\ol{\Q}}_p$. These determine the extensions $L_{x,\infty}/L_x$ and $\Q_{p,\infty}/\Q_p$. Note that restriction gives a natural isomorphism
\[\Gamma:= \Gal(L_{x,\infty}/\Q_{p,\infty})\cong \Gal(L_x/\Q_p).\]
We will consider the embeddings $W(\F_p)[\![v]\!]\inj W(\F_x)[\![u]\!]\inj A_\Inf=W(\cO_{\C^{\flat}})$ given by $v\mapsto u^e$ and $u\mapsto [\varpi_x^\flat]$. The action of $\Gamma$ on $A_\Inf$ preserves $W(\F_x)[\![u]\!]$, whose fixed subring is $W(\F_p)[\![v]\!]$.

Let $R$ be a $p$-adically complete topologically finite type $\cO$-algebra. Recall that a {\it Breuil-Kisin-Fargues} $G_{L_x}$-module with $\lsup LG$-structure and $R$-coefficients is a Tannakian functor $\mf{M}_{\Inf}$ from the category $^f\Rep_{\lsup LG}$ of finite dimensional algebraic representations of $\lsup LG$ to the category of Breuil-Kisin-Fargues $G_{L_x}$ modules over $A_{\Inf}$ with $R$-coefficients, as in \cite[\S 4.2]{EG23}. Equivalently, $\mf{M}_{\Inf}$ is an $\lsup LG$ torsor over $\Spec(A_\Inf\widehat{\otimes} R)$ with an appropriate $\varphi_0$-structure and $G_{L_x}$-action. We also have the similar notion of a {\it Breuil-Kisin} module with $\lsup LG$-structure and $R$-coefficient, which is an $\lsup LG$ torsor over $\Spec (W(\F_x)\otimes R)[\![u]\!]$ with an appropriate $\varphi_0$-equivariant structure. 
We will omit the phrase ``with $\lsup LG$'' in what follows.

When $R$ is a finite flat $\cO$-algebra (or a complete local Noetherian $\cO$-algebra) and $\mf{M}_x$ is a Breuil-Kisin module with $R$-coefficient, one gets an induced \'{e}tale $\varphi_0$-module
\[M_x=\mf{M}_x\otimes_{W(\F_x)[\![u]\!]} (A_\Inf)^{G_{L_{x,\infty}}},\]
corresponds to a parameter $V(M_x): G_{L_{x,\infty}}\to \lsup LG(R)$.

Now let $\Lambda$ be the ring of integers in a finite extension of $E$ and $V\in \cX^{\le \lambda,\tau}(\Lambda)$. Then by \cite[Proposition 2.2.6, Proposition 2.3.3]{Lin23b},
$V$ admits a unique Breuil-Kisin-Fargues $G_{L_x}$-lattice with descent data of type $\tau$, that is, a Breuil-Kisin-Fargues $G_{L_x}$-module $\mf{M}_\Inf$ with $\lsup LG$-structure and $\Lambda$-coefficients with an identification
\[\mf{M}_{\Inf}\otimes_{A_{\Inf}}W(\C^\flat)=V\otimes_{\Z_p}W(\C^\flat)=:M_\Inf\]
that respects $\varphi_0$-structure and $G_{L_x}$-action. Furthermore, such $\mf{M}_\Inf$ admits all descent, which by \cite[Definition 4.2.4]{EG23} and Tannakian formalism implies
there exists a unique Breuil-Kisin module $\mf M_x$ with $\lsup LG$-structure 
over $(W(\F_x)\otimes \Lambda)[\![u]\!]$
such that $\mf M_\Inf =\mf M_x\otimes_{W(\F_x)[\![u]\!]}A_\Inf$.

\begin{lemma}
(1) Under the embedding
\[\mf M_x \inj (\mf M_x \otimes_{W(\F_x)[\![u]\!]}W(\C^\flat))^{G_{L_{x,\infty}}}\cong M_\Inf^{G_{L_{x,\infty}}}\]
$\mf M_x$ is stable under the induced $\Gal(L_{x,\infty}/\Q_{p,\infty})$-action on the target.

(2)
The induced $\Gal(L_{x,\infty}/\Q_{p,\infty})$-action
on $\mf M_x \times^{\lsup LG} \lsup L\{*\}\cong \lsup L\{*\} \otimes_{\Z_p} W(\F_x)[\![u]\!]$
is the base change along $\lsup LG \to \Spec \Z_p$
of the standard $\Gal(L_{x,\infty}/\Q_{p,\infty})$-action on  $W(\F_x[\![u]\!])$.
Here $\lsup L\{*\}=\{*\}\rtimes \langle \sigma \rangle
=\langle \sigma \rangle$ is the $L$-group
of the singleton group.
\label{lem: Galois stability of Mx}
\end{lemma}

\begin{proof}

(1) The $\Gal(L_x/\Q_p)$-action on $M^{G_{L_{x,\infty}}}_\Inf$ permutes
its Breuil-Kisin lattices. The stability of $\mf M_x$ follows from the uniqueness of Breuil-Kisin lattices, which follows from the analogous statement for $G=\GL_n$ and the Tannakian formalism.

(2) $\mf M_x \times^{\lsup LG} \lsup L\{*\}$
corresponds to (via Fontaine's functors)
the unique parameter $G_{\Q_p} \to \lsup L\{*\}(\ol\Q_p)$
for the singleton group,
which is crystalline.
The statement is standard for crystalline representations.

\end{proof}
By \cite[\S 4.6]{EG23} and the Tannakian formalism, the meaning of $V$ having Weil-Deligne type $\tau$ is that the induced action of inertia $I(L_x/\Q_p)$ on $\mf M_x/u \mf M_x[1/p]$ is isomorphic to $\tau\otimes_{\cO}\Lambda[\frac{1}{p}]$. Since $\tau$ is tame, the isomorphism already exists without inverting $p$.

\begin{lemma}
Let $\Lambda=\cO'$, the ring of integer in some finite extension $E'/E$.
There exists a finite flat $\Lambda$-algebra $\Lambda'$ and $Y\in L^{\le \lambda}\cG_x(\Lambda')$ such that 
$M_\Inf\otimes_\Lambda\Lambda' \cong \iota_x(Y)\otimes_{\Z_p(\!(v)\!)}W(\C^\flat)$ as $(\varphi_0,G_{\Q_{p, \infty}})$-modules.

\label{lemma: existence of Kisin lattice}
\end{lemma}



\begin{proof}
After possibly passing to a finite flat $\Lambda$-algebra, we may assume $\mf M_x$ is a trivial $\lsup LG$-torsor. Fix a frame $\beta:\mf M_x\cong \bG(W(\F_x)\otimes \Lambda[\![u]\!])\rtimes \langle \sigma \rangle$.

For each $\gamma\in \Gal(L_x/\Q_p)$ and $y\in \bG(W(\F_x)\otimes \Lambda[\![u]\!])\rtimes \langle \sigma \rangle$, we have
\[\gamma\cdot y = \gamma(y)C_\gamma^{-1}\]
where $C_\gamma\in \bG(W(\F_x)\otimes \Lambda[\![u]\!])$.
Direct computation shows
$[C_\gamma]\in H^1(\Gamma, \bG(W(\F_x)\otimes \Lambda[\![u]\!]))$,
where $\Gamma$ acts solely on the coefficient ring $W(\F_x)\otimes \Lambda[\![u]\!]$.
On the other hand, the Frobenius structure can be written as
$\varphi_{\mf M_x}(y)=\varphi_0(y)(1\rtimes \sigma)X$
for some $X\in \bG(W(\F_x)\otimes \Lambda[\![u]\!][\frac{1}{u^e+p}])$
    and all $y\in \lsup LG(W(\F_x)\otimes \Lambda[\![u]\!])$.
The relation $\varphi_{M_x}(\gamma\cdot h)=\gamma\cdot \varphi_{M_x}(h)$
is equivalent to the equation
\[\varphi(C_\gamma)\gamma(X)C_\gamma^{-1}
=X\]
Hence $[C_\gamma]\in H^1(\Gamma, \bG(W(\F_x)\otimes \Lambda[\![u]\!]))^{\varphi=1}$.

The fact that $I(L_x/\Q_p)$ has order prime to $p$ and Shapiro's Lemma gives a chain of isomorphisms  
\begin{align*}
H^1(\Gal(L_x/\Q_p), \bG(W(\F_x)\otimes \Lambda[\![u]\!]))^{\varphi=1}
&\xrightarrow{\cong}
H^1(\Gal(L_x/\Q_p), \bG(W(\F_x)\otimes \Lambda))^{\varphi=1}\\
&\xrightarrow{\cong}
H^1(I(L_x/\Q_p), \bG(\Lambda))^{\varphi=1}
\end{align*}
where the last isomorphism is induced by the projection to the factor corresponding to our preferred embedding $\iota_\psi:W(\F_x)\inj \cO$.
The condition that $\mf M_x$ has inertial type $\tau$ is exactly the condition that the class $C_\gamma$ correspond to $\tau$. 

Recall that we fixed a presentation $\wt{w}(\tau)=v^\mu w$, which gives rise to $x\in X_*(\bT)\otimes \frac{1}{e}\Z/\Z$ such that $\varphi(x)=\mu+wx$. We identify $\wt{w}(\tau)$ with its Tits representative in $\bG(W(\F_x)\otimes \Lambda(\!(u)\!))^{\Gamma}=\bG(\Lambda(\!(v)\!))=L\cG_x(\Lambda)$.
Unravelling the above chain of isomorphisms and the relationship between $x$ and $\tau$, we see that    
$C_\gamma = C^{-1}v^{-x}\gamma(v^xC)$ for some $C\in \bG(W(\F_x)\otimes \Lambda[\![u]\!])$, and that
\[\gamma(\varphi(v^xC) X C^{-1} v^{-x})
=\varphi(v^xC) X C^{-1} v^{-x}.\]

Set 
\[Y:=\wt{w}(\tau)^{-1}\varphi(v^xC) XC^{-1}v^{-x} \in L\cG_x(\Lambda),\]
so that $\iota_x(Y) =(1\rtimes \sigma)\wt{w}(\tau)Y= \varphi_0(v^xC)(1\rtimes \sigma)X(v^xC)^{-1}$. Thus $\iota_x(Y)\otimes_{\Z_p(\!(v)\!)}W(\C^\flat)\cong M_\Inf$ as $(\varphi_0,G_{\Q_{p,\infty}})$-modules.

It remains to show $Y\in L^{\le \lambda}\cG_x(\Lambda)$. By Lemma \ref{lemma: Lang variant} below, we can
choose $b\in \bG(W(\F_x)\otimes\Lambda)$ such that
$w=\varphi(b)b^{-1}$. Setting $D=b^{-1} C\in  \bG(W(\F_x)\otimes\Lambda)$, we get
\[Y:=\wt{w}(\tau)^{-1}\varphi(v^xC) X C^{-1} v^{-x}=v^xb\varphi(D) X D^{-1} b^{-1}v^{-x}. \]

Since the Hodge-Tate weights of $\mf M_x$ is $\leq \lambda$, the above formula for $Y$ and the fact that $(-\frac{v}{p})^x\in \bG(\Lambda[\frac{1}{p}][\![v+p]\!])$ shows that it belongs to $L\cG_x^{\leq \lambda}(\Lambda[\frac{1}{p}])$, and hence $Y\in L\cG_x^{\leq \lambda}(\Lambda)$.

\end{proof}
\begin{remark}(Eigenframes and renormalized Frobenius) \label{rmk:eigenframes}The element $\wt{z}=v^xb$ has the property that
\[\varphi(\wt{z})\wt{z}^{-1}=\wt{w}(\tau)\]
If $\mf M_x$ is a Breuil--Kisin module with $\Gamma$-action of type $\tau$ with $R$-coefficient for some $p$-adically complete $\cO$-algebra $R$, we can always fppf locally choose a frame $\beta$ of $\mf M_x$ such that the matrix $D=1$.
We call such a trivialization an \emph{eigenframe} of $\mf M_x$. For such a frame we have the simple relationship
\[Y=\Ad(\wt{z})X\]
and we think of $Y$ as the renormalized matrix of Frobenius with respect to $\beta$.
\end{remark}

\begin{lemma}
\label{lemma: Lang variant}
Let $w\in \bG(\cO)= \bG(\Z_{p^f}\otimes_{\Z_p} \cO)^{\varphi_0=1}$ be an element
such that there exists an integer $f$ such that
$$
w\pi^{-1}(w)\pi^{-2}(w)\dots\pi^{-f+1}(w)=1.
$$
If $\Z_{p^f}\subset \cO$ and $G$ splits over $\Z_{p^f}$,
then
the equation
$w=c^{-1}\varphi(c)$
has a solution $c\in \bG(\Z_{p^f}\otimes_{\Z_p} \cO)$.

(We remind the reader that $\varphi=\varphi_0\pi^{-1}=\pi^{-1}\varphi_0$ where $\varphi_0$ acts on $\Z_{p^f}\otimes \cO$ by $\varphi_0(a\otimes b)=\Frob_p(a)\otimes b$).
\end{lemma}

\noindent When $\pi^{-1}=1$, this lemma is \cite[Remark 2.5.20]{Ei}.
Since the Tits representative $w$ of
Weyl group elements
$[w]\in W(\bT, \bG)$
generate a finite group
(which is called the extended Weyl group by Tits)
which is stable under the action of $\pi^{-1}$,
the lemma applies to all Tits representatives.

\begin{proof}

We can choose an isomorphism
$\bG(\Z_{p^f}\otimes \cO)=\bG(\cO^{\oplus f})\cong \bG(\cO)^{ \Z/f\Z}$ so that $\varphi_0$ correspond to the shifting automorphism $(g_i)\mapsto (g_{i+1})$. Then the given element $w$ is of the form $(g,g,\cdots g)$ where 
\[g\pi^{-1}(g)\cdots \pi^{-f+1}(g)=1\]
We need to find $c=(c_i)$ such that
\[g=c_i^{-1}\pi^{-1}(c_{i+1}).\]
But the condition on $g$ shows that setting 
\[c_0=1, c_{i+1}=\pi(c_i)g\]
gives a solution.

\end{proof}

\begin{prop} \label{prop:map to BK space}Assume that $G$ has connected center and $x$ is $(h_\lambda+1)$-generic.
The morphism
\[
[(L\cG^{\le \lambda}_{x})^{\wedge_p}/_{\varphi, x}(L^+\cG_{x})^{\wedge_p}]
\underset{\cR_{\lsup LG}}{\times}\cX^{\le \lambda, \tau}
\to \cX^{\le \lambda, \tau}
\]
is an isomorphism.

In particular, we obtain a morphism
\[r:\cX^{\le \lambda, \tau}\to [(L\cG^{\le \lambda}_{x})^{\wedge_p}/_{\varphi, x}(L^+\cG_{x})^{\wedge_p}]\]
\end{prop}

\begin{proof}
The proof is identical to that of \cite[Proposition 7.2.7]{LLLM22}, using Lemma \ref{lemma: existence of Kisin lattice}
and Proposition \ref{prop: monomorphism} as the input.
\end{proof}

Under slightly stronger genericity assumptions, one has 
\begin{prop}\label{prop: restriction on Galois is mono} Assume that $x$ is $(h_\lambda+2)$-generic. Then
\[\cX^{\le \lambda, \tau}\to \cR_{\lsup LG}\]
is a monomorphism, and hence the induced map 
\[r:\cX^{\le \lambda, \tau}\to [L\cG^{\le \lambda}_{x}/_{\varphi, x}L^+\cG_{x}]^{\wedge_p}\]
is a monomorphism.
\end{prop}

We need two preparatory lemmas:

\begin{lemma} \label{lemma: cyclo free pre-estimates}

Assume $x$ is $(h_\lambda+2)$-generic. Let $\ol \rho: G_{\Q_p}\to \lsup LG(\ol \F)$ be an $L$-parameter coming from  $\cX^{\le \lambda, \tau}(\ol \F)$.

\begin{enumerate}
\item
The semisimplification $\ol\rho^{\semis}$
also belongs to $X_{\lsup LG}^{\le \lambda, \tau}(\F_p)$.
\item The restriction map 
\[R\Gamma(G_{\Q_p},\ad(\ol \rho))\to R\Gamma(G_{\Q_{p,\infty}},\ad(\ol \rho))\]
induces an injection on $H^1$ and an isomorphism on $H^0$. The same statement holds after replacing $\Q_p$ with any unramified extension $\Q_{p^k}$.

\item Assume that $\rho|_{G_{\Q_{p^f}}}$ factors through a parabolic $\bP\subset \bG$ with unipotent radical $\bU$. Set $\Q_{p^f,\infty}=\Q_{p,\infty}\Q_{p^f}$.
Let $c\in Z^1(G_{\Q_{p^f}},\bU(\ol \F))$ be a cocycle whose restriction to $Z^1(G_{\Q_{p^f,\infty}},\bU(\ol \F))$ is trivial. Then $c$ is trivial.

\end{enumerate}
\end{lemma}

\begin{proof}
\begin{enumerate}
\item This follows from the standard continuity argument (c.f. \cite[Proposition 5.5.9]{LLLM22}).
\item 
We have $\ol \rho^{\semis}: G_{\Q_p}\to \lsup LS(\ol \F)$ for an unramified maximal torus $S\subset G$
(\cite{Lin23a}). As in \cite[Proposition 5.5.7]{LLLM22}, \cite[Proposition 3.1.2]{LLL19}, the restriction $\rho^{\semis}|_{I_{\Q_p}}$ is a tame inertial $L$-parameter admitting a presentation in
\[ \wt{w}(\tau)\Adm(\lambda)\]
Indeed, by unramified base change this reduces to the case $S=T$ is maximally split. In this case, the centralizer of $\ol\rho$ must fix the unique $\mf M_{\ol\rho}\in [L\cG^{\le \lambda}_{x}/_{\varphi, x}L^+\cG_{x}](\ol \F)$ giving rise to  $\ol \rho|_{G_{\Q_{p,\infty}}}$, hence $\mf M_{\ol\rho}$ can be represented by an element in $\chT(\ol \F)\Adm(\lambda)$.

But now the hypothesis that $\wt{w}(\tau)$ is $(h_\lambda+2)$-generic shows that $\ad(\ol\rho)$ is cyclotomic-free in the sense of \cite[Definition 7.2.9]{LLLM22}, and the result now follows from \cite[Lemma 3.10, 3.11]{LLLMold}.
\item Let $\bU_k=1\subset \bU_{k-1}=Z(\bU)\subset  \dots \subset \bU_0=\bU$ be the central series for $\bU$, whose associated graded $\bU_i(\ol \F)/\bU_{i+1}(\ol\F)$ is an $\ol \F$-vector space, stable under $G_{\Q_{p^f}}$ and cyclotomic-free.
We show that $c$ takes values in $\bU_k(\ol \F)$ for each $k$ by induction on $k$. Indeed, if this holds for some $k$, then \cite[Proposition 3.12]{LLLMold} shows that the restriction map
\[Z^1(G_{\Q_{p^f}},\bU_k(\ol \F)/\bU_{k+1}(\ol \F))\to Z^1(G_{\Q_{p^f,\infty}},\bU_k(\ol \F)/\bU_{k+1}(\ol \F))\]
is injective. But this implies $c$ takes values in $\ol U_{k+1}(\ol \F)$.
\end{enumerate}

\end{proof}

\begin{lemma}\label{lem:cyclo free estimates} Assume $x$ is $(h_\lambda+2)$-generic. Then the restriction functor $\ol\rho \mapsto \ol\rho|_{G_{\Q_{p, \infty}}}$
is a fully faithful embedding of the groupoid
of $L$-parameters $\ol\rho:G_{\Q_p}\to \lsup LG(\F_p)$
coming from $\cX^{\le \lambda, \tau}(\ol \F)$
to the groupoid
of $L$-parameters $\ol\rho:G_{\Q_{p,\infty}}\to \lsup LG(\F_p)$
coming from $\cX^{\le \lambda, \tau}(\ol \F)$.
\end{lemma}

\begin{proof}
Set $\ol\rho_\infty:=\ol\rho|_{G_{\Q_{p, \infty}}}$
for ease of notation.
It suffices to prove the following: if $\ol \rho, \ol \rho'$ are two parameters such that $\ol \rho_\infty=\ol \rho'_\infty$, then $\ol \rho=\ol \rho'$.

Since $G_{\Q_p}=G_{\Q_{p,\infty}}I_{\Q_p}$, it suffices to show $\ol \rho=\ol \rho'$ when restricted to $G_{\Q_{p^f}}$.
Hence for the rest of the proof we will replace $\ol \rho$, $\ol \rho'$ by their restriction to $G_{\Q_{p^f}}$, which allows to work with $\bG(\ol \F)$-valued parameters instead of $\lsup LG(\ol \F)$-valued parameters.

 Let $\bP,\bP'$ be two minimal parabolics of $\bG$ such that $\Img(\ol\rho),
 \Img(\ol\rho')$ factor through $\bP$, $\bP'$ respectively. 
For any Levi decomposition $\bP=\bU\rtimes \bM$, we can write
\[\ol \rho=\ol c \ol\rho^{\semis}\]
where the semisimplifcation $\ol \rho^{\semis}$ is obtained by the composition of 
$\ol \rho$ with  $\bP(\ol \F)\surj \bM(\ol \F)\inj \bG(\ol \F)$, and $\ol c\in Z^1(G_{\Q_{p^f}},\bU(\ol \F))$ is an $1$-cocycle for the action of $G_{\Q_{p^f}}$ on $\bU(\ol \F)$ via conjugation by $\ol \rho^{\semis}$.
We have the similar objects $(\ol \rho')^{\semis}, \ol c'$ for $\ol \rho'$.
By the minimality of $\bP$, $\ol\rho^{\semis}$ is $\bM$-irreducible. Since $\ol\rho^{\semis}(G_{\Q_{p^f}})$ factors the tame quotient of $G_{\Q_{p^f}}$, which coincides with the tame quotient of  $G_{\Q_{p^f,\infty}}$, $\ol\rho_\infty^{\semis}$ is also $\bM$-irreducible, and is also the a semisimplification of $\ol \rho_\infty$. We have the same statement for $\bP'$, $\ol \rho'$. But since $\ol \rho_\infty=\ol \rho'_\infty$, by the uniqueness of semisimplification up to $\bG$-conjugacy (see \cite[Proposition 3.3]{Ser05} and its proof), we see that $\bP$, $\bP'$ share a common Levi and $\bM$, $\bM'$ are $\bG$-conjugate. Thus by changing the choice of Levis, we can assume that $\bM=\bM'$, $\bP'=\Ad(w)\bP$ and $w\in N_{\bG}(\bM)$. 

Set $\bU_w=\bU\cap \Ad(w)\bU'$. Then $\Img(\ol\rho_\infty)\subset \bU_w(\ol \F)\bM(\ol \F)$. Hence the restriction $\ol c_\infty=\ol c|_{G_{\Q_{p^f,\infty}}}$ takes values in $\bU_w(\ol \F)$. 
The same argument as in the proof of Lemma \ref{lemma: cyclo free pre-estimates}(3) shows that $\ol c$ takes values in $U_w(\ol \F)$. Hence $\ol \rho$ takes values in $P\cap \Ad(w)P$. Symmetrically, $\ol \rho'$ takes values in $P\cap \Ad(w)P$.
But now setting $d=c'c^{-1}$, then $\ol \rho=d\ol\rho'$, and $d\in Z^1(G_{\Q_{p^f}},\bU(\ol \F))$ is a cocycle for the Galois action on $\bU(\ol \F)$ induced by conjugation of $\ol \rho$. Our hypothesis is that the restriction of $d$ to $Z^1(G_{\Q_{p^f,\infty}},\bU(\ol \F))$ is trivial, hence Lemma \ref{lemma: cyclo free pre-estimates}(3) shows that $d$ is trivial, and thus $\ol \rho=\ol \rho'$.

\end{proof}

\begin{proof}[Proof of Proposition \ref{prop: restriction on Galois is mono}] It suffices to show that for any Artinian $\cO/\varpi^m$-algebra $A$ and objects $y_1,y_2\in  \cX^{\le \lambda, \tau}(A)$ with images $z_1,z_2$ in $ \cR_{\lsup LG}(A)$, the induced map on Hom-spaces is an isomorphism. Standard deformation theory reduces this to the case $A=\ol \F$ (which holds by Lemma \ref{lem:cyclo free estimates}) and a tangent space computation, which hold by Lemma \ref{lemma: cyclo free pre-estimates}.
\end{proof}

\subsection{Monodromy condition}

\subsubsection{The $\dlog$ map}
Recall our standing assumption that $p>h$, the Coxeter number of $\bG$.
There exists a morphism 
$$
\dlog: L\cG^{\le \kappa}_x \to \Lie L\cG_x
$$
satisfying the relation
\begin{equation}
\label{eqn:dlog}
\dlog(fg) = \dlog(f) + \Ad(f)\dlog(g).
\end{equation}
The map $\dlog$ can be defined by $\dlog(g)=\frac{dg}{dv}g^{-1}$ via an arbitrary embedding
of $\cG_x$ to $\GL_N$
as smooth affine group schemes over $\cO[v]$.
It remains to show
the image of $\dlog$
is contained in $\Lie L\cG_x$.
Since $L\cG^{\le \kappa}_x$ is reduced,
it suffices to check on field-valued points;
and by the Bruhat decomposition $\bG=\bB W \bB$,
we only need to show $\dlog$ maps affine root groups
of $L\cG_x$
to $\Lie L\cG_x$.
Since $p>h$,
affine root groups of $L\cG_x$
are of the form $t\mapsto \exp_{\le h}(v^m\alpha~t)$
where $\alpha$ is a root vector and $\exp_{\le h}$
is the truncated exponential map.
Direct computation that $\dlog$ maps the 
affine root group
to the line spanned by $v^{m-1}\alpha$.

\begin{lemma}
\label{lemma: dlog height}
We have
$
\dlog (L\cG^{\le \kappa}_x) \subset (v+p)^{-\max(h_\kappa,1)}\Lie L^+\bG.
$
\end{lemma}

\begin{proof}
It is clear by the discussion above that
$\dlog (L^+\cG_x) \subset \Lie L\cG_x.$
Again, we only need to check on field-valued points.
For characteristic $p$ points,
we use the affine Bruhat decomposition and
for characteristic $0$ points, we use
the Cartan decomposition.
Combining Eqn. (\ref{eqn:dlog}) and Lemma \ref{lem:amplitude bound},
we are reduced to check 
$\dlog (L\cG^{\le \kappa}_x\cap L\bT) \subset (v+p)^{-1}\Lie L^+\cG_x$.
We can arrange it so that $L\bT$ is mapped to diagonal matrices of $L\GL_N$, the lemma follows from
$\frac{d (v+p)^m}{dv} = m(v+p)^{m-1}$, $m\in \Z$.
\end{proof}

\subsubsection{The canonical connection of a Breuil--Kisin modules}
Let $R$ be a $p$-adically complete, topologically of finite type flat $\cO$-algebra. Recall $\cO^{\mathrm{rig}}_{x,R}=\varprojlim (W(\F_x)\otimes_{\Z_p}R)[\![u,\frac{u^n}{p}]\!][\frac{1}{p}]$, the ring of rigid analytic functions on the open disk over $(\Spf R\otimes_{\Z_p} W(\F_x))^{\mathrm{rig}}$. It is endowed with the following structures:
\begin{itemize}
\item A differential operator $\partial=u\frac{d}{du}=ev\frac{d}{dv}$;
\item An $R$-linear Frobenius $\varphi_0$ which acts as the arithmetic Frobenius on $W(\F_x)$ and sends $u$ to $u^p$;
\item An action of $\Gamma=G(L_x/\Q_p)$ extending that on $W(\F_x)[\![u]\!]$. 
\end{itemize}
We have the element $\lambda=\prod_{i=0}^{\infty} \varphi_0^i(\frac{u^e+p}{p})\in \cO^{\mathrm{rig}}_{x,R}$.
\begin{prop} \label{prop:canonical derivation} Let $\mf M_x$ be a Breuil--Kisin module with $\lsup LG$-structure over $(W(\F_x)\otimes_{\Z_p} R)[\![u]\!]$. 
\begin{enumerate}
\item There is a unique connection $\nabla\in \Ad(\mf M_x) \otimes_{W(\F_x)[\![u]\!]\otimes R} \cO_{x,R}^{\mathrm{rig}}[\frac{1}{\lambda}]du$ which commutes with $\varphi=\varphi_{\mf M_x}\otimes \varphi_0$.
\item If $\beta$ is a frame of $\mf M_x$ and $Y$ the matrix of $\varphi_{\mf M_x}$ with respect to $\beta$, the differential operator $\nabla_{\partial}$ correspond to the element $N_{\partial}\in u\chfg(\cO_{x,R}^{\mathrm{rig}}[\frac{1}{\lambda}])$ given by
\[(1-p\Ad(X^{-1})\varphi)^{-1}\big(u\dlog_u(X^{-1})\big)=\sum_{i=0}^{\infty} p^i(\Ad(X^{-1})\varphi)^i(u\dlog_u (X^{-1}))\]
(here the symbol $\dlog_u$ stands for the logarithmic derivative with respect to $u$)
\item If $R=\Lambda$ is finite free over $\cO$, then $N_\partial$ has logarithmic poles along the locus $u^e+p=0$ if and only if $V(\mf M_x)$ correspond to a crystalline parameter $\rho:G_{L_{x}}\to \bG(\Lambda[\frac{1}{p}])$ such that $\rho(G_{L_{x,\infty}})\subset \bG(\Lambda)$. If furthermore $\mf M_x$ has Hodge type $\leq \kappa$ with $2h_{\kappa}<p-1$ and $N_{\partial}\in \frac{u}{p}\Ad(\mf M_x)\otimes_{W(\F_x)[\![u]\!]} W(\F_x)[\![u,\frac{u^e}{p}]\!]$ then $\rho(G_{L_x})\subset \bG(\Lambda)$.
\end{enumerate}
\end{prop}
\begin{proof} Choosing a frame $\beta$, the condition that $\nabla$ commutes with $\varphi$ is equivalent to the commutation relation
\[p\Ad(X^{-1}(1\rtimes \sigma)^{-1})\varphi_0(N_{\partial})=N_{\partial}-u\dlog_u(X^{-1})\]
or equivalently
\[(1-p\Ad(X^{-1})\varphi))(N_{\partial})=u\dlog_u(X^{-1})\]
The uniqueness and existence of such $N_{\partial}$ as well as the convergence of the series defining it thus follows from the case $\bG=\GL_n$, which is well-known. The same reasoning also shows the first half of the third item.

We turn to the last claim of the third item. Composing with the adjoint representation of $\bG$, we see that $\Ad(\mf M_x)$ is a Breuil--Kisin module with Hodge type concentrated in an interval of length $p-1$ and whose corresponding canonical derivation $N_{\partial}$ preserves the induced quasi-strongly divisible lattice in the sense of \cite[\S 3.4]{Liu08}. Hence by \cite[Proposition 3.5.1]{Liu08}, $\Ad\circ \rho(G_{L_x})\subset \GL(\bG(\Lambda))$. Composing with one-dimensional representations $\chi$ of $\bG$ instead, we see that $\chi(\rho(G_{L_x}))\subset \Lambda^\times$. Let $H=\bG/C$ be the image of $\bG$ in the sum of the adjoint and all one-dimensional representations of $\bG$, so $C$ is a finite central subgroup of $\bG$ of order prime to $p$ (since $p>2h_{\bG}$). We thus have $\rho(G_{L_x})\subset \bG(\Lambda)C(\Lambda[\frac{1}{p}])$. Let $G'=\rho^{-1}(\bG(\Lambda))\subset G_{L_x}$, then for any element $\sigma\in G_{L_x}$, $\sigma^m\in G'$ for some $m$ coprime to $p$. This shows that $G_{L_x(\zeta_{p^{\infty}})}\subset G'$, and hence $G_{L_x}=G_{L_{x,\infty}} G_{L_x(\zeta_{p^{\infty}})}\subset G'$.
 
\end{proof}
\begin{remark}(Monodromy operator with respect to eigenframes)\label{rmk:monodromy with descent data} Let $\mf M_x$ be a Breuil--Kisin module with $\lsup LG$-structure with $\Gamma$-action of type $\tau$. Let $\beta$ be an eigenframe and $\wt{z}\in \wt{W}$ be the element such that $\varphi(\wt{z})\wt{z}^{-1}=\wt{w}(\tau)$ as in Remark \ref {rmk:eigenframes}.
Then, using that $u\dlog_u=e v\dlog$ and recalling that $\varphi_x=\Ad(\wt{w}(\tau)^{-1})\varphi$, the formula for $N_{\partial}$ can be written in terms of the renormalized matrix of Frobenius $Y$ as
\begin{align*}
&\Ad(\wt{z}^{-1})\big((1-p\Ad(Y^{-1}\wt{w}(\tau)^{-1})\varphi))(ev\dlog(Y^{-1})+\Ad(Y^{-1})ex-ex )\big)\nonumber\\
=&\Ad(\wt{z}^{-1})\big(\sum_{i=0}^\infty (p\Ad(Y^{-1})\varphi_x)^i(ev\dlog(Y^{-1})+\Ad(Y^{-1})ex-ex)   \big)\nonumber
\end{align*} 
\end{remark}

\begin{defn} (Monodromy condition) \label{defn:monodromy}Let $R$ be a topologically finite type flat $\cO$-algebra, and $Y\in L\cG^{\leq \kappa}_x(R)$ and set
\begin{align*}
\cN_{\infty}(Y)&=(v+p)^{h_\kappa}\big(\sum_{i=0}^\infty (p\Ad(Y^{-1})\varphi_x)^i(v\dlog(Y^{-1})+\Ad(Y^{-1})x-x) \\
&=\big(\sum_{i=0}^\infty \frac{p^i}{(p+v^p)^{h_{\kappa}}\cdots (p+v^{p^i})^{h_{\kappa}}}((v+p)^{h_{\kappa}}\Ad(Y^{-1})\varphi_x)^i(v\dlog(Y^{-1})+\Ad(Y^{-1})x-x), 
\end{align*}
which belongs to $\Lie L\cG_x(R[\frac{1}{p}])$.

The \emph{monodromy condition} on $Y$ is the condition that  
\[\cN_\infty(Y)\in (v+p)^{h_{\kappa}-1}\Lie L^+\cG_x(R[\frac{1}{p}]),  \]
or more precisely the condition that the Taylor coefficients
\[(\frac{d}{dv})^t|_{v=-p}(\cN_\infty(Y))=0\]
\[(\frac{d}{dv})^t|_{v=-p}(v^{-1}\cN_\infty(Y)_{\alpha})=0\]
for $0\leq t \leq h_{\kappa}-2$ and $\alpha>0$ (here our convention is that for $X\in \chfg$, $X_\alpha$ is its component in the $\alpha$-root space).
These conditions determine a unique $\cO$-flat closed subfunctor $L\cG_x^{\leq \kappa,\nabla_\infty}$ of the $p$-adic completion $(L\cG^{\leq \kappa}_x)^{\wedge_p}$.
\end{defn}

\begin{prop} \label{prop: monodromy vs EG}Assume $x$ is $(h_\kappa+2)$-generic. Then we have a commutative diagram
\[
\begin{tikzcd}
 \cX^{\le \kappa, \tau} \ar[r, "{\cong}"]  \ar[dr,hookrightarrow] 
 & {\left[{L}\cG^{\leq \kappa,\nabla_\infty}_x {/}_{\varphi,x} ({L}^+\cG_x)^{\wedge_p}\right]} \ar[d,hookrightarrow, "{\iota}"] \\
&  \cR_{\lsup LG}
\end{tikzcd}
\]
\end{prop}
\begin{proof} The existence of the map horizontal map $r:\cX^{\le \kappa, \tau}\to [(L\cG^{\le \kappa}_{x})^{\wedge_p}/_{\varphi, x}(L^+\cG_{x})^{\wedge_p}]$ making the diagram commute follows from Proposition \ref{prop:map to BK space}, and the fact that all the maps are monomorphisms follows from Proposition \ref{prop: monomorphism} and Proposition \ref{prop: restriction on Galois is mono}. To show that $r$ identifies $\cX^{\le \kappa, \tau}$ with the substack $\left[{L}\cG^{\leq \kappa,\nabla_\infty}_x {/}_{\varphi,x} ({L}^+\cG_x)^{\wedge_p}\right]$, by \cite[Proposition 4.8.6]{EG23} it suffices to show both substacks have the same $\Lambda$-points for any finite flat $\cO$-algebra $\Lambda$. A $\Lambda$-point of $[L\cG^{\le \kappa,\nabla_\infty}_{x}/_{\varphi, x}(L^+\cG_{x})^{\wedge_p}]$ correspond to a Breuil--Kisin module $\mf M_x$ over $W(\F_x)\otimes \Lambda[\![u]\!]$ with $\Gamma$-action of type $\tau$ and Hodge type $\leq \kappa$, and whose canonical differential operator $N_\partial$ has logarithmic poles along $v=-p$. By Remark \ref{rmk:monodromy with descent data} and Lemma \ref{lem:monodromy estimates} below, $N_{\partial}\in \Ad(\mf M_x)\otimes_{W(\F_x)[\![u]\!]} \frac{u}{p}W(\F_x)[\![u,\frac{u^e}{p}]\!]$.
By Proposition \ref{prop:canonical derivation}, $\mf M_x$ correspond to a parameter $\rho:G_{\Q_{p,\infty}}\to \lsup LG(\Lambda)$ whose restriction to $G_{L_{x,\infty}}$ extends to a crystalline parameter $\rho:G_{L_x}\to \bG(\Lambda)$. Thus $\rho$ extends to a potentially crystalline parameter of $G_{\Q_p}$, and whose inertial type and Hodge--Tate weights are those of $\mf M_x$, so it has inertial parameter $\tau$ and Hodge--Tate weight $\leq \kappa$. But this is exactly a $\Lambda$-point of $\cX^{\le \kappa, \tau}$.
\end{proof}

\section{The mod $p$ local model\\ by B. Le Hung}\label{app: B}
We continue with the setup of Appendix \ref{app: EG stacks}. In particular we have an unramified pinned group $(G, B, T, \{X_{\alpha}\})$ over $\Q_p$ whose Langlands dual pinned group is $(\bG, \bB, \bT, \{Y_{\alpha^\vee}\})$. This gives rise to the Frobenius automorphism $\pi^{-1}$ on the root datum, and the endomorphism $\varphi=p\pi^{-1}$ on $X^*(T)$. 

We fix a tame inertial $L$-parameter $\tau$ and a dominant Hodge-Tate cocharacter $\lambda\in X_*(\chT)$. This determines the potentially crystalline stack $\cX^{\le \lambda,\tau}$ and $\cX^{\lambda,\tau}$.
We also put ourselves in Situation \ref{running assumption}, that is, we assume $\tau$ has a lowest alcove presentation $(w,\mu)$, which we fix. This choice is equivalent to choosing the element $x\in \frac{1}{p^f-1}X_*(\chT)=X^*(T)$ (for some sufficiently large $f$) uniquely characterized by the condition $\varphi(x)=\wt{w}(\tau)(x)=\mu+w(x)$. In particular we have $x\equiv -w^{-1}\mu$ modulo $p$. Set $e=p^f-1$. 

We work over a coefficient field $E$ with ring of integers $\cO$ and residue field $\F$. Assume $E$ is sufficiently large to contain a Galois closure of
$L_x=\Q_{p^{f}}( (-p)^{\frac{1}{e}})$.



\subsection{Twisted loop groups}
As in section \ref{subsect: PZ group schemes}, the choice of $x$ gives rise to the group schemes $L\cG_x$, and the endomorphism $\varphi_x$ on its $p$-adic completion $L\cG_x^{\wedge_p}$. Let $L^+_1\cG_x$ be the kernel of the map $L^+\cG_x\surj \bB\surj \bT$ where the first map is reduction mod $v$.

By construction $(L^+\cG_x)_\F$ is the Iwahori group scheme of $L^+\bG_\F$ corresponding to $\bB_\F$ with respect to the loop variable $v$, and $(L_1^+\cG_x)_\F$ is its pro-$v$ Iwahori subgroup. Thus $(L^+\bG)_\F=\chI$ in the notation of the main text (up to changing the loop variable from $v$ to $t$), and $(L_1^+\cG_x)_\F= \chI_1$.

As before, we have the twisted affine flag variety $\Gr_{\cG_x}=L\cG_x/L^+\cG_x$.

\subsection{Approximating the monodromy condition}
Let $R$ be a topologically finite type flat $\cO$-algebra and $Y\in L\cG_x^{\leq \lambda}(R)$. 
We recall from Definition \ref{defn:monodromy} the operator
\[\cN_{\infty}(Y)=\sum_{i=0}^\infty\frac{p^i}{(p+v^p)^{h_{\lambda}}\cdots (p+v^{p^i})^{h_{\lambda}}} (\Ad^*(Y)\varphi_x)^i\Big((v+p)^{h_{\lambda}}(v\dlog(Y^{-1})+\Ad(Y^{-1})x-x)\Big) \]
in $\Lie L^+\cG_x(R[\frac{1}{p}])$, where we abbreviate $\Ad^*(Y):=(v+p)^{h_{\lambda}}\Ad(Y^{-1})$. Note that $\Ad^*(Y)$ preserves both $\Lie L^+\cG_x(R)$ and $\Lie L^+\bG(R)=\chfg(R[\![v+p]\!])$.
We set 
\[\cN_0(Y)=(v+p)^{h_\lambda}(v\dlog(Y^{-1})+\Ad(Y^{-1})x-x)\in \Lie L_1^+\cG_x(R)\]
and
\[f_\infty(Y)=\sum_{i=1}^{\infty}\frac{p^i}{(p+v^p)^{h_{\lambda}}\cdots (p+v^{p^i})^{h_{\lambda}}}  \varphi_x((\Ad^*(Y)\varphi_x)^{i-1}\cN_0(Y)),\]
so that $\cN_\infty(Y)=\cN_0(Y)+(v+p)^{h_\lambda}\Ad(Y^{-1})(f_\infty(Y))$.

\begin{lemma}\label{lem:monodromy estimates} Let $R$ be a flat topologically finite type $\cO$-algebra and $Y\in L\cG_x^{\le \lambda}(R)$. Assume that $x$ is $m$-generic. 
Then for $i\geq 1$, 
\[\varphi_x(\Ad^*(Y)\varphi_x)^{i-1}(\cN_0(Y))\in v^{1+m\frac{p^i-1}{p-1}}\chfg(R[\![v+p]\!]).\]
Furthermore, if $Y'=AY$ where $A\in\ker( L^+\bG(R)\to L^+\bG(R/p^k))$ for some $k>0$, then 
\[\varphi_x(\Ad^*(Y')\varphi_x)^{i-1}(\cN_0(Y'))-\varphi_x(\Ad^*(Y)\varphi_x)^{i-1}(\cN_0(Y))\in p^kv^{1+m\frac{p^i-1}{p-1}}\chfg(R[\![v+p]\!]).\]
\end{lemma}
\begin{proof} 

By Lemma \ref{lemma: dlog height}
\[\cN_0(Y)=(v+p)^{h_\lambda}(v\dlog(Y^{-1})+\Ad(Y^{-1})x-x)\in \Lie L_1^+\cG_x(R).\]
Since $\Ad^*(Y)$ preserves $\chfg(R[\![v+p]\!])$ and 
\[\varphi_x(\Lie L_1^+\cG_x(R))=\Ad(w^{-1}v^{-\mu})\varphi(\Lie L_1^+\cG_x(R))\subset v^{1+m}\chfg(R[\![v+p]\!])\]
by the $m$-genericity of $x$, the lemma holds for $i=1$. 
Now the fact that 
\[\varphi_x(v^k\chfg(R[\![v+p]\!]))=\Ad(w^{-1}v^{-\mu})\varphi(v^k\chfg(R[\![v+p]\!]))\subset v^{p(k-1)+m+1}\chfg(R[\![v+p]\!])\]
gives the statement for $i\geq 2$ by induction.

The last assertion follows from the above and 
\[(\Ad^*(Y')-\Ad^*(Y))(\chfg(R[\![v+p]\!]))\in p^k\chfg(R[\![v+p]\!])\]
\[\cN_0(Y')-\cN_0(Y)=\Ad(Y^{-1})(\cN_0(A))\in p^k\chfg(R[\![v+p]\!]).\]
\end{proof}

\begin{prop} \label{prop: naive ideal}Assume that $x$ is $m$-generic for $m>2h_\lambda-3$. Let $R$ be a $p$-adically complete flat $\cO$-algebra, $N\in \{0,\infty\}$, and $Y\in L\cG_x^{\le \lambda}(R)$. 
Then the components of 
$(\frac{d}{dv})^t|_{v=-p}(\cN_N(Y))$ and $(\frac{d}{dv})^t|_{v=-p}(v^{-1}\cN_N(Y)_{\alpha})$ for $0\leq t \leq h_{\lambda}-2$ and $\alpha>0$ belong to $R\subset R[\frac{1}{p}]$.
The ideal generated by these elements form a functorial assignment
\[(R,Y) \mapsto I_{\nabla_{N,\nv}}(Y)\subset R\]
\end{prop}
\begin{proof}The statement for $N=0$ is obvious so we only need to deal with the case $N=\infty$.
Since
\[\frac{p^{i}}{(p+v^p)^{h_{\lambda}}\cdots (p+v^{p^i})^{h_{\lambda}}}\in \frac{1}{p^{(h_\lambda-1)i}}\Z_p[\![v,\frac{v^p}{p}]\!],\]
its $t$-th Taylor coefficient at $v=-p$ belongs to $p^{p-t-1-(h_\lambda-1)i}\Z_p$ for $t>0$.
It follows from Lemma \ref{lem:monodromy estimates}
that the $t$-th Taylor coefficient at $v=-p$ of 
\[\frac{p^i}{(p+v^p)^{h_{\lambda}}\cdots (p+v^{p^i})^{h_{\lambda}}} (\Ad^*(Y)\varphi_x)^i\big(\cN_0(Y)\big)\]
belongs to $p^{m\frac{p^i-1}{p-1}-i(h_\lambda-1)}R$.
Hence if $m>2h_\lambda-3$, the Taylor coefficients
\[(\frac{d}{dv})^t|_{v=-p}(\cN_\infty(Y))\]
\[(\frac{d}{dv})^t|_{v=-p}(v^{-1}\cN_\infty(Y)_{\alpha})\]
for $0\leq t \leq h_{\lambda}-2$ and $\alpha>0$ belong to $\chfg(R)$.
\end{proof}

This justifies the following
\begin{defn} Assume that $x$ is $(2h_\lambda-2)$-generic. For $N\in\{0,\infty\}$, the subfunctor $L\cG_x^{\le \lambda,\nabla_{N,\nv}}$ is the unique closed subfunctor of the $p$-adic completion $L\cG_x^{\le \lambda,\wedge_p}$ such that for all $p$-adically complete $\cO$-flat algebras $R$, $Y\in L\cG_x^{\le \lambda,\nabla_{N,\nv}}(R)$ if and only if $I_{\nabla_N,\nv}(Y)=0\subset R$.

We define $L\cG_x^{\le \lambda, \nabla_N}$ to be the $p$-saturation of $L\cG_x^{\le \lambda,\nabla_{N,\nv}}$.
\end{defn}

\begin{remark} (Naive monodromy vs deformed affine Springer fibers) \label{rmk:naive monodromy=AFS}
\begin{enumerate}
\item $L\cG_x^{\le \lambda,\nabla_{0,\nv}}$ can already be defined as a closed subfunctor of $L\cG_x^{\le \lambda}$ without needing to pass to $p$-adic completion.
\item Since $h_\lambda<p$, for $F\in R[\![v+p]\!]$, the condition that $(\frac{d}{dv})^t|_{v=-p}\big(F\big)=0$ for $0\leq t\leq h_{\lambda}-2$ is equivalent to $F\in (v+p)^{h_\lambda-1}R[\![v+p]\!]$. Thus $L\cG_x^{\le \lambda,\nabla_{0,\nv}}$ is cut out inside $L\cG_x^{\le \lambda}$ by the condition
\[(v+p)^{h_\lambda}(v\dlog(Y^{-1})+\Ad(Y^{-1})x-x)\in (v+p)^{h_\lambda-1}\Lie L^+\cG_x,\]
or equivalently
\[v(v+p)\dlog(Y^{-1})+\Ad(Y^{-1})((v+p)x)\in \Lie L^+\cG_x\]
In other words $[L\cG_x^{\le \lambda,\nabla_{0,\nv}}/L^+\cG_x]$ is exactly the deformed affine Springer fiber $\cX_{\gamma}^{\varepsilon=1}\cap \cS(\lambda)$ (where $\cS(\lambda)$ is the closure of the affine Schubert variety $S(\lambda)$ in the generic fiber of $\Gr_{\cG_x}$)) for the element $\gamma=\gamma(w,\mu)=-(v+p)x=(v+p)(w-\varphi)^{-1}(\mu)\in \chft[v]$ in Section \ref{ssec: deformed ASF}. In particular, its special fiber is $\rY_{\gamma}^{\varepsilon=1}(\le \lambda)$.
\item A simple computation shows that $L\cG_x^{\le \lambda,\nabla_{0,\nv}}$ is invariant under right translation by $L^+\cG_x$. A similar computation shows that $L\cG_x^{\le \lambda,\nabla_{\infty,\nv}}$ is invariant under the $x$-twisted $\varphi$-conjugation action of $(L^+\cG_x)^{\wedge_p}$ in Definition \ref{defn:twisted conjugation}. 
\end{enumerate} 
\end{remark}
\begin{prop} \label{prop:special fiber of approximations}Assume that $x$ is $m$-generic with $m>2h_\lambda-3$.
Then 
\[(L\cG_x^{\le \lambda,\nabla_{0,\nv}})_{\cO/p^{m-2h_\lambda+3}}=(L\cG_x^{\le \lambda,\nabla_{\infty,\nv}})_{\cO/p^{m-2h_\lambda+3}}\]
 as subfunctors of $(L\cG_x^{\le \lambda})_{\cO/p^{m-2h_\lambda+3}}$. 
\end{prop}
\begin{proof} 
The proof of Proposition \ref{prop: naive ideal} shows that the appropriate Taylor coefficients at $v=-p$ for all the terms occurring in $\cN_\infty(Y)$ other than $\cN_0(Y)$ are divisible by $p^{m-2h_\lambda+3}$.
\end{proof}

\begin{remark} (Fixing elementary divisors)
By Lemma \ref{lem: X lambda irred e neq 0} 
\[(L\cG_x^{\le \lambda,\nabla_{0,\nv}})_E=\coprod_{\lambda'\le\lambda, \lambda'\in X^*(T)^+} (L\cG_x^{\le \lambda,\nabla_{0,\nv}})_E\cap L\cG_x^{\lambda'}.\]
We define $L\cG_x^{\lambda',\nabla_0}$ to be the closure in $L\cG_x^{\le \lambda}$ of the factor corresponding to $\lambda'$.
This induces the decomposition of $[L\cG_x^{\le \lambda,\nabla_{0}}/L^+\cG_x]$ (which is $\cX^{\varepsilon=1}_{\gamma}(\le \lambda)$) as the scheme theoretic union of the closures $\cX^{\varepsilon=1}_{\gamma}(\lambda')$ of $\rX^{\varepsilon=1}_{\gamma}(\lambda')$.

We will see that this picture parallels the decomposition of $\cX^{\le \lambda,\tau}=[L\cG_x^{\le \lambda,\nabla_{\infty}}/_{\varphi,x} (L^+\cG_x)^{\wedge_p}]$ as the scheme theoretic union of the $\cX^{\lambda',\tau}$.

\end{remark}

\subsection{The homological model theorem}
Recall that we fixed $\lambda$, $x$, and that $\gamma=\gamma(w,\mu)=-(v+p)x$. We assume now that $x$ is $m$-generic for some $m>h_\lambda$.

By Proposition \ref{prop:special fiber of approximations}, $(L\cG_x^{\leq \lambda,\nabla_{0,\nv}})_{\F}=(L\cG_x^{\leq \lambda,\nabla_{\infty,\nv}})_{\F}$, and hence this common space is invariant under both the right action and the $x$-twisted $\varphi$-conjugation action of $(L^+\cG_x)_{\F}=\chI$. By \cite[Lemma 5.3.2]{Ei}, the $x$-twisted $\varphi$-conjugation action of $\chI_1$ and the right translation action are equivalent on $(L\cG_x^{\le \lambda})_\F$. The common quotient stack by either restricted action is then exactly the $\bT_{\F}$-torsor $\wt{\rY}_{\gamma}^{\varepsilon=1}(\le \lambda)$ over the deformed affine Springer fiber
$\rY_{\gamma}^{\varepsilon=1}(\le \lambda)=(\Gr^{\leq \lambda,\nabla_{0,\nv}}_{\cG_x})_{\F}$ obtained by pulling back along $\Fl_1=L\bG_\F/\chI_1 \rightarrow (\Gr_{\cG_x})_\F=\Fl$. 

Recall from section \ref{ssec: top-irred-comp generic} that we have the closed subschemes $\cX^{\varepsilon=1}_{\gamma}(\le \lambda)_\F$ and $\cX^{\varepsilon=1}_{\gamma}(\lambda')_\F$ of $\rY_{\gamma}^{\varepsilon=1}(\le \lambda)$ (for $\lambda'\le \lambda$ dominant).
 
\begin{thm} \label{thm:mod p model} Assume $x$ is $m$-generic for $m>2h_\lambda-2$. We have an equality
\[(L\cG_x^{\le \lambda,\nabla_0})_\F/\chI_1=[(L\cG_x^{\le \lambda,\nabla_\infty})_\F/_{\varphi,x}\chI_1]\]
as closed subschemes of $(L\cG_x^{\le \lambda})_\F/\chI_1$. In particular we have a local model diagram
\[
\begin{tikzcd}
& {(L\cG_x^{\le\lambda,\nabla_0})_\F/\chI_1 } \ar[dl, "{\bT_\F}"'] \ar[dr, "{\bT_\F}"] & \\
\cX^{\varepsilon=1}_{\gamma}(\le \lambda)_\F && \cX^{\le \lambda,\tau}_{\F}
\end{tikzcd}
\]
where the left arrow is the quotient by the right translation action of $\bT_\F$ and the right arrow is the quotient by the $\varphi$-twisted conjugation action of $\bT_\F$ (which coincides with the $\pi^{-1}$-twisted conjugation action).

Under this diagram, for each dominant $\lambda'\le \lambda$, $\cX^{\varepsilon=1}_{\gamma}(\lambda')_\F$ corresponds to $\cX^{\lambda',\tau}_\F$.
\end{thm}
The rest of this subsection is devoted to the proof of Theorem \ref{thm:mod p model}.

\begin{lemma}\label{lem:solving differential equation} Let $R$ be an $E$-algebra, and $B,C\in \chfg(R[\![v+p]\!])$. Then 
\begin{enumerate}
\item There is a unique solution $X\in \ker(\bG(R[\![v+p]\!])\surj \bG(R))$ to the equation
\[v\dlog(X^{-1})+\Ad(X^{-1})(B)=C\]
The assignment $(B,C)\mapsto X$ is a regular function on $(\Lie L^+\bG_E)^2$.
\item Suppose we have a subring $\mathring{R}\subset R$. Set $S_{\mathring{R}}\subset R[\![v+p]\!]$ to be the subring consisting of expressions
\[\sum_{i=0}^{\infty} a_i\frac{(v+p)^i}{i!p^i}\]
where $a_i\in \mathring{R}$.

Assume that $B,C\in \chfg(S_{\mathring{R}})$ and $B-C\in p^k \chfg(S_{\mathring{R}})$. Then $X\in \ker(\bG(S_{\mathring{R}})\to \bG(S_{\mathring{R}}/p^k))$.
\end{enumerate}
\end{lemma}
\begin{proof} The statement immediately reduces to the case $\bG=\GL_n$.
Writing $B=\sum B_i (v+p)^i, C=\sum C_i (v+p)^i$ and $X=1+\sum_{i=1}^{\infty}(v+p)^iX_i$ where $B_i,C_i,X_i\in \Mat_n(R[\![v+p]\!])$, our equation reduces to the recursion
\[p(i+1)X_{i+1}=iX_i+\sum_{j=0}^i X_jB_{i-j}-C_{i-j}X_{j}\]
for $i\geq 0$.
This uniquely solves for $X_{i+1}$ as a (non-commutative) polynomial $X_{i+1}(B,C)$ in terms of the coefficients $B_j, C_j$ for $j\leq i$, giving the first item. The second item follows from the observation that $X_i(B,C)$ belongs to the $2$-sided ideal generated by $B_j-C_j$ and an easy induction.
\end{proof}

\begin{lemma}\label{lem:monodromy congruent mod p} Assume $x$ is $m$-generic for some $m>2h_\lambda -2$. Let $\Lambda$ be a finite free $\cO$-algebra. 

\begin{enumerate}
\item Suppose $Y\in L\cG_x^{\leq \lambda,\nabla_{\infty}}(\Lambda)$. Then there exists $Y'\in L\cG_x^{\leq \lambda,\nabla_{0}}(\Lambda)$ such that $Y,Y'$ reduce to the same point in $L\cG_x^{\leq \lambda}(\Lambda/p^{m+2-2h_\lambda})$.
\item Conversely, suppose $Y\in L\cG_x^{\leq \lambda,\nabla_{0}}(\Lambda)$. Then there exists $Y'\in L\cG_x^{\leq \lambda,\nabla_{\infty}}(\Lambda)$ such that $Y,Y'$ reduce to the same point in $L\cG_x^{\leq \lambda}(\Lambda/p^{m+2-2h_\lambda})$.
\end{enumerate}
\end{lemma}

\begin{proof} 
\begin{enumerate}
\item 
Recall from Lemma \ref{lem:solving differential equation} the ring $S_{\Lambda}\subset \Lambda[\frac{1}{p}][\![v+p]\!]$ consisting of expressions
\[\sum_{i=0}^{\infty} a_i\frac{(v+p)^i}{i!p^i}\]
where $a_i\in \Lambda$. Note that $v\in pS_{\Lambda}$ and that
\[\frac{p^i}{(p+v^p)^{h_{\lambda}}\cdots (p+v^{p^i})^{h_{\lambda}}} \in \frac{1}{p^{(h_\lambda-1)i}}S_{\Lambda}^\times\]
We have
\[\cN_\infty(Y)=\cN_0(Y)+(v+p)^{h_\lambda}\Ad(Y^{-1})(Z) \]
where 
\[Z=f_\infty(Y)=\sum_{i=1}^{\infty}\frac{p^i}{(p+v^p)^{h_{\lambda}}\cdots (p+v^{p^i})^{h_{\lambda}}}  \varphi_x((\Ad^*(Y)\varphi_x)^{i-1}\cN_0(Y)). \]
By Lemma \ref{lem:monodromy estimates} and the fact that $\varphi_x(\cN_0(Y))\in v^{1+m}\chfg(\Lambda[\![v]\!])$, $Z\in p^{m+2-h_\lambda}\chfg(S_\Lambda)$.
Applying Lemma \ref{lem:solving differential equation} for $B=x$ and $C=x+Z$, we can find $X\in \ker(\bG(S_\Lambda)\to \bG(S_\Lambda/p^{m+2-h_\lambda})\times \bG(S_\Lambda/(v+p)))$ so that
\[v\dlog(X^{-1})+\Ad(X^{-1})x=x+Z,\]
i.e. $\cN_0(X)=(v+p)^{h_\lambda}Z$.
Hence 
\[\cN_0(XY)=\cN_0(Y)+\Ad(Y^{-1})\cN_0(X)=\cN_\infty(Y)\]
Set $\wt{X}=(X,1)\in L^+\cG_x(\Lambda[\frac{1}{p}])\subset \bG(\Lambda[\frac{1}{p}][\![v+p]\!])\times \bG(\Lambda[\![v]\!])$.
Since $Y\in L\cG_x^{\le \lambda,\nabla_\infty}(\Lambda)$, we learn that $\wt{X}Y\in L\cG_x^{\le \lambda,\nabla_0}(\Lambda[\frac{1}{p}])$.

We have a natural inclusion $\Lambda[\![v]\!]\to S_\Lambda\times \Lambda[\![v]\!]$ by expanding with respect to $v+p$ and $v$ respectively. Set $\bar S_{h_\lambda}=S_\Lambda/(S_\Lambda\cap (v+p)^{h_{\lambda}}\Lambda[\frac{1}{p}][\![v+p]\!])$, which (since $h_\lambda<p$) identifies with the subring $\{\sum_{i=0}^{h_\lambda-1}a_i\frac{(v+p)^i}{p^i}| a_i\in \Lambda\}\subset \Lambda[\frac{1}{p}][\![v+p]\!]/(v+p)^{h_\lambda}$. Then we have an inclusion
\begin{align*}
\Lambda[\![v]\!]/v(v+p)^{h_\lambda}&\inj \bar S_{h_\lambda}\times \Lambda[\![v]\!]/v\\
\sum_{i=0}^{h_\lambda} a_i(v+p)^i &\mapsto (\sum_{i=0}^{h_\lambda-1} a_i(v+p)^i,\sum_{i=0}^{h_\lambda} a_ip^i)
\end{align*}
We claim that the image $\bar X$ of $\wt{X}\in \bG(S_\Lambda)\times \bG(\Lambda[\![v]\!])$ in $\bG(\bar S_{h_\lambda}\times \Lambda[\![v]\!]/v)$ belongs to $\bG(\Lambda[\![v]\!]/v(v+p)^{h_\lambda})$ and reduces to $1$ mod $p^{m+2-2h_\lambda}$.
Indeed, choosing a faithful representation $\bG\inj \GL_n$, $\bar X$ pushes to a pair
\[(1+\sum_{i=1}^{h_\lambda-1} B_i\frac{(v+p)^i}{p^i},1)\subset \GL_n(\bar S_{h_\lambda})\times \GL_n(\Lambda)\]
with $B_i, C\in \Mat_n(\Lambda)$. But since $X$ pushes to $1$ in $\bG(\bar S_h/p^{m+2-h_\lambda})$, we see that $B_i\in p^{m+2-h_\lambda}\Lambda$.
Using that $m+2-h_\lambda> h_\lambda$, we see that the above pair is the image of
\[(1+ \sum_{i=1}^{h_\lambda-1} \frac{B_i}{p^i}(v+p)^i-\frac{\sum_{i=1}^{h_\lambda-1} B_i}{p^{h_\lambda}}(v+p)^{h_\lambda},1),\]
i.e. it comes from an element of $\ker (\GL_n(\Lambda[\![v]\!]/v(v+p)^{h_\lambda})\to \GL_n(\Lambda[\![v]\!]/(p^{m+2-2h_\lambda},v(v+p)^{h_\lambda})))$.

We can now lift $\bar X\in \ker(\bG(\Lambda[\![v]\!]/v(v+p)^{h_\lambda})\to \bG(\Lambda[\![v]\!]/(p^{m+2-2h_\lambda},v(v+p)^{h_\lambda})))$ to $A\in \ker(\bG(\Lambda[\![v]\!])\to \bG(\Lambda[\![v]\!]/p^{m+2-2h_\lambda}v)\subset L\cG_x(\Lambda)$, and obtain a factorization
\[\wt{X}Y=AKY\]
where $K\in L^+\cG_x(\Lambda[\frac{1}{p}])$ and $K\equiv 1$ mod $(v+p)^{h_\lambda}$. In particular $\Ad(Y^{-1})(K)\in L^+\cG_x(\Lambda[\frac{1}{p}])$. 
But since $AKY=AY\Ad(Y^{-1})(K)\in L\cG_x^{\le \lambda, \nabla_0}(\Lambda[\frac{1}{p}])$ and the property of being in $L\cG_x^{\le \lambda, \nabla_0}(\Lambda[\frac{1}{p}])$ is preserved under right translation by $\bG(\Lambda[\frac{1}{p}][\![v+p]\!])$, we see that $AY\in L\cG_x^{\le \lambda, \nabla_0}(\Lambda)$, and hence also belongs to $L\cG_x^{\le \lambda, \nabla_0}(\Lambda)$. We are now done since $AY\equiv Y$ mod $p^{m+2-2h_\lambda}$.
\item Let $Y\in L\cG_x^{\leq \lambda,\nabla_{0}}(\Lambda)$. Set $Y_1=Y$. 
We claim that we can recursively define a sequence $Y_n\in L\cG_x^{\le \lambda}(\Lambda)$, $A_n\in \ker (L^+\cG_x(\Lambda)\to L^+\cG_x(\Lambda/p^{n+m+1-2h_\lambda}))$ for $n\geq 1$ such that the following conditions are satisfied (where we interpret $f_\infty(Y_0)=0$): 
\begin{itemize}
\item $Y_{n+1}=A_nY_n$.
\item $\cN_0(Y_n)+\Ad(Y_n^{-1})(f_\infty(Y_{n-1}))\in (v+p)^{h_\lambda-1}\Lie L^+\cG_x(\Lambda[\frac{1}{p}])$.
\item $f_\infty(Y_{n})-f_\infty(Y_{n-1})\in p^{m+n+1-h_\lambda}\chfg(S_\Lambda)$.
\end{itemize}
This will finish the proof, since the sequence $Y_{n}$ converges in $L\cG_x^{\le \lambda}(\Lambda)$ to some $Y_\infty$ congruent to $Y$ mod $p^{m+2-2h_\lambda}$, and satisfies the condition
\[\cN_\infty(Y_\infty)=\cN_0(Y_\infty)+(v+p)^{h_\lambda}\Ad(Y_\infty^{-1})(f_\infty(Y_\infty))\in (v+p)^{h_\lambda-1}\Lie L^+\cG_x(\Lambda[\frac{1}{p}]),\]
so that $Y_\infty\in L\cG_x^{\le \lambda,\nabla_\infty}(\Lambda)$.

It remains to construct the $A_n$. Suppose we constructed $Y_n$. Then by Lemma \ref{lem:solving differential equation} and the third item in the inductive hypothesis we have, we can find $X_n\in \ker(\bG(S_\Lambda)\to \bG(S_\Lambda/p^{m+n+1-h_\lambda})\times \bG(S_\Lambda/(v+p)))$ so that
\[v\dlog(X_n^{-1})+\Ad(X_n^{-1})(x+f_\infty(Y_n))=x+f_\infty(Y_{n-1}),\]
i.e. $\cN_0(X_n)=(v+p)^{h_\lambda}f_\infty(Y_{n-1})-(v+p)^{h_\lambda}\Ad(X_n^{-1})(f_\infty(Y_{n}))$.

Now
\[\cN_0(X_nY_n)+(v+p)^{h_\lambda}\Ad((X_nY_n)^{-1})(f_\infty(Y_n))\]
\[=\cN_0(Y_n)+\Ad(Y_n^{-1})\big(\cN_0(X_n)+(v+p)^{h_\lambda}\Ad(X_n^{-1})(f_\infty(Y_n))\big)\]
\[=\cN_0(Y_n)+(v+p)^{h_\lambda}\Ad(Y_n^{-1})f_\infty(Y_{n-1})\in (v+p)^{h_\lambda-1}\Lie L^+\cG_x(\Lambda[\frac{1}{p}]).\]

Setting $\wt{X}_n=(X_n,1)$, the same argument as in the proof of the first part shows that the image of $\wt{X}_n$ in $\bG(\bar S_{h_\lambda}\times \Lambda[\![v]\!]/v)$ belongs to $\bG(\Lambda[\![v]\!]/v(v+p)^{h_\lambda})$ and reduces to $1$ modulo $p^{n+m+1-2h_\lambda}$. Hence we can find some $A_n\in \ker(\bG(\Lambda[\![v]\!])\to \bG(\Lambda[\![v]\!]/p^{n+m+1-2h_\lambda}v)\subset L\cG_x(\Lambda)$, and obtain a factorization
\[\wt{X}_nY_n=A_nKY_n=A_nY_n\Ad(Y_n^{-1})(K)\]
with $K\in L^+\cG_x(\Lambda[\frac{1}{p}])$ and $K\equiv 1$ mod $(v+p)^{h_\lambda}$. This guarantees that $Y_{n+1}=A_nY_n$ satisfies the first two of the above three item. In turn, the third item follows from Lemma \ref{lem:monodromy estimates}.
\end{enumerate}
\end{proof}
\begin{proof}[Proof of Theorem \ref{thm:mod p model}]
Note that $[L\cG_x^{\le \lambda,\nabla_{\infty}}/_{\varphi,x}(L_1^+\cG_x)^{\wedge_p}]$ and $[L\cG_x^{\le \lambda,\nabla_{0}}/L_1^+\cG_x]^{\wedge_p}$ are both topologically finite type $p$-adic formal stacks and are both analytically unramified. Indeed, for the first stack, these properties hold because they hold for the potentially crystalline stack $\cX^{\le \lambda,\tau}$ by \cite[Theorem 2.6.3]{Lin23b} and the corresponding property for deformation rings \cite{Ba}, while for the second stack it follows from Section \ref{ssec: top-irred-comp generic}.
By \cite[Corollary 15.2]{Bart23}, to show $[(L\cG_x^{\le \lambda,\nabla_{\infty}})_\F/_{\varphi,x}\chI_1]=(L\cG_x^{\le \lambda,\nabla_{0}})_\F/\chI_1$, it suffices to check they have the same Artinian points inside $(L\cG_x^{\le \lambda})_\F/\chI_1$.
By the argument in the proof of \cite[Lemma 7.2.6]{LLLM22}, for an Artinian $\cO$-algebra $A$, an $A$-point of either stack comes from a $\Lambda$-point for some finite free $\cO$-algebra $\Lambda$. But Lemma \ref{lem:monodromy congruent mod p} shows that we can find a $\Lambda$-point for the remaining stack which produce the same $\Lambda/p$-point of $(L\cG_x^{\le \lambda})_\F/\chI_1$, and hence the same $A$-point.

\end{proof}
\begin{remark} The argument in the proof in fact shows that $(L\cG^{\le \lambda,\nabla_0}_x)_{\cO/p^{m+2-2h_\lambda}}=(L\cG^{\le \lambda,\nabla_\infty}_x)_{\cO/p^{m+2-2h_\lambda}}$ and $(L\cG^{\lambda,\nabla_0}_x)_{\cO/p^{m+2-2h_\lambda}}=(L\cG^{\lambda,\nabla_\infty}_x)_{\cO/p^{m+2-2h_\lambda}}$ as subfunctors of $(L\cG_x)_{\cO/p^{m+2-2h_\lambda}}$. 
\end{remark}

Recall from Proposition \ref{prop: monomorphism} that we have a closed immersion
\[\iota_x: [(L\cG_x^{\le\lambda})_\F/_{\varphi,x}\chI]\inj \cR_{\lsup LG,\F}\]
where the Frobenius structure of the $\varphi_0$-module corresponding to $Y\in L\cG_x^{\le \lambda}$ is given by $\iota_x(Y)=(1\rtimes \sigma)\wt{w}(\tau)Y$.
On the other hand, Theorem \ref{thm:mod p model} shows that the restriction of $\iota_x$ to the substack $ [(L\cG_x^{\le \lambda,\nabla_0}/L^+_1\cG_x)_\F/_{\varphi,x}\bT_\F]$ factors through $\cX^{\le \lambda,\tau}_{\F}$, and in partcular through the map $\cX^{\EG}\to \cR_{\lsup LG}$. The following Lemma shows this factorization holds on the slightly larger locus where we replace $L\cG_x^{\le \lambda,\nabla_0}$ with $L\cG_x^{\le \lambda,\nabla_{0,\nv}}$:

\begin{prop}\label{prop:integrating Galois action} Assume that $x$ is $m$-generic with $m>h_{\lambda}$. Then there is a natural factorization
\[\iota_x: [\wt{\rY}^{\varepsilon=1}(\le \lambda)/_{\varphi,x} \bT_\F]\inj \cX^{\EG}\to \cR_{\lsup LG}\]
 
\end{prop}
\begin{proof} We use notations from subsection \ref{subsect:BKF modules}.
Suppose we have an $\F$-algebra $R$ and $Y\in L\cG_x^{\le\lambda,\nabla_{0,\nv}}(R)$ with $\iota(Y)$ giving rise to the \'{e}tale $\lsup LG$-torsor $\cM$ over $R(\!(v)\!)$. Recall the embedding ${R(\!(v)\!)}\inj \C^\flat_{R}:=\C^\flat\widehat{\otimes}_{\F_p} R$ given by sending $v$ to $(-p)^\flat$. We need to show that the obvious (i.e. only acting on the second factor) $G_{\Q_{p,\infty}}$-action on $\cM_\Inf:=\cM\otimes_{R(\!(v)\!)}\C_R^\flat$ extends (necessarily uniquely) to an action $G_{\Q_p}$ commuting with the Frobenius structure. To do this, it suffices to construct a natural Frobenius commuting semi-linear action of a generator $\tau$ of $\Gal(\Q_{p}^{\mathrm{ur}}\Q_{p,\infty}/\Q_p^{\mathrm{ur}}(\zeta_{p^\infty}))\cong \Z_p$.

We fix $\varepsilon\in C^{\flat}$ a compatible sequence of $p^\infty$-roots of $1$. We remind the reader that $C^{\flat}$ is a (non-discrete) valuation field with valuation ring $\cO_{C^{\flat}}$, and that $|\varepsilon-1|=|v|^{\frac{p}{p-1}}$. We choose the generator $\tau$ such that
\[\tau(v)=\varepsilon v, \tau(\varepsilon)=\varepsilon\]
Unwinding the commutation relation as in Lemma \ref{lemma: existence of Kisin lattice}, we need to construct an element $C_\tau\in \bG(\C^\flat_R)$ such that
\[C_\tau^{-1}=\tau(Y^{-1}\wt{w}(\tau)^{-1})\varphi(C_\tau)^{-1}\wt{w}(\tau)Y\]
Setting $C_\tau^{-1}=\varepsilon^xX$, using the fact $\tau(\varepsilon)=\varepsilon$ and $(\varphi-w)(x)=\mu$, the above relation is equivalent to
\[X=\tau(\varepsilon^{-x}Y^{-1}\varepsilon^x)\varphi_x(X)Y.\]
As in the proof of \cite[Proposition 18.1]{Bart23}, to construct such an $X$ it suffices to check
\begin{itemize}
\item Set $C$ to be the kernel of $\bG\to \GL(\chfg)\times \prod_{\chi:\bG\to \GL_1} \GL_1$, which is a finite central subgroup of $\bG$ of order prime to $p$. Set $\chI^\flat(R)$ to be the subgroup of $\bG(\cO_{\C^\flat,R})$ which is congruent to $1$ modulo $(\varepsilon-1)^{1/p}$ and congruent to the unipotent radical $\bN$ of $\bB$ modulo $\varepsilon-1$. 
Then
\[\tau(\varepsilon^{-x}Y^{-1}\varepsilon^x)Y\in C(\cO_{\C^\flat,R})\chI^\flat(R)\]
\item The twisted adjoint action $Z\mapsto \tau(\varepsilon^{-x}Y^{-1}\varepsilon^x)\varphi_x(Z)Y$ preserves and topologically nilpotent on $v^2\chfg(\cO_{\C^\flat,R})$. 
\end{itemize}
Indeed, the first item implies if we set $X_0=\tau(\varepsilon^{-x}Y^{-1}\varepsilon^x)Y$, then $\varphi_x(X_0)\in C(\cO_{\C^\flat,R})\ker(\bG(\cO_{\C^\flat,R})\to \bG(\cO_{\C^\flat,R}/v^{h_\lambda+2}))$. By replacing $\tau$ and $\varepsilon$ with a $\tau^m$, $\varepsilon^m$ with $m$ prime to $p$ depending only on $C$, we can get rid of the $C(\cO_{\C^\flat,R})$ factor, and hence assume that $\varphi_x(X_0)\in \ker(\bG(\cO_{\C^\flat,R})\to \bG(\cO_{\C^\flat,R}/v^{h_\lambda+2}))$. Thus
\[X_1=\tau(\varepsilon^{-x}Y^{-1}\varepsilon^x)\varphi_x(X_0)Y\in \ker(\bG(\cO_{\C^\flat,R})\to \bG(\cO_{\C^\flat,R}/v^{2})).\]
The second item now shows that if we continue to define inductively $X_{i+1}=\tau(\varepsilon^{-x}Y^{-1}\varepsilon^x)\varphi_x(X_i)Y$, then $X_i$ converges to some $X_\infty$ which satisfies the desired commutation relation. Furthermore, note that $X_\infty$ is the unique possible such elements that belong to $\chI^\flat(R)$, and hence its construction is functorial in $R$.

It remains to check the above two items. The second item is immediate from the fact that 
\[
\tau(\varepsilon^{-x}Y^{-1}\varepsilon^x)\varphi_x(v^k\chfg(\cO_{\C^\flat,R}))Y\subset v^{p(k-1)+m+1-h_\lambda}\chfg(\cO_{\C^\flat,R})
\]
and the fact that $m>h_\lambda$. 

For the first item, we first reduce to checking it after pushing $Y$ along the adjoint representation $\bG\to \GL(\chfg)$, since the statement is clearly true after pushing along any one-dimensional representation $\chi:\bG\to \mathbb{G}_m$ (due to the fact that $\tau(v^k)v^{-k}=\varepsilon^k\equiv 1$ mod $\varepsilon-1$). We can choose a standard choice of Borel $B_{\GL}$ and maximal torus $T_{\GL}$ of $\GL(\chfg)=\GL_d$ which is compatible with $\bB$ and $\bT$. This gives a choice of Iwahori subgroup $\chI_{\GL}$ and a pro-$v$ Iwahori subgroup $\chI_{1,\GL}$.
Let $Z\in \GL_d(R(\!(v)\!))$ be the image of $Y^{-1}$ in the adjoint representation. Let $D=v\frac{d}{dv}-[x,]$ be the differential operator on $\Mat_d( R[\![v]\!])$ which preserves $\Lie \chI_{\GL}(R)$ (where we abusively still write $x$ for its image under the adjoint representation).
The hypothesis that $Y\in L\cG_x^{\le \lambda,\nabla_{0,\nv}}(R)$ implies that
\begin{itemize}
\item $v^{h_\lambda}Z^{\pm 1}\in \Mat_d( R[\![v]\!])$.
\item $D(Z)Z^{-1}=J\in\frac{1}{v}\Lie \chI_{1,\GL}(R)$.
\end{itemize}
Set $Z'=v^{h_\lambda}Z\in \Mat_d( R[\![v]\!])$.
Now the image of $\tau(\varepsilon^{-x}Y^{-1}\varepsilon^x)Y$ is 
\[\tau(\varepsilon^{-x}v^{-h_\lambda}Z'\varepsilon^x)v^{h_\lambda}(Z')^{-1}=\varepsilon^{-h_\lambda}\tau(\varepsilon^{-x}Z'\varepsilon^x)(Z')^{-1}\]
Let $\log_p\varepsilon=\sum_{k=1}^{p-1}(-1)^{k-1} \frac{(\varepsilon-1)^k}{k}$ be the truncated logarithm. Then
\[\tau(\varepsilon^{-x}Z'\varepsilon^x)\in\sum_{k=0}^{p-1} \frac{1}{k!}D^k(Z')\log_p(\varepsilon)^k+(\varepsilon-1)^p\Mat_d(\cO_{C^{\flat},R})\]
An easy induction shows that $D^k(Z')\in \frac{1}{v^k} \Lie \chI_{\GL}(R)Z'$, hence
\[\tau(\varepsilon^{-x}Z'\varepsilon^x)(Z')^{-1}\in 1+\sum_{k=1}^{p-1}(\frac{\log_p\varepsilon}{v})^k\Lie \chI_{1,\GL}(R)+\Mat_d(\cO_{C^{\flat},R})(\varepsilon-1)^pZ'^{-1}\]
Since $p>2h_\lambda+2$ due to the existence of an $m$-generic $x$, the last term belongs to $(\varepsilon-1)\Mat_d(\cO_{C^{\flat},R})$. Since $\frac{\log_p\varepsilon}{v}\cO_{\C^\flat}=(\varepsilon-1)^{1/p}\cO_{\C^\flat}$, we have shown that the image of $\tau(\varepsilon^{-x}Y^{-1}\varepsilon^x)Y$ reduces to $1$ modulo $(\varepsilon-1)^{1/p}$ and reduces to $U_{\GL}$ modulo $\varepsilon-1$, which is what we needed.


\end{proof}

\subsection{Proof of Theorem \ref{thm: companion}}
We now apply the results in the previous section with $\lambda$ replaced by $\lambda+\rho$. Recall that 
\[\dim \rY^{\varepsilon=1}_\gamma(\le \lambda+\rho)=\dim \cX^{\lambda+\rho,\tau}_{\F}=\dim \cX_{\red}^{\EG}=\dim \bG/\bB\]
The pullback maps induce a chain of isomorphisms
\[ \topCh(\rY_\gamma^{\varepsilon=1}(\leq \lambda+\rho))\xrightarrow{\sim}  \topCh(\wt{\rY}_\gamma^{\varepsilon=1}(\leq \lambda+\rho))\xleftarrow{\sim}  \topCh([\wt{\rY}_\gamma^{\varepsilon=1}(\leq \lambda+\rho)/_{\varphi,x} \bT_\F])\]
Since $d$ is the top degree for all of these spaces, we have $\topCh(-) = \dZ(-)$ for each. On the other hand, Proposition \ref{prop:integrating Galois action}, shows that $\iota_x$ induces a natural map
\[(\iota_x)_*: \dZ([\wt{\rY}_\gamma^{\varepsilon=1}(\leq \lambda+\rho)/_{\varphi,x} \bT_\F])\inj \dZ(\cX_{\lsup LG,\red}).\]
Composing this with the above chain of isomorphisms gives the desired map
\[\transfer_\gamma \co  \dZ(\ul\rY_\gamma^{\varepsilon=1}(\leq \lambda+\rho)) \inj \dZ(\cX^{\EG}_{\red})\]
By Theorem \ref{thm:mod p model}, for each $\lambda'\le \lambda+\rho$,
\[ \transfer_\gamma([\cX^{\varepsilon=1}_{\gamma}(\lambda')_\F])=[\cX^{\lambda',\tau}_\F].\]
This gives the first part of Theorem \ref{thm: companion}, since $[\cX^{\varepsilon=1}_{\gamma}(\lambda')_\F]$ represents $\fsp_{p\to 0}[X_{\gamma}^{\varepsilon=1}(\lambda')]$ in $\topCh(\ul\rY_\gamma^{\varepsilon=1}(\leq \lambda+\rho))$ (note that both are $0$ if $\lambda'$ is not regular).

We now establish the second part of Theorem \ref{thm: companion}. Thus we are given a pair $(\wt{u},\wt{v})\in \ul{\wt{W}}_1\times \ul{\wt{W}}^+$ with $\wt{w}=\wt{v}^{-1}w_0\wt{u}\in \Adm^{\reg}(\lambda+\rho)$. Write $\wt{u}=v^{\rho_u}u$ and $\nu=\wt{v}^{-1}(0)$, so that $\wt{v}\in W t_{-\nu}$.
By Lemma \ref{lem:transformation of components}, for every $z\in W$, there is a suitable $\gamma_z$ such that
\[\rY^{\varepsilon=1}_\gamma(\wt{w})=\wt{v}^{-1}z\rY^{\varepsilon=1}_{\gamma_z}(w_0\wt{u}),\]
and that the right-hand side has a $\bT_\F$-stable open dense subset belonging to $ \wt{v}^{-1}zv^{w_0\rho_u} L\bN_\F w_0u \chI/\chI$.
Equivalently, after adjusting $z$ and using $\wt{u}=v^{\rho_u}u$, the pullback $\wt{\rY}^{\varepsilon=1}_\gamma(\wt{w})$ of $\rY^{\varepsilon=1}_\gamma(\wt{w})$ to $\Fl_1$ has an open dense subset $U_z$ belonging to
\[\bT_\F v^{\nu}zw_0 L\bN_\F w_0 \wt{u} \chI_1/\chI_1.\]
Thus the \'{e}tale $\varphi_0$-modules in $\iota_x(U_z)$ have Frobenius structure of the form
\[(1\rtimes \sigma)\wt{w}(\tau)\bT_\F v^{\nu}zw_0 L\bN_\F w_0 \wt{u}\]
which in turn is equivalent to 
\[(1\rtimes \sigma)\varphi(\wt{u})\wt{w}(\tau) v^{\nu}zw_0 \bT_\F L\bN_\F w_0\]
We now choose $z$ so that $\varphi(\wt{u})\wt{w}(\tau) v^{\nu}z=v^{\kappa}$, then
\[\kappa=\pi^{-1}(\wt{u})\bup(t^{\mu}w \wt{v}^{-1}(0)-\rho)+\rho,\]
Thus a $\varphi_0$-module in $\iota_x(U_z)$ has Frobenius structure of the form
\[(1\rtimes \sigma)v^{\kappa}w_0 \bT_\F L\bN_\F w_0,\]
and hence has corresponding parameter $G_{\Q_{p,\infty}}\to \lsup LB$ whose composition with the projection $\lsup LB \to \lsup  LT$ and restriction to inertia is the parameter with inertial weight $\kappa$.

On the other hand, by Proposition \ref{prop:integrating Galois action}, $\iota_x(\wt{\rY}^{\varepsilon=1}_\gamma(\wt{w}))$ comes from  $\cX^{\EG}_{\lsup LG, \F}$, and contains an open dense substack $\iota_x(U_z)$ of dimension $d=\dim \chG/\chB$ consisting of $\lsup LB$-valued parameters whose semisimplification has inertial weight $\kappa$. On the other hand, since $\kappa$ is $2$-generic, for $\sigma=F(\kappa-\rho)$, the closed substack $\cC_\sigma$ is defined as the scheme-theoretic image of a monomorphism
\[\cX_{\lsup LB,\F}^{\kappa}\to \cX_{\lsup LG,\F}^{\EG}\] 
where $\cX_{\lsup LB,\F}^{\kappa}$ is the stack of $\lsup LB$-valued parameters whose semisimplification has inertial weight $\kappa$. Since $\kappa$ is $2$-generic, $\cX_{\lsup LB,\F}^{\kappa}$ is irreducible and smooth of dimension $d$. Since $\iota_x(U_z)$ lies in the image of $\cX_{\lsup LB,\F}^{\kappa}$ and have the same dimension $d$, we conclude that $\iota_x(\wt{\rY}^{\varepsilon=1}_\gamma(\wt{w}))$ must coincide with $\cC_\sigma$.
This shows that
\[\transfer_\gamma \left( [\ulY^{\varepsilon=1}_{\gamma}(\wt{w})] \right)  =  [\cC_{\sigma}]\]
with $\sigma=F(\kappa-\rho)=F(\pi^{-1}(\wt{u})\bup(t^{\mu}w \wt{v}^{-1}(0)-\rho))$ as desired.

\bibliographystyle{amsalpha}
\bibliography{Bibliography}

\newcommand{\etalchar}[1]{$^{#1}$}
\providecommand{\bysame}{\leavevmode\hbox to3em{\hrulefill}\thinspace}
\providecommand{\MR}{\relax\ifhmode\unskip\space\fi MR }
\providecommand{\MRhref}[2]{%
  \href{http://www.ams.org/mathscinet-getitem?mr=#1}{#2}
}
\providecommand{\href}[2]{#2}
\begin{thebibliography}{BBAMY25}

\bibitem[ABG04]{ABG04}
Sergey Arkhipov, Roman Bezrukavnikov, and Victor Ginzburg, \emph{Quantum
  groups, the loop {G}rassmannian, and the {S}pringer resolution}, J. Amer.
  Math. Soc. \textbf{17} (2004), no.~3, 595--678.

\bibitem[AKL{\etalchar{+}}22]{AKLPR}
Dhyan Aranha, Adeel~A. Khan, Alexei Latyntsev, Hyeonjun Park, and Charanya
  Ravi, \emph{Localization theorems for algebraic stacks}, 2022.

\bibitem[Bal12]{Ba}
Sundeep Balaji, \emph{G-valued potentially semi-stable deformation rings},
  ProQuest LLC, Ann Arbor, MI, 2012, Thesis (Ph.D.)--The University of Chicago.

\bibitem[BAL24]{BL21}
Pablo Boixeda~Alvarez and Ivan Losev, \emph{Affine {S}pringer fibers, {P}rocesi
  bundles, and {C}herednik algebras}, Duke Math. J. \textbf{173} (2024), no.~5,
  807--872, With an appendix by Boixeda Alvarez, Losev and O. Kivinen.

\bibitem[Bar23]{Bart23}
Robin Bartlett, \emph{Cycles relations in the affine grassmannian and
  applications to {B}reuil--{M}ézard for {G}-crystalline representations},
  2023.

\bibitem[BBAMY22]{BBMY}
Roman Bezrukavnikov, Pablo Boixeda-Alvarez, Michael McBreen, and Zhiwei Yun,
  \emph{Non-abelian {H}odge moduli spaces and homogeneous affine {S}pringer
  fibers}, 2022.

\bibitem[BBAMY25]{BBMY2}
\bysame, \emph{Affine {S}pringer fiber and the small quantum group}, 2025.

\bibitem[BBASV22]{BBSV}
Roman Bezrukavnikov, Pablo Boixeda-Alvarez, Peng Shan, and Eric Vasserot,
  \emph{A geometric realization of the center of the small quantum group},
  2022.

\bibitem[Bez16]{Bez16}
Roman Bezrukavnikov, \emph{On two geometric realizations of an affine {H}ecke
  algebra}, Publ. Math. Inst. Hautes \'{E}tudes Sci. \textbf{123} (2016),
  1--67.

\bibitem[BF08]{BF08}
Roman Bezrukavnikov and Michael Finkelberg, \emph{Equivariant {S}atake category
  and {K}ostant-{W}hittaker reduction}, Mosc. Math. J. \textbf{8} (2008),
  no.~1, 39--72, 183.

\bibitem[BG14]{BG14}
Kevin Buzzard and Toby Gee, \emph{The conjectural connections between
  automorphic representations and {G}alois representations}, Automorphic forms
  and {G}alois representations. {V}ol. 1, London Math. Soc. Lecture Note Ser.,
  vol. 414, Cambridge Univ. Press, Cambridge, 2014, pp.~135--187.

\bibitem[BHS19]{BHS19}
Christophe Breuil, Eugen Hellmann, and Benjamin Schraen, \emph{A local model
  for the trianguline variety and applications}, Publ. Math. Inst. Hautes
  \'{E}tudes Sci. \textbf{130} (2019), 299--412.

\bibitem[BL23]{BL23}
Jeremy Booher and Brandon Levin, \emph{{$G$}-valued crystalline deformation
  rings in the {F}ontaine-{L}affaille range}, Compos. Math. \textbf{159}
  (2023), no.~8, 1791--1832.

\bibitem[BM02]{BM02}
Christophe Breuil and Ariane M\'{e}zard, \emph{Multiplicit\'{e}s modulaires et
  repr\'{e}sentations de {${\rm GL}_2({\bf Z}_p)$} et de {${\rm
  Gal}(\overline{\bf Q}_p/{\bf Q}_p)$} en {$l=p$}}, Duke Math. J. \textbf{115}
  (2002), no.~2, 205--310, With an appendix by Guy Henniart.

\bibitem[BM13]{BM13}
Roman Bezrukavnikov and Ivan Mirkovi\'{c}, \emph{Representations of semisimple
  {L}ie algebras in prime characteristic and the noncommutative {S}pringer
  resolution}, Ann. of Math. (2) \textbf{178} (2013), no.~3, 835--919, With an
  appendix by Eric Sommers.

\bibitem[BMR06]{BMR06}
Roman Bezrukavnikov, Ivan Mirkovi\'{c}, and Dmitriy Rumynin, \emph{Singular
  localization and intertwining functors for reductive {L}ie algebras in prime
  characteristic}, Nagoya Math. J. \textbf{184} (2006), 1--55.

\bibitem[BMR08]{BMR08}
\bysame, \emph{Localization of modules for a semisimple {L}ie algebra in prime
  characteristic}, Ann. of Math. (2) \textbf{167} (2008), no.~3, 945--991, With
  an appendix by Bezrukavnikov and Simon Riche.

\bibitem[CEGS25]{CEGS}
Ana Caraiani, Matthew Emerton, Toby Gee, and David Savitt, \emph{The geometric
  {B}reuil-{M}\'{e}zard conjecture for two-dimensional potentially
  {B}arsotti-{T}ate {G}alois representations}, Algebra Number Theory
  \textbf{19} (2025), no.~2, 287--312. \MR{4859067}

\bibitem[CG10]{CG10}
Neil Chriss and Victor Ginzburg, \emph{Representation theory and complex
  geometry}, Modern Birkh\"{a}user Classics, Birkh\"{a}user Boston, Ltd.,
  Boston, MA, 2010, Reprint of the 1997 edition.

\bibitem[DJK21]{DJK21}
Fr\'{e}d\'{e}ric D\'{e}glise, Fangzhou Jin, and Adeel~A. Khan,
  \emph{Fundamental classes in motivic homotopy theory}, J. Eur. Math. Soc.
  (JEMS) \textbf{23} (2021), no.~12, 3935--3993.

\bibitem[DL76]{DL76}
P.~Deligne and G.~Lusztig, \emph{Representations of reductive groups over
  finite fields}, Ann. of Math. (2) \textbf{103} (1976), no.~1, 103--161.

\bibitem[EG23]{EG23}
Matthew Emerton and Toby Gee, \emph{Moduli stacks of \'{e}tale ({$\varphi,
  \Gamma$})-modules and the existence of crystalline lifts}, Annals of
  Mathematics Studies, vol. 215, Princeton University Press, Princeton, NJ,
  [2023] \copyright 2023.

\bibitem[EG24]{EGSurvey}
\bysame, \emph{Moduli stacks of {$(\varphi,\Gamma)$}-modules: a survey},
  Arithmetic geometry, Tata Inst. Fundam. Res. Stud. Math., vol.~41, Tata Inst.
  Fund. Res., Mumbai, [2024] \copyright 2024, pp.~229--306. \MR{4812705}

\bibitem[FH25]{FH}
Tony Feng and Michael Harris, \emph{Derived structures in the {L}anglands
  {C}orrespondence}, 2025, To appear in the Proceedings of the 2022 IHES Summer
  School on the Langlands program.

\bibitem[Fre07]{F07}
Edward Frenkel, \emph{Langlands correspondence for loop groups}, Cambridge
  Studies in Advanced Mathematics, vol. 103, Cambridge University Press,
  Cambridge, 2007.

\bibitem[Ful98]{Ful98}
William Fulton, \emph{Intersection theory}, second ed., Ergebnisse der
  Mathematik und ihrer Grenzgebiete. 3. Folge. A Series of Modern Surveys in
  Mathematics [Results in Mathematics and Related Areas. 3rd Series. A Series
  of Modern Surveys in Mathematics], vol.~2, Springer-Verlag, Berlin, 1998.

\bibitem[FYZ]{FYZ3}
Tony Feng, Zhiwei Yun, and Wei Zhang, \emph{Modularity of {H}igher theta series
  {I}: {C}ohomology of the generic fiber}.

\bibitem[FZ10]{FZ10}
Edward Frenkel and Xinwen Zhu, \emph{Any flat bundle on a punctured disc has an
  oper structure}, Math. Res. Lett. \textbf{17} (2010), no.~1, 27--37.

\bibitem[GHS18]{GHS18}
Toby Gee, Florian Herzig, and David Savitt, \emph{General {S}erre weight
  conjectures}, J. Eur. Math. Soc. (JEMS) \textbf{20} (2018), no.~12,
  2859--2949.

\bibitem[GK14]{GK14}
Toby Gee and Mark Kisin, \emph{The {B}reuil-{M}\'{e}zard conjecture for
  potentially {B}arsotti-{T}ate representations}, Forum Math. Pi \textbf{2}
  (2014), e1, 56.

\bibitem[GKM04]{GKM04}
Mark Goresky, Robert Kottwitz, and Robert Macpherson, \emph{Homology of affine
  {S}pringer fibers in the unramified case}, Duke Math. J. \textbf{121} (2004),
  no.~3, 509--561.

\bibitem[GKM06]{GKM06}
Mark Goresky, Robert Kottwitz, and Robert MacPherson, \emph{Purity of
  equivalued affine {S}pringer fibers}, Represent. Theory \textbf{10} (2006),
  130--146.

\bibitem[GLS15]{GLS15}
Toby Gee, Tong Liu, and David Savitt, \emph{The weight part of {S}erre's
  conjecture for {$\mathrm{GL}(2)$}}, Forum Math. Pi \textbf{3} (2015), e2, 52.

\bibitem[Hje24]{Ei}
Eivind~Otto Hjelle, \emph{The {M}oduli {S}tack of {B}reuil-{K}isin {M}odules
  with {D}escent {D}ata for {R}eductive {G}roups}, ProQuest LLC, Ann Arbor, MI,
  2024, Thesis (Ph.D.)--Northwestern University.

\bibitem[Jan03]{Jan03}
Jens~Carsten Jantzen, \emph{Representations of algebraic groups}, second ed.,
  Mathematical Surveys and Monographs, vol. 107, American Mathematical Society,
  Providence, RI, 2003.

\bibitem[Jan04]{Jan04}
\bysame, \emph{Representations of {L}ie algebras in positive characteristic},
  Representation theory of algebraic groups and quantum groups, Adv. Stud. Pure
  Math., vol.~40, Math. Soc. Japan, Tokyo, 2004, pp.~175--218.

\bibitem[Kis06]{Kis06}
Mark Kisin, \emph{Crystalline representations and {$F$}-crystals}, Algebraic
  geometry and number theory, Progr. Math., vol. 253, Birkh\"{a}user Boston,
  Boston, MA, 2006, pp.~459--496.

\bibitem[Kis08]{Kis08}
\bysame, \emph{Potentially semi-stable deformation rings}, J. Amer. Math. Soc.
  \textbf{21} (2008), no.~2, 513--546.

\bibitem[Kis09]{Kis09}
\bysame, \emph{The {F}ontaine-{M}azur conjecture for {${\rm GL}_2$}}, J. Amer.
  Math. Soc. \textbf{22} (2009), no.~3, 641--690.

\bibitem[KL88]{KL88}
D.~Kazhdan and G.~Lusztig, \emph{Fixed point varieties on affine flag
  manifolds}, Israel J. Math. \textbf{62} (1988), no.~2, 129--168.

\bibitem[Lee23]{Lee23}
Heejong Lee, \emph{Emerton--{G}ee stacks, {S}erre weights, and
  {B}reuil--{M}ézard conjectures for $\mathrm{GSp}_4$}, 2023.

\bibitem[LH]{companion}
Bao~V. Le~Hung, \emph{Local models for non-generic potentially crystalline
  {E}merton-{G}ee stacks}, in preparation.

\bibitem[LHLM23a]{LLLMextremal}
Daniel Le, Bao V.~Le Hung, Brandon Levin, and Stefano Morra, \emph{Extremal
  weights and a tameness criterion for mod $p$ {G}alois representations}, 2023.

\bibitem[LHLM23b]{LLLM22}
\bysame, \emph{Local models for {G}alois deformation rings and applications},
  Invent. Math. \textbf{231} (2023), no.~3, 1277--1488.

\bibitem[Lin23a]{Lin23c}
Zhongyipan Lin, \emph{A {D}eligne-{L}usztig type correspondence for tame
  $p$-adic groups}, 2023.

\bibitem[Lin23b]{Lin23a}
\bysame, \emph{The {E}merton-{G}ee stacks for tame groups}, 2023.

\bibitem[Lin23c]{Lin23b}
\bysame, \emph{The {E}merton-{G}ee stacks for tame $p$-adic groups {I}{I}},
  2023.

\bibitem[Liu08]{Liu08}
Tong Liu, \emph{On lattices in semi-stable representations: a proof of a
  conjecture of {B}reuil}, Compos. Math. \textbf{144} (2008), no.~1, 61--88.

\bibitem[LLHL19]{LLL19}
Daniel Le, Bao~V. Le~Hung, and Brandon Levin, \emph{Weight elimination in
  {S}erre-type conjectures}, Duke Math. J. \textbf{168} (2019), no.~13,
  2433--2506.

\bibitem[LLHLM18]{LLLMold}
Daniel Le, Bao~V. Le~Hung, Brandon Levin, and Stefano Morra, \emph{Potentially
  crystalline deformation rings and {S}erre weight conjectures: shapes and
  shadows}, Invent. Math. \textbf{212} (2018), no.~1, 1--107.

\bibitem[LLHLM20]{LLLM20}
\bysame, \emph{Serre weights and {B}reuil's lattice conjecture in dimension
  three}, Forum Math. Pi \textbf{8} (2020), e5, 135.

\bibitem[Lus96]{Lus96}
George Lusztig, \emph{Affine {W}eyl groups and conjugacy classes in {W}eyl
  groups}, Transform. Groups \textbf{1} (1996), no.~1-2, 83--97.

\bibitem[Min24]{Min24}
Yu~Min, \emph{Classicality of derived {E}merton--{G}ee stack {I}{I}:
  generalised reductive groups}, 2024.

\bibitem[MR99]{MR99}
Ivan Mirkovi\'{c} and Dmitriy Rumynin, \emph{Centers of reduced enveloping
  algebras}, Math. Z. \textbf{231} (1999), no.~1, 123--132.

\bibitem[Ngo10]{Ngo10}
Bao~Chau Ngo, \emph{Le lemme fondamental pour les alg\`ebres de {L}ie}, Publ.
  Math. Inst. Hautes \'{E}tudes Sci. (2010), no.~111, 1--169.

\bibitem[NS20]{NS20}
David Nadler and Vivek Shende, \emph{Sheaf quantization in {W}einstein
  symplectic manifolds}, 2020.

\bibitem[Pas15]{Pas15}
Vytautas Paskunas, \emph{On the {B}reuil-{M}\'{e}zard conjecture}, Duke Math.
  J. \textbf{164} (2015), no.~2, 297--359.

\bibitem[PZ13]{PZ13}
G.~Pappas and X.~Zhu, \emph{Local models of {S}himura varieties and a
  conjecture of {K}ottwitz}, Invent. Math. \textbf{194} (2013), no.~1,
  147--254.

\bibitem[Ric08]{Ric08}
Simon Riche, \emph{Geometric braid group action on derived categories of
  coherent sheaves}, Represent. Theory \textbf{12} (2008), 131--169, With a
  joint appendix with Roman Bezrukavnikov.

\bibitem[Sag97]{Sa97}
D.~S. Sage, \emph{A construction of representations of affine {W}eyl groups},
  Compositio Math. \textbf{108} (1997), no.~3, 241--245.

\bibitem[Ser87]{Ser87}
Jean-Pierre Serre, \emph{Sur les repr\'{e}sentations modulaires de degr\'{e}
  {$2$} de {$\mathrm{Gal}(\overline{\mathbb{Q}}/\mathbb{Q})$}}, Duke Math. J.
  \textbf{54} (1987), no.~1, 179--230.

\bibitem[Ser05]{Ser05}
\bysame, \emph{Compl\`ete r\'{e}ductibilit\'{e}}, no. 299, 2005, S\'{e}minaire
  Bourbaki. Vol. 2003/2004, pp.~Exp. No. 932, viii, 195--217.

\bibitem[Yun14]{Yun14}
Zhiwei Yun, \emph{The spherical part of the local and global {S}pringer
  actions}, Math. Ann. \textbf{359} (2014), no.~3-4, 557--594.

\bibitem[Yun17]{Yun17}
\bysame, \emph{Lectures on {S}pringer theories and orbital integrals}, Geometry
  of moduli spaces and representation theory, IAS/Park City Math. Ser.,
  vol.~24, Amer. Math. Soc., Providence, RI, 2017, pp.~155--215.

\end{thebibliography}




\end{document}